\documentclass[a4paper]{amsart}

\usepackage{amsmath}
\usepackage{amsthm}
\usepackage{mathtools}
\usepackage{xspace}
\usepackage[psamsfonts]{amssymb}
\usepackage[utf8]{inputenc}
\usepackage{float}\restylefloat{figure}
\usepackage{color}
\usepackage{tikz}\usetikzlibrary{arrows}
\graphicspath{{img/}}
\usepackage{enumitem}
\usepackage{hyperref}
\usepackage{cleveref}
\usepackage{tikz}
\usepackage{tikz-cd} 

\theoremstyle{plain}
\newtheorem{lemma}{Lemma}
\newtheorem{theorem}[lemma]{Theorem}
\newtheorem{proposition}[lemma]{Proposition}
\newtheorem{definition}[lemma]{Definition}
\newtheorem{corollary}[lemma]{Corollary}
\theoremstyle{remark}
\newtheorem{remarkn}[lemma]{Remark}
\newtheorem*{remark}{Remark}

\renewcommand{\epsilon}{\varepsilon}
\renewcommand{\emptyset}{\varnothing}
\renewcommand{\Re}{{\operatorname{Re}\,}}
\renewcommand{\Im}{{\operatorname{Im}\,}}

\newcommand{\bEA}{\begin{eqnarray*}}
\newcommand{\eEA}{\end{eqnarray*}}
\newcommand{\bEAn}{\begin{eqnarray}}
\newcommand{\eEAn}{\end{eqnarray}}

\newcommand{\ds}{\displaystyle}
\newcommand{\tend}{\longrightarrow}

\newcommand{\wt}{\widetilde}
\newcommand{\wh}{\widehat}
\newcommand{\on}{\operatorname}

\newcommand{\C}{{\mathbb C}}
\newcommand{\D}{{\mathbb D}}
\renewcommand{\H}{{\mathbb H}} % \H is by default "long Hungarian umlaut"
\newcommand{\N}{{\mathbb N}}

\newcommand{\R}{{\mathbb R}}
\newcommand{\Z}{{\mathbb Z}}

\newcommand{\eps}{\epsilon}

\newcommand{\asin}{\operatorname{asin}}

\newcommand{\cst}{\operatorname{cst}}
\newcommand{\wtends}{\xrightharpoonup{\phantom{m}}}
\newcommand{\la}{\left\langle}
\newcommand{\ra}{\right\rangle}
\newcommand{\loc}{{\operatorname{loc}}}
\newcommand{\de}{\partial}
\newcommand{\bd}{\bar\partial}
\newcommand{\Leb}{\operatorname{Leb}}
\newcommand{\avg}{\operatorname{avg}}
\newcommand{\of}{{\,:\,}}
\newcommand{\prA}{\pi}
\newcommand{\prB}{F}
\newcommand{\res}{\operatorname{res}}
\newcommand{\Aff}{\operatorname{Aff}}
\newcommand{\Sconf}{\mathrm{S\text{-}Conf}}
\newcommand{\Conf}{\mathrm{Conf}}
\newcommand{\Sk}{\operatorname{\mathit{Sk}}}
\newcommand{\Int}{\operatorname{Int}}

\newcommand{\Per}{\operatorname{Per}}
\newcommand{\Glu}{\operatorname{Glu}}
\newcommand{\Eff}{\mathrm{Eff}}

\newcommand{\HU}{\mathrm{H\tilde Z}}

\newcommand{\cal}[1]{\mathcal{#1}}
\newcommand{\ov}[1]{\overline{#1}} 
\newcommand{\setof}[2]{\big\{{#1}\,\big|\,{#2}\big\}}
\newcommand{\bigsetof}[2]{\left\{{#1}\,\middle|\,{#2}\right\}}
\newcommand{\pp}[2]{\frac{\partial{#1}}{\partial{#2}}}
\newcommand{\loctitle}[1]{\par\medskip\noindent\textit{#1}
\nopagebreak\par\smallskip}
\newcommand{\avgx}[1]{\operatorname{\underset{#1}{\avg}}}

\title[Similarity surfaces, Connections and the MRMT]{Similarity surfaces, connections and the measurable Riemann mapping theorem} 

\author{Arnaud Chéritat}
\address[Arnaud Chéritat]{C.N.R.S./Institut Mathématique de Toulouse, UMR 5219, Université de Toulouse, France}
\email{arnaud.cheritat@math.univ-toulouse.fr}
\author{Guillaume Tahar}
\address[Guillaume Tahar]{Beijing Institute of Mathematical Sciences and Applications, Huairou District, Beijing, China}
\email{guillaume.tahar@bimsa.cn}

\begin{document}

\begin{abstract}
    This article studies a particular process that approximates solutions of the Beltrami equation (straightening of ellipse fields, a.k.a.\ measurable Riemann mapping theorem) on $\C$.
    Introducing a sequence of similarity surfaces constructed by gluing polygons, we explain the relation between their conformal uniformization and the Schwarz-Christoffel formula.
    Numerical aspects, in particular the efficiency of the process, are not studied, but we draw interesting theoretical consequences.
    First, we give an independent proof of the analytic dependence, on the data (the Beltrami form), of the solution of the Beltrami equation (Ahlfors-Bers theorem).
    For this we prove, without using the Ahlfors-Bers theorem, the holomorphic dependence, with respect to the polygons, of the Christoffel symbol appearing in the Schwarz-Christoffel formula.
    Second, these Christoffel symbols define a sequence of parallel transports on the range, and in the case of a data that is $C^2$ with compact support, we prove that it converges to the parallel transport associated to a particular affine connection, which we identify.
\end{abstract}


\maketitle

\setcounter{tocdepth}{1}
\tableofcontents

%
% ----------------------------------------------------------
%

\section*{Introduction}

More than fifty years later, Chapter~II of \cite{A} remains an excellent introduction to the theory of quasiconformal mappings and we highly recommend it.\footnote{Or the second edition \cite{Abis}, which was issued in 2006.}
This book also includes in Chapter~V a proof of the measurable Riemann mapping theorem (see \Cref{sec:mrmt} here for a statement). This proof makes use of integral operators with singular kernels (Ahlfors-Beurling operator), $L^p$ spaces for $p>2$, letting $p\tend 2$, and an inequality due to Calderón and Zygmund.
There has been prior and posterior proofs with different approaches in special cases or in the general case, depending on the data $\mu$ in the statement (the article \cite{Glu} gives a short historical overview and references):
\begin{itemize}
\item Gauss (1825) for an $\R$-analytic $\mu$, by complexifying $\R^2$ into $\C^2$ and using a clever trick that does not extend to the $C^\infty$ class;
\item Korn (1914) and Lichtenstein (1916) for Hölder-continuous $\mu$; Lichtenstein's method already involves integral operators. Korn uses a fixed point method involving solving the Laplacian;
\item Lavrentiev (1935) for a continuous $\mu$, by constructing approximations to the solution using conformal geometric methods (we sum-up the method in \Cref{part:beltrami}) \cite{Lav1,Lav2};
\item Morrey (1936) for a general, measurable, $\mu$, using a density argument to reduce to solving the case where $\mu$ is analytic, \cite{Mo};
\end{itemize}
to cite only a few.
According to \cite{Bo1} (translated in \cite{Bo2}), the idea to introduce the Ahlfors-Beurling operator (it did not bear that name at that time) for this problem is due to Vekua \cite{V}.
Concerning the holomorphic dependence on the data $\mu$, all the proofs that the authors of the present article know use the Ahlfors-Beurling operator. These works include:
\begin{itemize}
\item Ahlfors and Bers, via $L^p$ spaces
\item Buff-Douady, via $L^2$ spaces only and the Fourier transformation
\item Glutsyuk, via a limit of a version for the torus, which is proved using Fourier series.
\end{itemize}
\Cref{part:beltrami} in the present article is an addition to this list, with the difference that we do not use the Ahlfors-Beurling operator.
Our proof uses distributions and $L^2$ spaces for the definition of quasiconformal maps, a classical compactness argument via an $L^2$ estimate, a new ingredient that is similarity surfaces, and the Poincaré-Koebe theorem.\footnote{\ldots\ Of which there are several proofs that do not use Beltrami forms. See for instance Chapter~10 of \cite{Ah4}, or \cite{Hu}.}. It bears resemblance with the approach of Lavrentiev and in fact the article of Lavrentiev avoids the use of the Poincaré-Koebe theorem so it is not impossible that the proof we present may be adapted too to avoid this use, but this seems not obvious.

\loctitle{Structure and content of the article}

\Cref{part:simsurf} is an introduction to the subject of similarity surfaces with a focus on ones that are conformally equivalent to punctured Riemann spheres and obtained by gluing finitely many convex polygons, possibly unbounded (\Cref{sec:simsurfpolygon}).
In \Cref{sec:sc} we draw a link with a generalization of the Schwarz-Christoffel formula: it expresses similarity charts. We believe this fact was known before but we do not have references for this.
In the formula a particular rational map appears, which is the Christoffel symbol (expression in a chart) of a particular conformal connection.
\Cref{sec:holodep} is devoted to a proof of holomorphic dependence of the rational map appearing as a Christoffel symbol, when the polygons that are glued are modified by letting the vertices vary holomorphically.
A key point is to completely avoid using the measurable Riemann mapping theorem in this proof. Such an insistance induces complications that are dealt with in \Cref{app:pf:prop:hsol}.

\medskip

\Cref{part:beltrami} states the measurable Riemann mapping theorem (solution of the Beltrami equation associated to a Beltrami form $\mu$) and details our proof. A crucial point is the holomorphic dependence, which relies on the holomorphic dependence of the Christoffel symbol proved in \Cref{part:simsurf}.
We first describe an (already known) density method with some amount of generality. 
Then we apply the results of the previous section to finalize the proof.

The density method does not involve the Ahlfors-Beurling operator nor $L^p$ spaces beyond $L^1$, $L^2$.
It involves the notion of distribution and the Sobolev space $W^{1,2}_\loc$, which serve in particular for the definition of quasiconformal maps.
It also involves compactness statements for quasiconformal maps: one about normal families, which is proved in \cite{A} by conformal geometry techniques; one about $L^2$ bounds, which is proved in \cite{A} by proving differentiability almost everywhere and bounding the differential with the Jacobian.

The presentation of our proof is not completely self contained: part of the book of Ahlfors, \cite{A,Abis}, is considered as a basic reference.
In particular several lemmas of this book are used here without proof, the interested reader will have to check such proofs there.\footnote{None of these proofs uses the Ahlfors-Beurling integral operator.}
Details about others aspects, proofs, or prerequisite knowledge are given in appendices here.

\medskip

While \Cref{part:beltrami} presents an original approach to known results,  \Cref{part:cnx} is prospective and contains new results.
In \Cref{part:beltrami}, the solution $f$ of the Beltrami equation in the measurable Riemann mapping theorem is obtained as a limit of approximations $f_n$.
Each approximation is associated to a particular abstract similarity surface $\cal S_n$, with $(2n^2+1)^2+1$ singularities.
Uniformizing the underlying Riemann surfaces to the Riemann sphere gives rise to a meromorphic conformal connection on $\C$, whose associated Christoffel symbol $\zeta_n$ is a rational map with $(2n^2+1)^2$ poles.
A natural question arises: is there a limit to these similarity surfaces?
In the particular case where the Beltrami form $\mu$ is $C^2$ with compact support, we state and prove limit theorems concerning the sequence of connections $\zeta_n$.
In particular, it tends in some sense to a conformal connection, but whose Christoffel symbol $\zeta$, which we characterize, is not anymore a holomorphic function.
In some sense, we gave an answer to the above question: the limit of the similarity surfaces $\cal S_n$ is the plane endowed with the connection $\zeta$.

\loctitle{Acknowledgements}

The authors would like to thank the anonymous referee for valuable remarks and the editors of the S.M.F.\ for their support.
Research by both authors is supported by the ANR project TIGerS (ANR-24-CE40-3604).
Research by G.T. is also supported by the Beijing Natural Science Foundation (Grant IS23005).
A.C.\ thanks several colleagues for useful discussions, in particular Xavier Buff, Christian Henriksen, John Hamal Hubbard and Michel Zinsmeister.
A.C.\ got introduced to the Ahlfors-Bers theorem by Adrien Douady, and got interested in the link with affine surfaces through the study \cite{C} that was sparkled after seeing a drawing of by Douady and Touzot that was pinned in Douady's office.

%
% ----------------------------------------------------------
%

\part{Similarity surfaces}\label{part:simsurf}

Similarity surfaces are also known under other names. In \cite{T} they are called \emph{affine surfaces} (but for other authors the term affine manifold may refer to other notions), and also \emph{$G$-manifolds} (of dimension $2$) for $G=$ the group of similarities of the Euclidean plane (this is a different notion from the $G$-structures of Elie Cartan\footnote{See \cite{Ste} Definition~2.1 Chapter~VII page~310.})
The article \cite{Ve}, considered as foundational, calls them \emph{affine complex surfaces}. In \cite{M} they are called \emph{branched affine surfaces} (with an emphasis on some type of singularities). They also have an interpretation in terms of flat symmetric conformal connections (holomorphic/meromorphic connections) on Riemann surfaces, see \Cref{sec:affco}.

%
% ----------------------------------------------------------
%

\section{Definitions and basic properties}\label{sec:simsurfdef}

Here we define similarity surfaces and an associated notion of monodromy.

\subsection{Definition}

A similarity surface is a two dimensional topological manifold together with an atlas whose transition maps are locally $\C$-affine maps $z\mapsto az+b$, i.e.\ on each connected component of the domain of a transition map there is a map $z\mapsto az+b$ that coincides with it.

Since $\C$-affine maps are holomorphic, it follows that similarity surfaces are special cases of one dimensional analytic manifolds, a.k.a.\ Riemann surfaces.\footnote{Though some authors require Riemann surfaces to be connected, in which case either one can require similarity surfaces to be connected too, or only the connected similarity surfaces will be Riemann surfaces.\label{foo:bar}}

Translation surfaces (transition maps are translations) are example of similarity surfaces, as are half translation surfaces (transition maps are of the form $z\mapsto a\pm z$). See \cite{Z}.

The best illustration of a similarity surface is M.C.\ Escher's lithograph \emph{Print gallery} where $\C^*$ is endowed with the atlas that consists of the branches of $z\mapsto z^\alpha$ with $\alpha = (2\pi i + \log 256)/(2\pi i)$, see \cite{dSLDS}.

\subsection{Morphisms}

A map between open subsets of similarity surfaces is called \emph{affine} if its expression in charts is locally $\C$-affine.
The set $\C$ is a similarity surface endowed with the \emph{canonical atlas}, which consists in a single map: the identity map of $\C$.
A map from a similarity surface to $\C$ is called affine if it is affine for the canonical atlas.

\subsection{Monodromy}

Let a \emph{germ} of chart at a point $M\in \cal S$ be an equivalence class in $E/\sim$, where $E$ is the set of non-constant affine maps $s: V\to\C$ defined in neighbourhoods $V$ of $M$, and $s_1\sim s_2$ whenever there is a neighbourhood of $M$ on which they are equal.

Given a path $\gamma:[0,1]\to\cal S$ and a germ of chart $\tilde s_0$ at $\gamma(0)$, there is a unique way to follow the germ along $\gamma$ so that for all $t$, $\tilde s_t$ is a germ of chart at $\gamma(t)$ and so that the germs locally match,\footnote{This can be expressed by saying that $t\mapsto \tilde s_t$ is continuous, for an appropriate topology on the space of germs.} i.e.\ for all $t$, $\exists\eta>0$ and a representative $s_t$ of $\tilde s_t$ such that for $|t-t'|<\eta$, $\gamma(t')$ is in the domain of $s_t$ and $\tilde s_{t'}$ is the germ of $s_t$ at $\gamma(t')$.
Existence can be proved by using an open cover of the image of $\gamma$ by charts and patching restrictions thereof appropriately along finitely many pieces of $\gamma$.
Uniqueness can be proved by an open-closed argument on $[0,1]$.

If the path is a loop ($\gamma(0)=\gamma(1)=M$) then $s_1\circ s_0^{-1}$ coincides with some $\C$-affine self map $\phi:\C\to\C$ in a neighbourhood of $s_0(M)\in\C$. The map $\phi$ and the point $s_0(M)$ are independent of the choice of the representatives $s_0$ and $s_1$ of $\tilde s_0$ and $\tilde s_1$ of the germ. The \emph{monodromy} of $\gamma$ is defined as the conjugacy class of $\phi$ in the group of non-constant affine self-maps of $\C$.
This class depends only on the free\footnote{Without basepoint.} homotopy class of $\gamma$.
By abuse of language, we will say that the monodromy is $z\mapsto z$ or $z\mapsto z+1$ in the translation case and $z\mapsto \lambda z$ in the other cases.
The \emph{monodromy factor} is defined as $1$ in the translation case and $\lambda$ in the other cases.

\begin{remark} 
A similar point of view, which we will not develop here, is to endow the universal cover $U\cal S$ of $\cal S$ based on $M$ with a similarity atlas compatible with that of $\cal S$, and to extend the initial germ $s_0$ to a global affine map $U\cal S\to\C$.
This allows to define a finer monodromy invariant: a conjugacy class in the set of representatiton of the fundamental group in the affine group.
\end{remark}

%
% ----------------------------------------------------------
%

\section{Polygons}\label{sec:simsurfpolygon}

In this section we explain how gluing bounded or unbounded polygons naturally defines similarity surfaces on the complement of the vertices, explore the conformal nature (puncture or hole) of their singularities at the vertices, and in the puncture case compute the associated monodromy.

\subsection{Gluing bounded convex polygons}

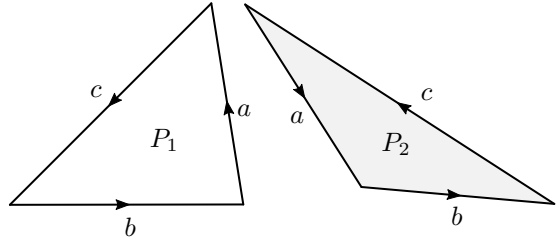
\begin{figure}
\begin{tikzpicture}
  \node at (0,-0.05) {\includegraphics{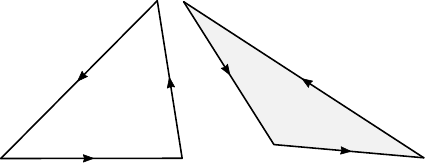}};
  \node at (-1.55,-0.5) {$P_1$};
  \node at (-0.5,-0.1) {$a$};
  \node at (-2,-1.65) {$b$};
  \node at (-2.45,0.15) {$c$};
  \node at (1.5,-0.55) {$P_2$};
  \node at (0.2,-0.2) {$a$};
  \node at (2.3,-1.5) {$b$};
  \node at (1.9,0.1) {$c$};
\end{tikzpicture}
\caption{An example of polygon gluing pattern. The letters indicate the edge pairing and the gluing reverses the edge orientation. Here the quotient is homeomorphic to a sphere and has three vertices and three edges.}\label{fig:a}
\end{figure}

Assume we are given a finite collection $(P_j)_{j\in J}$ of bounded convex polygons. Each polygon must have finitely many sides (edges) and their vertices are the endpoints of these sides. As subsets of $\C$, we will choose these polygons including their sides, i.e.\ they are closed for the topology of $\C$. 
We allow flat angles ($\pi$ radians) at a vertex, so keep in mind that the data of the polygon as a subset of $\C$ does not necessarily characterize the set of its vertices.
Each edge is naturally oriented by following the polygon boundary anticlockwise.
See \Cref{fig:a,fig:c}.

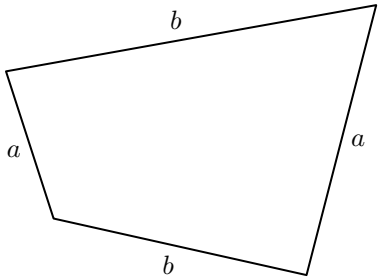
\begin{figure}
\begin{tikzpicture}
  \node at (0,0) {\includegraphics{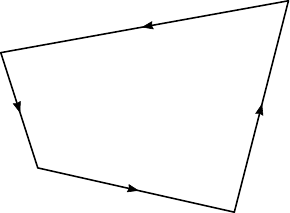}};
  \node at (-2.33,-0.15) {$a$};
  \node at (2.2,0) {$a$};
  \node at (-0.2,1.6) {$b$};
  \node at (-0.3,-1.65) {$b$};
\end{tikzpicture}
\caption{Another example. This is the classical torus fundamental domain, but with a random quadrilateral instead of a square. Here there is only one polygon and only one vertex after projection to the quotient.}\label{fig:c}
\end{figure}

Consider the collection of all the sides $e$ with their polygons $P_j$ and denote them $e\of j$. Assume we are given a pairing\footnote{A pairing is a partition into sets of cardinality $2$.} of these sides $e\of j$.
From this, one can define a compact oriented topological manifold $\cal S$ of dimension $2$, as follows: for each pair, choose an ordering of the pair as $(e_1\of j_1,e_2\of j_2)$; note that the case $j_1=j_2$ is allowed by the pairing, but not the case $e_1\of j_1=e_2\of j_2$; there is a unique complex affine map $s$ sending $e_1$ to $e_2$ and reversing their orientations.
Let $\cal P$ be the disjoint union of all polygons, whose elements will be denoted $z\of j$ to distinguish two elements of differently indexed polygons with the same affix, which may happen.
Let $\cal S$ be $\cal P$ quotiented by $z\of j_1\sim s(z)\of j_2$ for all the paired edges as above. Call $\prA: \cal P \to \cal S$ the quotient map.\footnote{We use the same notation for Archimedes' constant, context should prevent confusion.}

Interior points of polygons are not identified with any other points. An interior point $z:j$ of an edge is identified with exactly one other point $s(z):j'$, which sits in the interior of the matching edge.
Depending on the situation, a vertex can be matched with any finite number of other vertices.

\subsection{Local model at vertices}

Let us look at what the quotient looks like at a vertex $\prA(v_0:j_0)$.
Consider all the $v:j$ such that $\prA(v:j)=\prA(v_0:j_0)$: these are the vertices of the $P_j$ that project to the same point. (Note that it is quite possible for a given $P_j$ to have several vertices with the same image by $\prA$.) 
Near these vertices the polygon $P_j$ looks like a sector with some opening angle and the point $v$ belongs to two edges, one termed ``before'' and the other one ``after'' so that the sequence before$\ \rightarrow\ $sector$\ \rightarrow\ $after follows the anticlockwise order with respect to the vertex.\footnote{Be attentive to the fact that this is designed to study a vertex neighbourhood. As a side effect, the segment termed ``before'' comes after the segment termed ``after'' if we follow the anticlockwise orientation of the boundary $P_j$.}
Consider the data $(v,e,j)$, that we will call \emph{flag}, where $e$ is the edge that is termed ``before''.
Denote $e'$ the edge termed ``after'': it is glued via some affine map $s$ to some edge of some other polygon $P_{j'}$ and we call \emph{successor} of $(v,e,j)$ the flag $(s(v),s(e'),j')$.
Following the consecutive gluings, we get from flag to flag and eventually back to the initial flag since there are finitely many polygons.
By explicit transformations, we can map homeomorphically each sector to sectors of the same opening with apex of affix $0$ and so that the gluings, except maybe one, become the identity. This shows that the quotient is indeed a manifold at the vertices.

\subsection{Removing vertices and getting a similarity surface}

Let $\cal V \subset \cal S$ denote the set of all vertices (after passing to the quotient).
On $\cal S' = \cal S-\cal V$ one can define an atlas of a similarity surface as follows: inside each polygon $\prA(\Int P_j)$, take $\prA^{-1}$ as a chart.
Near an edge $e = \prA(e_1\of j_1)=\prA(e_2\of j_2)$ with endpoints removed (call this $\Int e$), we consider two cases.
If $P_1\neq P_2$ then take the neighbourhood $V=\prA(\Int P_1)\cup \Int e \cup \prA(\Int P_2)$ and the chart $\phi: V\to s(\Int P_1) \cup \Int e_2 \cup \Int P_2$ mapping $\prA(z:P_1)$ to $s(z)$ and $\prA(z:P_2)$ to $z$ (this coincides in $\Int e$).
If $P_1=P_2$ we do the same but using only non-overlapping neighbourhoods of $\Int e_1$ and $\Int e_2$ in $P_1$ instead of the full polygon.
In both cases, the sets $s(\Int P_1)$ and $\Int P_2$ cannot overlap because we took convex polygons and $s$ reverses the orientation of edges.
See \Cref{fig:b}.

\begin{figure}
\begin{tikzpicture}
  \node at (-2.8,-1.5) {\includegraphics{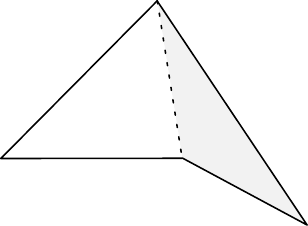}};
  \node at (1,0) {\includegraphics{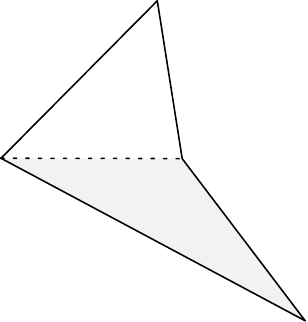}};
  \node at (-5.2,1.5) {\includegraphics{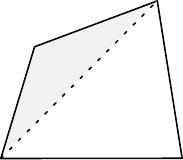}};
\end{tikzpicture}
\caption{Continuation of \Cref{fig:a}; by gluing the polygons along their paired edges with appropriate $\C$-affine maps, we get local chart of a similarity surface in neighbourhoods of the edges minus their endpoints.}\label{fig:b}
\end{figure}

\subsection{Conformal erasability of the singularities}

\begin{figure}
\begin{tikzpicture}
  \node at (0,0) {\includegraphics[scale=0.8]{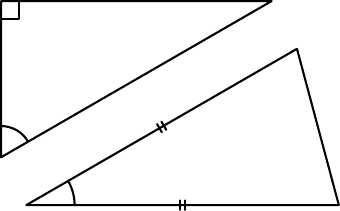}};
  \node at (0.3,-0.3) {$a$};
  \node at (-0.55,0.55) {$a$};
  \node at (0.5,-1.6) {$b$};
  \node at (-2.5,0.5) {$b$};
  \node at (-2.4,-0.95) {$p$};
  \node at (-2.15,-1.4) {$p$};
  \node at (-1.85,-0.05) {$\tau/6$};
  \node at (-0.8,-1.1) {$\tau/12$};
  \node at (0,-1.8) [anchor=north] {\begin{minipage}{5cm}Assume the vertex $p$ is shared by these two triangles and no other.\end{minipage}};
    
  \node at (6,0.1) {\includegraphics[scale=0.8]{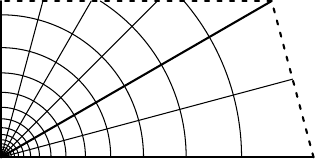}};
  \node at (6.0,-1.25) {$b$};
  \node at (5.9,0.0) {$a$};
  \node at (3.6,0.15) {$b$};
  \node at (3.3,-1.5) [anchor=north west] {\begin{minipage}{6cm}Here we glued the side pair labelled $a$ together and added a exp-polar grid. We disregard what the dashed sides are attached to.\end{minipage}};
  
  \node at (0,-4.3) {\includegraphics[scale=0.7]{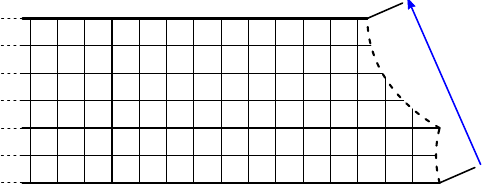}};
  \node at (-0.7,-3.15) {$b$};
  \node at (0,-5.6) {$b$};
  \node at (2.6,-4.0) {$u$};
  \node at (0,-5.8) [anchor=north] {\begin{minipage}{5.5cm}Image by $z\mapsto \log z$. The vector has affix $u=-\log 2 +i\tau/4$ and gluing of sides $b$ becomes translation by $u$.\end{minipage}};

  \node at (6.0,-6) {\includegraphics[scale=0.85]{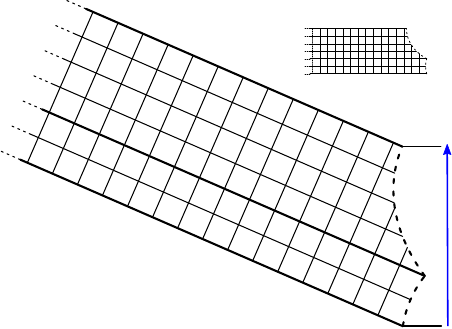}};
  \node at (8.95,-6.7) {$i\tau$}; 
   \node at (9.4,-8.3) [anchor=north east] {\begin{minipage}{5.7cm}Applying $z\mapsto \frac{i\tau}{u} z + \on{const}$ turns this set (the small copy) into a bigger one and the vector $u$ into $i\tau$.\end{minipage}};

  \node at (0.4,-10.1) {\includegraphics[scale=0.8]{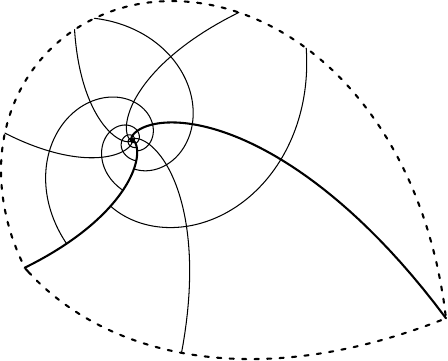}};
  \node at (2.1,-10.9) {$a$};
  \node at (-1.4,-11.1) {$b$};
  \node at (9.4,-10.4) [anchor=north east] {\begin{minipage}{5.2cm}Last we apply the exponential, which performs the gluing of the side pair $b$, and gives us a Riemann chart near the vertex $p$.\end{minipage}};
\end{tikzpicture}
\caption{Riemann chart near a vertex for a specific example. We have $\lambda = 1/2$ and $\theta = \tau/4$, the monodromy is $z\mapsto \frac{i}{2} z$. (For convenience we used $\tau=2\pi$ in this figure.)}\label{fig:gluex1}
\end{figure}

Seeing $\cal S'$ as a Riemann surface,\footnote{Or a finite union thereof, see \cref{foo:bar}.} the vertices could be either \emph{punctures} (the Riemann surface can be extended by adding the point and an appropriate chart) or \emph{holes} (no such extension can be done; this happens iff the vertex has a neighbourhood in $\cal S'$ conformally isomorphic to a ring of finite modulus: $1<|z|<1+\eps$, see for instance \cite{Ne}).
Let us prove that we are in the first case by completing the Riemann surface $\cal S'$ at the vertices, thus promoting $\cal S$ to a Riemann surface.
Let $p\in \cal S$ be a vertex.
Consider the circularly ordered sequence of flags $(v,e,j)$ associated to $p$ as described earlier in this section, and their associated angular sectors.
Let $\theta>0$ denote the sum of the angles of the sectors.
When we pass from a flag $(v,e,j)$ to its successor $(s(v),s(e'),j')$, there is an affine map $s(z)=az+b$ performing the transition. Let $\lambda$ denote the \emph{inverse} of the product of the dilation factors $|a|$. This product is \emph{not} necessarily $1$.
The monodromy\footnote{Defined at the end of \Cref{sec:simsurfdef}.} of the similarity surface for an anticlockwise loop winding once around the vertex is equal to
\[z\mapsto {\lambda}{e^{i\theta}} z.\]
\begin{remark}
It is easy to get confused and believe that the monodromy would be the inverse of this function.
If we scale an rotate the different sectors to attach them successively, putting the vertex at $0\in\C$, we get a big sector (living in the universal cover of $\C^*$ if $\theta>2\pi$) of which there remains only one pair of sides to glue together by some similitude.
Call them ``before'' and ``after'' so that what we successively meet in the anticlockwise order is ``before'', sector, ``after''.
Then the monodromy is the similitude that sends ``before'' to ``after'', not the other way round.
\end{remark}

\Cref{fig:gluex1} illustrates an example and also illustrates the construction of a Riemann chart that we now explain.
Map each sector to a horizontal band, infinite on the left, by a branch of $z\mapsto w=\log(z-v)$. The gluing map from a band to its successor is then a translation.
Let us translate them around in $\C$ and stack them according to these gluings.
We obtain a band $B$ of height $\theta$, on which the only remaining gluing goes from the top line to the bottom line and is the translation by the complex vector $-u$ where $u=\log(\lambda) +i\theta$ has positive imaginary part.
Now, for a well-chosen complex number $\alpha$, the similarity $w\mapsto \alpha w$ send the vector $u$ to $2\pi i$ and the band to a band that is not necessarily horizontal any more, yet is still infinite on the left, and whose two sides are glued by the translation by $-2\pi i$ from top to bottom. A uniformizing map for this quotient is just the exponential map $\exp$ from $\alpha B$ to a neighbourhood $U$ of $0$ in $\C^*$.
(If we would have sent $u$ to $-2\pi i$ instead of $2\pi i$, the band would be infinite towards the right, and its exponential would be a neighbourhood of infinity.)
Adding $\{0\}$ to $U$ we get an analytic chart near the vertex $p$. We have turned $\cal S$ into a compact Riemann surface.\footnote{Or a finite union thereof, see \cref{foo:bar}.}

\subsection{Unbounded polygons}\label{page:unbdd}

\begin{figure}
\begin{tikzpicture}
\node at(-0.15,-0.3) {\includegraphics{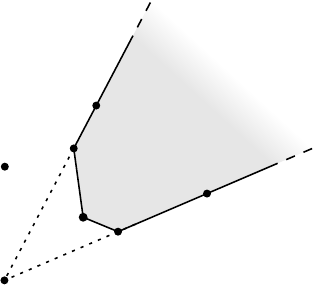}};
\node at(-2.45,-0.72) {$c$};
\node at(-3,-2.6) {$f$};
\node at(0.05,-1.65) {$a$};
\node at(-1.52,0.02) {$a$};
\end{tikzpicture}
\caption{Example of an unbounded polygon with the two unbounded sides glued together by a similitude whose centre $c$ is different from the focus $f$ of the unbounded sector.}\label{fig:unbddpoly}
\end{figure}
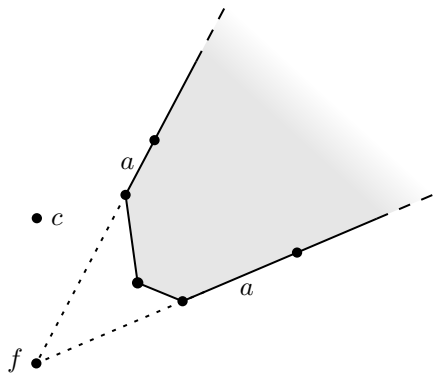

We will need to let unbounded polygons $P$ enter the game. As above, we start by including only convex ones.
They are assumed closed as subsets of $\hat{\C}$, i.e.\ they include the vertex $\infty$.
We call \emph{finite} vertex a vertex that is in $\C$, i.e.\ different from $\infty$.
We only consider unbounded polygons that have exactly two edges reaching infinity, we call them \emph{unbounded} edges, and require that these edges are half-lines, i.e.\ that each reach to a finite vertex, possibly the same. Our requirement rules out in particular the case where an edge is a whole straight line, but still allows the two unbounded edges to have a union that is a whole straight line (iff.\ they end at the same finite vertex and make there an angle of $\pi$).

This assumption implies obviously that there are at least two edges, and less obviously that the complement has exactly one connected component.
So it rules out the case $P=\C$, which has no edge, and also rules out a half plane with no finite vertex, and also the case where $P$ is a strip.
However, we can recover all these cases by pasting appropriate polygons satisfying our assumptions (see also the last paragraph of the present section).

With these conditions each unbounded polygon has only one infinite vertex and each of its two unbounded edges is a half line.
Recall that we assume it is convex.
We now consider the previous constructions, but allowing the presence of unbounded polygons that satisfy the conditions.
We have to assume that the pairing between edges satisfies that unbounded edges are paired with unbounded edges, and so that the infinite vertex is on the first vertex of one edge and the second vertex of the other, where the polygon boundary is still oriented counterclockwise. 
Moreover, the orientation reversing $\C$-affine map sending an unbounded $e_1:j_1$ to an unbounded $e_2:j_2$ exists but is not anymore unique.
So the data that must be given prior the construction has to include the choice of this map for each such pair (for a visual way to do this, one can mark unbounded edges with supplementary non-vertex points, that must match under the gluings, as in \Cref{sub:hd:p}).

Then the construction can be carried out and we still get a topological surface $\cal S$, but the set of vertices $\cal V \subset\cal S$ now includes a new kind of element: \emph{infinite} vertices.
These are the vertices which come from the vertex at infinity of an unbounded polygon.
As for vertices of bounded polygons we define $\theta\geq 0$ as the sum of opening angles at infinity of the unbounded polygons around the vertex.

Again something new may happen: conformally, \emph{infinite vertices are not anymore necessarily punctures}. An example is given in \Cref{app:hole}.
To the infinite vertex of an unbounded convex polygon satisfying our restrictions, one can associate the angle between the two edges reaching to infinity.
It is a non-negative number.
In the example of \Cref{app:hole} this angle is $0$.

To a infinite vertex $v$ of the quotient let us again associate the sum of the angles at infinity of the unbounded polygons whose infinite vertex projects to $v$. Then a sufficient condition for an infinite vertex to be a puncture is that this sum is positive.
The uniformizing map is obtained in a similar way as for finite vertices, with a few complications. 

An infinite (closed) sector is of the form $a+[R,+\infty)\cdot e^{i[\theta_1,\theta_2]}$ for some $R\geq 0$ and some $a\in\C$ that we will call its \emph{focus}. A half infinite strip is the image of the set $[0,1]\times [0,+\infty)\subset \R^2 \simeq\C$ by a complex affine map.

\begin{remark}
In the case of a finite vertex, we changed variable using the complex logarithm; it was natural to put the vertex at $0$ and we saw that this allows to give a nice local Riemann chart.
If we were to try to do the same for an infinite vertex, we would need to decide which translation to perform before the logarithm, i.e.\ where to put the origin in the picture.
Here is a list of some complications that may occur.
The angle of an unbounded polygon could be $0$. The total angle at an infinite vertex could be $0$.
The focus of an infinite sector extracted from a polygon does not necessarily coincide with one of the finite vertices of its unbounded edges.
The focuses of the infinite sectors do not necessarily match under the gluing maps between unbounded edges.
Taking a logarithm based on a focus would not conjugate the gluing (a similarity) to a map as simple as a translation (and if the two finite vertices of the two unbounded edges do not coincide, we of course cannot base the log on both).
\end{remark}

Instead of working as in the remark above, we take subsets of sectors or strips that are neighbourhoods of infinity in the cycle of unbounded polygons around an infinite vertex, and glue them in a universal cover $U\C^*$ of $\C^*$ into one big set sitting in finitely many sheets of the cover (one sheet is enough near $\infty$ if $\theta<2\pi$).
A neighbourhood of infinity in this set is a sector or a strip. Its two sides must be glued together according to some $\C$-affine map (the vertices are now excluded from of the picture, except the one at infinity) whose factor is the monodromy factor or its inverse depending on conventions.
To fix the conventions,  we order the two sides of any infinite sector/strip in the clockwise order: it seems to be the opposite of the ordering we chose for finite vertices, but note that under an inversion $z\mapsto 1/z$ it is actually the same anticlockwise convention.
The monodromy sends the second edge of the aforementioned big set to the first one, as in the finite case.
The monodromy factor must be of the form
\[\lambda e^{-i\theta}\]
for some $\lambda>0$. Now there are several cases

\begin{enumerate}
\item $\theta>0$.
\begin{enumerate}
\item Either the affine map $s$ has a unique fixed point.
Then we can take a logarithm based on this fixed point.
The image of the two sides will be almost horizontal curves\footnote{Because we work close to infinity. If bothered by that, one can make these curves actual horizontals by shifting the original edges, which can, in a neighbourhood of infinity, be achieved by cutting and pasting pieces of the sectors.} that extend infinitely to the right and that we must glue by a translation of vector $v = \log\lambda - i\theta$.
Multiplication by $2\pi i / v$ followed by $\exp$ realizes this gluing.

\item Or the affine map $s$ is the identity.
Then the total angle is a multiple of $2\pi$, and the two curves to be glued together have the same projection by $\phi : U \C^* \to \C^*$.
The quotient is a $k$-fold cover of a punctured neighbourhood of infinity.
A uniformization is given by a determination of $z\mapsto z^{-1/k}$.

\item Or the affine map $s$ is a translation by a non-zero vector $a$.
In this case $\lambda=1$, and $\theta$ must be a multiple of $2\pi$:
\[\theta = k2\pi.\]
By composing everything with $z\mapsto a^{-1}z$ we can assume that $a = 1$.
By cutting and pasting we can assume that the edges to glue have their infinite direction of argument $0$.
We provide explicitly a formula defining a function from a punctured neighbourhood of $0$ to $\C$ \emph{whose inverse} has a branch realizing the gluing:\footnote{In this case the map $\phi$ is also the straightening coordinate of the vector field $\dot z = \frac{z^{k+1}}{-k+2\pi i z^k}$, which provides another interpretation.}
\[\phi(z) = z^{-k} + \frac{1}{2\pi i}\log z.\]
We omit the details.
\end{enumerate}

\item $\theta=0$.
\begin{enumerate}
\item Either the gluing map is a translation, which can be taken as $z\mapsto z+2\pi i$ by an appropriate $\C$-affine change of coordinate. Then a uniformization is provided by $z\mapsto \exp(z)$.
\item Or the gluing map is of the form $z\mapsto az+b$ for some $a\in\R$ with $a>0$ and $a\neq 0$. Then the vertex at infinity is not a puncture but a hole. 
\end{enumerate}

\end{enumerate}

\subsection{Signed angle and monodromy factor}

To unify the formulas for the monodromy we define the \emph{signed angle} $\sigma\in\R$ of a vertex as $\sigma = \theta$ if the vertex is finite and $\sigma = -\theta$ if the vertex is infinite, where $\theta\geq 0$ is the total angle defined in the previous paragraphs. Then the factor of the monodromy of a small loop winding anticlockwise around the singularity is in both cases of the form
\[\lambda e^{i\sigma}\]
for some $\lambda>0$. Note that the monodromy factor determines $\lambda$ completely but only determines $\sigma$ modulo $2\pi i$.

\subsection{Even fancier polygons}

In fact we can use polygons that are non-convex, or not simply connected, and without restrictions on the type of unbounded edges or their number, still finite.
Indeed such polygons can be cut further, possibly introducing auxiliary finite vertices (in which case, obviously, the resulting similarity surface structure will have an erasable singularity at these vertices), so as to respect the previous restrictions and yield the same objects $\cal S$ and $\cal S'$ once we put back the auxiliary vertices.
Then one thing will have to be kept in mind: that if a polygon $P\subset\C$ has $n$ unbounded complementary components in $\C$, then it must be considered as having $n$ distinct infinite vertices (before the gluing; after, this number may get reduced).

One could also use polygons that overlap themselves, i.e.\ are spread over $\C$ in an non-injective way, but we will not develop this aspect.

\subsection{Geodesics}

In a similarity surface, there is a well-defined notion of parametrized geodesics: these are maps from an interval in $\R$ to the surface, such that in each chart, the map is locally $\R$-affine and non-constant.
Equivalently, the curve is differentiable and its derivative is locally constant and non-zero in charts.
Note that changes of similarity charts preserve this property.
For each initial time $t_0$, starting point $p$ in the surface and initial tangent vector $v\neq 0$ at this point, there is a unique maximal geodesic $\phi$ with $\phi(t_0)=p$ and $\phi'(t_0)=v$.

%
% ----------------------------------------------------------
%

\section{The Schwarz-Christoffel formula}\label{sec:sc}

In this section we explain how to associate a function, denoted $\zeta$ in this article, to conformal charts of similarity surfaces.
We call it the Christoffel symbol, as it turns out to be the coefficient of a meromorphic connection, but we will not need the notion  of connection in the present part nor in \Cref{part:beltrami}, but only in \Cref{part:cnx}.
We give basic properties of $\zeta$, and for puncture type singularity of a similarity surface $\cal S'$ obtained by gluing polygons and removing vertices as in the previous section, we identify the function $\zeta$ in an appropriate conformal chart of the puncture.
Moreover, in the case the affine surface is conformally isomorphic to the Riemann sphere minus a finite number of points, we give an explicit formula for $\zeta$, which is a generalization, already known, of the Schwarz-Christoffel formula.

\subsection{Christoffel symbol}

Consider any similarity surface $\cal S'$. Recall that it is in particular a Riemann surface. Assume that we are given a similarity chart $s : U\to\C$ and a Riemann chart $r: U\to\C$ on the same\footnote{We can always restrict a chart to an open subset of its domain and add it to the atlas, this gives the same structure on $\cal S'$.} open subset $U\subset \cal S'$.
Consider the map $\phi = s\circ r^{-1}$: it expresses a similarity chart in a Riemann chart.
Note that $\phi$ is holomorphic, by definition of a Riemann surface atlas.
Now if we are given two similarity charts $s_1$ and $s_2$ on $U$ and if $U$ is connected then by definition of a similarity atlas, there must exist $a\in\C^*$ and $b\in\C$ such that $s_2=a\,s_1+b$.
(In fact the definition ensures the existence of $a$ and $b$ locally, and analytic continuation ensures they are constant on the connected set $U$.)
It follows that
\[\phi_2 = a\,\phi_1 + b\]
hence
\[\phi_2' = a\,\phi_1'\]
thus
\[\frac{\phi''_1}{\phi'_1} = \frac{\phi''_2}{\phi'_2}.\]
For any holomorphic function on an open subset of $\C$ we now use the notation
\[N\phi =N(\phi) := \frac{\phi''}{\phi'}.\]

Assume we are given a subset $U$ of $\cal S'$.
Note that for every point $M\in U$ there are similarity charts defined in a neighbourhood of $M$ but there may fail to exist a similarity chart that would be defined on the entirety of $U$, even when there is a Riemann chart $r$ defined on $U$.

If we are given a Riemann chart $r:U\subset \cal S' \to \C$, we can use it to endow $r(U)$ with the similarity atlas $\{\phi=s\circ r^{-1}\}$ where $s$ varies in the similarity atlas of $\cal S'$.
Note that if we denote $V$ the domain of $s$ then the domain of $\phi$ is $r(U\cap V)$.

We will thus momentarily study open subsets $O$ of $\C$ with a similarity surface atlas whose charts $\phi$ are holomorphic.
By the above computation, all similarity charts have the same quantity $N(\phi)$ at a given $z\in O$. So even if there is not necessarily a single chart $\phi$ on $O$, there is always a well-defined holomorphic function $\zeta : O\to\C$ such that for all similarity charts $s: V\subset O\to\C$, then
\[N\phi=\zeta\]
on $V$.
It is easy to recover $\phi$ from $\zeta$, as the equation $\phi''/\phi' = \zeta$ is equivalent to $\log \phi' = \int \zeta$ (locally), so 
\[\phi = \int \exp\left(\int \zeta\right)\]
holds locally.
The integration constants make $\phi$ known locally only up to post-composition by $\C$-affine maps, which is coherent with the fact that we have a similarity atlas:
$\phi = b+\int \exp \left( c+\int \zeta\right) = b+a \int \exp \left( \int \zeta\right)$ with $a=e^c$.
Note also that, for a loop $\gamma$ contained in a single chart $O$ the monodromy factor of a loop will be equal to the exponential of the integral of $\zeta$ along this loop:
\[ \lambda = \exp\left(\int_\gamma \zeta(z)dz\right)
.\]

Let us come back to abstract similarity surfaces $\cal S'$. Up to now we only expressed similarity charts in \emph{one} Riemann chart, but what happens if we change the Riemann chart?
Assume $r_1$ and $r_2$ are both defined on $U\subset \cal S'$ and let $\psi = r_1\circ r_2^{-1}$.
Then the expressions $\phi^{(i)} = s\circ r_i^{-1}$ of a similarity chart $s$ via the maps $r_i$ are related by
$\phi^{(2)} = \phi^{(1)}\circ \psi$.
It is easy, for holomorphic functions, to compute $N(\phi\circ\psi)$:
\begin{equation}\label{eq:Nchvar}
N(\phi\circ\psi) = \psi' \times N(\phi)\circ\psi + N(\psi).
\end{equation}
For our situation let us write this as
\begin{equation}\label{eq:zetachvar}
 \zeta^{(2)} = \psi' \times \zeta^{(1)} \circ \psi + \frac{\psi''}{\psi'}.
\end{equation}
where $\zeta^{(i)}$ is the function of the previous paragraph, for the similarity atlas $\phi^{(i)}$ on the subset $r_i(U)$ of $\C$.

\begin{remark}
\Cref{eq:Nchvar} is reminiscent of the formula of a change of variable for a $1$-form.
In fact $N\phi$ does not express a $1$-form on the Riemann surface $\cal S'$ but a \emph{Christoffel symbol} in complex dimension one.
Christoffel symbols are the coefficients that appear in the infinitesimal expression of affine connections in charts, see \cite{BG} or \Cref{sub:affcoconf}.
More precisely, on the similarity surface there is a
standard flat connection: a parallel transport of a vector along a path is given by the trivial parallel transport in similarity charts.
Consider now an arbitrary holomorphic chart with coordinate $z$. 
Let $\phi(z)$ be a similarity chart, and let
$\zeta(z) = \phi''(z)/\phi'(z)$.
Then the connection in question takes the expression $\nabla Y = \partial_z Y + \zeta (z)Y$ in coordinate $z$. 
The parallel transport equation is $\partial_z Y + \zeta (z)Y=0$.
The family of vectors $\frac{1}{\phi'(z)}\frac{\partial}{\partial z}$  is parallel along every path for this connection.
\Cref{eq:Nchvar} is also reminiscent of a similar formula for the \emph{Schwarzian derivative} $S\phi$. In fact if we look at surfaces with atlases whose change of charts are homographies instead of affine maps, the operator $S$ will play the role of the operator $N$.
\end{remark}

If we have an isolated point in $\C- U$, then by \cref{eq:zetachvar} the polar part of $\zeta$ behaves under a change of variable exactly like the polar part of a meromorphic $1$-form. In particular the residue ``$\res $'' is an invariant. This is coherent with the fact that the monodromy of the similarity surface around this singularity is independent of any choice of Riemann surface charts and is a conjugacy class of affine map whose linear factor is equal to $\exp(2\pi i\res )$.

\subsection{Polygon gluing and residue of the Christoffel symbol near singularities with non-zero angle sum}\label{sub:presc}

We place ourselves here in the situation of \Cref{sec:simsurfpolygon} where a topological surface $\cal S$ and a similarity surface $\cal S'=\cal S-\cal V$ were constructed by gluing bounded or unbounded polygons, where $\cal V$ is the set of vertices, and we make the assumption that all infinite vertices have an angle sum $\theta$ that is non-zero.
We saw that in this case the atlas associated to $\cal S'$ extends to a Riemann surface atlas of $\cal S$.
Let us study more closely the Christoffel symbol near the vertices.

For any vertex which is a finite vertex we constructed a Riemann chart whose image is a neighbourhood of $0$ and in which the similarity charts are branches of $z\mapsto z^\alpha$ for the complex number $\alpha = (\log \lambda + i \theta)/2\pi i$ where $\theta$ is the sum of the angles of the polygons at this vertex and $\lambda>0$ is such that the monodromy factor at the vertex is $\lambda e^{i\theta}$.
For any branch of $z\mapsto z^\alpha$, a direct computation gives
\[N z^\alpha = \frac{\alpha-1}{z}.\]
By \cref{eq:zetachvar} and the consequence on the polar part explained after it, it follows that in any other Riemann chart the Christoffel symbol $\zeta$ will have a simple pole at the vertex, of residue 
\begin{equation}\label{eq:resfin}
\res  = \alpha-1 = \frac{\log (\lambda) + i \theta}{2\pi i} -1.
\end{equation}
In particular
\[\Re\res >-1.\]

A similar construction was done for infinite vertices for which $\theta>0$.
The conclusion was similar: the monodromy has factor $\lambda e^{-i\theta}$ for some $\lambda>0$ and there is a Riemann chart mapping a neighbourhood of this vertex to a neighbourhood of $0$ and for which similarity charts are given by branches of $z\mapsto z^\alpha$ with $\alpha = (\log\lambda -i\theta)/2\pi i$ or, if the monodromy is a non-identity translation (hence $\lambda=1$ and $\theta = k2\pi$), by branches of $z\mapsto z^{-k}+\frac{1}{2\pi i}\log z$.
The Christoffel symbol has in both cases a \emph{simple pole} at the vertex, of residue
\begin{equation}\label{eq:resinf}
\res  = \alpha-1 = \frac{\log (\lambda) - i \theta}{2\pi i} -1
\end{equation}
in this chart and thus in any other Riemann chart.
In particular
\[\Re\res <-1.\]

\begin{remark}
The notion of signed angle $\sigma$ in \Cref{sec:simsurfpolygon} allows to unify \Cref{eq:resfin,eq:resinf}: in both cases
\[\res  = \frac{\log (\lambda) + i \sigma}{2\pi i} -1.
\]
The monodromy factor $\lambda e^{i\sigma}$ has a unified expression in terms of $\res$:
\[\lambda e^{i\sigma} = \exp(2\pi i \res).\]
\end{remark}

\subsection{Polygons gluing to a sphere and the Schwarz-Christoffel formula}

We now make the supplementary assumption that the topological surface $\cal S$ is homeomorphic to a sphere.
(We still place ourselves in the situation of \Cref{sub:presc}, which we recall: $\cal S$ is obtained by gluing together finitely many bounded or unbounded polygons; all infinite vertices are moreover required to have an angle sum $\theta$ that is non-zero, so that in particular all singularities are punctures and the Riemann surface atlas of $\cal S'$ extends to $\cal S$.)
Then by the \emph{Poincaré-Koebe} theorem, $\cal S$ is isomorphic as a Riemann surface to the Riemann sphere $\hat \C$. Let
\[\prB : \cal S\to\hat\C\]
be such an isomorphism.
Let $\{z_1,\ldots,z_{m}\}\subset\hat\C$ be the image by $\prB$ of the vertices.
The point $\infty\in\hat\C$ may or may not be one of them.
For each vertex $v_k$, let $\theta_k$ be the total angle, $\sigma_k=\pm\theta_k$ the signed angle, and $\lambda_k$ the dilation ratio of this vertex.
Note that $\prB$ gives us a \emph{global chart} of the Riemann surface $\cal S-\{\prB^{-1}(\infty)\}$. Removing vertices, $\prB$ is a global Riemann chart for $\cal S'-\{\prB^{-1}(\infty)\}$.
The associated Christoffel symbol $\zeta$ is a holomorphic function from $\C - \{z_1,\ldots,z_{m}\}$ to $\C$.

Let $\on{res_k}$ denote the residue of $\zeta$ at the vertex $z_k$.

\begin{theorem}\label{thm:sc} We have
\[\zeta(z) = \sum_{\stackrel{k=1}{z_k\neq\infty}}^m\frac{\res_k}{z-z_k}.\]
Moreover
\[\sum_{k=1}^m \res _k = -2.\]
\end{theorem}
\begin{proof}
By the above discussion, $\zeta$ has at most a simple pole at each vertex $z_k\neq \infty$ and we know the corresponding residue $\res _k$.
To analyse what happens at infinity we use the change of variable $w=\psi(z)=1/z$, call $\tilde\zeta$ the expression of the Christoffel symbol in the coordinate $w$ and use \Cref{eq:zetachvar}:
\[\zeta(z) = -z^{-2} \tilde\zeta(z^{-1}) -2z^{-1}.\]
If $\infty$ is not a vertex then $\zeta(z) = -2z^{-1} + \cal O(z^{-2})$ as $z\to\infty$. If $\infty$ is a vertex $z_k$ then $\tilde\zeta(w) = \frac{\res _k}{w}+ \cal O(1)$ as $w\to 0$, from which:
\[\zeta(z) = (-\res _k-2)z^{-1} + \cal O(z^{-2})\]
as $z\to\infty$. In all cases: $\zeta$ tends to $0$ to at $\infty$.

Let $f(z)$ denote the sum of $\res _k/(z-z_k)$ as in the statement of the theorem.
The difference $\zeta-f$ is a holomorphic function on $\C$ minus finitely many points, whose singularities are erasable and which tends to $0$ at infinity.
It follows that $\zeta-f=0$.
This proves the first claim.
Moreover, from the above it follows that $\zeta(z)-f(z) = \left(-2-\sum  \res _k\right) z^{-1}+\cal O(z^{-2})$ as $z\to\infty$, whence the second claim.
\end{proof}

\begin{remark}
That the sum of residues is $-2$ (in the case of $\hat \C$) and not $0$ is another difference between the Christoffel symbol and $1$-forms.
It can be proved more generally for any Christoffel symbol on $\hat\C$ with finitely many singularities (polar or essential) by using \Cref{eq:zetachvar} as above with the change of variable $z\mapsto 1/z$ and expressing $\int \zeta(z)dz$ on a big circle in two ways, using the residues on each side.
Another (equivalent) way is to use that the difference of two Christoffel symbols is a $1$-form and to use that the sum of residues of a $1$-form is always $0$ to reduce the computation to a particular symbol (for instance $\zeta = 0$).
Last, let us remark that this formula generalizes to genus $g$ compact Riemann surfaces: the sum of the residues is then equal to $2g-2$.
This can be seen as a generalization of the Gauss-Bonnet formula: indeed any compact surface with a flat metric having finitely many conical singularities gives rise to a connection whose residues are all real and the Gaussian curvature is concentrated at these singularities as a Dirac mass with value $2\pi-\alpha$ where $\alpha$ is the angle of the cone.
\end{remark}

From \Cref{thm:sc}, the similarity charts $\phi$ satisfy locally that
\[\phi = a + b\int \prod_{\stackrel{k=1}{z_k\neq\infty}}^m\left(z-z_k\right)^{\res_k}.\]

\begin{remark} This is a generalization of the \emph{Schwarz-Christoffel formula}, which expresses the conformal mapping from a half plane to a simply connected polygon. See \cite{C} for more details.
\end{remark}

\subsection{Example}\label{sub:ex}

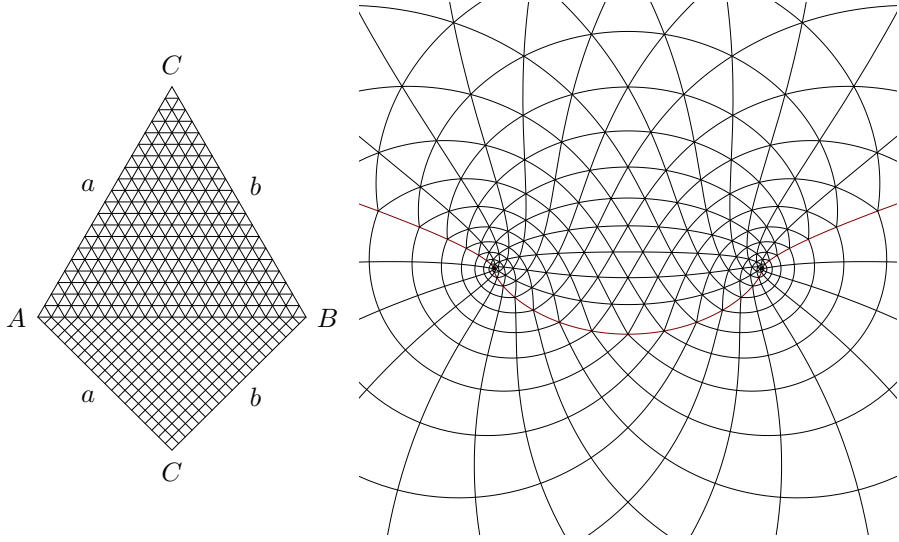
\begin{figure}
\begin{tikzpicture}
  \node at (0,0) {\includegraphics[scale=0.5]{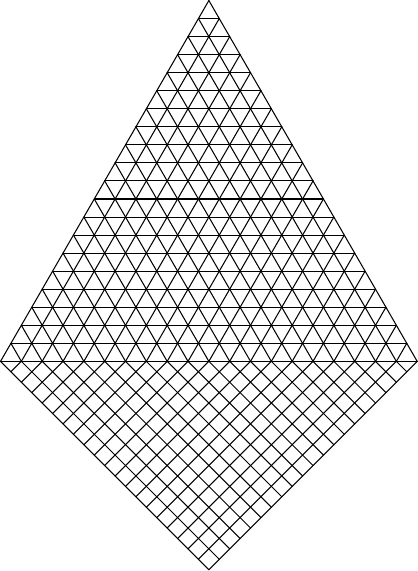}};
  \node at (6,0) {\includegraphics[scale=0.5]{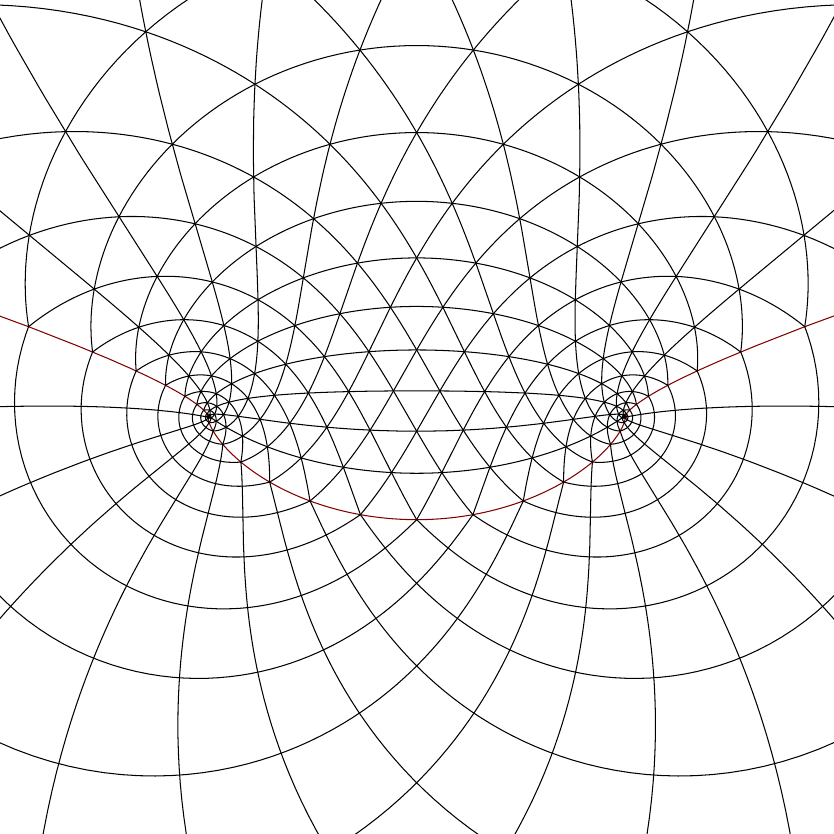}};
  \node at (-1.1,1.1) {$a$};
  \node at (1.1,1.1) {$b$};
  \node at (-1.1,-1.7) {$a$};
  \node at (1.1,-1.7) {$b$};
  \node at (-2.05,-0.65) {$A$};
  \node at (2.05,-0.65) {$B$};
  \node at (0,2.7) {$C$};
  \node at (0,-2.7) {$C$};
\end{tikzpicture}
\caption{The example of \Cref{sub:ex}. Left: the edges with the same label are glued together. Right: image under a uniformization.}\label{fig:ex}
\end{figure}

In \Cref{fig:ex} we illustrate the following example: one polygon is an equilateral triangle and the other is a isocele rectangle triangle.
The three sides of the first are glued to the three sides of the second as on the figure on the left, so that the quotient is homeomorphic to the sphere.
It has three vertices $A$, $B$ and $C$, three edges and of course two polygons.
We get the following, where $\theta$ is the total angle and $\lambda$ the monodromy factor (using $\tau = 2\pi$ on the left to clarify the computations):
\begin{align*}
A: && \theta &= \frac{\tau}{6} + \frac{\tau}{8} = \frac{7}{24}\tau & \lambda &= \sqrt{2} & \res &= \frac{7}{24}-1 + \frac{\log\sqrt{2}}{2\pi i} \\
B: && \theta &= \frac{\tau}{6} + \frac{\tau}{8} = \frac{7}{24}\tau & \lambda &= 1/\sqrt{2} & \res &= \frac{7}{24}-1 - \frac{\log\sqrt{2}}{2\pi i} \\
C: && \theta &= \frac{\tau}{6} + \frac{\tau}{4} = \frac{5}{12}\tau & \lambda &= 1 & \res &= \frac{5}{12}-1 \\
\end{align*}
In the global chart $\hat\C$ where $C$ is at $\infty$, $A$ at $-1$ and $B$ at $1$ we get
\[\zeta = \frac{\res_A}{z+1} + \frac{\res_B}{z-1}.\]
On \Cref{fig:ex} we see on the left the two triangles with two different hatchings (think of them as a kind of coordinate systems); on the right the corresponding Riemann surface $\cal S^0$ has been mapped to $\hat\C$ as above and we drew the image of the hatchings using the methods of \cite{C}, Section~4.3.2 to follow geodesics associated to $\zeta$. 

\subsection{Any Christoffel symbol}\label{sub:any}

\subsubsection{From holomorphic Christoffel symbol to similarity surface}

Given an open subset $U$ of $\C$ and a holomorphic function $\zeta:U\to\C$, we can define an atlas of similarity surface whose Christoffel symbol is $\zeta$, by solving locally $N\phi = \zeta$ and taking the maps $\phi$ as charts.
The maps $\phi$ are also called \emph{straightening maps} of the Christoffel symbol $\zeta$.

\subsubsection{Simple poles}

Assume that $\zeta$ has an isolated singularity $z_0\in\C$ which is a simple pole of residue $\res$.
\begin{lemma}\label{lem:lds}
Under the conditions above, if $\res\notin \{-1,-2,\ldots\}$ then there is a holomorphic change of coordinate $z\mapsto w$ sending $z_0 \to 0$ such that similarity charts are expressed as branches of $w\mapsto w^{\res+1}$, up to affine maps, and hence the Christoffel symbol takes the expression
\[\tilde \zeta(w) = \frac{\res}{w}.\]
\end{lemma}
\begin{proof}
By a translation we can assume that $z_0=0$.
We have the power series expansion as $z\to 0$:
\[\zeta(z) = \frac{\res}{z} + a_0 + a_1z+ a_2 z^2 +\ldots\]
of which one antiderivative satisfies
\[\int\zeta = \res \times \log z + o(1)\]
hence
\[\exp\left(\int\zeta\right) = z^{\res} \exp(o(1)) = \sum_{n=0}^{+\infty} b_n z^{\res + n}\]
for some sequence $b_n\in\C$ with $b_0= 1$.
Since $z^{\res+n}$ has several values, we need to specify a meaning:
we temporarily work in the universal cover of $\C^*$, let $\log z$ be a global branch and define $z^{\res+n} = \exp((\res+n)\log z) = \exp(\res\times\log z) z^n$.
This generalized power series expansion has complex exponents but positive radius of convergence and we can integrate term by term to get
\[\int \exp\left(\int\zeta\right) = \sum_{n=0}^{+\infty} b_n\frac{z^{\res+n+1}}{\res+n+1}
=z^{\res+1}\exp(R(z))\]
for some holomorphic function $R$ of $z$ (i.e.\ without singularity at $0$). 
Then \[w = z\exp\frac{R(z)}{\res+1}\]
satisfies the claims. 
\end{proof}
Note that we did not claim local uniqueness of the coordinate $w$ that satisfies the conclusion of the lemma. We do not wish to study the uniqueness problem in this article.

\begin{remark} Let us state how to deal with the remaining cases (we give them without proofs because we will not use these normal forms).
In the case where $\res = -1$ we can also get $\tilde\zeta(w) = -1/w$, which has branches of $w\mapsto \log w$ as similarity charts up to affine maps.
In the case $\res = -k \in\{-2,-3,\ldots\}$, we can reduce to either of the normal forms $\tilde\zeta(w) = -k/w$ or $\tilde \zeta(w) = -k/w + w^{k-2}$, and these two cases are mutually exclusive.
The first one has $w\mapsto 1/w^{1-k}$ as similarity charts but the similarity charts of the second one has a complicated expression.
We may prefer an alternative normal form, for example one which gives branches of $w\mapsto\frac{-1}{w^{k-1}}+\log w$ as similarity charts up to affine maps, in which case $\ds \tilde \zeta(w) = \frac{-k}{w} + \frac{w^{k-2}}{1+\frac{w^{k-1}}{(k-1)}}$.
\end{remark}

An isolated singularity of a Christoffel symbol with a simple pole with
\[\Re( \res) >-1\]
is called \emph{bounded}.
It is called
\emph{unbounded} if
\[\Re (\res) <-1.\]
We do not give a particular name to the case $\Re(\res) = -1$.

According to the lemma above, near any \emph{bounded} isolated singularity of a similarity surface, the latter is isomorphic to a bounded sector of apex $0$ and opening $\theta = 2\pi (\Re(\res) +1)$ (if $\theta>2\pi$, the sector sits in the universal cover of $\C^*$) endowed with the canonical atlas of $\C$, and whose two sides are glued by the similarity $z\mapsto e^{2\pi i \res}z$ from $\arg z = \theta_0$ to $\arg z=\theta_0+\theta$.

Near any \emph{unbounded} singularity of \emph{non-integer} residue, this is the same but the sector is a neighbourhood of infinity and has opening $\theta = -2\pi (\Re(\res) +1)$, and the similarity is still $z\mapsto e^{2\pi i \res}z$ and maps $\arg z = \theta_0$ to $\arg z=\theta_0-\theta$.

We leave to the reader the task to determine local models for the remaining cases if they wish. Such models will not be used here.

\medskip

\subsubsection{A convenient change of variables for the unbounded type singularities}

Assume that $0$ is a simple pole of $\zeta$ of unbounded type ($\Re \res <-1$) and consider the Laurent series expansion
\[ \zeta(z) = \frac{\res}{z} + a_0 + a_1 z+ a_2 z^2 +\ldots
.\]
Close to $0$ the term $\frac{\res}{z}$ is dominant, so it is natural to try and compare the similarity charts $\phi$ of $\zeta$ to the similarity charts $z\mapsto z^{1+\res}$ of the ``ideal'' Christoffel symbol $\frac{\res}{z}$.
For this we make the change of variable 
\[\tilde z = z^{1+\res}.\]
Note that the correspondence $z\longleftrightarrow \tilde z$ is a bijection from the universal cover $U\C^*$ of $\C^*$ to $U\C^*$, via
\[ \log \tilde z = (1+\res)\times\log z
\]
where $\log : U\C^* \to \C$ is a well-defined bijection.
From now and up to the end of the present~\ref{sub:any}, we assume that $z$ and $\tilde z$ live in $U\C^*$.

\subsubsection{Tending to $0$}

If $\res\in\R$ then $|z|\tend 0$ $\iff$ $|\tilde z|\tend+\infty$.
But if $\res\notin\R$ then the set $|z|<r$, where $r>0$, corresponds to a set of $\tilde z\in U\C^*$ whose boundary spirals inward to $0$ as $(\Im\res)\times( \arg \tilde z) \to -\infty$ and outward to $\infty$ as $(\Im\res)\times( \arg \tilde z) \to +\infty$: indeed the set of values of $\log \tilde z$ is the half plane
\begin{equation}\label{eq:HR:def}
\HU(r) := \setof{\tilde z\in U\C^*}{\Re\frac{\log \tilde z}{1+\res} < \log r},
\end{equation}
it is delimited by a straight line that is not vertical and this half plane contains a neighbourhood of $+\infty$ in $\R$, because $\Re (1+\res) <0$.
In any cases, we have $|z|\tend 0$ iff $\tilde z$ enters and stays in every $\HU(r)$.

\subsubsection{$\tilde z$-sectors}

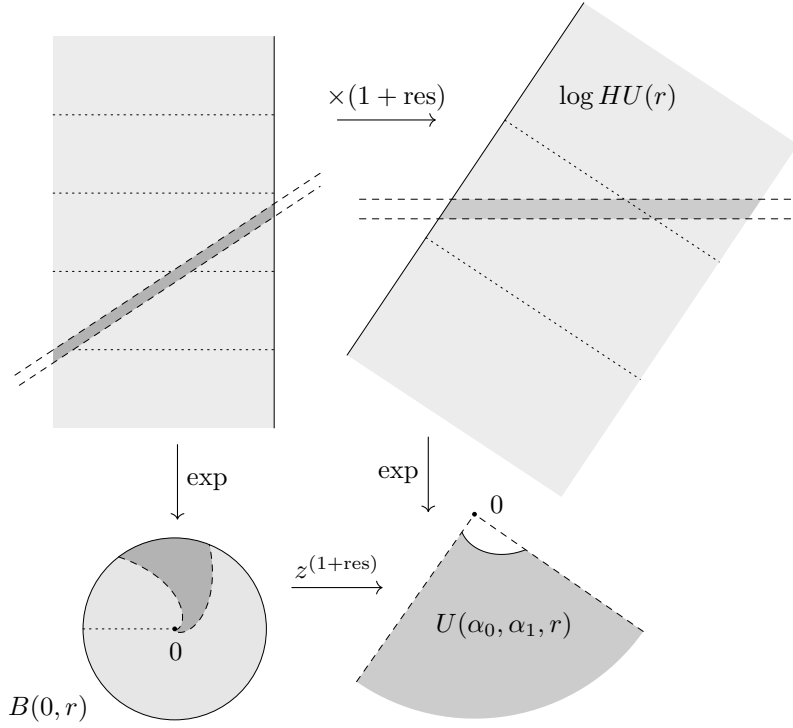
\begin{figure}
\begin{tikzpicture}
\node at (-0.3,0) {\includegraphics[scale=0.3]{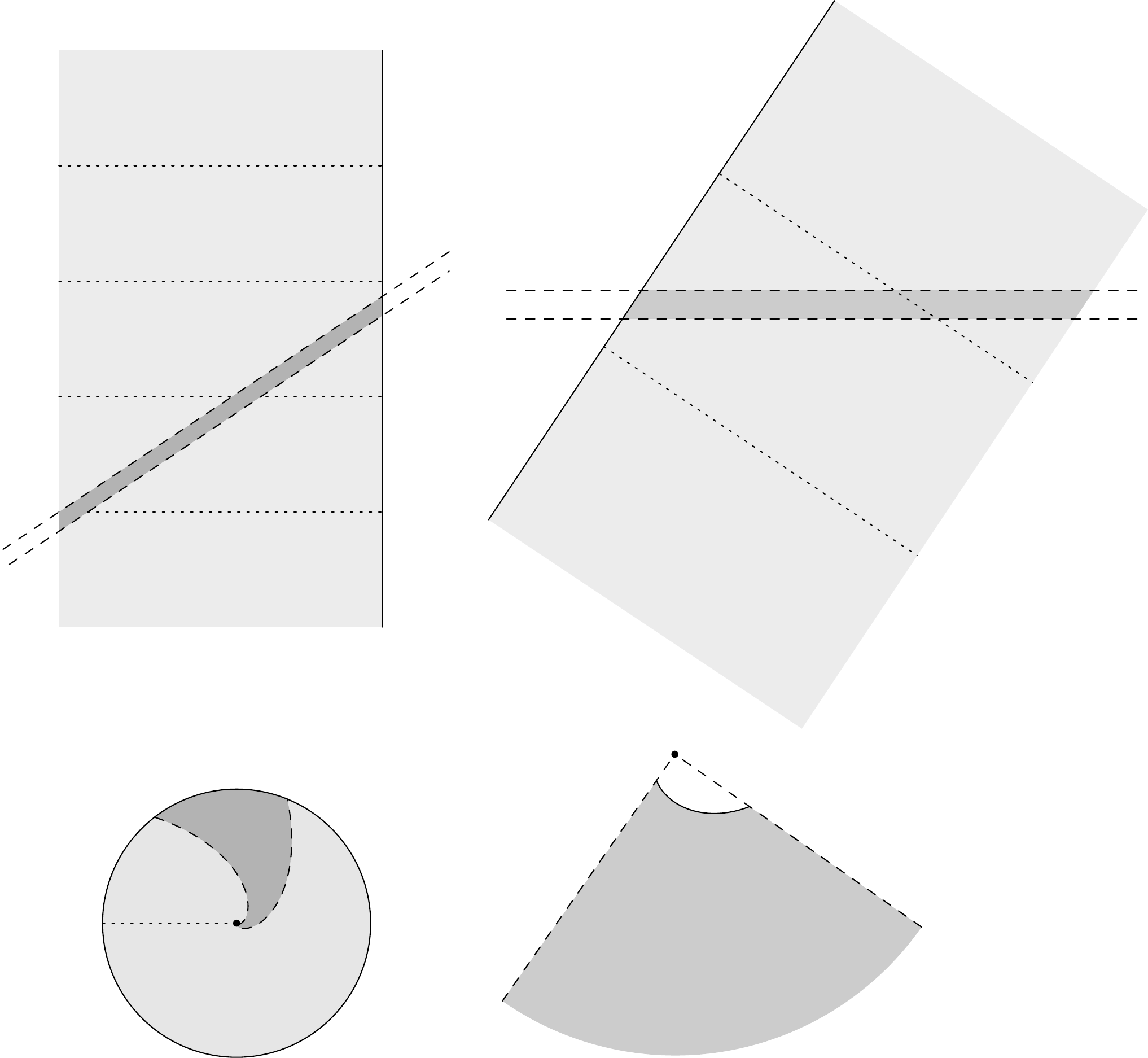}};
\node at (-3.32,-3.85) {$0$};
\node at (0.9,-1.9) {$0$};
\node at (-5,-4.6) {$B(0,r)$};
\draw[->] (-1.2,3.2) -- node[above] {$\times (1+\res)$} (0.1,3.2);
\draw[->] (-3.3,-1.1) -- node[right] {$\exp$} (-3.3,-2.1);
\draw[->] (0,-1) -- node[left] {$\exp$} (0,-2.0);
\draw[->] (-1.8,-3.0) -- node[above] {$z^{(1+\res)}$} (-0.6,-3.0);
\node at (1,-3.5) {$U(\alpha_0,\alpha_1,r)$};
\node at (2.5,3.5) {$\log \HU(r)$};
\end{tikzpicture}
\caption{Example of domain $U(\alpha_0,\alpha_1,r)=\setof{\tilde z}{\arg \tilde z \in [-\alpha_0,\alpha_1]\text{ and } |z|<r}$ in the case $\Re \res <-1$. Dots correspond to $\arg z \equiv \pi \bmod 2\pi$, dashes to $\arg \tilde z =\alpha_0$ or $\alpha_1$.}\label{fig:spi-5}
\end{figure}

Since $\tilde z\in U\C^*$, there is a well defined function $\arg \tilde z$ taking values in $\R$.
The sector in $\tilde z$ coordinates defined by $\arg \tilde z \in (\alpha_0,\alpha_1)$ corresponds in $\log \tilde z$ coordinate to the strip $\Im \log \tilde z \in (\alpha_0,\alpha_1)$.
In this strip, the subset for which the corresponding $z$ belongs to $B(0,r)$ is the intersection with $\log \HU(r)$.
In $\tilde z$ coordinate, this intersection is a subset of the sector with a spiralling neighbourhood of $0$ removed, see \Cref{fig:spi-5}.
We call this set
\[U(\alpha_0,\alpha_1,r) = \setof{\tilde z\in U\C^*}{\alpha_0<\arg \tilde z<\alpha_1,\ |z|<r} \subset U\C^*
\]
where $\arg \tilde z =\Im \log \tilde z$ is well-defined on $U\C^*$ and $\tilde z\mapsto z$ also.
For any given $U=U(\alpha_0,\alpha_1,r)$ we have
\[\log |z| = \Re\left(\frac{1}{1+\res}\right)\log |\tilde z| + \cal O(1)\] where $\cal O(1)$ is quantity bounded over all $U$ with a bound that depends continuously on $(\alpha_0,\alpha_1,r)$ for $\alpha_0<\alpha_1$ and $0<r$.
In particular, under the constraint $\tilde z\in U(\alpha_0,\alpha_1,r)$, we have $|z|\tend 0 \iff |\tilde z|\tend+\infty$.
Finally, note that $U$ contains the sector defined by $\arg \tilde z\in (\alpha_0,\alpha_1)$ and $|\tilde z|> R$ where $R=R(\alpha_0,\alpha_1,r,\res)$ can be computed easily.

\begin{lemma}\label{lem:uch}
  Assume that $\Re \res <-1$.
  There is some $r>0$ that depends on $\zeta$ and such that the following holds.
  Let $\phi$ be any solution\footnote{They exists since $B(0,r)-\{0\}$ lifts to a simply connected subset of $U\C^*$} on the part of $U\C^*$ on which $|z|<r$, of the equation $\phi''(z)/\phi'(z) = \zeta(z)$. The quantity $\phi$ can be interchangeably considered as depending on $z$ or $\log z$ or $\tilde z$, etc.
  Then there exists $a\in\C^*$ such that:
  \begin{enumerate}
  \item there exists a holomorphic function $z\mapsto r_1(z)$ on $B(0,r)$ such that
  \[\frac{\partial \phi}{\partial \tilde z} = a \times (1+r_1(z))\]
  and $r_1(0)=0$;
  \item for all $\alpha_0,\alpha_1\in\R$ with $\alpha_0<\alpha_1$, there exists $C,C'\in\R$ and a holomorphic function $r_2 : U(\alpha_0,\alpha_1,r) \to \C$ such that we have, $\forall \tilde z \in  U(\alpha_0,\alpha_1,r)$, 
  \[\phi = a\tilde z + r_2(\tilde z)\]
  with
  \[ |r_2(\tilde z)|\leq \max\left(|z|\!\times\!|\tilde z|,1\right)(C'+C\log |\tilde z|)
  .\]
  In particular
  \[\phi \underset{|\tilde z|\to \infty}\sim a \tilde z.\]
  \end{enumerate}
\end{lemma}
\begin{proof}
We restrict $\zeta$ to some $B(0,r)$ so that
\[ \zeta(z) = \frac{\res}{z} +h(z)\]
with $h$ bounded, say $|h|\leq M$.
Let $v = \log z$ and let us express the Christoffel symbol $\hat \zeta$ in the $v$ coordinate: since $z=e^v$, formula~\eqref{eq:zetachvar} gives
\[\hat\zeta(v) = e^v \zeta(e^v) + 1\]
and using the development of $\zeta$ above we get
\[\hat\zeta(v) = 1+\res + e^v h(e^v)
.\]
Let $\hat\phi(v) = \phi(z)$. Then $\hat\phi$ is a solution of $\hat\phi'' / \hat\phi' = \hat \zeta$. Integrating on the half plane $\Re v < \log r $ we get that
\[\log (\hat\phi'(v)) = c_0 + (1+\res) v + g(e^v)\]
where $g$ is the antiderivative of $h$ mapping $0$ to $0$ and $c_0\in\C$ is a constant.
Recall $\log \tilde z = (1+\res) v$.
Hence
\[\frac{\partial\phi }{\partial \log\tilde z} = \frac{1}{1+\res}\hat\phi'(v)  = a \tilde z (1+ r_1(z))\]
for $a=e^{c_0}\in\C^*$ and for the holomorphic bounded function $r_1:B(0,r)\to\C$ mapping $0$ to $0$ and defined by $r_1(z) = \exp(g(z))-1$.
Since $\partial \phi/\partial \log \tilde z = \tilde z\, \partial \phi/\partial \tilde z$, this proves the first point.

We have $|r_1(z)|\leq c|z|$ for some $c>0$ and thus
\begin{equation}\label{eq:uv}
\left|\frac{\partial\phi}{\partial \log \tilde z} - a \tilde z\right| \leq  c' |\tilde z|\times|z|
\end{equation}
with $c' = c|a|$.
Let us abbreviate $U=U(\alpha_0,\alpha_1,r)$.
The set of $\tilde z\in \partial U$ such that $|z|=r$ is a compact set for which $\phi$ is hence bounded.
Any element $\tilde z\in U$ satisfies $\log \tilde z = t+\log \tilde z_0$ for a $t>0$ and a $\tilde z_0$ as above.
We have that $|\tilde z|$ and $e^t$ are comparable in the sense that their quotient is bounded away from $0$ and $\infty$ by factors depending only on $U$.
Similarly $|z|$ and $e^{-\alpha t}$ are comparable where
\[\alpha = -\Re(1/(1+\res))>0.\]
Consider the path $t\mapsto \tilde z(t)=e^t \tilde z_0$ and the corresponding path $\log \tilde z(t) = \log \tilde z_0 + t$.
Let us integrate \cref{eq:uv} along the latter path with respect to $dt$.
We have $\partial \phi/\partial \log \tilde z = \frac{\partial}{\partial t} \phi$.
Note that an antiderivative of $\tilde z$ with respect to the variable $t$ is $\tilde z$.
On the right hand side of \cref{eq:uv}, note that $|\tilde z(t)| = |\tilde z_0|  e^t$ and $|z(t)|= |z_0| e^{-\alpha t}$.
Hence
\[ \left| \phi(z(t)) - a \tilde z(t) - \phi(z(0)) + a \tilde z_0 \right| \leq  c' |\tilde z_0| |z_0| \int_0^t e^s e^{-\alpha s} ds \]
If $\alpha\geq 1$ then $e^s e^{-\alpha s} \leq 1$ whence $\int_0^t e^s e^{-\alpha s} ds \leq t$.
If $\alpha\leq 1$ then $e^s e^{-\alpha s}\leq e^{(1-\alpha)t}$  whence $\int_0^t e^s e^{-\alpha s} ds \leq  te^{(1-\alpha)t}$.
\end{proof}

Recall that $U(\alpha_0,\alpha_1,r)$ contains the sector defined by $\arg \tilde z\in (\alpha_0,\alpha_1)$ and $|\tilde z|> R$ where $R=R(\alpha_0,\alpha_1,r,\res)$ can be computed easily.

\begin{corollary}\label{cor:uch2}
  In the conditions of \Cref{lem:uch}, for any $\alpha_0<\alpha_0'<\alpha_1'<\alpha_1$ there exists $0<r''<r'\leq r$ such that for any $\theta\in\R$, if one denotes $S$ the (simply connected) set of $z\in U\C^*$ for which $\tilde z \in U(\alpha_0,\alpha_1,r')$, then the restriction of $\phi$ to $S$ avoids $0$, and has a lift to $U\C^*$ that is analytic, injective and whose image contains $a \times U(\alpha_0',\alpha_1',r'')$, where $a$ is the factor of \Cref{lem:uch}.
\end{corollary}
\begin{proof}
It follows from the definitions and from the estimate of the second point in \Cref{lem:uch} using, for instance, the argument principle on $\phi$ if $\alpha_1-\alpha_0<2\pi$ or for the general case the argument principle on $\log \phi$.\footnote{To avoid the $\log$ one can instead invoke homological arguments.}
\end{proof}

For the next lemma, we call bounded (resp.\ unbounded) sector of apex $0$ the subsets of $\C$ of the form
\[ \setof{z\in\C^*}{\alpha_0<\arg z<\alpha_1\text{ and } |z|<r} \]
resp.\ 
\[ \setof{z\in\C^*}{\alpha_0<\arg z<\alpha_1\text{ and } |z|>r} \]
with $\alpha_1-\alpha_0\in(0,2\pi]$ and $r>0$.

We recall that a geodesic is a curve whose image in affine charts are locally $\R$-affine maps and non-constant.
Given a meromorphic function $\zeta$ on an open subset of $\C$, we saw that we can associate to it an atlas of similarity surface on the complement of the poles, hence an associated notion of geodesic.\footnote{It is not hard to see that geodesics are the non-constant solutions $\gamma$ of the O.D.E.\ $\gamma'' = - (\gamma')^2 \zeta \circ \gamma$, which we will not use here.}

\begin{lemma}\label{lem:ucs}
For any geodesic parametrized by $t$ and tending to a bounded resp.\ unbounded singularity $z_0$ of $\zeta$ (this excludes the case $\Re\res = -1$):
\begin{enumerate}
\item The geodesic takes a finite time to reach the singularity in the bounded case, and an infinite time in the unbounded case.
\item There exists $r>0$ and a simply connected subset $U$ of $B(z_0,r)-\{z_0\}$, and an (injective) affine chart $\phi:U \to \C$ whose image is a bounded resp.\  unbounded sector with apex $0$, and such that $\phi$ sends the geodesic, for $t$ large enough, into the bisecting line of the sector, tending to $0$ for a bounded singularity and to $\infty$ for an unbounded one.
\end{enumerate}
\end{lemma}

\begin{proof}
We distinguish two overlapping cases: $\res\notin\{-2,-3,\ldots\}$ and $\Re\res <-1$.

In the case $\res\notin\{-2,-3,\ldots\}$ the lemma follows from the local model deduced from \Cref{lem:lds}, both in the bounded and unbounded case.
In the bounded case, this model is the quotient of a sector (with the canonical atlas of $\C$) by a similarity fixing its apex, and a geodesic is a locally straight line in the atlas, hence one tending to the singularity must end up being a ray through the apex, as can be seen by the same trick as for the fact that a billiard trajectory in a sector can only undergo finitely many bounces: if we have a a straight line geodesic that is not aimed at the apex then patching finitely many copies of the sector by the similarity covers a sector of angle $>\pi$ that will contain the geodesic until it escapes.
In the unbounded case, tending to $0$ in the $z$ coordinate corresponds to tending to infinity in the affine model, hence the time has to tend to infinity for $z(t)$ to reach $0$, and similarly to the previous case, one see that the geodesic eventually remains in one fundamental domain.
Let us apply a translation to the chart so that the geodesic radiates from $0$ in this chart.
Then the fundamental domain still contains an unbounded sector with apex $0$.

In the case $\Re\res <-1$ we first focus on proving that the time tends to infinity. For this we recall the first point of the conclusion of \Cref{lem:uch}: $\frac{\partial \phi}{\partial  \tilde z} = a \times (1+r_1(z))$, which is equivalent to
\[\frac{\partial \phi}{\partial \log \tilde z} = a \times \tilde z\times (1+r_1(z))\]
where $r_1$ is holomorphic and bounded in $z\in B(0,r)$ and $r_1(0)=0$.
By taking $r$ smaller we can assume that $r_1'(z)$ is bounded too.
We have
\[\arg \frac{\partial \phi}{\partial \log \tilde z} = \arg a + \Im \log\tilde z + \Im \log(1+r_1(z)).\]
By taking $r$ small enough we can assume that the last summand has modulus $<1/10$ for all $\tilde z \in \HU(r)$ (we recall that $\tilde z\in  \HU(r) \iff |z|<r$ and that both $z$ and $\tilde z$ live in $U\C^*$).
We have
\[ \frac{\partial\log(1+r_1(z))}{\partial \log \tilde z} = \frac{r'_1(z)}{1+r_1(z)} \times \frac{z}{1+\res}
.\]
By taking $r$ small enough we can assume that this quantity also has modulus $<1/10$ for all $\tilde z \in \HU(r)$.
We let $t_0$ for which for all $t\geq t_0$ then $|z(t)|<r$ and we now only consider times $t\geq t_0$.
As $t$ varies and $z$ varies accordingly along the geodesic, the variables $\tilde z$, and $\phi$ can be considered as depending on $t$.
The geodesic has the property that $\frac{\partial \phi}{\partial t}$ remains constant, in particular the argument of this quantity remains constant.
Since
$\frac{\partial \log \tilde z}{\partial t} = \frac{\partial \phi}{\partial t} \left/ \frac{\partial \phi}{\partial \log \tilde z} \right.
$
we have that, as a function of $t$, the quantity $\log \tilde z$ follows a field line of a vector field whose direction is given by
\begin{equation}\label{eq:theta}
\theta(\log \tilde z) = c -\,\Im \log \tilde z - \Im \log (1+r_1(z))
\end{equation}
for some constant $c\in\R$ that depends on $a$ and on the direction of the geodesic in the $\phi$ coordinate.
Isoclines of this vector field are given by the level lines of $\theta \bmod 2\pi$.
By the estimates on $\log(1+r_1(z))$ and its derivative with respect to $\log \tilde z$ we have that:
\begin{itemize}
\item the isocline of angle $\beta$ is a union of curves $I_{\beta+2k\pi}$, $k\in\Z$, each contained within distance $<1/10$ of the horizontal line of ordinate $c-\beta$;
\item the tangent to $I_\beta$ at every point makes an angle $<\asin \frac{1}{10}$ with the horizontal direction;
\item each of these curves extends from $\partial \HU(r)$ to infinity;
\item each is asymptotic at infinity to the horizontal curve of imaginary part $c-\beta$.
\end{itemize}
We have that $\theta(\log \tilde z)$ strictly decreases when $\log \tilde z$ moves upwards along a vertical line; we will say that $\log \tilde z$ is \emph{above} $I_\beta$ if $\theta(\log \tilde z)<\beta$ and \emph{below} if $\theta(\log \tilde z)>\beta$ (note the inversion).
We have $\asin \frac{1}{10}<\pi/4$ which, together with the above points, implies the following:
Consider the isocline $I_\beta$ with $\beta = 2k\pi-\pi/4$.
If $\log \tilde z(t_0)$ is below $I_\beta$ then for $t>t_0$ the curve $\log \tilde z(t)$ is disjoint from $I(\beta)$ so $\log \tilde z(t)$ stays below $I_\beta$.
Similarly if $\log \tilde z(t_0)$ is above $I_\beta$ with $\beta = 2k\pi+\pi/4$ then after $t_0$, $\log \tilde z$ stays above $I_\beta$.
Let $\beta_*$ denote the isocline on which the geodesic starts.
We take the smallest $\beta_0\in 2\pi\Z -\pi/4$ such that $\beta_0\geq\beta_*$ and the biggest $\beta_1 \in 2\pi\Z +\pi/4$ such that $\beta_1\leq\beta_*$ then $\log \tilde z(t)$ has to remain between the isoclines of index $\beta_0$ and $\beta_1$.
Note that the value of $\beta_1-\beta_0$ is either $2\frac{\pi}{4}$ or $2\pi+2\frac{\pi}{4}$, depending on $\alpha_*$.
In particular, $\log\tilde z(t)$ stays within a horizontal strip of height $\leq 2\pi + 2\frac{\pi}{4}+\frac{2}{10}$.
Now the only way for $z(t)$ to tend to $0$ while $\log \tilde z(t)$ has bounded imaginary part is that $\Re \log\tilde z(t) \tend +\infty$.
By the second point of \Cref{lem:uch}, $\phi(t)\sim a \tilde z(t)$ so  $\phi(t)$ tends to infinity, hence the geodesic is defined for an infinite amount of time.
This proves the first claim.

Since geodesics are affine maps of $t$ in affine charts, we have $\phi(t) = a't+b$ for some $a'\in\C^*$ and $b\in \C$.
The previous analyis applies to any affine chart so we may sutract $b$ from $\phi$ and hence assume that $b=0$.
Note that we only need to consider $t>0$. So $\phi(t)\in a'\R_+$ and
$\arg \phi(t)$ is constant.
Since $\tilde z(t) \sim \phi(t)/a$, it follows that $\log \tilde z(t)$ is asymptotic to a horizontal line, let $\alpha_\infty$ denote its imaginary part.
Let $\alpha_0<\alpha_\infty <\alpha_1$.
We choose any $\alpha'_0$ and $\alpha'_1$ such that $\alpha_0<\alpha'_0<\alpha_\infty<\alpha'_1<\alpha_1$ and apply \Cref{cor:uch2}: the set  $\tilde z \in U(\alpha_0,\alpha_1,r')$ corresponds bijectively to a set of values of $\phi(z) \in U\C^*$ that contains $a\times U(\alpha'_0,\alpha'_1,r'')$.
This set contains a sector with apex $0$ and central axis $a'\R_+$, possibly with smaller angle.
\end{proof}

\subsection{Resting place and paths to a bounded singularity}\label{sub:rest}

As explained in \Cref{sec:simsurfdef}, given a continuous path $\gamma:[0,1]\to\cal S'$ and a germ of chart $[\phi_0]$ (the square brackets refer to the fact that the germ is an equivalence class) at the beginning of the path $\gamma(0)$, there is a way to follow this germ along the path $\gamma(t)$ as a family of germs $[\phi_t]$.
The value $\phi_t(\gamma(t))$ is well-defined and the function $\tilde\gamma: t\mapsto \phi_t(\gamma(t))$ is called the \emph{development} of the path $\gamma$ with respect to $[\phi_0]$.
Let us call the point $\phi_1(\gamma(1))$ the \emph{resting place} of the development (\emph{principal value} could be a nice alternative denomination).

Can this be done for a path tending to a bounded singularity?
The germ will in general degenerate, but the development will have a limit, provided we take some precaution.
To see this, let us work in a Riemann chart, as we saw there exists, for which the singularity is at $0$ and where for some $\alpha\in\C$ with $\Re \alpha>0$, the similarity charts include the branches of $z\mapsto z^\alpha$.
Then if $\Im\alpha\neq 0$ and $z$ winds along a circle around $0$ many times in the sense opposite to the sign of $\Im\alpha$, then $z^\alpha$ tends to infinity.
Or if $z$ spirals down to $0$ too slowly, winding this way, then $z^\alpha$ still tends to $\infty$.
On the other hand, if $z$ tends to $0$ along a well-defined tangent direction, then $z^\alpha$ tends to $0$, and this limit $0$ does not depend on the direction.
Hence the development $\tilde \gamma$ w.r.t.\ $[\phi_0]$ of a path $\gamma$ tending to a bounded singularity may or may not converge, but if we replace its last moments by a straight line (in any given Riemann chart), its development has a limit that is independent of the replacement done, and we can thus associate to $\gamma$ and $[\phi_0]$ a well-defined resting place.

Consider two paths $[0,1]\to\cal S$ that both start from the same point and end at the same bounded singularity.
Assume that they take their values in $\cal S'$ except at $t=1$ and
assume that they are homotopic in $\cal S$ by a homotopy rel.\ $\{0,1\}$ (i.e.\ that fixes the ends) and that takes value in $\cal S'$ between the ends.
Then their developments starting from the same germ $[\phi_0]$ have the same resting place (indeed the development can be followed under the homotopy and the resting place stays immobile).

%
% ----------------------------------------------------------
%

\section{Holomorphic dependence}\label{sec:holodep}

As the section title suggests, we investigate here various holomorphic dependence properties:
\begin{itemize}
\item as a function of the Christoffel symbol $\zeta$: those of the special coordinates (where the Christoffel symbol takes a simple expression) in \Cref{sub:hdsc}, of developing maps (\Cref{sec:hd:rp}), of saddle connections (\Cref{sub:hd:sc}),
\item as a function of the polygons when they are glued into a surface homeomorphic to the sphere and which, when the vertices are removed, is conformally isomorphic to the Riemann sphere $\hat\C$ minus finitely many points: of the position of these points on $\hat\C$ and associated residues (\Cref{sub:hd:p}).
\end{itemize}

\subsection{Holomorphic dependence of special coordinates.}\label{sub:hdsc}

We now state and prove a holomorphic dependence statement related to \Cref{lem:lds}.
We will restrict to parameter spaces of dimension one to simplify the presentation and since in the application in \Cref{part:beltrami} we only need that.

Let $\zeta_0$ be a particular map satisfying the conditions of \Cref{lem:lds}.
We assume that $r>0$ and that we are given complex numbers
$\res_\tau$ and $z_\tau$ that depend holomorphically on $\tau\in B(0,r)$.
We assume that $\eps>0$, that $z_\tau\in B(z_0,\eps)$ and that we have a family
$(\tau,z)\mapsto\zeta_t(z)$ for $\tau\in B(0,r)$ and $z\in B(z_0,\eps)-\{z_\tau\}$ given by
\[\zeta_\tau(z) = \frac{\res_\tau}{z-z_\tau} + h_\tau(z)\]
where the map $(\tau,z)\in B(0,r)\times B(z_0,\eps)\mapsto h_\tau(z)\in\C$ is analytic.

\begin{lemma}\label{lem:ldshol}
Under these conditions, there exists $\eta>0$ and $r'>0$ such that for $\tau\in B(0,r')$, one can choose the change of variable $(z\mapsto w)$ of \Cref{lem:lds} to be defined on $B(z_0,\eta)$ and depend holomorphically on $\tau$.
\end{lemma}
\begin{proof}
In \Cref{app:pf:lem:ldshol}.
\end{proof}

If for $\tau=0$ we have $\res_0\in \{-2,-3,\ldots\}$ then $\zeta_0$ is not in the situation of \Cref{lem:lds}. Instead we adapt \Cref{cor:uch2}. Recall the definition of the $\tilde z$-coordinate given before~\Cref{cor:uch2}: $\tilde z=z^{1+\res}$, which is well-defined if $z$ and $\tilde z$ are considered as elements of the universal cover $U\C^*$ of $\C^*$.
Recall also the definition of the sets $U(\alpha_0,\alpha_1,r) = \setof{\tilde z\in U\C^*}{\alpha_0<\arg \tilde z<\alpha_1,\ |z|<r}$.
Complete $U\C^*$ into a topological space $\wh U\C^* = U\C^*\cup\{\hat 0\}$ by adding a point $\hat 0$ whose neighbourhoods consists in the sets containing $V_\eps=\{\hat 0\}\cup$ the set of points whose projection to $\C^*$ has modulus $<\eps$.

\begin{lemma}\label{lem:geodeschol}
Assume that $\Re \res<-1$.
Denote $\tilde z = z^{1+\res_\tau} \in U\C^*$ where $z\in U\C^*$.
Then there exists $\eps>0$ and $0<r''<r'<r$ such that for any $\theta\in\R$ and any $\tau\in B(0,\eps)$, there exists a simply connected open subset $S_\tau\subset U\C^*$, $a_\tau\in\C^*$
and a holomorphic function $\phi_\tau : S_\tau \to\C$ such that
\begin{itemize}
\item $z_\tau+\pi(S_{\tau}) \subset B(z_0,r)$ where $\pi$ is the projection $U\C^*\to\C^*$;
\item $S_\tau$ depends continuously on $\tau$ by an isotopy of $\wh U\C^*$;
\item in $\tilde z$ coordinate, the image of $S_\tau$ contains $U(\theta-\pi/2,\theta+\pi/2,r')$;
\item $\phi_\tau''(z)/\phi_\tau'(z) = \zeta_\tau(z_\tau+\pi(z))$;
\item $\phi_\tau$ depends holomorphically on $\tau$;
\item $a_\tau$ depends holomorphically on $\tau$;
\item $\phi_\tau (z) \sim a_\tau \tilde z$ when $z\tend 0$ within $S_\tau$, uniformly\footnote{I.e.\ for all $\eta>0$ there is $\eta'>0$ such that for all $\tau\in B(0,\eps)$, and all $z\in B(0,\eta')\cap S_\tau$, $\Big|\frac{\phi_\tau(z)}{\tilde z} -a_\tau\Big|<\eta$.} over $\tau$;
\item the image of $\phi_\tau$ contains $a_\tau \times U(\theta-\pi/4,\theta+\pi/4,r'')$.
\end{itemize}
\end{lemma}
\begin{proof}
By the change of variable $z'=z-z_\tau$, possibly reducing $r$, it is enough to treat the case where $z_\tau=0$ for all $\tau$.
Recall that $\tilde z\mapsto z$ is well-defined from $U\C^*$ to $U\C^*$ and note that $z$ depends analytically on $(\tilde z,\res)$.
Let $\tilde \zeta_\tau$ be the pull-back of $\zeta_\tau$ in $\tilde z$-coordinate: it depends analytically on $(\tau,\tilde z)$ and hence there exists a global solution $\phi=\tilde\phi_\tau$ of $\phi''/\phi'=\tilde \zeta_\tau$ on $\HU(r)\subset U\C^*$ that depends holomorphically on $\tau$.

The variable $\log \tilde z$ depends analytically on $\tau$ and by looking at the proof of \Cref{lem:uch}, one sees that $r_1(z)$ and thus $a$ also depend analytically on $\tau$.
Moreover, in \Cref{lem:uch}, all estimates can be taken uniform provided $\tau$ is small enough and $\alpha_0$ and $\alpha_1$ are fixed.
The same holds for \Cref{cor:uch2}.
\end{proof}

\subsection{Holomorphic dependence: resting place}\label{sec:hd:rp}

Let $U$ be an open subset of $\C$.
Let $\cal G$ be the set of paths $\gamma:[0,1]\to U$ that satisfy $\forall s\in[0,1]$, $\gamma(s)\neq \gamma(1)$.
It is the disjoint union for $z_0\neq z_1$, both in $U$, of the sets $\cal G(z_0,z_1)$ of those elements $\gamma\in \cal G$ such that $\gamma(0)=z_0$ and $\gamma(1)=z_1$.
We consider the homotopy classes $[\gamma]$ in $\cal G(z_0,z_1)$, i.e.\ the path-connected components of this set, and denote $\sim$ the associated equivalence relation.
Let $\Gamma(z_0,z_1) = \cal G(z_0,z_1)/\sim$ and $\Gamma = \bigcup \Gamma(z_0,z_1)$ where the union is over all pairs $z_0,z_1\in U$  such that $z_0\neq z_1$.

To state the next lemma we endow the set $\Gamma$ with a complex manifold structure, whose charts are locally given by the position of the endpoints $z_0$, $z_1$.
For this we need to determine when two homotopy classes are close when their endpoints do not match, i.e.\ we need a topology on $\Gamma$.

%Consider two distinct points $z_0\neq z_1$.
%Consider the set $\cal G(z_0,z_1)$ of paths $\gamma:[0,1]\to U$ from $z_0$ to $z_1$ such that $\gamma(s)\neq z_1$ when $s\neq 1$.
%Consider a homotopy class $[\gamma]$ in the set above, i.e.\ a path-connected component of this set.
%To state the next lemma we endow the set $\Gamma$ of all $(z_0,z_1,[\gamma])$ with a complex manifold structure.
%For this we need to determine when two homotopy classes are close when their endpoints do not match, i.e.\ we need a topology on $\Gamma$.

Given two paths $\gamma$, $\tilde\gamma$ in $\C$, let \[d(\gamma,\tilde \gamma) = \sup_{s\in[0,1]} |\tilde\gamma(s)-\gamma(s)|
.\]
%Let $\cal G$ be the set of paths $\gamma$ in $U$ that satisfy $\forall s\in[0,1]$, $\gamma(s)\neq \gamma(1)$: it is the disjoint union of the $\cal G(z_0,z_1)$ for $z_0\neq z_1$, both in $U$.
Consider the quotient map
\[ Q: \gamma\in\cal G\mapsto [\gamma]\in \Gamma
\]
and endow $\Gamma$ with the quotient topology.
For all $\eps>0$ and all $\gamma\in\cal G$, let $V_{\gamma,\eps} = Q(B(\gamma,\eps))$, in other words it is the set of $[\tilde \gamma]$ where $\tilde \gamma\in\cal G$ and $d(\gamma,\tilde\gamma)<\eps$.

\begin{proposition}\label{prop:topo}
The following hold:
\begin{enumerate}
\item for all $\gamma\in\cal G$, the collection of $V_{\gamma,\eps}$ forms a basis of neighbourhoods of $[\gamma]$;%\footnote{This implies that the topology is unique, since a set is open iff.\ it is a neighbourhood of all its points.}
\item the topology on $\Gamma$ is Hausdorff separated;
\item the map $Q$ is continuous;
\item\label{item:topo:open} the map $Q$ is open (i.e.\ the sets $V_{\gamma,\eps}$ are open);
\item the map 
\[ [\gamma]\in\Gamma \mapsto (\gamma(0),\gamma(1))\in\C^2
\]
is a local homeomorphism.
\end{enumerate}
\end{proposition}

See \Cref{app:GammaTop} for a proof of these statements.
We can thus endow the set $\Gamma$ with a two dimensional complex manifold structure, locally given by the position of the points $z_0$ and $z_1$.

\begin{lemma}\label{lem:rphol}
Let $U$ be an open subset of $\C$ (or of a Riemann surface).
Consider a moving point $z_0(\tau)$ and a family of Christoffel symbols $\tau\in B(0,r)\mapsto \zeta_\tau$
of the form
\[\zeta_\tau(z) = h_\tau(z) + \frac{\res_\tau}{z-z_1(\tau)}\]
with $z_0(\tau)$, $z_1(\tau)$ and $\res_\tau$ holomorphic in $\tau$, $h_\tau(z)$ holomorphic in $z\in U$ and $\tau$, and for all $\tau$: $z_0(\tau)\neq z_1(\tau)$ and $\Re \res_\tau>-1$ (i.e.\ the singularity $z_1(\tau)$ of $\zeta_\tau$ is bounded, possibly removable).
Choose $a_\tau\in\C^*$ and $b_\tau\in \C$ holomorphic in $\tau$.
Consider the unique germ of similarity chart $\phi_\tau$ of $\zeta_\tau$ with $\phi_\tau(z_0(\tau))=b_\tau$ and $\phi_\tau'(z_0(\tau))=a_\tau$.
Assume that $X_\tau:=[\gamma_\tau]\in\Gamma$ varies continuously with $\tau$, and that $\gamma_\tau(0)=z_0(\tau)$ and $\gamma_\tau(1)=z_1(\tau)$, in particular $X_\tau$ varies holomorphically with $\tau$.
To all this data, associate the resting place $c_\tau$ (see \Cref{sub:rest}) of the development of $\phi_\tau$ along $\gamma_\tau$.
Then $c_\tau$ depends holomorphically on $\tau$.
\end{lemma}
\begin{proof}
See \Cref{app:pf:lem:rphol}.
\end{proof}

\subsection{Holomorphic dependence: following a saddle connection when varying $\zeta$}\label{sub:hd:sc}

Recall that a notion of \emph{geodesic} can be defined on similarity surfaces: as parametrized curves which in similarity charts have a constant and non-zero speed vector, a condition that is indeed independent of the chosen similarity chart.
Geodesics do not need to be injective.
A geodesic from a singularity to another one is called a \emph{saddle connection} by analogy with dynamical sytems.
We include the endpoint in its parametrization, so it is parametrized by a closed subinterval of the extended real line $[-\infty,+\infty]$.

We saw in \Cref{lem:ucs} that a geodesic reaching a bounded singularity does so in finite time.
A saddle connection between bounded singularities is hence parametrized by a bounded interval. We can reparametrize it by the segment $[0,1]$ if needed, with a real affine change of variable.

\begin{proposition}\label{prop:fseg}
Let $U$ be an open subset of $\C$ and $z^*_0$, $z^*_1$ two distinct points in $U$.
Assume $\zeta^*$ is a Christoffel symbol holomorphic on $U-\{z^*_0,z^*_1\}$ and that 
\[\zeta^*(z) = h^*(z) + \frac{\res^*_0}{z-z^*_0} + \frac{\res^*_1}{z-z^*_1}\]
with $h$ holomorphic and bounded on $U$.
Assume that $\Re \res^*_k>-1$ for $k\in\{0,1\}$ (i.e.\ the singularities $z_k$ are either bounded or better: removable).
Last, assume that a saddle connection $f^*:[0,1]\to U$ exists between the two singularities.
Then there exists $\eps>0$ such that for all Christoffel symbols $\zeta$
of the form
\[\zeta(z) = h(z) + \frac{\res_0}{z-z_0} + \frac{\res_1}{z-z_1}\]
with $h:U\to \C$ holomorphic, $\sup_{z\in U} |h(z)-h^*(z)|<\eps$, and for all $k\in\{0,1\}$, $|z_k-z_k^*|<\eps$ and $|\res_k-\res_k^*|<\eps$,
then there exists a saddle connection $f:[0,1]\to U-\{z_0,z_1\}$ from $z_0$ to $z_1$.

Under the same conditions, if $h$, $\res_k$, $z_k$ depend analytically on a complex parameter $\tau\in \Lambda$ for an open subset $\Lambda$ of $\C$, then $f=f_\tau$ depends analytically on $\tau$ in the following sense: there exists a open subset $V$ of $\C\times\Lambda \subset \C^2$, containing $(0,1)\times\Lambda$ (note that $(0,1)$ is the open interval) and for which $(t,\tau)\in (0,1)\times \Lambda\mapsto f_\tau(t)$ extends to $V$ into an analytic function of two complex variables.\footnote{Obviously, $f_\tau(0)=z_0(\tau)$ and $f_\tau(1)=z_1(\tau)$ also depend analytically on $\tau$, but we usually cannot extend $V$ to contain a neighbourhood of $\{0,1\}\times\Lambda$.}
Moreover $f_\tau$ depends continuously on $\tau$ for the uniform norm on the Riemann sphere.

If moreover $f^*$ is injective then for $\tau$ small enough, $f$ is injective.
\end{proposition}
\begin{proof}
In \Cref{app:pf:prop:fseg}.
\end{proof}

Similarly a saddle connection form a singularity to itself can be followed holomorphically, but we will not use that fact.

\medskip

We extend the notion of saddle connection to include geodesics from a bounded singularity to an unbounded one (or in the opposite direction), under the assumption that the unbounded singularity has positive total angle.
Such a geodesic can be parametrized by $t\in(0,+\infty)$, but there is no more uniqueness of the parametrization: it can be reparametrized by a linear change of the variable $t$.
Moreover, since the total angle of the infinite vertex is positive, it follows that the connection is not completely rigid: one can vary to some extent the direction, in a similarity chart, along which $\infty$ is reached.
One way to recover a form of uniqueness of the saddle connection is by specifying a marked point which it has to go through, and which must be the image of $t=1$ by the parametrization.
As in the bounded case, we include the endpoints in the saddle connection, so it is parametrized by $[0,+\infty]$.

\begin{proposition}\label{prop:fseg2}
Let $U$ be an open subset of $\C$ and $z^*_0$, $z^*_1$, $z^*_2$ three distinct points in $U$.
Assume $\zeta^*$ is a Christoffel symbol holomorphic on $U-\{z^*_0,z^*_2\}$ and that 
\[\zeta^*(z) = h^*(z) + \frac{\res^*_0}{z-z^*_0} + \frac{\res^*_2}{z-z^*_2}\]
with $h$ holomorphic and bounded on $U$.
Assume that $\Re \res^*_0>-1$ and $\Re \res^*_2<-1$ (i.e.\ $z^*_0$ is a bounded singularity and $z^*_2$ is an unbounded one).
Last, assume that a saddle connection $f^*:[0,+\infty]\to U$ exists from $z^*_0$ to $z^*_2$ with $f^*(1) = z^*_1$.
Then there exists $\eps>0$ such that for all Christoffel symbols $\zeta$
of the form
\[\zeta(z) = h(z) + \frac{\res_0}{z-z_0} + \frac{\res_1}{z-z_2}\]
with $h:U\to \C$ holomorphic, $\sup_{z\in U} |h(z)-h^*(z)|<\eps$, for all $k\in\{0,1,2\}$, $|z_k-z_k^*|<\eps$ and for all $k\in\{0,2\}$, $|\res_k-\res_k^*|<\eps$,
then there exists a saddle connection $f:[0,+\infty]\to U-\{z_0,z_2\}$ from $z_0$ to $z_2$, with $f(1)=z_1$.

Under the same conditions, if $h$, $\res_k$, $z_k$ depend analytically on a complex parameter $\tau\in \Lambda$ for an open subset $\Lambda$ of $\C$, then $f=f_\tau$ depends analytically on $\tau$ in the following sense: there exists a open subset $V$ of $\C\times\Lambda \subset \C^2$, containing $(0,+\infty)\times\Lambda$ (note that $(0,+\infty)$ is the open interval) and for which $(t,\tau)\in (0,+\infty)\times \Lambda\mapsto f_\tau(t)$ extends to $V$ into an analytic function of two complex variables.
Moreover $f_\tau$ depends continuously on $\tau$ for the uniform norm on the Riemann sphere.

If moreover $f^*$ is injective then for $\tau$ small enough, $f$ is injective.
\end{proposition}
\begin{proof}
See \Cref{app:pf:prop:fseg}.
\end{proof}

\subsection{Holomorphic dependence with respect to the polygons}\label{sub:hd:p}

Let us come back to our example of similarity surface $\cal S'$ constructed by gluing polygons, and its completion into a Riemann surface $\cal S$ that is isomorphic to the Riemann sphere with isomorphism
\[\prB : \cal S\to\hat\C.\]
The map $\prB$ is unique up to an automorphism of $\hat\C$, i.e.\ a homography.
We recall that the group of homographies is sharply $3$-transitive.

In this section we allow unbounded polygons but we assume that their angle at infinity is always $>0$, which rules out for instance vertices with residue of real part $-1$.

We can always subdivide the polygons further (adding vertices is allowed) to meet the following conditions:
\begin{enumerate}
\item All polygons are strictly convex,
\item\label{item:ubb} every unbounded polygon has exactly two unbounded edges, and their angle is $>0$,
\item $\cal S$ has at least three vertices.
\end{enumerate}

The space of bounded strictly convex polygons with $n$ indexed vertices is endowed with a complex structure, which is simply given by the affix of its vertices.
For (strictly convex) \emph{unbounded polygons} satisfying \eqref{item:ubb}, we need to add a supplementary information: a marked point in the interior of the unbounded edges, but this marked point is not considered as a vertex. For any paired unbounded sides, there is a unique orientation reversing affine map sending one side to the other and matching their marked points.
Reciprocally, for any orientation reversing affine map sending one such side to the other, we can choose (non unique) pairs of matching marked points.

\begin{remark}
The non-injectivity of this representation, due to this non-uniqueness, is not seen as a problem.
Injectivity can be recovered by replacing the marked points by the dilation factor $a\in\C^*$ of the affine map $z\mapsto az+b$ pairing the unbounded sides.
Note that, anyway, the space of polygons will be factored further by considering polygons up to affine transformations.
\end{remark}

In the two statements below, we fix the type of each polygon $P_j$---number of edges, boundedness thereof---but not the position of their vertices or marked points, and we fix the combinatorial data of which edges are paired together.

\begin{lemma}\label{lem:rh}
Under these conditions, the residues $\res _k$ depend holomorphically on the bounded polygons and marked unbounded polygons.
\end{lemma}
\begin{proof}
Recall \Cref{eq:resfin,eq:resinf} in their unified version:
\[\res  = \frac{\log \lambda + i \sigma}{2\pi i} -1\]
for some $\sigma\in\R$, and that the monodromy factor of a small loop winding anticlockwise around the vertex is $\lambda e^{i\sigma}$.

For a finite vertex, $\sigma=\theta$ is the total angle of the vertex. The monodromy is also a product of quotients $\ds\frac{a'-a}{b'-b}$ where $a$, $a'$, $b$, $b'$ are vertices.

For an infinite vertex $\sigma=-\theta$ varies continuously and the monodromy factor is also a product of quotients $\ds\frac{a'-a}{b'-b}$ where $a$, $b$ are vertices and $a'$, $b'$ are vertices or marked points. The claim follows.
\end{proof}

Choose three distinct vertices $v_1$, $v_2$, $v_3$ and choose the unique $\prB$ as above that sends them respectively to $\infty$, $0$ and $1$. Label the other vertices up to $v_p$ for some $p\geq 3$.
Let $z_k = \prB(v_k)$.

\begin{proposition}\label{prop:hsol}
Under these conditions and the ones stated before \Cref{lem:rh}, the points $z_k$ depend holomorphically on the bounded polygons and marked unbounded polygons.
\end{proposition}
\begin{proof}
In \Cref{app:pf:prop:hsol}.
\end{proof}

\section{Proofs}

The proofs of some of the statements made in the previous sections have been moved here, for a smoother reading.

\subsection{Holomorphic dependence of special coordinates: proof of \Cref{lem:ldshol}}\label{app:pf:lem:ldshol}

To make the reading easier, we will not provide every detail.
We assume that we have a change of variable for $\zeta_0$ as in \Cref{lem:lds} on $B(z_0,\eps)$.
For $\tau$ small enough, $z_\tau$ remains in $B(0,\eps)$ and $\res_\tau$ remains in $\C-\{-1,-2,\ldots\}$.
Let us note that a first change of variable $z\mapsto z-z_\tau$ obviously depends holomorphically on $\zeta$, which allows us to restrict to the case where $z_\tau$ does not change and is equal to $0$.
We go through the computations of \Cref{lem:lds} again:
\[\int\zeta_\tau = \res_\tau \log z + \int h_\tau \]
where $\int h_\tau$ is the antiderivative of $h_\tau$ that vanishes at the origin.
This map $\int h_\tau$ depends holomorphically on $\tau$ and $z$.
Then
\[\exp\left(\int \zeta_\tau\right) = z^{\res_\tau}\exp\left(\int h_\tau\right).\]

Denote $g_\tau = \exp \left(\int h_\tau\right)$ and let $g_\tau(z) = \sum b_n(\tau) z^n$ be its power series expansion with respect to $z$. We have
\[ \int z^{\res_\tau} g_\tau(z)dz  = z^{\res_\tau+1} \sum_{n=0}^{+\infty}\frac{b_n(\tau) z^n}{\res_\tau+n+1} = z^{\res_\tau +1} \frac{f_\tau(z)}{\res_\tau+1}
\]
with
\[f_\tau(z) = \sum_{n=0}^{+\infty} \frac{\res_\tau+1}{\res_\tau+n+1} b_n(\tau) z^n = \sum a_n(\tau) z^n 
.\]
This series has at least the radius of convergence of $g_\tau$.
By using Cauchy's estimates $|b_n(\tau)|\leq \eta^{-n}\sup_{B(0,\eta)}|g_\tau|$ and that $|g_\tau(z)|$ is uniformly bounded for $(z,\tau)$ small enough we get that $\tau\mapsto f_\tau(z)$ is for each fixed $z$ a convergent series of holomorphic functions of $\tau$, hence holomorphic.

Then we take $R_\tau(z) = \log \frac{f_\tau(z)}{\res_\tau+1}$. For $\eta$ and $r'$ small enough the map $f_\tau$ will be non-vanishing and we can choose the branch of the log to be analytically varying with $(z,\tau)$.
It follows that $w = z\exp(R_\tau(z)/(\res_\tau+1))$ is holomorphic in $\tau$.

For $\tau=0$, $z\mapsto w$ it is the restriction of a change of variable defined on $B(0,\eps)$. By a theorem of Hurwitz, it will still be injective on $B(0,\eta)$ for $|\tau|<r'$ provided $r'$ is small enough.

\subsection{Following the resting point: Proof of \Cref{lem:rphol}}\label{app:pf:lem:rphol}

Holomorphy is a local property and by a change of parameter it is enough to prove the claim on a subset $B(0,r')\subset B(0,r)$ of values of $\tau$.
Let $B=B(z_1(0),R)$ with $R$ small enough so that $\ov B\subset U -\{z_0(0)\}$.
Choose also $R$ small enough so that we can apply \Cref{lem:ldshol} to get a special coordinate $z\mapsto w$ on $B$ for the singularity $z_1(\tau)$ provided $\tau$ is small enough.
Choose\footnote{Project $\gamma_0(s)$ radially to the boundary of a slightly bigger disk, until $s$ is big enough so that the rest of the path is inside and then replace by a radial semgent: one gets a homotopic path.} a representative $\gamma_0$ of $[\gamma_0]$ that reaches $\partial B$ for the first time for some $s=s_1$ and then goes down to $z_1$ along a radial segment of $B$.
In particular the ball $B(z_1(0),R)$ is disjoint from $\gamma_0([0,s_1])$.

The resting place $c$ for small $\tau$ will be computed by patching a finite number of local solutions $f$ of $f''/f' = \zeta_\tau$ as follows.
Choose $\eta>0$.
For $N\geq 2$ big enough the intervals $[s_1 k/N,s_1 (k+1)/N]$, $0\leq k\leq N-1$, are mapped by $\gamma_0$ to subsets of the balls $B_k=B(u_k,\eta)$ where for convenience we denote
\[u_k = \gamma_0\left(s_1\frac{k}{N}\right).\]
For $\eta$ small enough these balls are contained in $U$.
Note that $[0,s_1]\subset[0,1)$.
Choose some $\eps>0$ so that $\eps < \min(\eta,\eps<R/2,\eps_0)$, where $\eps_0$ is the $\eps$ given by \Cref{lem:GammaTop3} in \Cref{app:GammaTop}, applied to $\gamma_0$.
For $\tau$ small enough, $X_\tau$ belongs to the neighbourhood $V_{\gamma,\eps}$ of $X_0$:  $z_0(\tau)\in B(z_0(0),\eps)$, $z_1(\tau)\in B(z_1(0),\eps)$, and by~\Cref{lem:GammaTop3}, $\gamma_\tau$ is homotopic to $\gamma=\delta\cdot(\psi\circ\gamma_0)$
for \emph{any} path $\delta$ from $\tilde z_0(\tau)$ to $z_0(0)$ within $B(z_0(0),\eps)$, and \emph{any} homeomorphism $\psi$ of $U$ that maps $z_1(0)$ to $z_1(\tau)$ and which is the identity outside $B(z_1(0),\eps)$.
We choose $\delta = \delta_\tau$ to be a straight segment and $\psi_\tau$ to send the rays $[z_1(0),u]$ of the disk $\ov B'=\ov B(z_1(0),\eps)$ to the segments $[z_1(\tau),u]$ and let
\[\gamma_\tau:=\delta_\tau\cdot(\psi_\tau\circ\gamma_0)
\]
In particular $\psi_\tau\circ\gamma(s)$ follows the straight segment $[\gamma_0(s_1),z_1(\tau)]$ for $s\in[s_1,1]$.
For a fixed $\tau$, for all $0\leq k\leq N-1$ let inductively $\phi_{k,\tau}$ denote the solution of $\phi''_{k,\tau} / \phi'_{k,\tau} =\zeta_\tau$ on $B_k$ such that
\begin{itemize}
\item the germ of $\phi_{0,\tau}$ at $z_0(\tau)$ is $\phi_\tau$,
\item for $k\geq 1$, $\phi_{k,\tau}$ and $\phi_{k-1,\tau}$ coincide in a neighbourhood of $u_k$.
\end{itemize}
The maps $\phi_{k,\tau}$, as solutions $f$ of $f''/f' = \zeta_\tau$ with initial conditions $f(u)=v$ and $f'(u)=w$, depend holomorphically on $\tau$, $u$, $v$ and $w$.
It follows by induction that the points $\phi_{k,\tau}(u_k)$, $\phi_{k,\tau}(u_{k+1})$ and the derivatives $\phi_{k,\tau}'(u_k)$, $\phi_{k,\tau}'(u_{k+1})$, for $0\leq k\leq N-1$, are holomorphic functions of the data.

For the last patch recall that $\psi_\tau\circ\gamma$ is a straight segment from $\gamma_0(s_1)$ to $z_1(\tau)$. We use a branch of $z\mapsto w^{\res_\tau+1}$ on a straight sector $S$ in $z$ coordinate with apex $z_1(\tau)$ and containing the segment $[u_N,z_1(\tau))$, post-composed with an appropriate $\C$-affine map $\sigma_\tau$ so as to match with $\phi_{N-1,\tau}$ in a neighbourhood of $u_N$.
The map $\sigma_\tau$ depends holomorphically on the data and $c=\sigma_\tau(0)$.

\subsection{Following saddle connections: proof of \Cref{prop:fseg,prop:fseg2}}\label{app:pf:prop:fseg}

\begin{figure}
\begin{tikzpicture}
\node at (0,-2) {\includegraphics[scale=0.66]{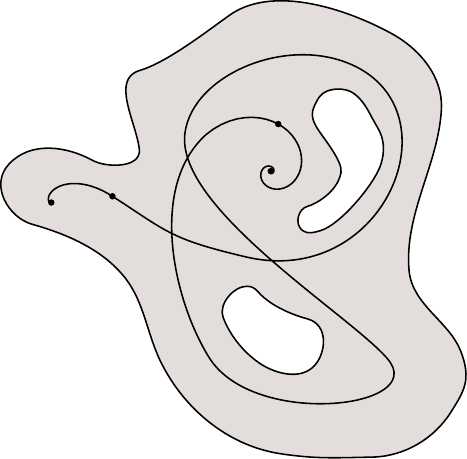}};
\draw[-] (-1.5,-3.5) node[anchor = north east]{$U$} -- (-1,-3);
\draw[-] (-2.2,-2.5) node[anchor = north]{$z^*_1$} -- (-2.05,-1.8);
\draw[-] (-0.8,0.3) node[anchor = south east]{$z^*_0$} -- (0.35,-1.2);
\draw[-] (0,-5) node[anchor = north]{saddle connection} -- (1,-4);
\node[text width=5cm] at (0,-6.2) {Sketch in Riemann coordinates for $\zeta^*$.};

\node at (6,0) {\includegraphics[scale=0.5]{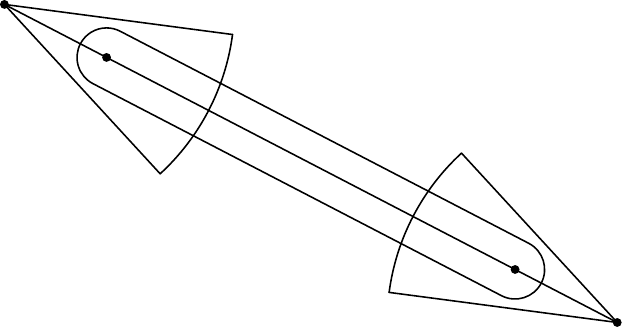}};
\node at (3.05,1.35) {$w_1^*$};
\node at (8.65,-1.65) {$w_0^*$};
\draw[-] (4.4,-0.6) node[anchor = east]{$W$} -- (5.8,-0.2);
\draw[-] (6.4,-0.9) node[anchor = east]{$L(t_0)$} -- (7.6,-0.9);
\draw[-] (5.6,0.9) node[anchor = west]{$L(t_1)$} -- (4.4,0.9);
\node[text width=5.5cm, anchor=north] at (6,-1.8) {Picture in developped similarity coordinates, before perturbation. $W$ and the two sectors are contained in the domain of $\psi^*$.};

\node at (6,-5.3) {\includegraphics[scale=0.5]{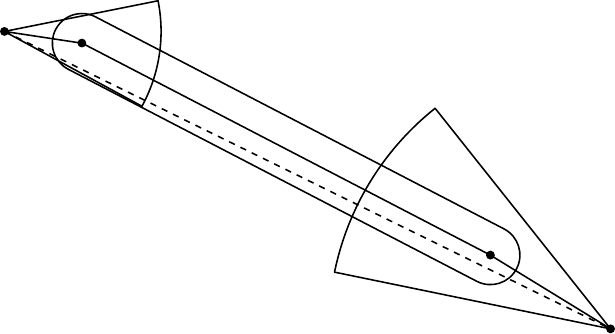}};
\node[text width=7.5cm] at (5,-7.9) {If the perturbation is small enough, there is still a segment between the sector apexes in the domain of $\psi$. The sectors have scaled and rotated but the central bar has not changed.};
\end{tikzpicture}
\caption{Objects in the proof of \Cref{prop:fseg}. The boundaries of a sector are not glued together (they are strict subsets of the sectors defining a neighbourhood of the vertices).}
\label{fig:seg}
\end{figure}

We will say that $\zeta$ depends analytically on $\tau$ whenever $h$, $\res_k$, $z_k$ do.

Let us call $\eps$-close (to $\zeta^*$) the maps $\zeta$ as in the statement of \Cref{prop:fseg}.
For $\eps$ small enough the two singularities remain of residue of real part $>-1$, i.e.\ bounded.

Since geodesics go straight with constant speed in similarity charts and since the extent of times on which $f$ is defined is finite, by pasting inverses of local charts, and of those of \Cref{lem:lds} at the extremities, one can prove that there exists an open segment $(w_0^*,w_1^*)\subset \C$ (one can choose $w_0^*=0$ and $w_1^*=1$), a neighbourhood $V$ of it, and a map $\psi^* : V
\to U -\{z^*_0,z^*_1\}$ such that
\begin{itemize}
\item for all $t\in(0,1)$, $\psi^*((1-t)w_0^* + t w_1^*) = f(t)$,
\item $\psi^*$ is analytic, locally invertible with local inverses being similarity charts for $\zeta^*$,
\item $V$ contains a sector based on each end of $[w_0^*, w_1^*]$ and having this segment as a central axis (see also \Cref{lem:ucs}).
\end{itemize}
Since $f$ is not necessarily injective, $\psi^*$ is not either.

Let
\[L(t) = (1-t)w_0^*+t w_1^*.\]

By a compactness argument, for any compactly contained open subset $W$ of $V$ containing the middle point $w_{1/2}^* = L(1/2) = (w_0^*+w_1^*)/2$, there is $\eps>0$ such that a map $\psi$ depending on $\zeta$ persists on $W$ for all $\eps$-close $\zeta$ such that,
\begin{itemize}
\item $\psi(w_{1/2}^*)=\psi^*(w_{1/2}^*)$, $\psi'(w_{1/2}^*)=(\psi^*)'(w_{1/2}^*)$,
\item $\psi$ is analytic, locally invertible with local inverses being similarity charts for $\zeta$.
\end{itemize}
Moreover, $\psi$ depends analytically on $\tau$ if $\zeta$ does. This comes from the expression of local similarity charts in terms of $\zeta$ (or as a variant: from the fact that $\psi$ is a solution of the complex O.D.E.\ $-\psi''/\psi'^2 = \zeta \circ \psi$.)

Near singularities, according to \Cref{lem:ldshol}, the formula for the holomorphic change of variable $z\mapsto w$ in \Cref{lem:lds} depends holomorphically on $\tau$ if $\zeta$ does.
Recall that branches of $w^{1+\res}$ are similarity charts.
We can thus choose for $\eps$ small enough some times $t_0$ close to $0$ and $t_1$ close to $1$, independent of $\zeta$, such that there are similarity chart inverses $\psi_0$, $\psi_1$ defined on the sector
\[\setof{u=r e^{i\theta}}{0<r<2,\ -\eta<\theta<\eta}\]
for some $\eta$ independent of $\zeta$, and such that, as $u\to 0$, $\psi_0(u) \tend z_0$, $\psi_1(u) \tend z_1$, $\psi_0(1) = \psi(L(t_0))$ and $\psi_1(1)=\psi(L(t_1))$.

Choose $W$ as above to contain the sub-segment $[L(t_0),L(t_1)]$ of $[w^*_0,w^*_1]$.
We can patch together $\psi$, $\psi_0$ and $\psi_1$ by 
replacing $\psi_i$ with $\tilde\psi_i: w\mapsto \psi_i(s_i(w))$ where $s_i$ is is the $\C$-affine map such that $\tilde \psi_i(L(t_i)) = \psi(L(t_i))$ (i.e.\ $s_i(L(t_i))=1$), $\tilde \psi_i'(L(t_i)) = \psi'(L(t_i))$, and possibly reducing $\eta$ and $W$.
We obtain an extended map $\tilde \psi$ defined on the union of a fixed neighbourhood $W$ of $[L(t_0),L(t_1)]$ and of two sectors centred on points $w_0$ and $w_1$ that depend on $\zeta$, of symmetry axis $[w_i,L(t_i)]$, of opening angle $\eta$ independent of $\zeta$, and of radius $2|w_i-L(t_i)|$.
See \Cref{fig:seg}.

Everything in this construction depends holomorphically on $\tau$ if $\zeta$ does. 
For $\eps$ small enough, the segment $[w_0,w_1]$ (which depends on $\zeta$) is contained in the domain of $\tilde\psi$ (which also depends on $\zeta$) and we can take for $t\in(0,1)$: $f(t) = \tilde\psi((1-t)w_0+tw_1)$ and $f(0)=z_0$, $f(1)=z_1$, as a map satisfying the conclusion of \Cref{prop:fseg}.
Continuity of $\tau\mapsto f$ for the uniform norm follows from the continuity of $\tau\mapsto \tilde\psi$, which itself is due to $\tilde\psi$ being the solution of an O.D.E.\ and near $w_0$ and $w_1$ to \Cref{lem:ldshol} and the construction.

\begin{figure}
\begin{tikzpicture}
\node at (0,0) {\includegraphics[scale=0.5]{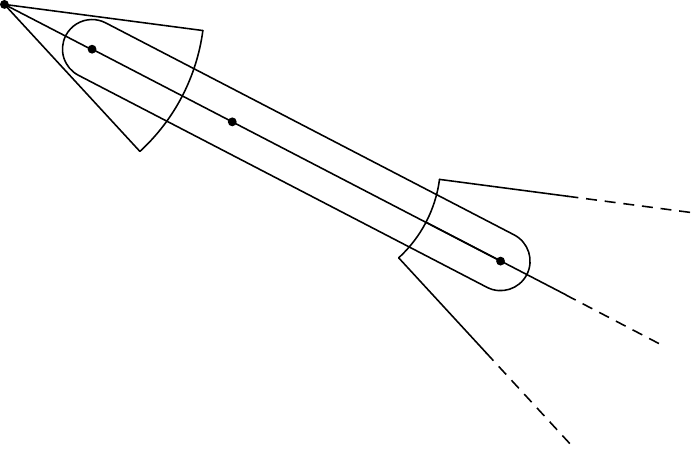}};
\node at (-2.8,2.1) {$w^*_0$};
\draw[-] (-0.5,1.5) node[anchor = south west]{$w^*_1$} -- (-0.93,0.95);
\node at (5.5,0) {\includegraphics[scale=0.5]{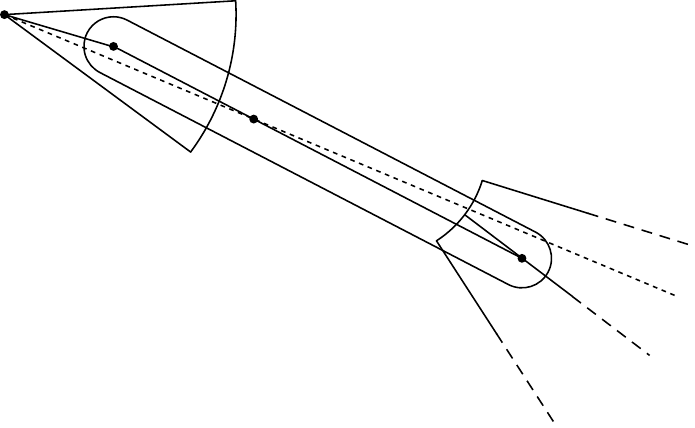}};
\draw[-] (-2.17,0.3) node[anchor = north]{$L(t_0)$} -- (-2.17,1.4);
\draw[-] (1.285,0.7) node[anchor = south]{$L(t_1)$} -- (1.285,-0.2);
\end{tikzpicture}
\caption{Adaptation of \Cref{fig:seg} for the proof of \Cref{prop:fseg2}}
\label{fig:seg2}
\end{figure}

Last assume the initial connection $f^*$ is injective. Let us prove that for $\eps$ small enough, $f$ is injective.
Assume by way of contradiction that there is a sequence of such $\zeta_n$ with $\eps_n\tend 0$ such that, denoting the connections $\gamma_n\tend\gamma$, we have $\gamma_n(t_n)=\gamma_n(t'_n)$ with $0\leq t_n< t'_n\leq 1$.
Necessarily $0<t_n$ and $t'_n<1$.
We could assume by extraction that $t_n\tend t$ and $t'_n\tend t'$.
By continuity $\gamma^1(t)=\gamma^1(t')$, so $t=t'$.
If $t$ is not an end of $[0,1]$, since $\gamma_\tau(t)$ depends smoothly on the pair ($t,\tau$), this would imply that $(\gamma^1)'(t)=0$ contradicting that $\gamma^*$ is a geodesic.
If $t=0$ or $1$, this contradict the fact that the geodesic is, in the affine coordinate, a straight segment to the apex of the sector near $w_0$ and $w_1$, and the uniform local injectivity of $\psi$ near $w_0$ and $w_1$, which follows from holomorphy in \Cref{lem:ldshol}.

The proof of \Cref{prop:fseg2} is very similar, with one of the sectors replaced by an unbounded sector, using \Cref{lem:ucs,lem:geodeschol} for this part. 
See \Cref{fig:seg2}.

%
% ----------------------------------------------------------
%

\subsection{Proof of holomorphic dependence (\Cref{prop:hsol})}\label{app:pf:prop:hsol}

This section introduces a lot of notations.
However, most, if not all, will only be used here.

\subsubsection{Adding the image of the marked points}

For each paired unbounded edges there was a marked point on the interior of each.
We extend the sequence $z_k$ by appending the image by $\prB:\cal S\to\hat\C$ of the marked points, to get an (injective) element of $\hat\C^m$ for some $m\geq p$
and we will prove that this extended sequence depends holomorphically on the polygons.

\subsubsection{Summary of the proof}

Before entering in a fully detailed argument, let us give an overview of the proof.

The construction of \Cref{prop:hsol} starts from a collection of polygons and gluings.
We will define a space $\Sconf$ of collections of strictly convex polygons up to $\C$-affine transform, which represents the possible geometric shapes that can take each of the polygons when they vary while remaining strictly convex.
(For unbounded polygons, marked points are taken into account.)
The space $\Sconf$ is a complex manifold. Given a point in this space, the construction can be carried out and yields in particular the points $z_k\in\hat\C$ and their associated residues $\res_k$.
We already know by \Cref{lem:rh} that the $\res_k$ depend holomorphically on the polygons and marked points.

In fact we will see that one can get from the construction not only the position of the points $z_k$ (up to an automorphism of $\hat\C$) but an element of the Teichmüller space $\cal T$ of the Riemann sphere $\hat \C$ with $m$ marked points.
So we consider the map
\[\Glu : \Sconf \to \cal{TR}\]
that sends a collection to the Teichmüller space element thus constructed, together with the collection of residues.
The name of the map reflects the fact that it is defined by a process that starts by \emph{gluing} polygons.
A priori we do not know anything on $\Glu$, and even its continuity should not be considered as a trivial statement.

To prove analytic dependence of the $z_k$ we will prove the stronger statement that $\Glu$ is analytic. For this, we will introduce sort of a reciprocal map
\[\Per : \cal{TR} \to \Conf\]
where $\Conf$ is a space containing $\Sconf$.
The name $\Per$ is an abbreviation of \emph{periods}.
It is defined as follows.
An element of $\cal{TR}$ consists in an element of $\cal T$, in particular points $z_k$, and the specification of complex numbers $\res_k$ associated to each $z_k$ with some compatibility conditions.
One can associate the unique Christoffel symbol having these residues at the $z_k$, and we thus get a similarity surface structure on $\hat\C-\{z_k\}_{k=1}^m$.
To the element of $\cal T$ corresponds a decomposition of the sphere into topological polygons, unique up to isotopy of the sphere fixing the marked points.
Local straightening maps of the Christoffel symbol can be extended to these topological polygons into holomorphic maps, which are not necessarily anymore injective.
However, one can look at the image of those topological polygons by these extended straightening maps, and the position of the vertices in these images, up to a $\C$-affine map, is independent of the isotopy.
This collection of positions, up to $\C$-affine maps, is what $\Per$ associates to the element of $\cal{TR}$ we started from.

The name \emph{periods} of the map has been chosen because the relative position can be computed by integrating an O.D.E. along paths between the $z_k$ (those that are bounded vertices or marked point) and the result only depends on the homotopy class of the path; this is akin to computing translation vectors between branch points of a translation surface by integration of a one-form, and those vectors are sometimes called periods by extension of the notion of period of a one-form along a closed loop.

One immediately gets that 
\[\Per\circ\Glu = \on{Id}_{\Sconf}.\]

The space $\Conf$ is not an analytic manifold but we identify a subset $\Conf^*$ that carries an analytic manifold structure and such that
\[\Sconf \subset \Conf^*\subset \Conf,\]
$\Sconf$ is an open subset of $\Conf^*$ and their analytic structures match.
We then check in \Cref{lem:thetaan} that $\Per$ is continuous on $\cal{TR}$, and analytic on $\Per^{-1}(\Conf^*)$.

Let \[\Eff:=\Glu(\Sconf)\subset \cal{TR}.\]
The name $\Eff$ is an abbreviation for \emph{effective}.
From $\Per\circ\Glu = \on{Id}_{\Sconf}$ it immediately follows that \[\Glu\circ \Per|_{\Eff} = \on{Id}_{\Eff}.\]
Intuition tells us that $\Eff$ must be an open subset of $\cal TR$ but this is a non-trivial fact. 
Recall that at this point of the proof, we do not even know that $\Glu$ is continuous, so  $\Eff$ could just be any kind of uncountable subset of $\cal TR$. On the other hand, it is immediate that the subset $\Per^{-1}(\Sconf)$ of $\cal{TR}$ is open, since $\Per$ is holomorphic hence continuous.

A key point (see below), which will require some effort, will be to prove that each point in $\Eff$ has a neighbourhood $W$ in $\Per^{-1}(\Sconf)$, on which the equality 
\[\Glu\circ \Per|_{W} = \on{Id}_{W}\]
holds.
This implies that $W\subset\Eff$ and proves that $\Eff$ is an open subset of $\cal{TR}$. 

To conclude, we use a theorem that if one has an analytic map $\phi:U\to V$ between open subsets of complex manifolds of dimensions $a$ and $b$, and if $\phi$ is a bijection, then $a=b$ and $\phi^{-1}$ is analytic too.
We apply this to $\phi=\Per|_{\Eff}: U=\Eff \to V=\Sconf$.

The aforementioned key point follows from \Cref{lem:ma} and the paragraph after it. An effective element $x$ of $\cal {TR}$ corresponds to a collection of strictly convex polygons $P_j$.
Nearby elements $x'$ of $\cal {TR}$ will map by $\Per$, which is continuous, to configurations of points for which there exist strictly convex polygons $P'_j$ having these points as vertices and marked points.
The difficulty is to prove that gluing these new polygons indeed gives back $x'$.
For this we follow, as a function of $x'$ the saddle connections of the Christoffel symbol associated to $x'$.
This defines a cellular decomposition of the sphere and we check that the cells actually map to the $P'_j$ by the straightening maps of $x'$, provided it is close to $x$.
Since the cellular decomposition for $x'$ is obtained by deformation of the one for $x$, we can check that gluing the $P'_j$ also gives the same isotopy class as the Teichmüller part of $x'$.

\subsubsection{Polygon space}

Consider strictly convex bounded polygons $P\subset\C$: they are assumed with non-empty interior and have at least three vertices. Let us orient the boundary of $P$ in the counter clockwise way, and index its vertices accordingly as follows: $w_1,\ldots, w_n$.
The set of strictly convex polygons with indexed vertices thus forms an open subset of $\C^n$.

For the unbounded polygons, recall that we assumed: that they are strictly convex; that there are exactly two unbounded edges; that they make a positive angle (this is a standing assumption made in \Cref{sub:hd:p}, in particular in \Cref{prop:hsol}); that we consider $\infty$ as a vertex; that we include a marked point in the interior of each of the two unbounded edges.
We still orient the boundary in the anti-clockwise way and index the vertices and marked points accordingly \emph{but omitting $\infty$}: $w_1,\ldots, w_n$. We still get an open subset of $\C^n$, where now $n$ is the number of edges plus one, so $n\geq 3$.

\subsubsection{Quotienting the polygon space}

Changing one polygon by mapping it under an affine map does not change the quotient surfaces $\cal S$, $\cal S'$ (there are trivial canonical isomorphisms), nor the sequence $(z_k)$, hence we can work in the space of (possibly marked) polygons $P_j$ up to the action of affine maps: this is still a complex manifold (it amounts to fixing the affix of two particular vertices/marked points of $P_j$) and the quotient map by this action is analytic and open.

In fact we will use later larger spaces. Let $\Aff\C$ denote the group of $\C$-affine maps $w\mapsto aw+b$ of $\C$ and let it act on $\C^n$ component-wise: for $(w_k)\in\C^n$, $s\cdot(w_k) := (s(w_k))$.
Denote
\[A\C^n =\C^n/\Aff\C.\]
The quotient map $\C^n\to A\C^n$ is open for the quotient topology.\footnote{I.e.\ the saturate of an open subset of $\C^n$ is open, which is easily checked here since affine maps are open.}
Let $A^*\C^n \subset A\C^n$ be the quotient of $\C^n-\Delta$ where $\Delta = \setof{(w,\ldots,w)}{w\in\C}$ i.e.\ $\Delta$ represents the case where all points have merged into a single point.
The space $A^*\C^n$ holds a complex manifold structure such that the quotient map $\C^n-\Delta \to A^*\C^n$ is analytic. Locally, charts of $A^*\C^n$ are given by fixing the affix of two distinct points.
The space $A\C^n$ has only one point more than $A^*\C^n$ and we do not seek to extend the complex structure to this point. For the quotient topology $A\C^n$ is not Hausdorff separated: indeed the special point has only one neighbourhood, which is the whole space $A\C^n$ (this is sometimes called a \emph{focal point}).

Recall that we indexed the collection of polygons $P_j$ by a finite set $J$.
Let $n_j$ be the number of vertices on $P_j$ (including marked points, excluding infinity).
Let
\bEA
&\Conf=\prod_{j\in J} A\C^{n_j}&
\\
&\Conf^*=\prod_{j\in J} A^*\C^{n_j}&
\eEA
For each $j$, the number $n_j$ is fixed and we also fix the set of $j$ for which the $w_k$ include marked points, whose indexes we fix and must be two integers that are consecutive modulo $n_j$.
Let $\Sconf\subset\Conf$ denote the elements such that for each $j$ without marked point, the $(w_k)$ are given by the anti-clockwise indexed vertices of a strictly convex bounded $P_j$, and in the case with marked point, the $(w_k)$ are given as explained in the paragraph above called ``Polygon space''.
Then the subset $\Sconf$ of $\Conf$ is open and 
\[\Sconf \subset \Conf^*\subset \Conf.\]
As a product of analytic manifolds, $\Conf^*$ is an analytic manifold and this makes $\Sconf$ an analytic manifold too.

\subsubsection{Moduli space and Teichmüller space of the sphere with marked points}

We fix a particular polygon configuration in $\Sconf$, which we will call $[\cal P^0]$.
Denote $\cal S^0$, $\cal V^0 = \setof{v_k^0}{1\leq k \leq m}$, and $(\cal S')^0 = \cal S^0-\cal V^0$ the respective Riemann surface, vertices, and similarity surface associated to $\cal P^0$.

A quick reminder on the moduli and Teichmüller spaces of the sphere with marked points can be found in \Cref{app:teich}. However, for the reader's convenience, we copy here some of the definitions given there.

The moduli space is the subset $\cal M \subset \hat\C^m$ of $m$-uplets $(z_k)\in\hat\C^n$ for which the $z_k$ are distinct and $z_1=0$, $z_2=1$, $z_3=\infty$. It is a complex manifold of dimension $m-3$.
To a configuration $\cal P\in \Sconf$ we associate the image in $\hat \C$ of its vertices and marked points: $Z=(z_k)=(\prB(v_k))\in\cal M$ for the associated conformal isomorphism $\prB:\cal S\to\hat \C$ sending $v_1$ to $0$, $v_2$ to $1$ and $v_3$ to $\infty$ (as defined before \Cref{prop:hsol}), and we let $\cal Z = \setof{z_k}{1\leq k \leq m}$.
To $\cal P$ we also associate the Christoffel symbol $\zeta_{\cal P}$ on $\hat\C - \cal Z$ that our construction yields.
Naturally, we denote $\prB^0$, $Z^0$ and $\cal Z^0$ the objects associated to $\cal P^0$. 

Let $\cal T$ be the Teichmüller space associated to $\cal S^0, \cal V^0$ defined (without the use of quasiconformal maps) as follows: 
\[\cal T = \cal F/ \cal H^0\]
where $\cal F$ is the set of orientation preserving homeomorphisms $f:\cal S^0 \to \hat\C$ sending respectively $v^0_1$, $v^0_2$ and $v^0_3$ to $0$, $1$ and $\infty$ and $\cal H^0$ is the set of orientation preserving self homeomorphisms of $\cal S^0$ that are isotopic to the identity rel.\ $\cal V^0$, i.e.\ by an isotopy fixing each vertex $v_k^0$.

If $f\in\cal F$, its equivalence class modulo right composition with $\cal H^0$ will be denoted $[f]\in \cal T$. Hence, $[f_1]=[f_2]$ iff $\exists \phi\in\cal H^0$ such that $f_2 = f_1\circ\phi$. But note that this is equivalent to: $f_2$ is isotopic to $f_1$ rel.\ $\cal Z = \{z_k=f_1(v_k^0)\}$.

Let us endow $\hat \C$ with a spherical metric\footnote{A metric for which the circles are geometric circles via a stereographic projection.} $d$ and the set $\cal F$ with the metric
\[d(f_1,f_2) = \sup\setof{d(f_1(x),f_2(x))}{x\in\cal S^0}.\]
The space $\cal T$ is endowed with the quotient topology.
It is separated (distinct points have some disjoint neighbourhoods).
The projection $\Pi: \cal F\to \cal T$ is open and has continuous local sections (it is a fibre bundle, see \Cref{app:teich}). 
Since $\Pi$ is open, a subset of $\cal T$ is a neighbourhood of $[f]\in\cal T$ iff for some $\eps>0$ it contains the set of $[f']$ for which $d(f,f')<\eps$.
We will use the following topological lemma.
Endow $\cal S^0$ with a distance $d$ inducing its topology and $\cal H^0$ with the distance 
\[d(\phi_1,\phi_2) = \sup\setof{d(\phi_1(x),\phi_2(x))}{x\in\cal S^0}.\]
Then we have the following lemma (see \Cref{app:teich}, \Cref{lem:topo0}):
\begin{lemma}\label{lem:topo}
There exists $\eps>0$ such that for every orientation preserving self-homeomorphism $\phi$ of $\cal S^0$, if $\phi$ fixes every point in $\cal V^0$ and satisfies $d(\phi,\on{id}_{\cal S^0})<\eps$ then $\phi\in\cal H^0$ (i.e.\ $\phi$ is isotopic to the identity rel.\ $\cal V^0$).
\end{lemma}

\subsubsection{Stating the objective in terms of analyticity of a function $\Glu$}

\begin{figure}
\begin{tikzpicture}
\node at (0,-0.5) {\includegraphics[scale=0.45]{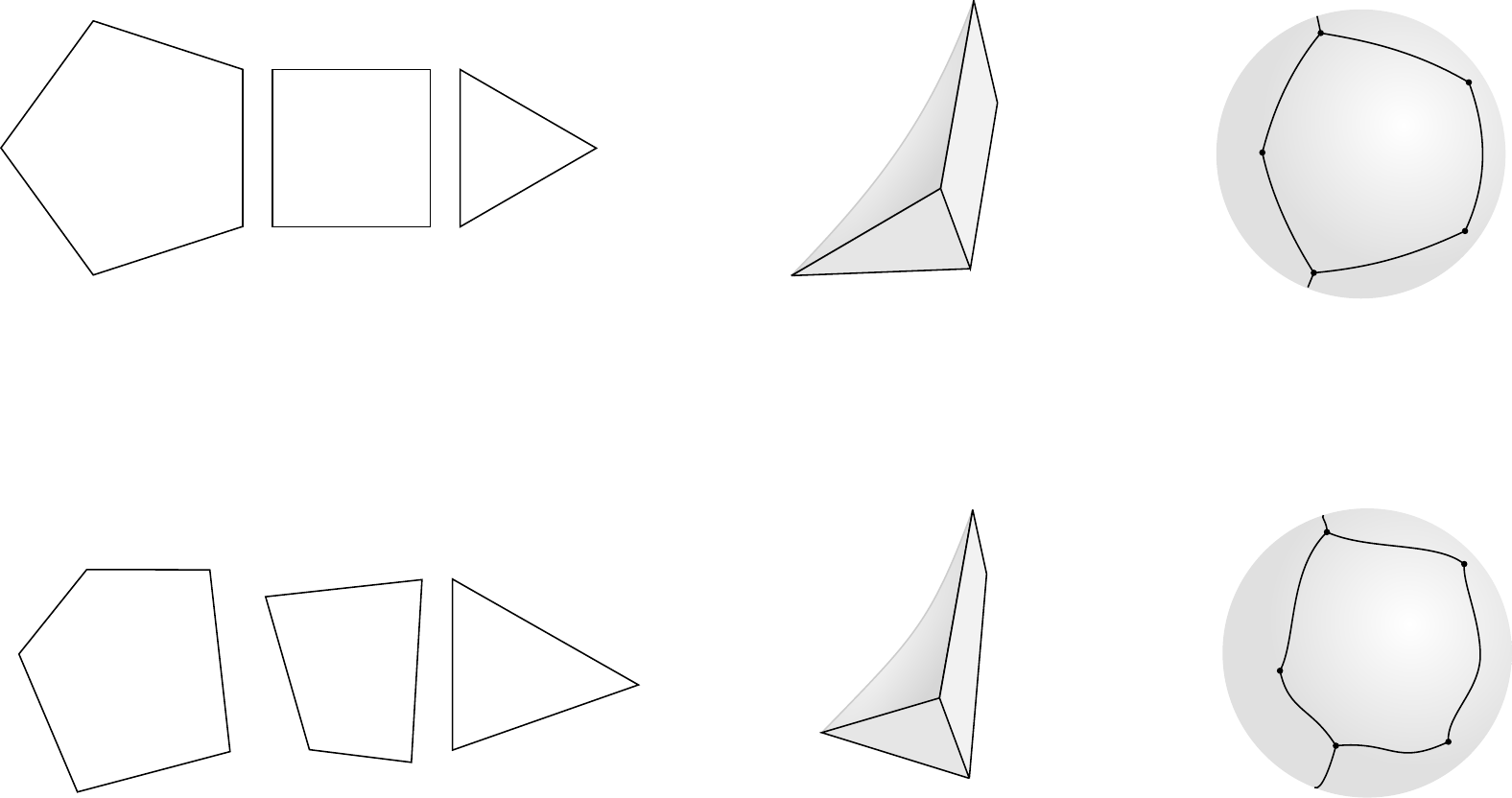}};

\node at (-4.9,1.5) {$P_1^0$};
\node at (-3.2,1.5) {$P_2^0$};
\node at (-1.95,1.5) {$P_3^0$};
\node at (-1.3,2.3) {$\cal P^0$};

\node at (-4.9,-2.7) {$P_1$};
\node at (-3.2,-2.7) {$P_2$};
\node at (-1.95,-2.7) {$P_3$};
\node at (-1.3,-1.9) {$\cal P$};

\node at (1,2) {$\cal S^0$};
\node at (1,-2) {$\cal S$};

\node at (4.6,0.8) {$z_k^0$};
\node at (4.8,-3) {$z_k$};

\draw[->] (-0.6,1.5) -- node[above] {$\pi^0$} (0.4,1.5);
\draw[->] (-0.6,-2.5) -- node[above] {$\pi$} (0.4,-2.5);
\draw[->] (2.2,1.5) -- node[above] {$\prB^0$} (3.2,1.5);
\draw[->] (2.2,-2.5) -- node[above] {$\prB$} (3.2,-2.5);

\draw[->] (-3.3,0) -- node[right] {$\phi$} (-3.3,-1);
\draw[->] (1.2,0) -- node[right] {$\psi$} (1.2,-1);
\draw[->] (2.5,0) -- node[above right] {$f$} (3.5,-1);
\end{tikzpicture}
\caption{This commuting diagram sums up the objects involved in the construction of $\Glu$. For a fixed $\cal P$, the descending maps are not unique but $\pi$, $F$ and $[f]$ are.}
\label{fig:glu}
\end{figure}

The different objects that will be introduced and their relation are summed up in \Cref{fig:glu}.

For a polygon configuration in $\Sconf$, let us explain why our gluing construction not only gives us an element of $\cal M$ but in fact an element of $\cal T$. 
For this, given a polygon configuration $[\cal P] \in \Sconf$ we build below a homeomorphism $\psi: \cal S^0\to\cal S$. It is not unique but its class modulo $\cal H^0$ will be unique.

Consider a collection $\phi = (\phi_j)$ of homeomorphisms $\phi_j: \partial P^0_j \to \partial P_j$ that send vertices and marked points of $P^0_j$ to their corresponding points in $P_j$ (in particular, they are orientation preserving) and such that the $\phi_j$ are compatible with the identification of paired edges:
if $(e\,{:}\,j)$ and $(e'\,{:}\,j')$ are paired in $\cal P^0$, and $s$ is the $\C$-affine map matching them, calling $(\tilde e\,{:}\,j)$, $(\tilde e'\,{:}\,j')$ and $\tilde s$ the corresponding objects in $\cal P$, we ask that $\phi_{j'} \circ s = \tilde s \circ \phi_j$.
Such a collection exists, but is not unique.
Then $\phi_j$ can be extended to an orientation preserving homeomorphism $\phi_j: P^0_j \to P_j$.
Let us join all these maps into a map $\phi : \bigsqcup_j P^0_j \to \bigsqcup_j P_j$.
This join descends to an orientation preserving homeomorphism $\psi:\cal S^0 \to \cal S$ that depends only on $\phi$ such that $\prA\circ\phi = \psi \circ \prA^0$, where
\[\pi:\bigsqcup_j P_j \to \cal S\text{ and }\pi^0 : \bigsqcup_j P^0_j \to \cal S^0\]
are the projection to the respective quotients.
This is summed up in the following commutative diagram:
\[
\begin{tikzcd}
\bigsqcup_j P^0_j  \arrow{r}{\phi} \arrow[swap]{d}{\pi^0} & \bigsqcup_j P_j  \arrow{d}{\pi} \\%
\cal S^0 \arrow{r}{\psi}& \cal S
\end{tikzcd}
\]
Different choices for $\phi$ yield different maps $\psi$ but they are all equivalent modulo $\cal H^0$.\footnote{This is directly seen to be equivalent to: the maps $\psi^{-1}$ are all isotopic rel. the vertices and marked points of $\cal S^0$. Actually by taking inverse this is also equivalent to stating that the maps $\psi$ are isotopic rel the vertices and marked points of $\cal S$.}
This follows from the following topological facts, via a homeomorphism from $P^0_j$ to $\ov\D$: the set of orientation preserving self-homeomorphisms of a segment is connected (it is convex!), an isotopy of the circle can be radially extended to an isotopy of the closed unit disk $\ov \D$ (immediate), the set of self-homeomorphisms of $\ov\D$ that are the identity on $\partial \D$ is connected (easy thanks to Alexander's trick \cite{Al}).

Recall that $\prB$ denotes a specific Riemann surface isomorphism from $\cal S$ to $\hat \C$.
Then the map
\begin{equation*}\label{eq:fFpsi}
f:=\prB \circ \psi
\end{equation*} 
is an element of $\cal F$.
By what we proved in the previous paragraph, its class modulo pre-composition by elements of $\cal H^0$ is independent of the choice of $\phi$ and is the element of $\cal T$ we will associate to $\cal P$. This independence is used in the proof of \Cref{lem:ma}.

We partition the index set $\{1,\ldots m\}$ into three subsets as follows: let $F\subset\{1,\ldots m\}$ denote the set of those indices $k$ for which the vertex $v^0_k$ of the base surface $\cal S^0$ is bounded,
$I\subset\{1,\ldots m\}$ those for which it is unbounded and $M\subset\{1,\ldots m\}$ be the marked points.
Let
\begin{equation}\label{eq:defR}
\cal R =\bigsetof {(r_k)\in\C^m}{\begin{array}{c}
\sum r_k = -2
\\ \forall k\in F,\ \Re r_k>-1
\\ \forall k\in I,\ \Re r_k<-1
\\ \forall k\in M,\ r_k=0
\end{array}}
\end{equation}
and
\[\cal {TR} = \cal T\times \cal R.\]
To a polygon configuration $[\cal P]\in \Sconf$ let us associate the element $[f]$ of $\cal T$ that the construction above yields. Let us also associate the set of residues $\res_k$ of the Christoffel symbol $\zeta_P$ at the $m$ vertices/marked point and temporarily denote
\[\res = (\res_k) \in \C^m.\]
We saw earlier that $\sum \res_k=-2$, and the types of the singularities (bounded or unbounded) of $\zeta$ are the same as for $\cal S^0$ according to \Cref{sub:presc} and for marked points there is no singularity, so $\res \in \cal R$.
We denote
\[
\Glu: \left.
\begin{array}{rcl}
\Sconf & \to & \cal {TR}
\\{} % leave {} here! this is to prevent interaction of \\ with [
[\cal P] & \mapsto & ([f],\res)
\end{array}
\right.
\]
where $f$ and $\res$ are defined above.
We call $\Glu$ the \emph{gluing} map.
We want to prove that $\Glu$ is analytic.
But we do not even know yet that it is continuous.

\subsubsection{An analytic inverse candidate $\Per$}

Consider any element $([f],\res)$ of $\cal {TR}$, $f:\cal S^0\to\hat\C$ and $(\res_k)\in\cal R$.
Let $z_k=f(v_k^0)$ be the affixes that any map in $[f]$ gives to the vertices and marked points and let $\cal Z=\setof{z_k}{1\leq k \leq m}$ denote their set.
Consider the unique Christoffel symbol $\zeta$ that is holomorphic on $\hat\C - \cal Z$ and whose singularities are at most simple poles and have residues $\res_k$ at $z_k$.
By definition of $\cal R$ if the type of $v_k^0$ is finite (F), infinite (I), marked (M) then the corresponding singularity $z_k$ of $\zeta$ is respectively bounded ($\Re \res_k >-1$), unbounded ($\Re \res_k <-1$), erasable ($\Re \res_k =0$).
Consider a polygon $P_j^0$ of $\cal P^0$, choose any point $w_*$ in the interior of $P^0_j$, and let $z_*=f\circ\prA^0(w_*)$ (recall that $\pi^0 :\bigsqcup_j P^0_j \to \cal S^0$).
Choose any germ  $\phi_0$ of similarity chart for $\zeta$ at $z_*$.
For all bounded vertices and all marked points $w$ of $P^0_j$, the germ can be developed along any path within $f\circ\prA^0(P^0_j)$ going from $z_*$ to $z=f\circ \pi^0(w)$ without hitting the vertices (except at the end) and the associated resting place (see \Cref{sub:rest}) is independent of the choice of the path since $P^0_j$ is simply connected.
Let 
\[R\in\C^{n_j}\]
be the indexed collection of all these resting places.
Recall that $A\C^n$ denotes the quotient of $\C^n$ by the action of the affine group $\Aff\C$ acting component-wise.

\begin{lemma}\label{lem:Rindep}
The projection of $R$ in $A\C^{n_j}$ is independent of the representative $f$ in $[f]$, of the choice of $w_*$ and of the choice of $\phi_0$.
\end{lemma}

Before proving it, let us gather all those projections for the different $j\in J$: we have hence defined a map
\[\Per: \cal{TR} \to \Conf\]
that we call the \emph{periods} map.

To prove \Cref{lem:Rindep}, and for further uses, we extend the definition of $R$.
Let $([f],\res)\in\cal{TR}$ and $\zeta$ be the associated Christoffel symbol. Assume that $z_*\in\hat\C-\cal Z$ is any point, that $\phi_0$ is a germ of affine chart at $z_*$ for $\zeta$, that $\gamma_i : [0,1]\to\hat\C$, $i\in\{1,\ldots,p\}$, $q\in \N$ are paths such that $\gamma_i(0)=z_*$, $\gamma_i([0,1))\subset\hat\C-\cal Z$ and $\gamma_i(1)\in\cal Z$ is a bounded singularity of $\zeta$ (i.e.\ corresponds to a $v_k^0$ that is either a finite vertex or a marked point). Develop $\phi_0$ along each $\gamma_i$ and denote $c_i\in\C$ its resting place. Let $R=(c_i)_{i=1}^q$, which depends on all the choices above.

\begin{lemma}\label{lem:Rindep2}
Consider two choices of $z_*$, $\phi_0$, $\gamma_i$ and $z_*'$, $\phi'_0$ and $\gamma'_i$ as above and assume that there exists a path $\delta$ from $z_*'$ to $z_*$ such that each path $\gamma'_i$ is homotopic to the concatenation $\delta\cdot\gamma_i$ by a homotopy leaving fixed the starting and endpoints, and avoiding $\cal Z$ except at the endpoint.
Then the projections of the corresponding $R$ and $R'$ in $A\C^{q}$ are identical.
\end{lemma}
\begin{proof}
Let us develop $\phi'_0$ along $\delta$. We obtain at the end a germ $\tilde\phi_0$ at $z_*$. This germ takes the form $\tilde\phi_0=s\circ\phi_0$ for some $\C$-affine map $s$.
The development of $s\circ\phi_0$ along $\gamma_i$ is the post-composition by $s$ of the development of $\phi_0$ along the same path $\gamma_i$ and its resting place is $s(c_i)$.
By the homotopy assumption and properties of resting places (see \Cref{sub:rest}) the resting places $c_i'$ composing $R'$ are the same as the resting places obtained from $z'_*$, $\phi'_0$, $\delta\cdot\gamma_i$.
The development of $\phi'_0$ along $\delta\cdot\gamma_i$ is the concatenation of the development of $\phi'_0$ along $\delta$ and the development of $\tilde\phi_0$ along $\gamma_i$.
Hence the $c'_i$ are the same as the resting places obtained from
$z_*$, $s\circ\phi_0$, $\gamma_i$, i.e.\ $c'_i=s(c_i)$.
\end{proof}

\begin{proof}[Proof of \Cref{lem:Rindep}]
Assume that $[f]=[f']$, that $w_{*}$ and $w_{*}'$ are in the interior of $\cal P^0_j$ and that $\phi_0$ and $\phi_0'$ are germs of straightening coordinates of $\zeta$ at these points.
Let $z_*=f(w_*)$ and $z'_*=f'(w_*')$.
Consider paths $\eta_i$, $i\in\{1,\ldots,n_j\}$, in $\cal P^0_j$ from $w_*$ to the finite vertices and marked points and a path $\beta$ from $w'_*$ to $w_*$ within $\cal P^0_j$.
Let $\gamma_i=f\circ\eta_i$.
Let $\gamma'_i = f'\circ(\beta\cdot\eta_i)= (f'\circ\beta)\cdot(f'\circ\eta_i)$.

Since $[f]=[f']$ there is an isotopy $s\in[0,1] \mapsto f_s$ from $f'$ to $f$ fixing the images of all vertices and marked points.\footnote{Write $f' = f\circ h$ for some $h\in \cal H^0$, let $h_s$ be a path from $h$ to $\on{Id}$ in $\cal H^0$ and let $f_s = f \circ h_s$.}
In particular $f_0=f'$ and $f_1=f$.
Let $\alpha:\, t\in[0,1] \mapsto \alpha(t) = f_t(w_*)$ which goes from $\alpha(0)=f'(w_*)$ to $\alpha(1)=f(w_*)=z_*$.
Let $\delta$ be the concatenation $(f'\circ\beta)\cdot\alpha$.
Then the conditions of \Cref{lem:Rindep2} are satisfied, with the following homotopy from $\delta\cdot \gamma_i$ to $\gamma'_i$:
first $f'\circ \eta_i$ is homotopic to $\alpha\cdot\gamma_i$ via
$(s,t)\in[0,1]^2\mapsto \alpha(2t)$ if $t\leq s/2$ and $f_{s}\left(\eta_i\left(\frac{t-\frac{s}{2}}{1-\frac{s}{2}}\right)\right)$ if $t\geq s/2$.
Then $\delta\cdot \gamma_i = ((f'\circ\beta)\cdot \alpha)\cdot\gamma_i$
is homotopic to $(f'\circ\beta)\cdot (\alpha\cdot\gamma_i)$ hence to $(f'\circ\beta)\cdot (f'\circ\eta_i) = \gamma'_i$.

The conclusion of \Cref{lem:Rindep2} proves \Cref{lem:Rindep}.
\end{proof}

The proof of the next statement contains no essential difficulty but we will have to deal with many details related to the definitions.
Recall that we have put a complex manifold structure on $\Conf^*$ but not on $\Conf$.

\begin{lemma}\label{lem:thetaan}
The subset $\Per^{-1}(\Conf^*)\subset\cal{TR}$ is open and on it, the map $\Per$ is analytic.
\end{lemma}
\begin{proof}
We will work locally so we consider an element $(T',\res')\in\cal {TR}=\cal T\times\cal R$ such that $\Per(T',\res')\in\Conf^*$,  and we will allow ourselves to restrict a finite number of times to smaller and smaller neighbourhoods of this element.
Let us focus on one polygon $P^0_j$ and the collection of resting places $R\in\C^{n_j}$ used in the definition of $\Per$.

Fix any $f'\in\cal F$ such that $T' = [f']$. Fix a point $w_*'\in P_j^0$ and let
\[z_*'=f'\circ\pi^0(w_*')\in\hat\C.\]

We recall that the collection $R\in\C^{n_j}$ is associated to $([f],\res)\in\cal{TR}$ and to a choice of point $w_*$ in the interior of $P^0_j$ and of an initial germ $\phi_0$ of similarity chart at $z_*=f\circ\pi^0(w_*)$.

We want to keep $z_*$ fixed and equal to $z_*'$, and hence we want to define $w_* =(f\circ \pi_0)^{-1}(z_*')$.
This poses several problems.
For one thing, $f$ is only defined up to isotopy rel.\ the vertices and marked points, and thus $w_*$ depends on the choice of the representative $f$ of $[f]$, and for some representatives, the point $w_*$ is not even in $P^0_j$.
We could prove that the construction can be carried out for such $w_*$ and work from that, but by simplicity we prefer here to use that the map $\cal F\to\cal T:\ f\mapsto [f]$ has local sections (\Cref{prop:Piloc} in \Cref{app:teich}).
So we restrict to a neighbourhood $V$ of $T'$ in $\cal F$ on which there is a continuous map $V\to \cal F: T\mapsto f_T$ with $[f_T] = T$.
We can moreover arrange so that $f_{T'} = f'$:
indeed, if this is not the case, just replace every $f_T$ by $f_T \circ f_{T'}^{-1}\circ f'$ and note that $f_{T'}^{-1}\circ f'\in\cal H^0$.
We then define
\[w_*(T) = (f_T\circ \pi_0)^{-1}(z_*').\]
Note that $w_*(T')=w_*'$ and that $T\mapsto w_*(T)$ is continuous.
By making $V$ smaller, we ensure that $w_*(T)$ remains in the interior of $P^0_j$.

We managed to fix $z_*=f_T\circ\pi_0(w_*(T))$ so that it stays equal to $z_*'$ for all $T\in V$.
We also want to keep fixed the value and derivative at $z_*$ of the initial germ used in the definition of $R$.
Hence given any point $(T,\res) \in \cal{TR}$ such that $T \in V$, we choose, as initial germ of straightening coordinate $\phi=\phi_T$ at $z_*$ for the Christoffel symbol associated to $(T,\res)$, the unique one satisfying $\phi_T(z_*)=\phi_0(z_*)$ and $\phi_T'(z_*)=\phi_0'(z_*)$.
Note that $z_* = z_*'\neq \infty$, because $\infty$ is one of the vertices or marked points by the convention that $z_3$ is always $\infty$.

Let us prove that the map $(T,\res)\in V\times \cal R\mapsto R(T,\res)$ is analytic (recall that $\cal {TR} = \cal T \times \cal R$ and that $\cal R$ is defined by \cref{eq:defR}).
We will take advantage of a theorem of Hartogs that asserts that a function $F$ from a finite dimensional complex manifold $X$ to another is analytic if and only if $F\circ \xi$ is analytic for all one complex variable analytic map $\xi$ taking values in $X$.
The space $\cal{TR}$ is finite-dimensional.
By Hartogs' theorem, it is enough to prove that $\tau\mapsto R(T_\tau,\res(\tau))$ is analytic for all analytic map $\tau\mapsto T_\tau\in V$ and all analytic map $\tau\mapsto \res(\tau)\in\cal R$, where $\tau$ is a complex number.

Let us show how this follows from \Cref{lem:rphol}. 
Since $R(T_\tau,\res(\tau))$ is an element of the cartesian product $\C^{n_j}$, we focus on one of its components $\C$.
It corresponds to some vertex or marked point $v\in P^0_j$, more precisely it is the resting place of any path $\gamma_\tau$ of the form $\gamma_\tau = f_{T_\tau} \circ \pi^0 \circ \delta_\tau$ where $\delta_\tau$ is a path in $P^0_j$ from $w_*(T_\tau)$ to $w_*'$ that avoids the vertices and marked points of $P^0_j$ except at the end.
We choose $\delta_\tau$ to depend continuously on $\tau$ (for instance a segment, linearly parametrized, works since $P^0_j$ is convex and $w_*(T)$ is in its interior).
Then $\gamma_\tau$ depends continously on $\tau$.

Let $z_k'=f'(v^0_k)$.
Let $z'_{k_0}$ be the one that corresponds to the specific vertex we called $v$.
The set $U$ on which we will apply \Cref{lem:rphol} is the complement in $\hat\C$ of the union of a small closed disk around each $z'_k$, except $z'_{k_0}$.

The holomorphy hypotheses of \Cref{lem:rphol} are immediately seen to be satisfied.
Continuity of $\tau\mapsto X_\tau=[\gamma_\tau]\in\Gamma$ follows from the continuity of $\tau\mapsto \gamma_\tau$.

We have thus proved that $R$, normalized as described earlier in this proof, depends analytically on $(T,\res)$ near $(T',\res')$.
Since we assumed that $R(T',\res')$, which is an element of $\C^{n_j}$, does not have all its components equal, this is still the case for nearby $(T,\res)$.
Hence $\Per^{-1}(\Conf^*)\subset\cal{TR}$ is open.

Analyticity of $\Per$ follows from the analyticity of $R$.
\end{proof}

\subsubsection{Left inverse property of $\Per$}

We have 
\begin{equation}\label{eq:txi}
\Per\circ\Glu = \on{Id}_{\Sconf}.
\end{equation}
Indeed, we defined $\Glu([\cal P]) = ([f],\res)$ for some, non-unique yet specific, $f\in\cal F$.
Using one of these specific $f$ in the definition of $\Per$ above immediately yields, from the way $f$ is defined, that for all $P=P^0_j\in\cal P_0$ there is a similarity chart on $f\circ \pi^0(P)$ whose image is exactly $P_j$ (minus the vertices).
The resting places hence correspond to the vertices/marked points of $P_j$, whence $\Per([f],\res) = [\cal P]$.

\subsubsection{Right inverse property of $\Per$ and conclusion}

\begin{figure}
\begin{tikzpicture}

\node at (-0.25,-2.9) {\includegraphics[scale=0.5]{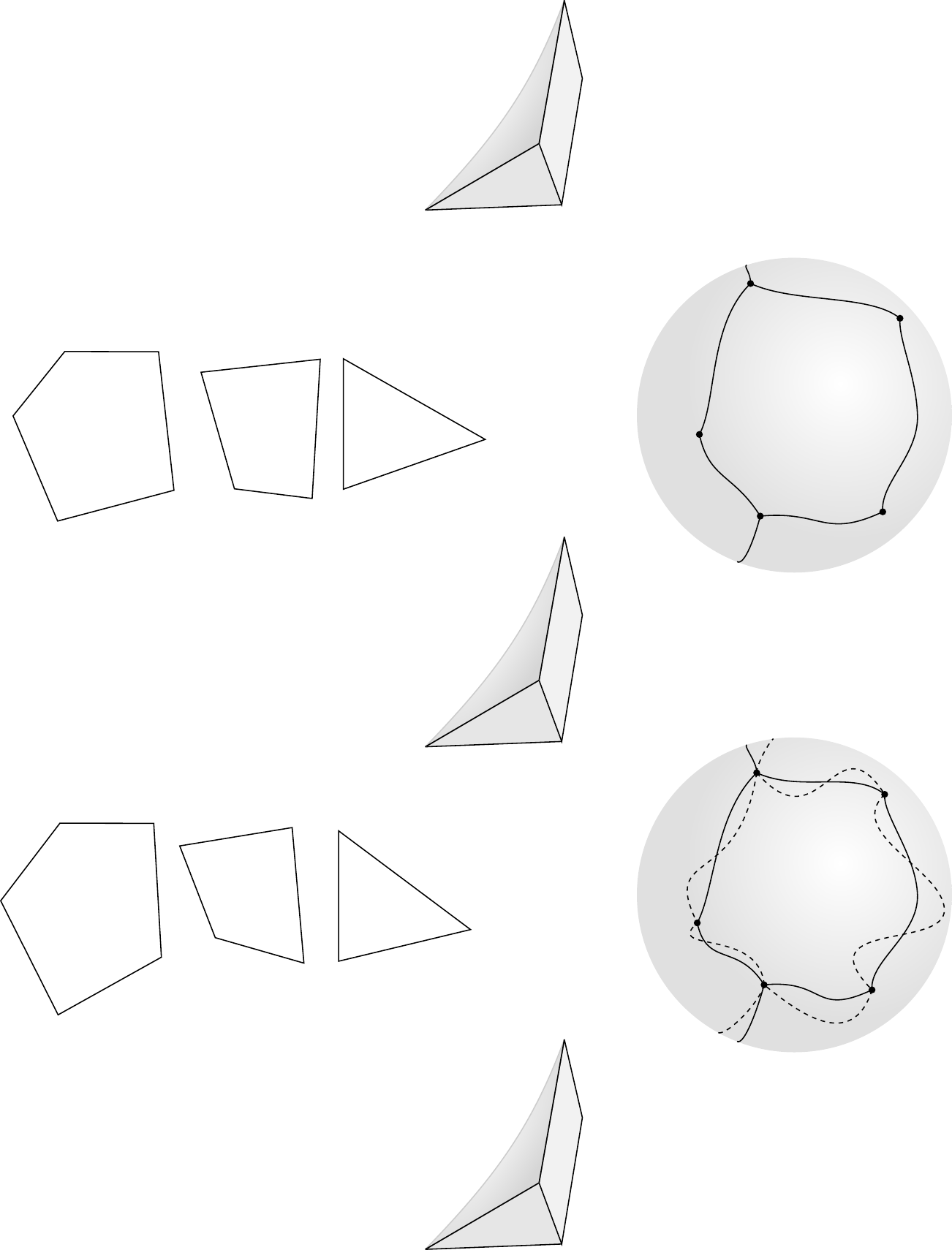}};

\node at (-0.1,3.4) {$\cal S^0$};
\node at (-0.1,-2.8) {$\cal S^0$};
\node at (-0.1,-8.65) {$\cal S^0$};

\node at (-4.5,-0.7) {$P^1_1$};
\node at (-2.7,-0.7) {$P^1_2$};
\node at (-1.35,-0.7) {$P^1_3$};
\node at (-0.5,0.2) {$\cal P^1$};

\draw[->] (0.3,-0.6) -- node[above] {$\prB^1\circ\pi^1$} (1.2,-0.6);

\draw[dashdotdotted, ->] (-0.2,0.3) -- node[above left] {$\Glu$} (0.8,0.8);
\draw[->] (1,1.6) -- node[below left] {$f^1$} (1.8,0.8);

\draw[dashed,->] (1.1,-3.8) -- node[below left] {$f$} (1.9,-4.6);
\draw[dashdotdotted, ->] (1.0,-4.65) -- node[above left] {$\Per$} (0.2,-5.7);
\node at (0,-6.0) {$\cal P$};
\draw[->] (1.2,-8.4) -- node[above left] {$\tilde f$} (2.0,-7.5);
\draw[dashdotdotted, ->] (0.2,-6.3) -- node[below left] {$\Glu$} (1.2,-7.3);

\end{tikzpicture}
\caption{Illustration of the proof of \Cref{lem:ma}.
Solid lines on the sphere depict geodesics while dashed lines are not necessarily geodesics.
We are given $f$ close to $f_1$ and from this we construct $\cal P=\Per([f])$ close to $\cal P^1$ and we must justify that we can choose a representative $\tilde f$ of $\Glu(\cal P)$ that is close to $f^1$ too. Since the vertices $z_k$ are at the same place on the bottom sphere, this will prove that $[f]=[\tilde f]=\Glu([\cal P])$.
The map $\tilde f$ is actually defined by following the saddle connections on the $1$-skeleton on and then by some interpolation in the $2$-cells, and we prove that $[\tilde f]=\Glu([\cal P])$.}
\label{fig:scene}
\end{figure}

Let
\[\Eff=\Glu(\Sconf)\]
where the notation $\Eff$ stands for \emph{effective}.

\begin{lemma}\label{lem:ma}
  Every point in $\Eff$ has an open neighbourhood $W$ in $\cal {TR}$ on which $\Glu\circ\Per$ is defined (i.e.\ $\Per(W)\subset\Sconf$) and equals the identity of $W$.
\end{lemma}
\begin{proof}
Let $\mathrm{eff}^1=\Glu([\cal P^1]) = ([f^1],\res^1)$ for some $f^1\in\cal F$, $\res^1\in\cal R$ and $\cal P^1 = (P^1_j)\in \Sconf$.
We recall that the homeomorphism $f^1:S^0\to\hat\C$ is obtained by mapping each polygon $P^0_j$ to $P^1_j$ and passing to the quotients, then composing with a uniformization of the Riemann surface $\cal S^1$ to $\hat\C$.
By \Cref{eq:txi} we have $\Per(\mathrm{eff}^1) = [\cal P^1]$.
Since $\Sconf$ is open and $\Per$ is continuous at $\mathrm{eff}^1$, there is some neighbourhood $W$ of $\mathrm{eff}^1$ such that $\Per(W)\subset \Sconf$, so $\Glu\circ\Per$ is defined in particular on $W$.

There remains to check that, for a possibly smaller $W$, for all $([f],\res)\in W$, denoting $[\cal P] = \Per([f],\res)$ we have $\Glu([\cal P]) = ([f],\res)$.
By taking $W$ small we can ensure that the strictly convex polygon configuration $[\cal P]$ is close to $[\cal P^1]$.
For any $\eps>0$, by taking $W$ small we also know that for any $([f],\res)\in W$, there is a representative $f$ of $[f]$ that is $\eps$-close to $f^1$ (see point~\ref{item:topo:open} of \Cref{prop:topo}). %which is also \Cref{prop:neigh} in \Cref{app:teich}).
 
Let $\zeta$ be the Christoffel symbol associated to $([f],\res)$ and $\zeta^1$ the one associated to $([f^1],\res^1)$.
The image by $\prB^1$ of an edge of $\cal S^1$ is a saddle connection of $\zeta^1$. 
Let us recall our notation:
\[\cal P^1= (P^1_j)_{j\in J}\text{ and }\cal P= (P_j)_{j\in J}.\]
Let $Q^1_j=\prB^1\circ\prA^1(P^1_j)$. 
This forms a cellular decomposition of $\hat\C$.
The two endpoints of each of these saddle connection must be distinct, since no polygon has only one vertex and all have at most one infinite vertex.
By \Cref{prop:fseg,prop:fseg2}, we can follow continuously (analytically) each of the saddle connections mentioned earlier as a function of $([f],\res)\in W$ provided $W$ is small enough.
The initial saddle connections are injective and thus the followed one remain injective if $W$ is small enough, by the end of \Cref{prop:fseg,prop:fseg2}.

Since the saddle connections are disjoint in the case of $[\cal P^1]$, except at the vertices, a continuity argument proves they are still disjoint away from any neighbourhood of the vertices provided $W$ is small enough.
Saddle connections tend to finite vertices along a straight line to the vertex in the similarity charts whose images are sectors based on the vertex and whose sides are glued by a similarity.
Hence distinct saddle connections remain disjoint near the finite vertices too.
The same is true near infinite vertices because the direction at which the different connections converge to infinity make a non-zero angle with respect to each other, because of the restriction we put on the unbounded polygons.
It follows that the saddle connections remain disjoint, except at their ends.

The followed connections hence define a cellular decomposition of $\hat\C$, with cells $Q_j$ corresponding to deformations of the cells $Q^1_j$.
More precisely let $\Sk\subset\hat\C$ denote 1-skeleton of this cell decomposition for $(Q_j)$ and $\Sk^1$ the same for $(Q^1_j)$.
Matching the saddle connections as parametrized curves defines a homeomorphism $h : \Sk^1\to \Sk$ that is $\eps$-close to the identity provided $W$ is small enough.
It extends to a homeomorphism $h$ of the sphere that is $\eta$-close to the identity, with $\eta\tend 0$ as $\eps\to 0$: this claim follows from a general theorem of topology.\footnote{For instance apply on the image of each polygon the isotopy extension theorem of a circle imbedding in the plane, see for instance Corollary~1.4 of \cite{EK}.}
We then let % $Q_j=Q_j^1$ and
\[\tilde f = h \circ f^1.\]
See the following footnote for an alternative approach.\footnote{Another way to proceed would be to use the regularity of the curves and a control on how they reach the vertices, to identify which region $Q_j$ correspond to $Q_j^1$, then to prove, in a way similar to what we do below, that $Q_j$ is mapped bijectively by any solution of $\phi''/\phi'=\zeta$ to a $\C$-affine image of $P_j$, then to exploit continuous dependence of $\phi$ w.r.t\ $\zeta$ to reduce to finding homeomorphisms close to the identity between $P^1_j$ and $P_j$.}
The map $\tilde f$ is an orientation preserving homeomorphism from $\cal S^0$ to $\hat \C$ and maps each cell $\pi(P^0_j)$ to $Q_j$.

Claim:
\[ [f]=[\tilde f] .\]
Indeed, let $\tilde h: \cal S^0\to\cal S^0$ be such that $\tilde f = f\circ \tilde h$, i.e.\ $\tilde h = f^{-1}\circ \tilde f = f^{-1} \circ h \circ f^1$.
The map $\tilde h$ fixes every vertex and marked point (but it does not send the edges and polygons of $\cal S^0$ to themselves), and is close to the identity.
By \Cref{lem:topo} we have $\tilde h\in\cal H^0$, which proves the claim.

Since our polygons are simply connected, there is a solution $\phi_j=\phi$ of $\phi''/\phi' = \zeta$ defined on $Q_j-\{z_k\}$.
This solution maps the edges to straight lines, because $Q_j$ is bounded by saddle connections, which are geodesics, and between two finite vertices or a vertex and a marked point, the image is a bounded open segment whose ends are the resting place of the paths used to define $\Per([\tilde f],\res)$. But recall that $[\tilde f] = [f]$ hence $\Per([\tilde f],\res) = \Per([f],\res) = [\cal P]$.
The image by $\phi_j$ of an edge of $Q_j$ from a marked point to an infinite vertex is a half-line tending to $\infty$.
Hence all these straight lines form the boundary the strictly convex, possibly marked, polygon $P_j$ composing the collection $[\cal P]$ we started from (up to a $\C$ affine map).
Replacing $\phi_j$ by its post-composition with an affine map and adding the vertices to its domain of definition, we may assume that $\phi_j(Q_j) = P_j$:
\[\phi_j : Q_j \to P_j.\] 
By the winding number theorem, the map $\phi_j$ is a bijection from $Q_j$ to $P_j$, holomorphic in the interior.

Hence $\tilde f$ allows to realize the gluing of the $P_j$. 
More precisely let $\cal S$ be the abstract surface constructed from the $P_j$ and $\pi : \coprod P_j \to \cal S$ the quotient map.
Let $\cal S'$ be $\cal S$ minus the vertices.
The maps $\phi_j\circ\tilde f$ send $\pi^0(P^0_j)$ to $P_j$ in a way that is compatible with the gluings of the $P_j$ (recall that $\tilde f = h \circ f^1$, and $h$ and $\phi_j$ respect the linear parametrization of geodesics) and define a map
\[\psi: \cal S^0 \to \cal S\]
such that $\psi = \pi \circ \phi_j \circ \tilde f$ on each $\pi^0(P^0_j)$.
Let now
\[F := \tilde f\circ\psi^{-1}:\cal S\to\hat \C\]
which satisfies
\[F\circ\pi|_{P_j} = \phi_j^{-1}.\]

Claim: The map $F:\cal S\to\hat \C$ is an analytic isomorphism.
We first prove it on $\cal S'$, where it coincides with $\phi_j^{-1}\circ\pi^{-1}$ on $\pi(P_j)$, which is immediately seen to be holomorphic in the interior of $\pi(P_j)$ since $\pi$ and $\phi_j$ are holomorphic.
On a an edge minus endpoints $e^*$, recall that charts are given as follows: let $\pi(P_j)$ and $\pi(P_k)$ be the cells adjacent to $e^*$ in $\cal S$.
There is a neighbourhood $V$ of $e^*$ and a similarity chart $\phi:V\to\C$ for which
$\phi(\pi(z))=z$ when $z\in P_j \cap \pi^{-1}(V)$ and $\phi(\pi(z))=s(z)$ when $z\in P_k\cap \pi^{-1}(V)$ where $s$ is the $\C$-affine map gluing the corresponding edges of $P_k$ and $P_j$.
The image of this chart is the open neighbourhood $W = W_1\cup W_2$ of the appropriate edge of $P_j$ where $W_1 = P_j\cap \pi^{-1}(V)$ and $W_2=s(P_k\cap \pi^{-1}(V))$.
Then the expression of $F$ from this chart is $F\circ\phi^{-1}$, which is equal to $\phi_j^{-1}$ on $W_1$ and to $\phi_k^{-1}\circ s^{-1}$ on $W_2$. 
These two maps are holomorphic and coincide on the (straight) edge of $P_j$, hence their join is holomorphic.
We have thus proved that the map $F = \tilde f\circ\psi^{-1}$ is an analytic isomorphism from $\cal S'$ to $\hat \C$ minus a finite number of points. These singularities are erasable since $F$ is continuous, hence actually $F$ is an analytic isomorphism from $\cal S$ to $\hat \C$.

The surface $\cal S'$ carries a similarity surface structure (an atlas) coming from the polygons under the quotient map $\pi$, and $F$ sends this structure to $\hat\C - \cal Z$ where $\cal Z$ is the set of images of the vertices.
To this structure corresponds a Christoffel symbol $\tilde \zeta$, and we claim that this symbol is exactly $\zeta$.
Indeed, it is given by $\tilde \zeta=\phi''/\phi'$ for any similarity chart $\phi$.
We can take $\phi=\phi_j$ and by definition $\phi_j''/\phi_j'=\zeta$ holds on $Q_j$.
We conclude by analytic continuation of identities that $\tilde \zeta = \zeta$.

Coming back to the definition of the map $\Glu$, all this implies that $\Glu([\cal P])=([\tilde f],\res)$.
We then conclude as follows: $\Glu\circ\Per([f],\res) = \Glu([\cal P]) = ([\tilde f],\res) = ([f],\res)$.
\end{proof}

The open sets $W$ in \Cref{lem:ma} are contained in $\Eff=\Glu(\Sconf)$ because $\Glu(\Per(W)) = W$ and $\Per(W)\subset\Sconf$. Since by the lemma, every point of $\Eff$ is contained in such a $W$, it follows that $\Eff$ is open. 
Another consequence of the lemma is that $\Glu\circ \Per$ is the identity on $\Eff$.

It follows that the restriction $\Per|_\Eff$ is an injective analytic map that satisfies $\Glu\circ \Per|_\Eff = \on{Id}_{\Eff}$.
By \Cref{eq:txi}, $\Per|_\Eff\circ\Glu = \on{Id}_\Sconf$.
Hence $\Per|_\Eff$ is an analytic bijection between open subsets of complex manifolds and $\Glu$ is its inverse.
It follows that these manifolds have the same dimension and that $\Per|_\Eff$ and $\Glu$ are analytic isomorphisms (\cite{Gu}, Theorem~I.11).

This proves \Cref{prop:hsol}.

\begin{remark}
As an interesting consequence the fact that $\Glu$ and $\Per$ are analytic isomorphisms tells us $\on{dim}_\C \Sconf = \on{dim}_\C \cal {TR}$.
We can also check this identity directly: let $a_i$ be the number of cells of dimension $i$ in the cellular decomposition of $\cal S$ induced by the polygons $P_j$ (not including the marked points).
Let $a'_2$ be the number of unbounded polygons: $a'_2\leq a_2$.
Let $a'_0$ be the number of marked points in $\cal S$.
Each unbounded polygon has two marked points but the marked points are paired under the gluing, hence $a'_0=a'_2$.
We have
$\on{dim} \Sconf = \sum_{j\in J} \on{dim} A^*\C^{n_j} = \sum (n_j-2)$
(note that each $n_j\geq 3$).
If $P_j$ is a bounded polygon then $n_j$ is its number of vertices hence its number of edges.
If $P_j$ is unbounded then $n_j$ is one plus its number of edges.
Adding up, each edge will be counted twice.
It follows that $\on{dim} \Sconf = 2 a_1 + a'_2 - 2 a_2$.
On the other hand, $\dim \cal{TR} = \dim \cal M + \dim \cal R = (a_0 + a'_0 - 3) + (a_0-1)$.
The equality $\on{dim}_\C \Sconf = \on{dim}_\C \cal {TR}$ then follows from cancellation of $a'_2=a'_0$ and from Euler's identity $a_0-a_1+a_2=2$. 
\end{remark}

\section{Appendix}

The reader will find here: an example of non-puncture type vertex obtained by gluing (unbounded) polygons; a proof of \Cref{prop:topo} concerning a topology on the space of paths, which is necessary to state a holomorphic dependence result for developing maps of \Cref{sec:hd:rp}; a presentation of the Teichmüller space of the sphere with $n$ marked points and some of its properties, with (essentially) topological methods so as to avoid the use of the measurable Riemann mapping theorem.

%
% ----------------------------------------------------------
%

\subsection{Example of a hole type infinite vertex}\label{app:hole}

\begin{figure}
\begin{tikzpicture}

\node at (-3,0) {\includegraphics[scale=0.25]{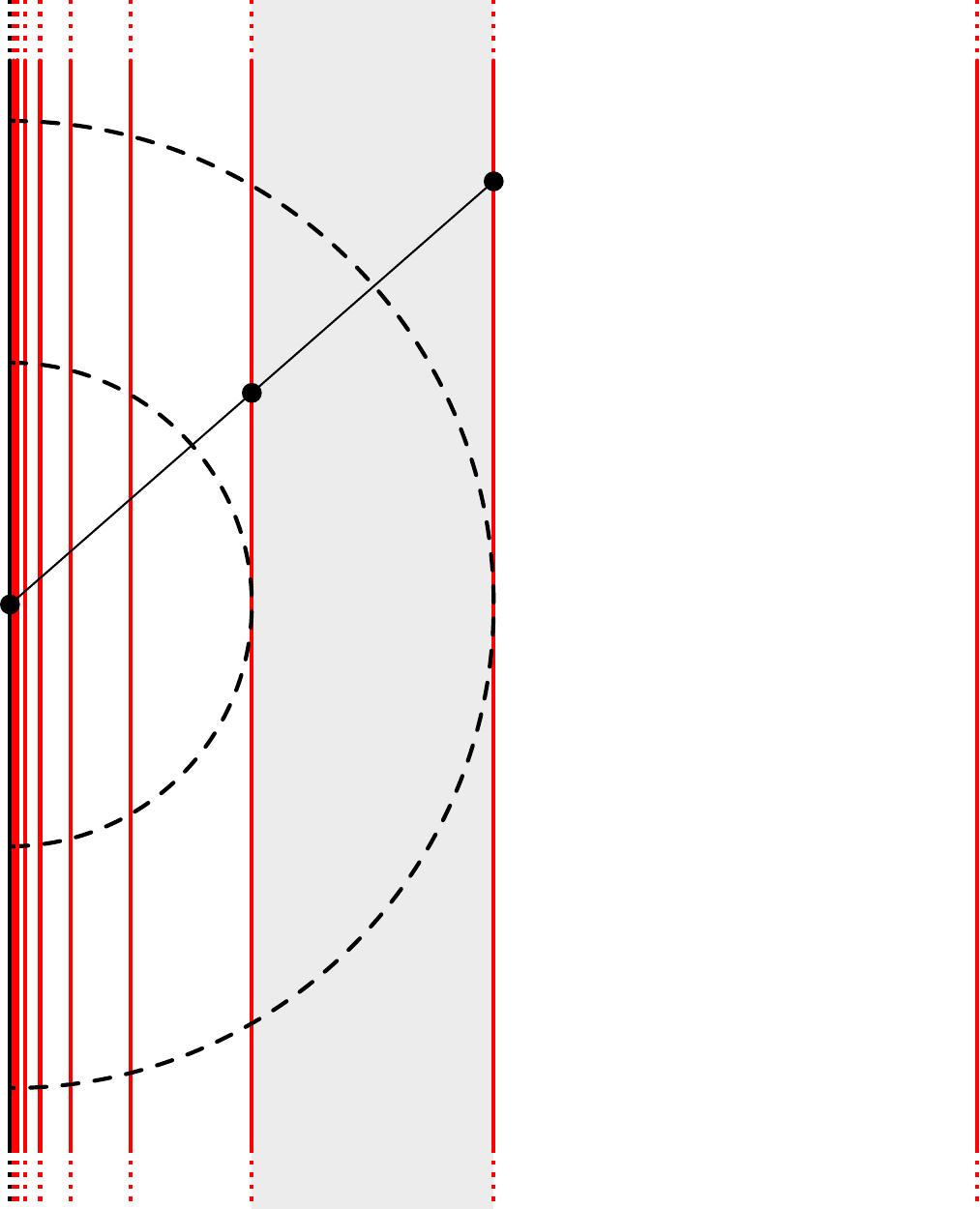}};
\node at (-4.2,1.2) {$z$};
\node at (-3.25,2.15) {$2z$};
\node at (-5,-3) {\ldots};
\node at (-4.55,-3) {$\frac12$};
\node at (-4,-3) {$1$};
\node at (-2.975,-3) {$2$};
\node at (-0.85,-3) {$4$};

\node at (0,0) {$\overset{\log}{\longrightarrow}$};

\node at (3,0) {\includegraphics{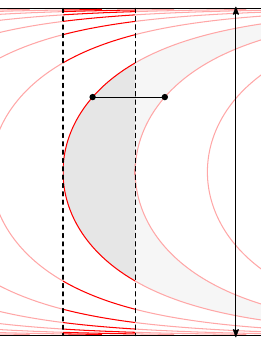}};
\node at (5,0) {$\pi$};
\node at (2.3,1.55) {$z$};
\node at (4.05,1.55) {$z+\log 2$};

\end{tikzpicture}
\caption{The construction of \Cref{app:hole}. The gray fundamental domain is glued via $s:z\mapsto 2z$. The quotient is isomorphic as a Riemann surface to $H/s$ and via the complex logarithm to a strip modulo $z\mapsto z+\log 2$.}
\label{fig:hole}
\end{figure}

Consider the strip $P:``\Re z \in [1,2]"$ and pair the left unbounded edge with the right one by sending $z\mapsto 2z$, i.e.\ $1+iy \to 2+2iy$.
(To respect the temporary restrictions of the section \emph{Unbounded polygons} on page~\pageref{page:unbdd} one can cut $P$ in two pieces along $\R$ and pair the two bounded edges together.) Then the quotient $\cal S$ is homeomorphic to a sphere and there are two vertices, both unbounded. It is then easy to uniformize the Riemann surface $\cal S'$: it is isomorphic to  $H/\sim$ where $H$ is the half plane ``$\Re z>0$'' and $z\sim 2z$. Mapping the situation by $z\mapsto \log z$ the set $H$ becomes the infinite band ``$|\Im z|<\pi/2$'' and the equivalence relation becomes $z\sim z+\log 2$. See \Cref{fig:hole}. Hence the Riemann surface $\cal S'$ is isomorphic to a cylinder of finite height, i.e.\ to a round annulus via the transformation $z\mapsto e^{2\pi i z / \log 2}$, and both vertices are holes.

\subsection{Topology on a space of paths}\label{app:GammaTop}

Let $U$ be an open subset of $\C$.
Denote $\cal G$ the set of paths $\gamma$ in $U$ that satisfy $\forall s\in[0,1]$, $\gamma(s)\neq \gamma(1)$.
For $z_0\neq z_1$, both in $U$, denote $\cal G(z_0,z_1)$ the subset of $\cal G$ of paths $\gamma$ with $\gamma(0)=z_0$ and $\gamma(1)=z_1$.
Given two paths $\gamma$, $\tilde\gamma$ in $\C$, let \[d(\gamma,\tilde \gamma) = \sup_{s\in[0,1]} |\tilde\gamma(s)-\gamma(s)|
.\]

As in \Cref{sec:hd:rp}, we define $\Gamma$ as the quotient space of $\cal G$ under the relation that $\gamma \sim \gamma'$ if and only if $\gamma$ and $\gamma'$ have the same endpoints $z_0$ and $z_1$ and belong to the same path-connected component of $\cal G(z_0,z_1)$.
We denote $[\gamma]$ the equivalence class of $\gamma$.

The distance $d$ on $\cal G$ defines a topology on $\cal G$ and we will define the topology on $\Gamma$ as the quotient under the natural map $\cal G\to \Gamma$. But before we need several lemmas.

\medskip

For all $\eps>0$, we denote $V_{\gamma,\eps}$ the image by the quotient map $\cal G\to \Gamma$ of the ball $B(\gamma,\eps)\subset \cal G$ for the distance $d$ introduced above. In other words, $V_{\gamma,\eps}$ is the set of $[\tilde \gamma]$ where $\tilde \gamma$ is a path in $U$ and $d(\gamma,\tilde\gamma)<\eps$.

We will use the notation
\[d(\gamma,\partial U):=\inf_{s\in[0,1],\, z\in\partial U}d(\gamma(s),z)
.\]
In all the statements below, we assume that all the paths take value in $U$. 

\begin{lemma}\label{lem:GammaTop2}
For all $\gamma$ there is $\eps>0$ such that if $\tilde\gamma(0)=\gamma(0)$, $\tilde\gamma(1)=\gamma(1)$ and $d(\gamma,\tilde \gamma)<\eps$ then $[\tilde\gamma]=[\gamma]$.
\end{lemma}
\begin{proof}
A linear interpolation works except near $z_1$: the condition that the path cannot hit $z_1=\gamma(1)$ at times $s\neq 1$ may fail to hold for the interpolation.
To solve this, let $r>0$ be such that $B(z_1,r)\subset U$.
Let $s_1\in(0,1)$ such that $\gamma([s_1,1])\subset B(z_1,r/2)$.
Let $\eps = \min(r/2, d(\gamma,\partial U), d(\gamma|_{[0,s_1]},\{z_1\}))$.

As in the statement we assume $d(\gamma,\tilde\gamma)<\eps$.
Let $\delta(s) = (1-s)\gamma(s_1)+s\tilde\gamma(s_1)$.
Since $\eps<r/2$, both points $\gamma(s_1)$ and $\tilde\gamma(s_1)$, and hence $\delta$, are  contained in $B(z_1,r)$.
The path $\gamma$ is homotopic within $\cal G(z_0,z_1)$ to $(\gamma|_{[0,s_1]}\cdot\delta)\cdot(\delta^{-1}\cdot\gamma|_{[s_1,1]})$ where the restrictions are appropriately reparametrized.

Dealing with $\gamma|_{[0,s_1]}\cdot\delta$: We have an endpoint fixing homotopy from $\gamma|_{[0,s_1]}\cdot\delta$ to $\tilde\gamma|_{[0,s_1]}$ 
by $t\in[0,1]\mapsto ((1-t)\gamma+t\tilde\gamma)|_{[0,s_1]} \cdot \delta|_{[t,1]}$ (with $\delta|_{[t,1]}$ appropriately reparametrized) ending on the concatenation of $\tilde\gamma|_{[0,s_1]}$ with a constant path, and this is homotopic to $\tilde\gamma|_{[0,s_1]}$.
The two homotopies avoid $z_1$ because $\eps<d(\gamma|_{[0,s_1]},\{z_1\})$ and stay in $U$ because $\eps<d(\gamma,\partial U)$.

Dealing with $\delta^{-1}\cdot\gamma|_{[s_1,1]}$: Denote $\gamma'$ this part and note that it is contained in $B(z_1,r)$.
Since $\eps<r/2$ the path $\tilde\gamma|_{[s_1,1]}$ is also contained in $B(z_1,r)$.
Both $\gamma'$ and $\tilde\gamma|_{[s_1,1]}$ are thus contained in $B(z_1,r)$, they have the same starting point and they avoid $z_1$ except at the end.
A linear interpolation in (lifted) polar coordinates provides a homotopy between them.
\end{proof}

\begin{lemma}\label{lem:GammaTop1}
Let $z_0=\gamma(0)$.
If $0<\eps<d(z_0,\partial U)$ and $\delta$ is any path in $B(z_0,\eps)$ that avoids $z_1=\gamma(1)$ and with $\delta(1)=z_0$ then $[\delta\cdot\gamma] \in V_{\gamma,\eps}$.
\end{lemma}
\begin{proof}
Let $t_0\in(0,1)$ and consider the path $\tilde\gamma$ that sends $s\in[0,t_0]$ to $\delta(s/t_0)$ and $s\in[t_0,1]$ to $\gamma\big(\frac{s-t_0}{1-t_0}\big)$.
Then $[\tilde\gamma]=[\delta\cdot\gamma]$.
For $t_0$ small enough we have $d(\tilde\gamma,\gamma)<\eps$.
\end{proof}

In the next statement we denote $z_0=\gamma(0)$, $z_1=\gamma(1)$, $\tilde z_0=\tilde\gamma(0)$, $\tilde z_1=\tilde\gamma(1)$.
\begin{lemma}\label{lem:GammaTop3}
Given $\gamma\in\cal G$, there exists $\eps>0$ with $\eps<d(\gamma,\partial U)$ and $\eps<|z_0-z_1|/2$ such that for all $\tilde\gamma\in\cal G$, if $d(\gamma,\tilde\gamma)<\eps$, if $\delta$ is any path within $B(z_0,\eps)$ from $\tilde z_0$ to $z_0$ and if $\psi$ is a self homeomorphisms of $U$ that is the identity outside $B(z_1,\eps)$ and such that $\psi(z_1) = \tilde z_1$ then $[\tilde \gamma] = [\delta\cdot(\psi\circ\gamma)]$.
\end{lemma}
\begin{proof}
Let $\eps_0$ be the $\eps$ given by \Cref{lem:GammaTop2}.
Choose any $\eps>0$ with $\eps<\min(\eps_0/4,|z_1-z_0|/2,d(\gamma,\partial U))$.
The relation to prove $[\tilde \gamma] = [\delta\cdot(\psi\circ\gamma)]$ is equivalent to 
$[\psi^{-1}\circ(\delta^{-1}\cdot\tilde\gamma)] = [\gamma]$.
We have $\psi^{-1}\circ(\delta^{-1}\cdot\tilde\gamma) = (\psi^{-1}\circ\delta^{-1})\cdot(\psi^{-1}\circ\tilde\gamma)$.
Since $\eps<|z_1-z_0|/2$, it follows that  $\psi^{-1}\circ\delta^{-1} = \delta^{-1}$.
By \Cref{lem:GammaTop1}, $[\delta^{-1} \cdot (\psi^{-1}\circ\tilde\gamma)]=[\alpha]$ for some path $\alpha$ with $d(\alpha,\psi^{-1}\circ\tilde\gamma)<\eps$.
We have $d(\psi^{-1}\circ\tilde \gamma,\tilde\gamma) < 2\eps$.
Hence $d(\alpha,\gamma)<4\eps<\eps_0$ so we can apply \Cref{lem:GammaTop2}: $[\alpha]=[\gamma]$.
\end{proof}

Note that given $\gamma$, $\tilde\gamma$ and $\eps$ as above, there always exist explicit objects $\delta$ and $\psi$ exist that satisfy the hypotheses.

\begin{lemma}\label{lem:GammaTop4}
If $[\gamma]=[\gamma']$ then for all $\eps'<d(\gamma',\partial U)$, there exists $\eps < d(\gamma,\partial U)$ such that $V_{\gamma,\eps}\subset V_{\gamma',\eps'}$.
\end{lemma}
\begin{proof}
Note that $d(\gamma,\gamma')$ is not assumed to be small.
Let $\eps_1$ be the $\eps$ given by \Cref{lem:GammaTop3}.
Let $\eps < \min(\eps'/3 , \eps_1)$. In particular $\eps <\min( d(\gamma,\partial U), |z_1-z_0|/2)$.
Consider any $[\tilde\gamma] \in V_{\gamma,\eps}$, for which by definition we can take a representative $\tilde\gamma$ with $d(\tilde \gamma,\gamma)<\eps$.
We want to prove that $[\tilde\gamma]\in V_{\gamma',\eps'}$, i.e.\ that $\tilde\gamma$ is equivalent to a path within distance $\eps'$ of $\gamma'$.
By \Cref{lem:GammaTop3}, we have $[\tilde\gamma] = [\delta\cdot(\psi\circ\gamma)]$ for some appropriate $\delta$ and $\psi$.
By post-composing a homotopy from $\gamma$ to $\gamma'$ with $\psi$ we get that $[\psi\circ\gamma] = [\psi\circ \gamma']$.
Hence $[\tilde\gamma] = [\delta\cdot(\psi\circ\gamma')]$.
By \Cref{lem:GammaTop1}, $[\delta\cdot(\psi\circ\gamma')]$ is homotopic in $\cal G(\tilde z_0,\tilde z_1)$ to a path $\alpha$ with $d(\alpha,\psi\circ\gamma')<\eps<\eps'/3$.
Moreover, $d(\psi\circ\gamma',\gamma')<2 \eps<2\eps'/3$.
Hence $d(\alpha,\gamma')<\eps'$.
\end{proof}

Consider the natural surjection
\[Q: \gamma \in \cal G \to [\gamma] \in \Gamma = \cal G/\sim
.\]
Recall that $\cal G$ is endowed with a metric $d$, in particular with a topology.

\begin{lemma}
For every open subset $O$ of $\cal G$, the saturate $Q^{-1}(Q(O))$ is open.
\end{lemma}
\begin{proof}
Let $\gamma\in Q^{-1}(Q(O))$: there exists $\gamma' \in O$ such that $Q(\gamma)=Q(\gamma')$, i.e.\ $[\gamma]=[\gamma']$.
Since $O$ is open, it contains $B(\gamma',\eps')$ for some $\eps'>0$, which we may assume $<d(\gamma,\partial U)$.
By \Cref{lem:GammaTop4}, there exists $\eps>0$ such that $V_{\gamma,\eps} \subset V_{\gamma',\eps'}$, i.e. $Q(B(\gamma,\eps))\subset Q(B(\gamma',\eps'))$. It follows that $B(\gamma,\eps)\subset Q^{-1}(Q(O))$.
\end{proof}

We now endow $\Gamma$ with the quotient topology.
The quotient map $Q$ is automatically continuous and the fact that saturate of open sets are open is equivalent to the fact that $Q$ is open.
Note that the set $V_{\gamma,\eps} = Q(B(\gamma,\eps))$ are in particular open and form a basis of the topology.

% It follows from \Cref{lem:GammaTop4} that if we define a ``neighbourhood'' of $[\gamma]$ as a set that contains $V(\gamma,\eps)$ for some $\eps>0$, then this definition is independent of the choice of representative $\gamma$ of $[\gamma]$.
% Declare a subset of $\Gamma$ ``open'' if for all of its points it contains a ``neighbourhood'' of it.

% Note that ``neighbourhoods'' of $[\gamma]$ exist, contain $[\gamma]$, the intersection of two ``neighbourhoods'' is still a ``neighbourhood'' (they both contain $V(\gamma,\eps)$ for the min of their $\eps$).
% It follows that the collection of ``open'' subsets of $\Gamma$ forms a topology. We can thus remove the quotes on the word: open.
% Moreover $V(\gamma,\eps)$ is open: all its elements are of the form $[\tilde \gamma]$ with $d(\tilde \gamma,\gamma)<\eps$ and thus $[\tilde\gamma] \in V(\tilde\gamma,\eps-d(\tilde \gamma,\gamma))\subset V(\gamma,\eps)$.
% It follows that our notion of ``neighbourhood'' coincides with the notion of neighbourhood associated to the topology: we can remove the quotes here too.

\begin{lemma}\label{lem:GammaTop5}
The topology on $\Gamma$ is Hausdorff separated.
\end{lemma}
\begin{proof}
If $(z'_0,z'_1)\neq(z_0,z_1)$ then it is enough to take $\eps<\frac{1}{2}\max(|z'_0-z_0|,|z'_1-z_1|)$ for $V_{\eps,\gamma}\cap V_{\eps,\gamma'}= \emptyset$ to hold.

If $(\tilde z_0,\tilde z_1)=(z_0,z_1)$ but $[\tilde\gamma]\neq[\gamma]$ then consider the value $\eps$ that \Cref{lem:GammaTop3} associates to $\gamma$ and the value $\eps'$ that it associates to $\gamma'$.
Let $\eps_0 = \min(\eps,\eps', |z_1-z_0|/2)$.
Then we claim that $V_{\eps_0,\gamma}\cap V_{\eps_0,\gamma'}= \emptyset$.
Otherwise there are two paths $\tilde \gamma$ and $\tilde \gamma'$ with $[\tilde \gamma] = [\tilde \gamma']$ and $d(\gamma,\tilde \gamma)<\eps_0$ and $d(\gamma',\tilde \gamma')<\eps_0$.
Let $\delta$ be a segment from $\tilde z_0=\tilde \gamma(0)=\tilde\gamma'(0)$ to $z_0=\gamma(0)=\gamma'(0)$ and $\psi$ be a self homeomorphism of $U$ that sends $z_1=\gamma(1)=\gamma'(1)$ to $\tilde z_1 = \tilde\gamma(1) = \tilde \gamma '(1)$ and is the identity outside $B(z_1,\eps)$.
Then by \Cref{lem:GammaTop3}, $[\tilde\gamma] = [\delta\cdot(\psi\circ \gamma)]$ from which it follows that $[\gamma]=[\delta^{-1}\cdot(\psi^{-1}\circ\tilde\gamma)]$.
Similarly
$[\gamma']=[\delta^{-1}\cdot(\psi^{-1}\circ\tilde\gamma')]$.
From the homotopy between $\tilde\gamma$ and $\tilde \gamma'$ one defines a homotopy between $[\delta^{-1}\cdot(\psi^{-1}\circ\tilde\gamma)]$ and $[\delta^{-1}\cdot(\psi^{-1}\circ\tilde\gamma')]$.
hence $[\gamma]=[\gamma']$, leading to a contradiction.
\end{proof}

Consider the map
\[\on{Ends}: [\gamma]\in\Gamma \mapsto (\gamma(0),\gamma(1))\in\C^2
.\]
It is continuous: indeed on $V(\gamma,\eps)$ it takes values in $B(z_0,\eps)\times B(z_1,\eps)$.

\begin{lemma}\label{lem:GammaTop6}
For the $\eps$ of \Cref{lem:GammaTop3} the map above is a homeomorphism from $V(\gamma,\eps)$ to $B(z_0,\eps)\times B(z_1,\eps)$.
\end{lemma}
\begin{proof}
We saw continuity.
Injectivity follows from \Cref{lem:GammaTop3}: for two $\tilde \gamma$ with the same extremities $\tilde z_0$ and $\tilde z_1$ we can choose the same $\delta$ and the same $\psi$.
Surjectivity can be achieved by letting $\tilde\gamma = \delta\cdot(\psi\circ\gamma)$, choosing any $\delta$ from $\tilde z_0$ to $z_0$ and choosing $\psi$ sending $z_1$ to $\tilde z_1$ so that it moves each point of $B(z_1,\eps)$ by at most $\eps$ (this is possible for the Euclidean distance\footnote{This may look specific to the distance we are working with on the plane. The reader may find more satisfying the weaker statement that for the $\eps$ of \Cref{lem:GammaTop3}, the map $\on{Ends}$ is a homeomorphism from $V(\gamma,\eps)\cap\on{Ends}^{-1}(B(z_0,\eps)\times B(z_1,\eps/2))$ to $B(z_0,\eps)\times B(z_1,\eps/2)$. Then any $\psi$ works and there is more flexibility on which distances allow for the proof to work.}): then $[\tilde \gamma]\in V_{\gamma,\eps}$ with a proof similar to that of \Cref{lem:GammaTop1}.
Continuity of the inverse is achieved by definining explicit $\delta$ (a straight segment) and $\psi$ that continuously depend on respectively $\tilde z_0$ and $\tilde z_1$ (i.e.\ we have explicit sections taking values in $V(\gamma,\eps)$).
\end{proof}

% Consider the natural surjective map
% \[Q: \cal G\to \Gamma\]
% that sends $\gamma$ to $(z_0,z_1,[\gamma])$.
% \begin{lemma}\label{lem:GammaTop6b}
% The map $Q$ is continuous and open.
% \end{lemma}
% \begin{proof}
% Let $B(\gamma,\eps)$ denote the ball of centre $\gamma$ and radius $\eps$ for the distance $d$ on $\cal G$.

% The map $Q$ is open: The elements $V(\gamma,\eps)$ form a basis of open neighbourhoods for $\Gamma$ and $V(\gamma,\eps)=Q(B(\gamma,\eps))$.

% The map $Q$ is continuous: It is equivalent to: for all $\gamma\in\cal G$ for all neighbourhood of $Q(\gamma)$, its preimage by $Q$ is a neighbourhood of $\gamma$. This is immediate as we saw that the sets $V(\gamma,\eps)$ forms a basis of neighbourhoods at $[\gamma]$ and that $B(\gamma,\eps)\subset Q^{-1}(V(\gamma,\eps))$.
% \end{proof}

% As a consequence we obtain the following alternative presentation of the topology.

% \begin{corollary}\label{cor:GammaTop7}
% The topology we defined on $\Gamma$ coincides with the quotient topology under $Q$.
% \end{corollary}

%
% ----------------------------------------------------------
%

\subsection{Teichmüller space}\label{app:teich}

We only focus here on the Teichmüller space of a sphere with marked points.

There are several points of view on Teichmüller spaces: 
sometimes they are seen as spaces of hyperbolic metrics, or spaces of conformal structures, spaces of deformations, spaces of marked points, etc.
Most references on Teichmüller space\footnote{There are a lot of books treating the subject. See for instance \cite{GL,Hu}.} study this space using quasiconformal maps and the measurable Riemann mapping theorem (MRMT), and its version with holomorphic dependence. But in the present article, we want to prove the MRMT's, so we cannot follow this approach.
Fortunately, for the study of the Teichmüller space of the sphere with marked points, the MRMT is not required to define and prove the basic properties that we need in the present article.
In particular, instead of defining the Teichmüller space via equivalence classes of quasiconformal maps, we use equivalence classes of homeomorphisms.\footnote{The proof of some classical theorems that we use about spaces of homeomorphisms may look hard. It is likely that the use of another space of maps, like smooth maps, would make the theory a bit easier, but we have not checked that.}

All the material here is classical, but it is hard to find references not involving the use of quasiconformal maps.

\medskip

Let $\cal S^0$ be a topological sphere and $v_k^0$ for $1\leq k\leq m$ be $m$ distinct points of $\cal S^0$.
We assume that
\[m\geq 3.\]
Let $V^0$ denote the $m$-uplet $(v_k^0)$ and $\cal V^0$ denote the set $\{v_k^0\}$.

Let $\cal F$ be the set of orientation preserving homeomorphisms $h:\cal S^0 \to \hat\C$ such that
\[h(v_1^0)=\infty,\ h(v_2^0)=0\text{ and }h(v_3^0)=1.\]
Let us endow $\hat \C$ with a spherical metric\footnote{A metric for which the circles are geometric circles via a stereographic projection to the Euclidean unit sphere in $\R^3$.} $d$ and the set $\cal F$ with the metric 
\[d(f_1,f_2) = \sup\setof{d(f_1(x),f_2(x))}{x\in\cal S^0}.\]
This defines a topology on $\cal F$.
The space $\cal F$ is path connected \cite{K}, in particular connected.
 
Let $\cal H^0$ be the group of orientation preserving self homeomorphisms of $\cal S^0$ that are isotopic to the identity rel.\ $\cal V^0$, i.e.\ by an isotopy fixing each marked point $v_k^0$.
It follows immediately from this definition that $\cal H^0$ is path connected.
The group $\cal H^0$ acts on the right on the set $\cal F$ by composition:
for $\phi\in\cal H^0$ and $f\in\cal F$ we let $f \cdot \phi = f\circ \phi$.
\begin{definition}
The Teichmüller space $\cal T$ associated to $\cal S^0, V^0$
is the quotient of $\cal F$ for this action:
\[\cal T = \cal F/\cal H^0.\]
\end{definition}
We endow $\cal T$ with the quotient topology. We denote
\[\Pi : \cal F \to \cal T\] the quotient map.

\subsubsection{A useful lemma}

Let us endow $\cal S^0$ with a metric $d$ compatible with its topology and the set $\cal H(\cal S^0)$ of self homeomorphisms $\phi$ of $\cal S^0$ with the metric
\[d(\phi_1,\phi_2) = \sup\setof{d(\phi_1(x),\phi_2(x))}{x\in\cal S^0}.\]
This endows $\cal H^0$ with a topology,\footnote{The metric $d$ is not complete on $\cal H^0$. A complete metric would be $d(\phi_1,\phi_2)+d(\phi_1^{-1},\phi_2^{-2})$. Yet, these two metrics induce the same topology on $\cal H^0$.} for which
it is a topological group.\footnote{By uniform continuity of $\phi\in\cal H^0$, left multiplication by $\phi$ is continuous.
Right multiplication by $\phi$ is an isometry.
It follows that multiplication is continuous.
Inversion is continuous at $\phi$ by uniform continuity of $\phi^{-1}$.}

We copy here for convenience \Cref{lem:topo}:
\begin{lemma}\label{lem:topo0}
There exists $\eps>0$ such that for all orientation preserving self homeomorphism $\phi$ of $\cal S^0$, if $\phi$ fixes every point in $\cal V^0$ and satisfies $d(\phi,\on{id}_{\cal S^0})<\eps$ then $\phi\in\cal H^0$ (i.e.\ $\phi$ is isotopic to the identity rel.\ $\cal V^0$).
\end{lemma}
It follows for instance from \cite{Ro}, which states in particular that the space of self-homeomorphisms of $S^2$ is locally arcwise connected, together with the use of an explicit correction to leave the points in $\cal V^0$ fixed, for instance using the tool described later here.

\subsubsection{Moduli space $\cal M$}

Let $\cal M$ be the set of $m$-uplets $(z_k)$ of pairwise distinct points in $\hat \C$ such that $z_1 = \infty$, $z_2=0$ and $z_3=1$.
It is a complex manifold of dimension $m-3$.
For $f\in\cal F$ we define
\[M : \begin{array}{rcl} \cal F & \to & \cal M \\
f & \mapsto & (f(v_k^0))\end{array}\]
which is a continuous map.
There are many ways\footnote{We leave it as an exercise to the reader. Igor Belegradek pointed to us an interesting generalization that can be found in \cite{AB}.} to see that the map $M$ is surjective, i.e.\ that for all $Z=(z_k)\in\cal M$ there exists an orientation preserving homeomorphism from $\cal S^0$ to $\hat\C$ sending each $v^0_k$ to $z_k$.
We also define
\[\pi:\cal T \to \cal M\]
as follows: to $[f]\in\cal T$ we associate $M(f)$, which does not depend on the representative $f$ of $[f]$.
So
\[\pi\circ\Pi = M\]
Since $M$ is surjective, $\pi$ is surjective too.

The map $M$ is open. Indeed, given $f\in\cal F$ and $Z=(z_k)\in\cal M$ close to $M(f)$, one constructs a homeomorphism close to $f$ and that send each $v^0_k$ to $z_k$, by post-composing $f$ with a homeomorphism that is the identity outside small disks around the points $f(v^0_k)$ and sends $f(v^0_k)$ to $z_k$.

\subsubsection{More topology}

By definition of the quotient topology, open sets are subsets of $\cal T$ whose preimages in $\cal F$ are open.
In our case, the saturate $\Pi^{-1}(\Pi(O))$ of open sets $O\subset \cal F$ are open: indeed we have a group action of $\cal H^0$ and for all $\phi\in\cal H^0$, the map $f\mapsto f\circ \phi$ is continuous, and its inverse too, so is a homeomorphism.
It follows that the the map $\Pi:\cal F\to\cal T$ is open, and open subsets of $\cal T$ are also the projection in $\cal T$ of open subsets of $\cal F$.
Hence:

\begin{proposition}\label{prop:neigh}
A basis of neighbourhoods of $[f]$ in $\cal T$ is given by the sets
\[\setof{[g]}{g\in\cal F,\ d(g,f)<\eps}\]
where $f$ is fixed and $\eps>0$ vary.
\end{proposition}

\begin{figure}
\begin{tikzpicture}
\node at (0,0) {\includegraphics{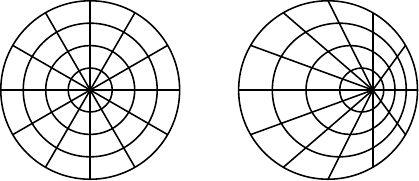}};
\end{tikzpicture}
\caption{Example of a homeomorphism $h_a$ of a Euclidean disk in the plane that is the identity on its boundary, sends its centre to $a$ and depends continuously on $a$. It is defined by $h_a(re^{i\theta})=(1-r)a+re^{i\theta}$. For a closed disk on the Euclidean sphere, conjugacy of $h_a$ by one's preferred homeomorphism to a Euclidean disk can be used.}
\label{fig:move}
\end{figure}

Note that these particular neighbourhoods are open since $\Pi$ is open. Note also that they depend on $(f,\eps)$ and not just on $([f],\eps)$.

\begin{lemma}
The quotient space $\cal T$ is separated.
\end{lemma}
\begin{proof}
Consider $f_1$ and $f_2$ in $\cal T$ and $\eps>0$.
Assume that the associated neighbourhoods defined in \Cref{prop:neigh} intersect.
Call them $V(f_1,\eps)$ and $V(f_2,\eps)$.
This means that $[g_1]=[g_2]$ for some $g_i$ in $\cal F$ with $d(g_i,f_i) <\eps$.
Hence $g_2 = g_1 \circ h_0$ with $h_0\in\cal H^0$.
In particular $g_2(v^0_k)=g_1(v^0_k)$.

Let us assume that $V(f_1,\eps)$ and $V(f_2,\eps)$ intersect for all $\eps>0$.
So for all $\eps>0$ we get maps $g_1$, $g_2$ and $h_0$ as above, that depend on $\eps$.
By passing $g_2(v^0_k)=g_1(v^0_k)$ to the limit as $\eps\tend 0$ we get $f_2(v^0_k) = f_1(v^0_k)$ for all $k$.
Now write $f_2 = f_1 \circ h$ with $h := f_1^{-1}\circ f_2$.
We have $h(v^0_k) = f_1^{-1}(f_2(v^0_k)) = v^0_k$.
Hence $h_0\circ h^{-1}$ fixes all $v^0_k$ too.
We also have $h_0\circ h^{-1} = (g_1^{-1} \circ f_1) \circ h \circ  (f_2^{-1} \circ g_2) \circ h^{-1}$.
So by continuity of composition and inversion, the map $h_0\circ h^{-1}$ is close to the identity when $\eps$ is small ($h$ does not depend on $\eps$).
It follows from \Cref{lem:topo0} that $h_0\circ h^{-1}\in\cal H^0$ (i.e.\ it is isotopic to the identity rel $\cal V^0$).
Since $h_0\in\cal H^0$ and $\cal H^0$ is a group, we get $h\in \cal H^0$.
Since $f_2 = f_1 \circ h$, we get $[f_2]=[f_1]$.
\end{proof}

\subsubsection{Complex structure on $\cal T$}

Recall that
\[\cal F \overset{\Pi}\longrightarrow \cal T \overset{\pi}\longrightarrow \cal M.\]

\begin{lemma}\label{lem:pilh}
The map $\pi$ is a local homeomorphism.
\end{lemma}
\begin{proof}
The map $\pi$ is continuous: indeed, given an open subset $U$ of $\cal M$, we have $\pi^{-1}(U) = \Pi\circ M^{-1}(U)$ ($\Pi$ is surjective), and since $\Pi$ is open and $M$ continuous, it follows that $\pi^{-1}(U)$ is open.

The map $\pi$ is also open: indeed, consider an open subset $O\subset \cal T$. By continuity of $\Pi$, the set $O'=\Pi^{-1}(O)$ is open. Since $\pi(O)=M(O')$ ($\Pi$ is surjective) and we saw that $M$ is open, it follows that $\pi(O)$ is open.

There remains to check that the map $\pi$ is locally injective.
That is, using the neighbourhoods of \Cref{prop:neigh}, to check that for every $f$, there is $\eps>0$ small enough such that for all $f_1$, $f_2$ with $d(f_i,f)<\eps$, if $\pi([f_1])=\pi([f_2])$ then $[f_1]=[f_2]$.

So we assume $\pi([f_1])=\pi([f_2])$, i.e.\ $M(f_1)=M(f_2)$.
Let us write $f_1=f_2\circ\psi$ where $\psi:=f_2\circ f_1^{-1}\in\cal H(\cal S^0)$.
Then $\psi$ fixes every point of $\cal V^0$.
If $\eps$ is close to $0$ then $\psi$ is close to the identity by continuity of composition and inversion,\footnote{This continuity can be proved directly, or we can conjugate by $f$ to reduce it to the already seen fact that $\cal H(\cal S^0)$ is a topological group.} hence by \Cref{lem:topo0} $\psi$ belongs to $\cal H^0$.
Since $f_1= f_2\circ\psi$, we get $[f_1]=[f_2]$.
\end{proof}

The fact that $\pi:\cal T\to\cal M$ is a local homeomorphism allows to endow $\cal T$ with a complex manifold atlas, of dimension $m-3$, for which $\pi$ is analytic.

Actually the map $\pi$ is a covering (better: a universal covering) but we do not use this fact in the present article.

\subsubsection{A useful tool}

Consider an indexed set of points $Z^*=(z_k^*)\in\cal M$ (so $z^*_1$, $z^*_2$ and $z^*_3$ are fixed to $\infty$, $0$ and $1$) and recall that we have put a spherical metric on $\hat\C$.
Choose $\eps>0$ so that the balls $B_k=B(z_k^*,\eps)$, for $k\geq 4$, are pairwise disjoint.
The product
\[ U=\{\infty\}\times\{0\}\times\{1\}\times\prod_{k\geq 4} B_k\]
is an open subset of $\cal M$.
  The set $B_k$ is a spherical disk on the Riemann sphere.
  It is easy to construct an explicit orientation preserving homeomorphism of the disk $B_k$ such that $z_k^*$ is mapped to a given $z_k\in B_k$ and that depends continuously on $z_k$. See for instance \Cref{fig:move}.
  By taking the join over $k$ and completing with the identity, we obtain an orientation preserving homeomorphism
  \[\psi_Z : \hat\C\to\hat \C\]
that depends continuously on $Z=(z_k)\in U$ and such that 
\begin{itemize}
\item $\psi_Z(z^*_k) = z_k$ for all $k$,
\item $\psi_{Z^*} = \on{Id}$,
\item $d(\psi_{Z},\on{Id})\leq 2\eps$ and $d(\psi_{Z}^{-1},\on{Id})\leq 2\eps$.
\end{itemize}
This family depends on $Z^*$ and $\eps$ but we chose to omit them in the notation $\psi_Z$.

\begin{proposition}\label{prop:Piloc}
 The map $\Pi:\cal F\to\cal T$ is a fibre bundle with typical fibre $\cal H^0$.
\end{proposition}
\begin{proof}
  Let $T^* = [f^*]\in\cal T$.
  We saw in \Cref{prop:neigh} that for any $\eta>0$, the set $V=\setof{[f']}{f'\in\cal F, d(f',f^*)<\eta}$ is an open neighbourhood of $T^*$.
  Let $Z^*=M(f^*)=\pi(T^*)$ and for $\eps$ small enough consider the set $U\subset \cal M$ and the family $\psi_Z$, $Z\in U$, introduced above.
  We choose $\eta=\eps/2$, so $\eta<\eps$.
  Then
  \[\forall T\in V,\ \pi(T)\in U.\]
   Consider the continuous map
  \[ \Theta :
  \begin{array}{rcl}
  V\times \cal H^0 & \to & \cal F
  \\
  (T,h) & \mapsto & \psi_{\pi(T)} \circ f^* \circ h
  \end{array}
  \]
  We claim that $\Theta$ is a local trivialization of $\Pi$ provided $\eps$ is small enough.
  
  We must first check that $\Pi\circ\Theta(T,h)=T$. 
  We have $\Pi\circ\Theta(T,h) = [\psi_{\pi(T)} \circ f^* \circ h]$ and $T = [f']$ for some $f'\in\cal F$ with $d(f',f^*)<\eta<\eps$.
  We must check that $(f')^{-1} \circ \psi_{\pi(T)} \circ f^* \in\cal H^0$.
  Note that this composition fixes every point in $\cal V^0$.
  By \Cref{lem:topo0} is is enough to prove that the composition is close to the identity.
  Now, $d(f',f^*)<\eta<\eps$ and $d(\psi_{\pi(T)},\on{Id})<\eps$ so by continuity of composition and inversion, if $\eps$ is small enough, $(f')^{-1} \circ \psi_{\pi(T)} \circ f^*$ is close to $(f^*)^{-1} \circ \on{Id}\circ f^*= \on{Id}$.
  
  Then we must check that the image of $\Theta$ is equal to $\Pi^{-1}(V)$ (the latter is automatically an open set since $\Pi$ is continuous).
  By the previous paragraph it is contained in $\Pi^{-1}(V)$.
  Let us prove the converse inclusion.
  Let $f\in\cal F$ with $\Pi(f)\in V$, i.e.\ there exists $f'\in \cal F$ and $h^0\in\cal H^0$ such that $d(f',f^*)<\eta$ and $f= f'\circ h^0$.
  Let $T = \Pi(f)$, so $\pi(T)=M(f)\in U$. 
  Then $f = \psi_{\pi(T)} \circ f^* \circ g$ for the map $g= (f^*)^{-1} \circ \psi_{M(f)}^{-1}\circ f = (f^*)^{-1} \circ \psi_{M(f)}^{-1}\circ f' \circ h^0$ and because $\cal H^0$ is a group, it is enough to check that
  \[(f^*)^{-1} \circ \psi_{M(f)}^{-1}\circ f'\in \cal H^0.\]
  This is checked using \Cref{lem:topo0} exactly as in the previous paragraph.

  Last, the reciprocal of $\Theta$ is given by
  \[\Xi : \begin{array}{rcl}
    \Pi^{-1}(V) & \to & V\times \cal H^0
    \\
    f & \mapsto & ([f],(f^*)^{-1}\circ \psi_{M(f)}^{-1}\circ f)
  \end{array}
  \]
  which is a continuous expression.
\end{proof}

Fibre bundles have in particular the homotopy lifting property and local sections. Let us make this last statement explicit.

\begin{proposition}\label{prop:localSections}
  For all $T^*\in\cal T$ and all representatives $f^*$ of $T$ (i.e.\ $T = [f^*]$), there exists $\eps>0$ such that for the open neighbourhood $V = \setof{[f']}{f'\in\cal F, d(f',f^*)<\eps}$ of $T$ and for the family $Z\mapsto \psi_Z$ described before \Cref{prop:Piloc},
  the following map is a local section of $\Pi$:
  \[
  \begin{array}{rcl}
  V & \to & \cal F
  \\
  T & \mapsto & \psi_{\pi(T)} \circ f^*
  \end{array}
  \]
\end{proposition}

%
% ----------------------------------------------------------
%

\clearpage

\part{The measurable Riemann mapping theorem}\label{part:beltrami}

The presentation here is not self contained.
In particular some lemmas of \cite{A} will be used without proof.
We will give precise reference to the second edition \cite{Abis}.
However, the arguments of the present article can be followed without these references: required material is recalled in the text and in \Cref{app:dl1,app:pf:lem:sssol}.
Complementary information on other aspects is given in three other appendices to the present part.

%
% ----------------------------------------------------------
%

\section{Statement}\label{sec:mrmt}

We give here definitions necessary to state the measurable Riemann mapping theorem (MRMT) and state it, and its version with holomorphic dependence (Ahlfors-Bers theorem).

Let $\D$ denote the unit disk in $\C$.
Let $\cal B$ denote the complex Banach space of complex-valued $L^\infty$ functions defined on $\C$, endowed with the essential supremum norm. Elements of $\cal B$ are not functions but classes of functions, for the equivalence relation of being equal almost everywhere.

\begin{definition}[Beltrami equation, quasiconformal map]\label{def:be}
Given $\mu$ in the unit ball of $\cal B$, and $z\mapsto\mu(z)$ a representative, a function $f:\C\to\C$ is called a \emph{solution of the  Beltrami equation associated to $\mu$} if:
\begin{itemize}
\item $f$ is a homeomorphism from $\C$ to $\C$,
\item the distribution derivative of $f$ is locally $L^2$,
\item the following holds almost everywhere:\footnote{The notation $\partial f$ and $\bar \partial f$ is explained below.}
\[\bd f(z) = \mu(z)\, \de f(z).\]
\end{itemize}
Such a map is called a \emph{quasiconformal homeomorphism of $\C$} and $\mu$ is called the \emph{Beltrami differential} of $f$. The map $f$ is also said to \emph{straighten} $\mu$. Sometimes, $\mu$ is also called a \emph{Beltrami form}.
\end{definition}

The second condition in the list means that the distribution partial derivatives $\de f/\de x$ and $\de f/\de y$ have locally $L^2$ representatives $f_x$ and $f_y$ and we let
\[df = f_x dx + f_y dy\]
which is a $1$-form (a current), but a priori not the differential of $f$ (see \Cref{rem:diff}).
There are too many good introductions and reference books for distributions to be cited here. Beyond the seminal \cite{Sc}, let us cite \cite{G} and \cite{LL}. Let us also cite \cite{St} for a pleasant introduction without going into technical details.
See \Cref{app:dl1} for information on distributions with locally $L^1$ or locally $L^2$ derivatives.

Bear in mind that we are asking that $\mu$ is in the unit ball of $\cal B$ and that this means that there exists $\eps>0$ such that $|\mu(z)|\leq 1-\eps$ for almost every $z$.

\begin{remark}
In his definition of a quasiconformal map, Definition B$'$ in %\cite{A} page~29
\cite{Abis} page~19,
Ahlfors only asks $df$ to be locally $L^1$.
He proves, and this is remarkable, %\cite{A} top of page~28,
(bottom of page~18) that $df$ is actually locally $L^2$, so \Cref{def:be} is equivalent.
However, we will not need that implication here and we will work directly with the $L^2$ version of the definition, which is also the most commonly used.
\end{remark}

In the third condition 
the complex numbers $\de f(z)$ and $\bd f(z)$ are defined by
\[df(z) = \de f(z)\cdot dz + \bd f(z)\cdot \overline{dz}\] 
where $\cdot$ is the multiplication of two complex numbers. This coincides with the classical operators $\de$ and $\bd$ if $f$ is $C^1$.
See \cite{Hu} or \cite{A,Abis} or \Cref{sub:prelim} for a geometric interpretation of the ratio $\bd f / \de f$ in terms of ellipse fields.
Here the above equation is to be taken for $L^2$ functions of $z$, so $\partial f$ and $\bar\partial f$ are complex valued $L^2$ functions.
The usual distribution derivatives of $f$ and the quantities above are related as follows:
\bEA
\pp{f}{x} &=& \phantom{i(}\de f + \bd f 
\\
\pp{f}{y} &=& i(\de f -\bd f)
\\
\de f &=& \frac{1}{2}\left(\pp{f}{x}-i\pp{f}{y}\right)
\\
\bd f &=& \frac{1}{2}\left(\pp{f}{x}+i\pp{f}{y}\right)
\eEA

\begin{remarkn}\label{rem:diff}
Note that we did \emph{not} claim that $f$ is differentiable almost everywhere. It turns out that this is true %(\cite{A}, Lemma~1 page~24).
(\cite{Abis}, Lemma~1 page~17).
We believe that this result is not needed to prove our main result, but that it is needed when one want to prove useful properties of quasiconformal maps in order to actually use them for deformation, surgeries, and other applications.
Note also that we did \emph{not} define $df$ as the differential of $f$ in the classical sense, but via the partial derivatives of $f$. Of course in the particular case where $f$ is $C^1$ then this coincides with the classical differential.
\end{remarkn}

\begin{theorem}[Measurable Riemann mapping theorem, existence]\label{thm:main1}
\hfill The Beltrami equation $\bd f = \mu\, \de f$ has a solution in the sense of \Cref{def:be}.
\end{theorem}

It is due to Gauss in the case $\mu$ is real-analytic, to Lavrentiev in the continuous case, to Morrey in the general case \cite{Mo}, with alternative proofs by Bers, Nirenberg, Boyarskii and others. These authors possibly use other formulations of the Beltrami equation. See the introduction to the present article and \cite{Glu} for more references.

\begin{definition}\label{def:nor}
A solution of the Beltrami equation is called \emph{normalized} if moreover
\begin{itemize}
\item $f(0)=0$ and $f(1)=1$.
\end{itemize}
\end{definition}
From any solution $f$, one easily defines a normalized solution
\[\tilde f(z) = \frac{f(z)-f(0)}{f(1)-f(0)}.\]

\begin{theorem}[Uniqueness, not proved here]\label{thm:mainuniq}
The normalized solutions of the equation $\bd f = \mu\, \de f$ are unique.
\end{theorem}
The proof of \Cref{thm:mainuniq} is out of the scope of the present article, see \cite{Abis}, Theorem~2 page~16 %\cite{A}, Theorem~2 page~23
together with the implication ``$B\implies A$'' page~20. 

\begin{definition}\label{def:holo}
We will say that a function $\mu: \tau\in\D \mapsto \mu_\tau\in \cal B$ is holomorphic if there is a power series expansion $\mu_\tau = \sum \tau^n a_n$, where $n\in\N$ and $a_n\in\cal B$, whose radius of convergence is at least one, i.e.\ $\limsup \|a_n\|^{1/n} \leq 1$.
\end{definition}
In this definition the sum is assumed to converge in the Banach space $\cal B$. For all $\tau\in\D$, $\mu_\tau$ is a well-defined element of $\cal B$, in particular it is an equivalence class of $L^\infty$ function from $\C$ to $\C$.
For the reader's interest, equivalent definitions of holomorphy are given in \Cref{app:weak}, but not used in this article, except in \Cref{sec:var-avg} where we study a variant of the construction.

\begin{theorem}[Holomorphic dependence]\label{thm:main2}
Under the holomorphy condition of \Cref{def:holo}, 
assuming moreover that for all $\tau\in\D$, $\|\mu_\tau\|_\infty<1$,
then there is a function $(\tau,z)\in\D\times\C \mapsto f_\tau(z)\in\C$ such that:
\begin{enumerate}
\item for every $\tau\in\D$, $f_\tau$ is a normalized solution of the Beltrami equation associated to $\mu_\tau$,
\item for every $z\in\C$, the function $\tau\in\D\mapsto f_\tau(z)\in\C$ is holomorphic.
\end{enumerate}
\end{theorem}

This is due to Ahlfors and Bers \cite{AhB}, and possibly others. They also studied other kind of smooth dependence. Of course by \Cref{thm:mainuniq}, $f_\tau$ is unique, so in fact the theorem can be stated as follows: the unique normalized solution $f_\tau$ is, pointwise in $z$, holomorphic in $\tau$.

\medskip

In this section we propose a proof of \Cref{thm:main1} and of \Cref{thm:main2} which uses the density argument and similarity surfaces.
A key point is that the obtained surfaces are isomorphic as Riemann surfaces to $\hat\C$ minus a finite set.
This is obtained by completing the atlas and making use of the Poincaré-Koebe theorem. 
In this approach where we cannot use the measurable Riemann mapping theorem (since we are proving it), we do not know how of an alternative method. 
The maps $f$ of \Cref{thm:main1} and $f_\tau$ of \Cref{thm:main2} will be obtained as extracted limits of maps $f_n$ straightening Beltrami forms that are piecewise constant with (finitely many) polygonal pieces.
The proof bears resemblance with Lavrentiev's approach \cite{Lav1,Lav2}, which we sum up in \Cref{sub:lavrentiev}.

\begin{remark}
Actually, using uniqueness (\Cref{thm:mainuniq}) one can prove that the full sequence $f_n$ converge to the solution of the Beltrami equation. See \Cref{app:beltrami:cv}.
\end{remark}

\begin{remark}
This hence defines an approximation scheme $f_n\tend f$ for solving the Beltrami equation.
However, we do not claim that this scheme is efficient.
\end{remark}

\section{Lavrentiev's approach}\label{sub:lavrentiev}

We sum up here \cite{Lav1}, whose translation in English can be found in \cite{Lav2}.

The Beltrami form $\mu$ is assumed \emph{continuous}.

A notion of $\eps$-near-solution to the Beltrami equation is defined with a parameter $\eps>0$, by asking that asymptotically near each point, small ellipses defined by $\mu$ are nearly mapped to small circles, up to $\eps$ (see the precise definition in Lavrentiev's article).
It is proved with relatively simple computations (including an equicontinuity statement on a class of $C^1$ diffeomorphisms that are quasiconformal) that if $\eps_n\tend 0$ then a sequence of $\eps_n$-near-solutions tends to a solution for $\mu$.

Then a sequence of near-solutions is constructed by
\begin{itemize}
\item taking simple $\eps$-near solutions defined on small squares: for instance by choosing $f: z\mapsto az+b\bar z$; this is one of the places where the continuity assumption on $\mu$ is used; 
\item proving that one can patch together these $\eps$-near solutions into a global $\eps$-near-solution; this is done by proving a sewing lemma (sewing two rectangles along an edge with a real and analytic sewing function) with an astute use of the Schwarz reflection principle.
\end{itemize}
In fact the second point could be proven by just invoking the Poincaré-Koebe uniformization theorem. What Lavrentiev does, besides the sewing lemma, is to argue (without giving details) that one can slightly modify the near-solutions before patching so that the sewing map is always analytic, while keeping a near-solution.

Moreover, if one chooses to take the $\R$-linear functions $a z+b\bar z$ on small squares, finding a global sewing is equivalent to uniformizing the surface $\cal S$ defined by gluing together linearly the many parallelograms that are the images of each square by the corresponding $\R$-linear function; though this is not what Lavrentiev does since he changes the maps slightly before each sewing.
The method that we are about to describe uses the uniformization of $\cal S$ and our proof of analytic dependence uses knowledge about its expression via the Schwarz-Christoffel formula.

\section{The use of similarity surfaces}\label{sub:sim}

Our proof of holomorphic dependence in the MRMT,  given in a subsequent section, is based on the use of an approximation scheme by a sequence of Beltrami forms $\mu_n$ that are piecewise constant on small square pieces.
The holomorphic dependence of the straightening of these approximations is a particular case of a more general statement that we prove here, based on the interpretation of the straightening of $\mu_n$ using conformal representation of affine surfaces glued from polygons, and using the results of \Cref{sub:hd:p}.

\medskip

An introduction to similarity surfaces built from polygons is given in \Cref{part:simsurf}.

\medskip

Assume that we are given a Beltrami form $\mu$ that is piecewise constant on $\C$, where the pieces are finitely many bounded or unbounded closed polygons $(Q_j)_{j\in J}$ without restriction, in particular they are not assumed convex nor simply connected, yet we assume that each has finitely many sides.
Each polygon comes with \emph{vertices} and its boundary is straight between vertices but we allow for consecutive edges to make a flat angle ($\pi$ radians), so it the subset $Q_j$ of $\C$ alone is not enough to determine the collection of its vertices.
We assume that if a point is a vertex for some $Q_j$ then it is also a vertex of any other polygon that contains it
(in a polygonal decomposition, there could be a corner of one polygon $Q_j$ belonging to the interior of the side of another $Q_k$; the required condition is easily achieved by just adding this point in the list of vertices of $Q_k$).
We make one more assumption: that no edge bounds the same polygon on its two sides.
This is easily achieved by removing any edge with this property.

To the unbounded edges $e$ we add marked points, and we mark the same point of $e$ for the two polygons on each side of $e$.

All our polygons are considered to be closed subsets of $\hat\C$.
We can associate to this a similarity surface as follows: for each polygon $Q_j$ choose an $\R$-affine map $a_j$ such that the Beltrami differential of $a_j$, which is constant, coincides with the value $\mu_j$ of $\mu$ on $Q_j$.
Now glue the polygons $P_j := a_j(Q_j)$ along their sides the same way that the polygons $Q_j$ were glued together, i.e.\ on each side with the affine map $a_k\circ a_{j}^{-1}$ for appropriate $(k,j)$ (the marked points are automatically matched: they are not needed here but later).

As in \Cref{sec:simsurfpolygon} we call $\cal S$ the topological surface thus obtained,
\[ \pi: \coprod P_j = \coprod a_j(Q_j)\to\cal S \]
the quotient map, $\cal V\subset\cal S$ the the set of vertices, and $\cal S'= \cal S-\cal V$.
To avoid the problem where two polygons $P_i$ and $P_k$ intersect as subset of $\C$ and $z$ belongs to this intersection, we use the notation $(z:i)$ and $(z:k)$ to distinguish them in the disjoint union $\coprod P_j$. This is particularly useful since it is likely that $\pi(z:i)$ and $\pi(z:k)$ are different.
 
\begin{figure}
\begin{tikzpicture}
\node at (-0.2,0) {\includegraphics[scale=0.5]{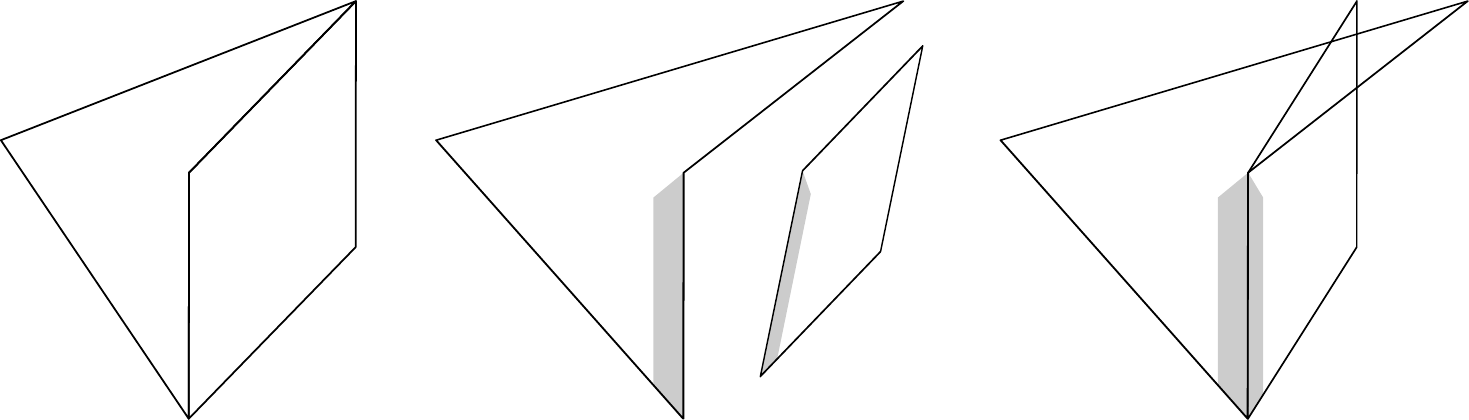}};
\node at (-5.4,0.3) {$Q_j$};
\node at (-3.9,0.3) {$Q_k$};
\node at (-1.5,0.3) {$P_j$};
\draw[-] (-2.2,-0.7) node[anchor=east]{$V_j$} -- (-0.8,-0.7) ;
\node at (-0.42,-0.7) {$e$};
\node at (0.85,0.2) {$P_k$};
\draw[-] (1.4,-0.7) node[anchor=west]{$V_k$} -- (0.25,-0.7) ;
\node at (3.4,0.3) {$P_j$};
\draw[-] (4.7,0.0) -- (5.3,-1.3) node[anchor=north]{$s^{-1}(P_k)$};
\end{tikzpicture}
\caption{Left: two polygons $Q_j$, $Q_k$. Middle: images $P_j$ and $P_k$ of $Q_j$ and $Q_k$ by the $\R$-affine maps $a_j$ and $a_k$, selection of an edge $e$. Right: the $\C$-affine map $s^{-1}$ allows to glue $P_k$ to $P_j$ along $e$; a neighbourhood of $e$ in shown in gray, that we use to provide a chart.}
\label{fig:gl2}
\end{figure}

Let us recall how we defined in \Cref{sec:simsurfpolygon} a similarity surface atlas on $\cal S'$ (see the part ``Removing vertices and getting a similarity surface''; there are slight differences due to the fact that we allow here more kinds of polygons than there, but the construction is equivalent).
On the interior of $\pi(P_j)$ we just use the identity.
Now consider an edge $e$ of $P_j$ and let $P_k$ be the polygon glued to $P_j$ along $e$ by $a_k\circ a_{j}^{-1}$. Denote $e^*$ the edge minus its endpoints.
The map $a_k\circ a_{j}^{-1}$ is a priori only $\R$-affine on $\C$.
There is a unique $\C$-affine map $s$ that extends the restriction of $a_k\circ a_{j}^{-1}$ to $e$.
Then a neighbourhood $V$ of $\pi(e^*:j)$ was identified,
of the form $V=\pi(V_k:k) \cup \pi(V_j:j)$ for some open sets $V_k\subset P_k$ and $V_j\subset P_j$ that intersect $\partial P_k$ resp.\ $\partial P_j$ exactly on the corresponding edges minus endpoints.
In \Cref{sec:simsurfpolygon} the polygons were convex so one could just take $V_j$ to be the interior of $P_j$ union $e^*$ and similarly for $V_k$. Here we will have to take smaller subsets.
A chart $\phi:V\to\C$ can be defined by $\phi(\pi(z:j)) = z$ if $z\in V_j$, and $\phi(\pi(z:k)) = s^{-1}(z)$ if $z\in V_k$, provided $V_j$ and $V_k$ are chosen so that $\phi(\pi(V_j:j))\cap \phi(\pi(V_k:k)) = \phi(\pi(e^*:j))$, which is possible.
Note that its image $\phi(V)$ is a neighbourhood of $e^*$ in $\C$.
See \Cref{fig:gl2}.

This atlas on $\cal S'$ is in particular a Riemann surface atlas. Using \Cref{sec:simsurfpolygon} we can see that all points in $\cal V$ are punctures, including $\infty$, which allows to extend the Riemann atlas to $\cal S$.
Indeed, near finite vertices of $Q_j$ we get what is called a bounded vertex in \Cref{sec:simsurfpolygon}, and all bounded vertices are punctures according to the part called ``Conformal erasability of the singularities''.
Near $\infty$, we have a positive total angle (negative signed angle), hence we also get a puncture as explained near the end of the part called ``Unbounded polygons''.
Since $\cal S$ is homeomorphic to a sphere, the \emph{Poincaré-Koebe theorem} implies that it is conformally isomorphic to the Riemann sphere $\hat\C$.
Call
\[F :\cal S\to\hat\C\]
such an isomorphism, chosen so that $\infty$ is mapped to $\infty$.
Let now $f:\C=\bigcup Q_j\to\C$ be defined as follows: on each polygon $Q_j$, $f(z)=F\circ\pi(a_j(z):j)$. This definition coincides on the shared edges and vertices.

\begin{lemma}\label{lem:sssol}
The map $f$ thus constructed is a solution of the Beltrami equation associated to $\mu$ in the sense of \Cref{def:be}.
\end{lemma}
\begin{proof}
See \Cref{app:pf:lem:sssol}.
\end{proof}

We can normalize $f$ (i.e.\ we can impose that $f(0)=0$ and $f(1)=0$) by choosing the isomorphism $F$ above appropriately.

\begin{lemma}\label{lem:holdep}
Fix the polygons $Q_j$ and assume that no unbounded polygon has an infinite vertex with angle $0$.
Then the normalized solution $f$ constructed above depends holomorphically on the family of values $(\mu_j)\in\D^J$ that $\mu$ takes on the $Q_j$ in the following sense: for all $z\in\C$, $\mu\mapsto f_\mu(z)$ is analytic.
\end{lemma}
\begin{proof}
It is enough to prove the claim locally, in the sense that $\mu$ can be restricted to be in a small neighbourhood of a given $\mu^0$.

We can refine the polygons so that the conditions of \Cref{sub:hd:p} are satisfied, which we recall here:
\begin{enumerate}
\item All polygons are strictly convex,
\item every unbounded polygon has exactly two unbounded edges, and their angle is $>0$,
\item $\cal S$ has at least three vertices.
\end{enumerate}

Let us pick three vertices $w_1=\infty$ and $w_2$ and $w_3$ in the initial polygonal decomposition $\bigcup Q_j$ of the plane $\C$ where the Beltrami form $\mu$ lives.
Number $w_4$ to $w_p$ the remaining ones.
Let $w'_1$ to $w'_q$ be the marked points and let $m=p+q$.
The points $w_k$ and $w'_k$ are independent of $\mu$.
Let $v_k$ and $v'_k$ be the corresponding points in $\cal S$.
They depend on $\mu$ since $\cal S$ does.
This gives a consistent labelling of the set of vertices $\cal V$ of $\cal S$.
The isomorphism $F  :\cal S\to \hat\C$ is normalized so that $f=F\circ\pi\circ \coprod_j a_j$ fixes $0$, $1$ and $\infty$.
We choose another isomorphism $\wh F:\cal S\to\hat\C$ by post-composing $F$ with a homography (fixing $\infty$, so actually an affine map): $\wh F = h\circ F$,
so that $\wh F$ sends $v_1$ to $\infty$, $v_2$ to $0$ and $v_3$ to $1$, or equivalently that $\hat f := \wh F \circ \pi\circ  \coprod_j a_j= h\circ f$ sends $w_1=\infty$ to $\infty$, $w_2$ to $0$ and $w_3$ to $1$.

It is enough to prove that:
\begin{itemize}
\item[(A)] For all $z\in\C$, $\hat f(z)$ depends holomorphically on $\mu$.
\end{itemize}
Indeed, $h(0)=\hat f(0)$ and $h(1)=\hat f(1)$ then vary holomorphically so the affine map $h$ and its inverse $h^{-1}$ have their coefficients that depend holomorphically on $\mu$, hence $f(z) = h^{-1} \circ \hat f(z)$ depends holomorphically on $\mu$ for a fixed $z$.

Let us prove claim (A). First, by \Cref{prop:hsol}, the points $\hat z_k=\hat f(w_k)$ depend holomorphically on $\mu$. For the other points $w\in\C-\{w_1,\ldots, w_p\}$, we use a trick:\footnote{We would like the anonymous referee for suggesting an initial simplification of the original proof that we pushed further to the present form.} we add $w$ to the set of vertices. For this, consider a finite polygonal decomposition of the polygons $Q_j$ containing $w$ (there are two such polygons if $w$ is on the side of a polygon $Q_j$, otherwise there is only one), and such that $w$ is a vertex.
There may be more than one vertex added, so we denote
\[w_{p+1} = w.\]

The key observation is that the corresponding map $\hat f_{\mathrm{new}}$ obtained using new polygonal decomposition is the same as the original map $\hat f_{\mathrm{orig}}$.
This is for instance an immediate consequence of the uniqueness of the solution of the Beltrami equation.
But it can also be deduced without this theorem in our particular situation.
Indeed the Riemann surface we build from gluing the new polygonal arrangement $P_j=a_j(Q_j)$ is naturally isomorphic to the Riemann surface before subdividing the polygons.
This translates into a conformal isomorphism $\phi:\hat\C\to\hat\C$ such that $\hat f_{\mathrm{new}} = \phi \circ \hat f_{\mathrm{orig}}$ that fixes in particular every $z_k$, hence it fixes $0$, $1$ and $\infty$ and such isomorphism is necessarily the identity.

Now the point $\hat f(w)=\hat f(w_{p+1})=\hat z_{p+1}$ is part of the set of singularities (it is an erasable singularity) for which \Cref{prop:hsol} proves holomorphic dependence in $\mu$.
\end{proof}

\begin{remark}
The conclusion of \Cref{lem:holdep} also holds in the presence of one or several unbounded polygons with zero angle at infinity, as will follow from the MRMT with parameter (\Cref{thm:main2}). But our approach required us to exclude this case.
If one wants to prove this without using \Cref{thm:main2} one can proceed by perturbing the polygons $Q_j$ so that all angles are strictly positive then taking limits of the solutions as some angles tend to $0$, and conclude by the fact that limits of holomorphic functions are holomorphic.
One has to prove bounds to justify that limits can be taken and that they are homeomorphisms straightening the intended Beltrami form.
Due to the simple form of the situation this should be doable by more direct and simpler methods than the proof of \Cref{thm:main2} (in particular no need to use \Cref{lem:nor,lem:L2}, the use of Koebe's distorsion theorems should suffice together with the explicit form of the maps $a_j$).
However, we will not develop this here.
\end{remark}

\section{The density argument, without parameter}\label{sub:MRMT}

We start by proving \Cref{thm:main1}, i.e.\ the MRMT without the parameter $\tau$, and explain in \Cref{sub:holodep} how to adapt the proof for holomorphic families $\mu_\tau$.

\begin{remark}
Density arguments like the one presented here have been known for a while.
A nice example is given in \cite{Hu}, Section~4.6, which presents a proof of the measurable Riemann mapping theorem, and Theorem~4.7.2.
Our presentation here is quite similar, with more details.
A slightly different type of density argument is performed in \cite{Lav1} (translated in \cite{Lav2}).
\end{remark}

\subsection{Approximation}

We are given some $\mu \in L^\infty(\C)$ with $\|\mu\|_\infty<1$.
We recall that the norm is the essential supremum of the representatives of $\mu$. 
Let $\mu_n$ be the following approximation of $\mu$:
\begin{itemize}
\item divide the square defined by $|\Re z|<n$, $|\Im z|<n$ into small squares of side $1/n$ (there will be $4n^4$ of them); let $\mu_n$ be constant in the interior of each of those squares, equal to the average\footnote{About the choice of averaging, see \Cref{sec:var-avg}.} of $\mu$ on this square.
\item elsewhere let $\mu_n(z)=0$.
\end{itemize}
Note that
\[\|\mu_n\|_\infty\leq \|\mu\|_\infty.\]

\begin{lemma}\label{lem:mnm}
For any compact subset $S$ of $\C$,
\[\|\mu_n-\mu\|_{L^1(S)} \tend 0.\]
\end{lemma}
\begin{proof}
We use here a classical argument.
We may as well suppose that $S$ is the square ``$|\Im z|<N$, $|\Re z|<N$'' with $N\in\N^*$.
The map $\Theta_n : \mu\mapsto\mu_n$ is linear and satisfies $\|\Theta_n(\mu)\|_{L^1(S)}\leq\|\mu\|_{L^1(S)}$.
Assume that for all $\eps>0$ one can find $\nu \in L^\infty(S)$ such that $\|\nu-\mu\|_{L^1(S)}\leq \eps$ and such that $\|\Theta_n(\nu)-\nu\|_{L^1(S)}\tend 0$ when $n\tend +\infty$.
Then $\|\Theta_n(\mu)-\mu\|_{L^1(S)} \leq \|\Theta_n(\mu)-\Theta_n(\nu)\|_{L^1(S)} + \|\Theta_n(\nu)-\nu\|_{L^1(S)}+\|\nu-\mu\|_{L^1(S)} \leq 3\eps$ when $n$ is big enough.

For $\nu$, we can for instance use that continuous maps are a dense subset of $L^1(S)$, and that the claim is easy to check for continuous maps because we have uniform convergence of $\Theta_n(\nu)$ to $\nu$ on $S$ by uniform continuity.

As an alternative choice, one can just check the claim for $\nu$ in the set of characteristic functions of product of intervals in $R^2 \cong \C$, then this holds for their linear combinations, and these linear combinations are also dense in $L^1(S)$.
\end{proof}

Usually local $L^2$ convergence is stronger than local $L^1$ convergence but here it is equivalent because the functions involved have uniformly bounded $L^\infty$ norm:

\begin{corollary}\label{cor:mnm}
For any compact subset $S$ of $\C$,
\[\|\mu_n-\mu\|_{L^2(S)} \tend 0.\]
\end{corollary}
\begin{proof}
Almost everywhere $|\mu|\leq 1$, so $|\mu_n|\leq 1$, and almost everywhere 
$|\mu_n-\mu|^2 \leq 2 |\mu_n-\mu|$.
\end{proof}

\begin{remark}
It is true that $\mu_n$ tends to $\mu$ for the weak topology (against continuous maps with compact support), and the proof is easy (uniform continuity of the test function). It also follows from local $L^1$ convergence. However, \Cref{app:contrex} gives an example that shows that weak convergence is not enough to extract a correct solution for $\mu$ from a sequence of solutions for $\mu_n$.
\end{remark}

Let $f_n$ be the normalized solution of the Beltrami equation for $\mu_n$ constructed in \Cref{sub:sim}.

\subsection{Compactness statements}

As announced in the introduction, we do not include here the proof of the following two statements, but indicate where they can be found in \cite{Abis} and \cite{Hu}.

To match the terminology of references we say that a quasiconformal map $f$ with $\|\mu\|_\infty\leq \kappa$ is $K$-quasiconformal with $K = \frac{1+\kappa}{1-\kappa}$.
The justification of this terminology comes from the geometric interpretations of quasiconformality, see \cref{sub:prelim} in the present article and \cite{A,Abis} for complements.

\begin{lemma}[Normal family]\label{lem:nor}
For a fixed real $K>1$, the set of normalized $K$-quasiconformal homeomorphisms from $\C$ to $\C$ forms a normal family in the sense below.
\end{lemma}
A normal family is a family such that from every sequence in this family, one can extract a subsequence that converges uniformly on every compact subset of $\C$.
In other words it is a compact family for the notion of uniform convergence on compact subsets of $\C$.
The lemma above is Theorem~2 page~33 in Chapter~III of \cite{Abis}. %page 51 in Chapter~III of \cite{A}.
For the closely related case of quasiconformal maps from $\D$ to $\D$ or between hyperbolic Riemann surfaces, see \cite{Hu} Section~4.4.

\begin{lemma}[$L^2$ bound for $df$]\label{lem:L2}
For any $K$-quasiconformal map $f$ on an open subset $U$ of $\C$,
\bEA
&& \int_U \left|\frac{\de f}{\de x}\right|^2 \leq K \Leb f(U),
\\
&& \int_U \left|\frac{\de f}{\de y}\right|^2 \leq K \Leb f(U).
\eEA
\end{lemma}
This is a direct consequence of \cite{Abis}, Chapter~II, last line of page~18 %page~27 
and the formulas $\de f/\de x = \de f + \bd f$ and $\de f/\de y = i(\de f - \bd f)$.
An alternative method can be found in \cite{Abis}, via Theorem~3 page~22 %page~33
(which uses Lemma~2 page~19 %page~28 
and the assumption that $f$ has $L^2$ distributional partial derivatives).
See also in \cite{Hu}, Proposition~4.2.4 and Corollary~4.2.6.

\subsection{Extraction of a limit}

From \Cref{lem:nor,lem:L2} it follows that we can extract a subsequence of $f_n$ (more generally from any subsequence of $f_n$) such that:
\begin{enumerate}
\item $f_n$ uniformly converges on every compact subset of $\C$ to a $K$-quasiconformal map $f$,
\item\label{item:wkv} $u_n=\de f_n/\de x$ and $v_n=\de f_n/\de y$ converge to some locally $L^2$ functions $u$ and $v$, for the weak convergence defined as follows: for all $L^2$ function $\phi$ with compact support, $\int (u_n-u)\phi\tend 0$ and $\int (v_n-v)\phi\tend 0$. 
\end{enumerate}
The first claim is an immediate consequence of \Cref{lem:nor} and the second claim is obtained by a further extraction justified as follows: for any (compact) rectangle $R\subset U$, by the first claim, the area of $f_n(R)$ is a bounded sequence. Hence the $L^2$ norm on $R$ of $\de f_n/\de x$ and $\de f_n/\de y$ are bounded sequences too by \Cref{lem:L2}. 
By the weak-star compactness theorem for $L^2(R)$ we can extract a subsequence such that the partial derivatives converge weakly in $L^2(R)$. By varying the rectangle and performing a diagonal extraction, the second claim follows.

\subsection{Final checks}

Let us first check that $u$ and $v$ are the distribution derivatives of the extracted limit $f$. For any test function $\phi \in C^\infty_c(\C)$, denote by $S$ its support
\[
\la\frac{\de f_n}{\de x},\phi\ra
= \la -f_n, \frac{\de \phi}{\de x}\ra
= \int_{S} -f_n \frac{\de \phi}{\de x}.
\]
On one hand by uniform convergence of $f_n$ to $f$ on $S$, as $n\to+\infty$
\[
\int_{S} -f_n \frac{\de \phi}{\de x}
\tend 
\int_{S} -f \frac{\de \phi}{\de x}
= D_x f(\phi)
\]
where $D_x f$ denotes the distribution partial derivative of $f$.
Hence
\[\la\frac{\de f_n}{\de x},\phi\ra \tend D_x f(\phi)
.\]
On the other hand by weak convergence on $S$ of $\de f_n/\de x$ to $u$:
\[
\la\frac{\de f_n}{\de x},\phi\ra =
\int_{S} \frac{\de f_n}{\de x} \,\phi
\tend \int_{S} u \,\phi.
\]
So we proved that for all test functions $\phi$,
\[D_x f(\phi)
 = \int_{S} u \,\phi.\]
In other words the distribution derivative $\de f/\de x$ is represented by the locally $L^2$ (hence locally $L^1$) function $u$. 
The proof for $\de f/\de y$ and $v$ is similar.

There remains to check that $f$ is a solution of the Beltrami equation associated to $\mu$.
Let $\phi\in C^\infty_c(\C)$. We know that
\[\la \de f_n,\phi\ra = \la \mu_n\bd f_n,\phi\ra\]
where by $\partial f_n$, $\bar\partial f_n$ we refer to the $L^2$ functions.
The fact that
\[\la \de f_n,\phi\ra \tend \la \de f,\phi\ra\]
follows from the formulas $\de f_n = (u_n-iv_n)/2$, $\de f = (u-iv)/2$ and the weak convergence of $u_n$, $v_n$ to $u$, $v$. We will be done if we can prove that
$\la \mu_n\bd f_n,\phi\ra\tend
\la \mu\,\bd f,\phi\ra$.
For this let us write
\[\la \mu_n\bd f_n - \mu\,\bd f,\phi\ra
=
\underbrace{\la (\mu_n-\mu) \bd f_n, \phi \ra}_{\ds A}
+
\underbrace{\la (\bd f_n- \bd f) \mu,\phi\ra}_{\ds B}
.\]
Let $S$ be the support of $\phi$.

The function $\mu\phi$ is $L^\infty$ with support contained in $S$, hence $L^2$ with support contained in $S$, so by the local weak $L^2$ convergence, we have $B\tend 0$.

Concerning
\[A = \int_S (\mu - \mu_n) \bd f_n \phi\]
by the Cauchy-Schwarz inequality
\[|A| \leq \left \|\mu-\mu_n\right\|_{L^2(S)} \times \|\bd f_n \phi\|_{L^2(S)}.\]
Now $\|\mu-\mu_n\|_{L^2(S)}\tend 0$ according to \Cref{cor:mnm}, 
$\|\bd f_n \phi\|_{L^2(S)} \leq \max|\phi| \times \|\bd f_n\|_{L^2(S)} $ and we have seen that $\|\bd f_n\|_{L^2(S)}$ is a bounded sequence (\Cref{lem:L2}).
It follows that $A\tend 0$ as $n\to+\infty$.

Note that if we only had extracted a weakly convergent sequence in $L^1_\loc$, we could still argue that $B\tend 0$ but not $A$, because instead of Cauchy-Schwarz we would use $|A| \leq \left \|\mu-\mu_n\right\|_{L^\infty(S)} \times \|\bd f_n \phi\|_{L^1(S)}$ but
$\|\mu-\mu_n\|_{L^\infty(S)}$ is unlikely to tend to $0$.

\subsection{Uniqueness and convergence}\label{app:beltrami:cv}

It is not the point in the present article to prove uniqueness of the solution of the Beltrami equation (\Cref{thm:mainuniq}), so we admit it (see the references below the statement of that theorem).
Let us explain how one deduces from it that the full sequence $f_n$ tends to this unique solution $f$ (uniformly on compact subsets of $\C$).

Since $\|\mu_n\|_\infty \leq \|\mu\|_\infty$, the sequence $f_n$ is $K$-quasiconformal for some common $K$.
By \Cref{lem:nor}, from any subsequence of $f_n$ one can extract a  sub-subsequence that converges uniformly on every compact subset of $\C$ to some function $f$.
The procedure of \Cref{sub:MRMT} can then be applied to this sub-subsequence, which prove that $f$ is a solution of the Beltrami equation.
Since the solution is unique, this means that all extracted limits $f$ of the sequence $f_n$ are identical.
Now in a Hausdorff-separated topological space, if a sequence $f_n$ has all its subsequences that contain sub-subsequence converging to $f$, then the full sequence $f_n$ tends to $f$.

\section{Holomorphic dependence}\label{sub:holodep}

We now explain how to adapt the proofs in \Cref{sub:MRMT} to get \Cref{thm:main2}, i.e.\ the MRMT with holomorphic dependence (Ahlfors-Bers theorem).

%We denote $\|\mu\|_\infty$ 
We start from the hypothesis that $\forall \tau\in\D$, $\| \mu_\tau\|_\infty<1$ and improve this inequality as follows:
\begin{lemma}\label{lem:muk}
$\forall \eps\in(0,1)$, $\exists \kappa(\eps)<1$ such that
\[\forall \tau\in B(0,1-\eps),\ \|\mu_\tau\|_\infty \leq \kappa(\eps).\]  
\end{lemma}
\begin{proof}
By hypothesis $\mu_\tau=\sum \tau^k a_k$ with $\limsup \|a_k\|_\infty^{1/k}\leq 1$.
It follows that $\tau\mapsto \|\mu_\tau\|_\infty$ is continuous on $\D$.
It thus reaches a maximum on $\ov B(0,1-\eps)$, and at this point its value is $<1$ by hypothesis.
\end{proof}

Let the Beltrami form $\mu_n(\tau)$ be defined from $\mu(\tau)$ exactly as in \Cref{sub:MRMT}.
Denote the square $S_n:$ ``$|\Re z|<n$, $|\Im z|<n$'' and call \emph{little square} the squares of side $1/n$ into which we have cut it.
In a little square $S'$ we have 
\[\mu_n(\tau) = \frac{\int_{S'} \mu_\tau}{\Leb S'} = \int_{S'} \sum _{k=0}^{+\infty} \tau^k \frac{a_k}{\Leb S'}= \sum_{k=0}^{+\infty} \tau^k \frac{\int_{S'} a_k}{ \Leb S' }\]
and the convergence of the series is normal\footnote{A series $\sum f_n$ of functions is said to converge normally if $\sum \|f_n\|_\infty<+\infty$.} over $\tau\in B(0,1-\eps)$ for all $\eps>0$ since $\limsup\|a_k\|_\infty^{1/k}\leq 1$.
This quantity is a holomorphic function of $\tau\in\D$ for all little squares $S'$.

By \Cref{lem:holdep} there are normalized solutions $z\mapsto f_n(\tau,z)$ of the Beltrami equation for $\mu_n(\tau)$, which vary holomorphically with $\tau$ for all fixed $n,z$.

% \begin{corollary}\label{cor:eqc}
% For all $z\in\C$, the sequence of holomorphic functions $t\in \D\mapsto f_n(t,z)$ is equicontinuous.
% \end{corollary}
% \begin{proof}
% Fix $\eps\in (0,1)$.
% For $t\in B(0,1-\eps)$, we have $\|\mu_n(t)\|_\infty \leq \|\mu_t\|_\infty\leq \kappa(\eps)<1$
% by \Cref{lem:muk}, so for all $t\in B(0,1-\eps)$ and all $n\in\N$, the functions $z\mapsto f_n(t,z)$ are $K$-quasiconformal for the same $K=(1+\kappa(\eps))/(1-\kappa(\eps))$.
% By \Cref{lem:nor} they form a normal family with respect to the variable $z$.
% In particular for any compact subset $C$ of $\C$, there is a bound $M>0$ such that $\forall n\in \N$, $\forall t\in B(0,1-|\eps|)$, and $\forall z\in C$, $|f_n(t,z)|\leq M$.
% In particular, for all $z\in\C$, the sequence of holomorphic functions $t\in B(0,1-\eps)\mapsto f_n(t,z)$ is bounded, hence equicontinuous.
% Since this is valid for all $\eps\in(0,1)$, we get the result.
% \end{proof}

\begin{corollary}\label{cor:eqc2}
The functions $f_n: (\tau,z)\in \D\times \C\mapsto f_n(\tau,z)$ are continous and this sequence of functions is normal.
\end{corollary}
\begin{proof}
Fix $\eps\in (0,1)$.
For $\tau\in \ov B(0,1-\eps)$, we have $\|\mu_n(\tau)\|_\infty \leq \|\mu_\tau\|_\infty\leq \kappa(\eps)<1$
by \Cref{lem:muk}, so for all $\tau\in \ov B(0,1-\eps)$ and all $n\in\N$, the homeomorphisms $F_{n,\tau} : z\in\C \mapsto f_n(\tau,z)$ are $K$-quasiconformal for the same $K=(1+\kappa(\eps))/(1-\kappa(\eps))$.
They are also normalized so by \Cref{lem:nor} they form a normal family with respect to the variable $z$.
In particular for any compact subset $C$ of $\C$, there is a bound $M>0$ such that $\forall n\in \N$, $\forall \tau\in \ov B(0,1-|\eps|)$, and $\forall z\in C$, $|f_n(\tau,z)|\leq M$.
In particular, the collection of holomorphic functions $G_{n,z}: \tau\in B(0,1-\eps)\mapsto f_n(\tau,z)$ for $n\in\N$ and $z\in C$ is bounded, hence uniformly equicontinuous on $\ov B(0,1-\eps')$ for any $\eps'>\eps$ by classical Cauchy estimates.
Since the maps $F_{n,\tau}$ are uniformly equicontinuous on $C$ for $\tau\in \ov B(0,1-\eps)$ hence in particular for $\tau\in \ov B(0,1-\eps')$, it follows that $f_n: (\tau,z)\in \ov B(0,1-\eps') \times C \to f_n (\tau,z)$ is a uniformly equicontinuous sequence of functions.
We saw that it is bounded too, so the sequence $f_n$ is normal on $\ov B(0,1-\eps') \times C$ by the Arzelà-Ascoli theorem.
Since this holds for any compact $C\subset\C$ and any pair $0<\eps<\eps'<1$, the claim follows.
\end{proof}

We will give two slightly different ways to conclude.
The first one avoids using uniqueness of the normalized solution of the Beltrami equation (\Cref{thm:mainuniq}).

\loctitle{Method~1.}

By \Cref{cor:eqc2}, one can extract a subsequence of $f_n$ that converges uniformly on every compact subset of $\D\times\C$.
Let $f$ denote this limit.
Since the uniform limit of holomorphic functions is holomorphic, the function $\tau\mapsto f(\tau,z)$ is holomorphic for each fixed $z$.

By \Cref{sub:MRMT} for each fixed $\tau$, there is a sub-subsequence such that $F_{n,\tau} :z\mapsto f_n(\tau,z)$ tends uniformly on compact subsets of $\C$ to a solution of the Beltrami equation for $\mu_\tau$.
It follows that $z\mapsto f(\tau,z)$ is actually a solution of the Beltrami equation for $\mu_t$

\loctitle{Method~2.}

By \Cref{app:beltrami:cv},
for each fixed $\tau$, the whole sequence (without extraction) of functions $F_{n,\tau}: z\mapsto f_n(t,z)$ converges uniformly on compact subset of $\C$ to a normalized solution that we denote $z\mapsto f(t,z)$ of the Beltrami equation for $\mu_t$.
The function $(t,z)\mapsto f(t,z)$ in is in particular the simple limit of the functions $(t,z)\mapsto f_n(t,z)$.
For each $z$, the function $t\in\D \mapsto f_n(t,z)$ is holomorphic.
The function $t\mapsto f(t,z)$ is hence the simple limit of a sequence of holomorphic functions, and this sequence is equicontinuous by \Cref{cor:eqc2}. Hence the convergence is uniform on compact subsets of $\D$, and the limit is holomorphic.

\medskip

\noindent\textbf{Note}: See \Cref{sec:var-avg} for a study of what happens if we average $\mu$ through a (reasonable) function of $\mu$ instead of directly taking $\mu_n = \int_S\mu /\int_S 1$ on each little square $S$.

%
% ----------------------------------------------------------
%

\section{Appendix}

%
% ----------------------------------------------------------
%

This appendix collects complements for the present part: several equivalent definitions of being holomorphic for functions from the unit disk to a Banach space; a counter-example to continuity of the straightening if the only the weak topology is used on the data $\mu$; a collection of facts about distribution with $L^1_\loc$ derivatives.
It also includes the postponed proof of the fact (\Cref{lem:sssol}) that for Beltrami forms that are piecewise constant on (finitely many) polygonal pieces, the construction proposed in \Cref{sub:sim} indeed yields a solution of the Beltrami equation.

\subsection{Equivalent formulations of the holomorphy condition}\label{app:weak}

In this section we recall equivalent definitions of being holomorphic, for functions from the unit disk to a Banach space.
The statements in the present section are not used in the proof of the main results of the article, we only use \Cref{prop:muhol} in the proof of \Cref{lem:fhmh} in \Cref{sec:var-avg} where we study variants of the construction.

\begin{proposition}\label{prop:holos}
    Let $E$ be a Banach space, and $f:\D\to E$ a function.
    Then the following are equivalent.
    \begin{enumerate}
        \item\label{bw1} (power series expansion) There is a power series expansion $f(\tau) = \sum \tau^n a_n$ over $n\in\N$ where $a_n\in E$, whose radius of convergence is at least one, i.e.\ $\limsup \|a_n\|^{1/n} \leq 1$. 
        \item\label{bw2} (differentiable) The function $f$ is differentiable, i.e.\ for all $t\in\D$, $\frac{1}{h}(f(\tau+h) - f(\tau))$ has a limit in $E$ as the complex number $h\to 0$.
        \item\label{bw3} (weak) Denote $E'$ the dual space to $E$. For any $\phi\in E'$, the function $\phi\circ f : \D\to\C$ is holomorphic.  
        \item\label{bw4} (weak*) Assume that $E$ is the dual to some Banach space $V$. For any $v\in V$, the function $\tau\in\D \mapsto f(\tau)(v) \in\C$ is holomorphic.  
        \item\label{bw5}(very weak) The function $f$ is bounded and there is a separating subset $X\subset E'$ (separating means $\forall e\in E-\{0\}$, $\exists \phi\in X$ $\phi(e)\neq 0$) such that $\forall \phi\in X$, the function $\phi\circ f$ is holomorphic.  
    \end{enumerate}
  In either case, $f$ is called holomorphic.
\end{proposition}

The implications \eqref{bw1} $\implies$ \eqref{bw2} $\implies$ \eqref{bw3} $\implies$ \eqref{bw4} and \eqref{bw3} $\implies$ \eqref{bw5} are immediate or easy.
That \eqref{bw3} implies \eqref{bw1} is proved for instance in \cite{Muj}, Lemma~8.13. In fact the proof adapts to give \eqref{bw4} $\implies$ \eqref{bw1} because the principle of uniform boundedness can still be applied.
That \eqref{bw5} implies \eqref{bw1} is proved in \cite{AN}, Theorem~3.1 and is an extension of a result whose proof can be found in \cite{Ka} (Remark~1.38, page~139) for which $X$ is assumed with dense span in $E'$.

Let us specialize this to the Banach space $\cal B$ of complex-valued $L^\infty$ functions on $\C$, endowed with the essential supremum norm, which we will denote $\|\cdot\|$, and add one characterization. First note that \Cref{def:holo} matches definition~\eqref{bw1} above.

\begin{proposition}\label{prop:muhol}
Let $\mu : t\in\D\mapsto\mu_t\in\cal B$ be a function.

Then $\tau\mapsto \mu_\tau$ being holomorphic in any of the senses of \Cref{prop:holos} is equivalent to any of the following criteria:

\begin{enumerate}
\item\label{aw4} For every $L^1$ function $\phi:\C\to\C$, the function $\tau\in\D\mapsto \int_\C \mu_\tau\phi$ is holomorphic.  
\item\label{aw5} The function $\tau\mapsto \|\mu_\tau\|_\infty$ is locally bounded and for every $C^\infty$ function with compact support $\phi:\C\to\C$, the function $\tau\in\D\mapsto \int_\C \mu_\tau\phi$ is holomorphic.  
\item\label{aw1} The function $t\mapsto \|\mu_\tau\|_\infty$ is locally bounded and there is a function $(\tau,z)\mapsto\mu_\tau(z)\in\C$ defined for $\tau\in \D$ and $z\in\C$ such that for all $\tau\in\D$, $z\mapsto \mu_\tau(z)$ is a measurable function and is a representative of $\mu_\tau$ and for all $z\in\C$, the function $\tau\mapsto \mu_\tau(z)$ is holomorphic.
\end{enumerate}
\end{proposition}

\begin{proof}
Point \eqref{aw4} is the criterion~\eqref{bw4} of \Cref{prop:holos} specialized to $L^\infty(\C)$ being the dual space to $L^1(\C)$.
    
Point \eqref{aw5} is the criterion~\eqref{bw5} of \Cref{prop:holos} with the set of smooth function with compact suport known to be separating.
    
Let us prove that if criterion~\eqref{bw1} of \Cref{prop:holos} holds, then \eqref{aw1} holds.
Indeed, each $a_n\in L^\infty(\C)$ has a representative that we denote $z\mapsto A_n(z)$. This representative could take values of modulus greater than $\|a_n\|$. It can only happen on a set of measure $0$ and we modify $A_n$ to take value $0$ on this set, which does not change its class in $L^\infty$. Then for \emph{all} $z\in\C$, the power series $\sum A_n(z) \tau^n$ converges on $\D$, we call $\mu_\tau(z)$ its value. Let us check that $z\mapsto\mu_\tau(z)$ is a representative of $\mu_\tau$.
By the criterion~\eqref{bw1} of \Cref{prop:holos}, $\mu_\tau$ satisfies
\[\|\mu_\tau - \sum_{k=0}^n a_k \tau^k\|\underset{n\to\infty}\tend 0\]
and this means $\mu_\tau$ has a representative for which the partial sums $z\mapsto \sum_{k=0}^n A_k(z) \tau^k$ tend to it almost everywhere (if a sequence of functions has its essential sup tending to $0$ then it tends to $0$ almost everywhere). But the sum tends to $\mu_\tau(z)$ everywhere, which is hence representative of $\mu_\tau$.

Let us now prove that if \eqref{aw1} holds, then $\tau\mapsto\mu_\tau$ satisfies criterion~\eqref{bw1} of \Cref{prop:holos}.
By hypothesis, for each $z\in\C$, $\tau\mapsto \mu_\tau(z)$ is holomorphic, so it has a power series expansion
\[\mu_\tau(z) = \sum A_n(z)\tau^n,\]
and its radius of convergence is at least $1$.
Let $r=1/2$.
We have for all $d\in\N$ and $z\in\C$:
\[A_d(z) = \lim_{N\to +\infty} \frac{1}{N}\sum_{k=0}^{N-1} \frac{\mu_{re^{ik/N}}(z)}{(re^{ik/N})^d}.\]
As a simple limit of measurable functions, the map $A_d$ is thus measurable.
This implies that $(\tau,z)\mapsto\mu_\tau(z) = \sum A_n(z)\tau^n$ is itself measurable (note the passage to two variables).

Let $r<1$. Local boundedness implies that $\tau\mapsto \|\mu_\tau\|$ is bounded on the compact set $\ov B(0,r)\subset \D$. Let $M=M_r$ be a bound.
Consider the set $E \subset \ov B(0,r)\times \C$ of pairs $(\tau,z)$ for which $|\mu_\tau(z)|\geq M$.
By definition of $\|\cdot\|$, for each $\tau\in\ov B(0,r)$, the section $E_\tau = \setof{z\in\C}{(\tau,z)\in E}\subset \C$ has Lebesgue measure $0$.
Since $(\tau,z)\mapsto\mu_\tau(z)$ is measurable, $E$ is measurable.
It follows that if we denote $P$ the set of $z$ for which the section $E^z = \setof{\tau\in\D}{(\tau,z)\in E}\subset \D$ has positive Lebesgue measure, then the Lebesgue measure of $P$ is $0$. For every $z\in\C-P$, the function $\tau\mapsto \mu_\tau(z)$ takes value of modulus $\leq M$ for almost every $\tau$, hence for every $\tau$ since it is holomorphic. It follows that for all $z\in\C-P$, $|A_n(z)| \leq Mr^{-n}$ for all $n$ by Cauchy's integral formula (alternatively, to avoid integration, one can use the average formula given earlier).
Then for all $\tau\in\ov B(0,r)$ and all $z\in\C-P$,
\[|\mu_\tau(z)-\sum_{k=0}^n A_k(z)\tau^k| \leq M\frac{(|\tau|/r)^{n+1}}{1-|\tau|/r}.\] Since $P$ has measure $0$, it follows that, denoting $a_n$ the element of $\cal B$ represented by $A_n$: $\forall \tau\in\ov B(0,r)$,
\[\|\mu_\tau - \sum_{k=0}^n a_k \tau^k\| \leq M\frac{(|\tau|/r)^{n+1}}{1-|\tau|/r}.\] Hence the power series $\sum a_n \tau^n$ has a radius of convergence $\geq r$ and $\tau\mapsto \mu_\tau$ is is its sum on $B(0,r)$.
Since $r$ can be chosen arbitrary in $(0,1)$, this concludes.
\end{proof}

%
% ----------------------------------------------------------
%

\subsection{Counterexample to weak continuity}\label{app:contrex}

Fix $\kappa\in(0,1)$.
Let $\Int $ denote the integer part of a real number.
Let $\mu_1(z) = 0$ if $\Int \Re z$ is even and $\mu_1(z) = \kappa$ otherwise.
The ellipse associated to $\mu=\kappa$ has vertical major axis and ratio $K = h(\kappa)$ where
\[h : [0,1) \to [1,+\infty),\ \kappa\mapsto \frac{1+\kappa}{1-\kappa}.\]
Note that $h$ is non-linear (it is strictly convex).

\begin{figure}
\begin{tikzpicture}
\node at (0,0) {\includegraphics[scale=0.75]{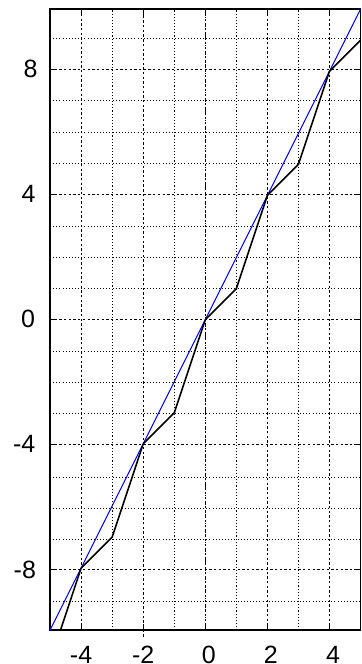}};
\end{tikzpicture}
\caption{Graph of the function $g$ in \Cref{app:contrex} for $\kappa=1/2$, hence $K=3$, in black. In blue the rescaling limit when conjugating by $x\mapsto nx$, i.e.\ $\lim g(nx)/n = K'x$ with $K'=2$. In this example the weak limit of $\mu_n$ is the constant $\mu=1/4$ and $K''=5/3 \neq K'$.}
\label{fig:gg}
\end{figure}

The normalized solution of the Beltrami equation for $\mu_1$ is given by $f(x+iy) = g(x)+iy$ where $g$ is the piecewise $C^1$ function with $g(0)=0$, $g'(x)=1$ whenever $\Int (x)$ is even and $g'(x) = K$ otherwise. See \Cref{fig:gg}.

Let $\mu_n(z) = \mu_1(nz)$.
Note that $\|\mu_n\|_\infty = \kappa$ for all $n$.
A normalized solution of the Beltrami equation for $\mu_n$ is given by $f_n(z) = f(nz)/f(n)$.
As $n\to\infty$, a computation gives
\[f_n(x+iy) \tend f(x) := x + iy/K'\]
with $K'=\frac{1+K}{2} = \frac{h(0)+h(\kappa)}{2}$.
On the other hand $\mu_n$ tends weakly to the constant $\mu = (0+\kappa)/2$, whose normalized straightening is $x+iy\mapsto x+iy/K''$ with $K'' = h(\frac{0+\kappa}2)$ so
\[K''\neq K'.\]

\begin{remark}
The dependence of $f$ on $\mu$ is non-linear and this is essentially the cause of the non-continuity with respect to the weak topology on $\mu$.
For instance, the simpler non-linear map $\mu \in L^\infty \mapsto \mu^2 \in L^\infty$ is not continuous if we take weak topology on both sides.
\end{remark}

%
% ----------------------------------------------------------
%

\subsection{Distributions with $L^1_\loc$ derivatives}\label{app:dl1}

Let $U$ be an open subset of $\R^n$.
The set of locally $L^1$ functions $f:U\to \R$ whose distribution derivatives are also locally $L^1$ is known as the set of \emph{locally Sobolev functions} and denoted $W^{1,1}_\loc(U)$.
Recall that such functions are equivalence classes up to modification on a set of Lebesgue measure $0$.

\begin{remark}
In fact, we may drop the first condition in the definition: a distribution whose partial derivatives are in $L^1_\loc$ is necessarily in $L^1_\loc$ (and this generalizes, see \cite{Maz}, the theorem of Section~1.1.2 page~7). However we will not need that fact. Moreover the functions $f$ we will consider in the end are maps that are assumed to be homeomorphisms, so they are continuous, which is better than $L^1_\loc$. So we stick with the original definition of $W^{1,1}_\loc$.
\end{remark}

More generally for $k\in\N$ and $p\in[1,+\infty]$, $W^{k,p}_\loc(U)$ denotes the Sobolev space of elements of $L^p_\loc$ whose distributional derivatives up to order $k$ are $L^p_\loc$.
In the version of the definition of quasiconformal maps $f$ that we chose (\Cref{def:be}) we assume that $f\in W^{1,2}_\loc(\C)$ and is continuous. In fact we assume still stronger: that $f$ is a homeomorphism from $\C$ to $\C$.

\begin{remark}
In one dimension, i.e.\ if $U$ is a subset of $\R$, then the elements of $W^{1,1}_\loc(U)$ have continuous representatives (in fact, absolutely continuous), and they are differentiable almost everywhere.\footnote{See for instance \cite{R}, Theorem~7.20 and the paragraph before~7.17.} 
If $n>1$, it is not true that any $W^{1,1}_\loc$ function has a continuous representative, nor an almost everywhere differentiable one (in the classical sense of being differentiable), not even that for almost every $x\in U$ there is a representative that is differentiable at $x$: a counterexample to these claims is given by the function 
\[f(x)=\sum_k \epsilon_k g(x-x_k)\]
where $x_k$ is a dense sequence in $U\subset \R^n$, $g\geq 0$ is a well chosen function described below and $\eps_k>0$ decreases sufficiently fast. To define $g$, let $r=\|x\|_2$ and if $n=2$ let $g(x) = \max(0,\log \frac{1}{r})$ or if $n>2$ let $g(x)=1/r^{n-2}$.
Then $f\in W^{1,1}_\loc(U)$ but the essential supremum of $f$ is infinite on every open set, because $f^{-1}((N,+\infty])$ is open and dense for all $N$, hence $f$ is differentiable nowhere even after modification on a set of measure $0$.
Yet, for any $W^{1,1}_\loc$ function $f$, on almost every line parallel to the main axes, $f$ is almost everywhere equal to an absolutely continuous function (and this remains true even after a $C^1$ change of variable $\psi$: $f\circ\psi$ is also $W^{1,1}_\loc$). 
This is called the ACL property, but we will not use it.
The interested reader may consult \cite{A,Abis}, Chapter~II\ B, and in particular Lemma~2.
\end{remark}

\begin{remark}
A similar counterexample for $W^{1,2}_\loc$ functions can also be built using $g(x)=\max(0,\log\frac{1}{r})$.
Actually for $k\in\N$, $p\in[1,+\infty)$, then all $W^{k,p}_\loc$ functions on $\R^n$ have continuous representatives if and only if $n< kp$ or ($p=1$ and $n=k$) where $n$ is the dimension number, see \cite{AF}, Theorem~4.12 page~85 and Sections~4.40 to~4.43.
\end{remark}

The next result is \cite{Abis}, Lemma~3 page~20, %page~19
we do not reprove it here. 
\begin{lemma}[Change of variable]\label{lem:chvar}
Consider a continuous function $f:U\to\C$ whose distribution derivatives $\de f/\de x_i$ are locally $L^1$.
Let $\psi : V\to U$ be a $C^2$ diffeomorphism from a subset $V$ of $\R^n$ to $U$.
Then the function $f\circ\psi$ has locally $L^1$ distribution derivatives, and the chain rule holds, i.e.\ 
\[\frac{\de (f\circ\psi)}{\de y_j} = \sum_{i=1}^n \frac{\de \psi_i}{\de y_j} \times \frac{\de f}{\de x_i} \circ \psi.\]
\end{lemma}

\begin{remark}
Actually $f$ does not need to be continuous: $f\in W^{1,1}_\loc$ is enough for the proof in \cite{A} to work.
Weaker conditions on the change of variables can also be made to work, for instance $C^1$, using different proofs.
Note that we will only need the lemma and its corollary below for $f$ continuous and a change of variable $\psi$ that is a rotation, see \Cref{app:pf:lem:sssol}.
\end{remark}

\begin{corollary}\label{cor:chvar}
In the previous lemma, if the distribution derivatives of $f$ are locally $L^2$ then the same holds for $f\circ \psi$.
\end{corollary}
\begin{proof}
Locally $L^2$ functions are in particular locally $L^1$ so we can apply the previous lemma.
The right hand side of the formula in the conclusion of this lemma adds up to a locally $L^2$ function.
\end{proof}

This also true with $L^2$ replaced by $L^p$, $p\in[1,+\infty]$.

%
% ----------------------------------------------------------
%

\subsection{Proof of \Cref{lem:sssol}}\label{app:pf:lem:sssol}

We already know that $f$ is a homeomorphism.
The main technical detail is to check that the distribution partial derivatives of $f$ are locally $L^2$.
Note that $f$ is smooth inside the polygons $Q_j$, which is even better.
There, the differential $da_j$ sends the ellipses $E$ associated to $\mu_j$ to circles while $F\circ\pi$, being holomorphic, sends circles to circles, so $df$ straightens the ellipses $E$ to circles, i.e.\ satisfies the Beltrami equation on $\Int Q_j$.
The rest of $\C$ consists points that are either interior points of edges or vertices.
Since this rest has Lebesgue measure $0$, once we know that $f\in W^{1,2}_\loc$ the above discussion immediately implies that $\bar \partial f = \mu\, \partial f$ holds almost everywhere.

There remains to check that $f\in W^{1,2}_\loc$. It is obvious inside the polygons: the map $f$ is smooth there, which is even better.
Near the edges and vertices we may use some erasability theorem available in the literature.\footnote{Like \cite{A2},  Theorem~4 page~9 (quasiconformal erasability of analytic arcs, here the arc is an open straight segment, so $f$ is q.c.\ at least on the complement of a finite set; a fortiori isolated points are erasable so $f$ is actually q.c.) or by repeated use of Rickman's lemma \cite{Ri}, Theorem~1 page~6 (near an inner point on an edge let $D$ be neighbourhood in $\C$ and $E$ the intersection of this neighbourhood with a closed polygon on one side of the edge; near a vertex use Rickman's lemma repeatedly to add the polygons one by one in trigonometric order, the last polygon added allows to include the vertex).}

Note that near an interior point of an edge one could also use the decomposition $f=F \circ \psi$ where $\psi$ is the join of $\pi\circ(a_j:j)$ and $\pi\circ(a_k:k)$ and work in the chart $\phi$ of $\cal S$ mentioned in \Cref{sub:sim},
for which $F\circ\phi^{-1}$ is holomorphic and $\phi\circ \psi$ is an explicty map that is $\R$-linear on both sides of a line.
Checking $W^{1,2}_\loc$ for such a map is easy.
It can thus be seen to be quasiconformal and the composition of a holomorphic bijection (or any quasiconformal map) with a quasiconformal map is still quasiconformal. This is not as simple near a vertex.

Alternatively, both near interior points of edges and near vertices, we can directly check that $f\in W^{1,2}_\loc$, as it is rather tractable in our case, and this is what we do in the rest of this section.

Near a point $z_0$ on an edge, if $z_0$ is not a vertex, we perform a change of variable $z$ in the domain, more precisely a rotation, so that this edge is vertical, which we assume now: this is enough by \Cref{cor:chvar}.
Denote $z_0=x_0+iy_0$.
Let $\phi$ be a test function whose support is contained in a square $S_\eps$ of equation $\max |\Re z-z_0|,|\Im z-z_0| < \eps$ for $\eps$ small enough. Then
\[
\la\pp{f}{y}, \phi\ra 
= \la f, -\pp \phi y \ra 
= \int_{-\eps}^\eps \left(da \times -\int_{-\eps}^\eps db \times \Big(f(z) \pp{\phi}{y}(z)\Big)\right)
\]
where in the inner integral: $z=(x_0+a)+i(y_0+b)$.

Claim: Let $e$ be an edge of $Q_j$ and $e^*$ the edge without its endpoints. Then $f$ is on $Q_j$ the restriction of a smooth function defined on $Q_j \cup W$ where $W$ is a neighbourhood of $e^*$.
Indeed, recall that for every edge $e$ of $P_j=a_j(Q_j)$ there is a chart $\psi:V\to\C$ on a neighbourhood $V$ of $\pi(e:j)$ that satisfies $\forall z\in P_j$, $\psi(\pi(z:j))=z$.
Then the map $\psi^{-1}\circ a_j$ is defined on $U=a_j^{-1}(\psi(V))$, its image is $V$, it is smooth and coincides with $\pi$ on $Q_j$. Since $F$ is holomorphic on $\cal S'$, the composition $F\circ\psi^{-1}\circ a_j$ is thus smooth and coincides with $f=F\circ\pi\circ a_j$ on  $U\cap Q_j$.

In particular we can integrate by parts without problem for \emph{all} $a\in(-\eps,\eps)$:
\[ - \int_{-\eps}^\eps db \times \Big(f(z) \pp{\phi}{y}(z)\Big)
= \int_{-\eps}^\eps db \times \Big(v(z) \phi(z)\Big)
\]
where $v$ is defined as follows: let $Q_j$ be the polygon on the right of the vertical segment, and $Q_k$ the one on the left. Let $f_j$ and $f_k$ be corresponding smooth extensions of $f$. If $a\geq 0$, let $v(z)=v_j(z)=\pp{f_j}{y}(z)$ and if $a\leq 0$ let $v(z)=v_k(z)=\pp{f_k}{y}(z)$. These two continuous maps coincide on the vertical line with the derivative of $f$ along this line, hence $v$ is continous.
The map $v$ is a representative of the distribution derivative of $f$ on $S_\eps$, since by Fubini's theorems we can finish the computation and get $\la f,-\pp \phi y\ra = \la v,\phi\ra$.
As a continuous map on $S_\eps$, $v$ is even better than locally $L^2$.

This is similar for $\pp{f}{x}$. Write
\[
\la\pp{f}{x}, \phi\ra 
= \la f, -\pp \phi x \ra 
= \int_{-\eps}^\eps \left(db \times -\int_{-\eps}^\eps da \times \Big(f(z) \pp{\phi}{x}(z)\Big)\right)
\]
with $z=(x_0+a)+i(y_0+b)$ and
\[ -\int_{-\eps}^\eps da \times \Big(f(z) \pp{\phi}{x}(z)\Big)
= -\int_{-\eps}^0 da \times [\cdots]
 -\int_{0}^\eps da \times [\cdots].
\]
Now
\[ \int_{0}^{\pm\eps} da \times \Big(f(z) \pp{\phi}{x}(z)\Big)
= f(z_0+ib)\phi(z_0+ib) -\int_{0}^{\pm\eps} da \times \Big(u(z) \phi(z)\Big)
\]
where
$u(z) = u_j(z) = \pp {f_j} x$ if $a>0$ and
$u(z) = u_k(z) = \pp {f_k} x$ if $a<0$.
This time $u_j$ and $u_k$ do not match anymore on the vertical axis so we cannot join them into a continuous function. Nevertheless, let $u=u_j$ or $u_k$ inside the polygons and $0$ on the vertical edge: then $u$ is measurable and bounded, hence locally $L^2$, and the formula $\la f,-\pp\phi x\ra = \la u,\phi\ra$ holds by Fubini and cancellation of the term $f(z_0+ib)\phi(z_0+ib)$. This cancellation can be seen as a consequence of the continuity of $f$.

That $f$ is also $W^{1,2}_\loc$ at vertices $z_0$ is proved similarly by decomposing the integral $\int f\pp \phi {x_i}$ over the square neighbourhood $S_\eps$ of $z_0$ with Fubini and cutting each horizontal or vertical segment on which the inner integral is taken into pieces cut by the edges (recall the edges are straight segments in the domain of $f$).
We can omit the vertical/horizontal segment going through $z_0$ since its contribution to the integral is $0$.
On each piece an integration by parts can be taken and cancellations of boundary terms will occur thanks to the continuity of $f$.
We thus get a function $u_i$ that is $0$ on the edges and is on $\Int Q_j$ the restriction of a function that is continuous on $Q_j-\{\text{vertices}\}$ and such that the distribution derivative $\partial f/\partial x_i$ is represented by $u_i$ by Fubini's theorems.
There remains to check that $u_i$ is locally $L^2$ at $z_0$ and in fact to use Fubini above it was already necessary to check that it is locally $L^1$, as a condition for Fubini's theorem is that the integrand has an absolute value of finite integral.
To check this, recall that as explained in \Cref{sec:simsurfpolygon}, a local conformal coordinate is given by a composition of $a_j$ with the complex logarithm $\log(a_j(z-z_0))$ followed by well-chosen translations followed by multiplication by a complex constant $\alpha$ independent of $j$ followed by the exponential.
Another way to express it is branches of $z\mapsto c_j\times a_j(z-z_0)^{\alpha_j}$ where $c_j\in\C^*$.
Moreover $\alpha_j = \frac{2\pi i}{i \theta + \log \lambda}$ for some $\theta>0$ and $\lambda \in (0,+\infty)$, in particular
\[\Re(\alpha_j)>0.\]
The map $f$ near $z_0$ is the join of such maps on finitely many sectors followed by holomorphic map $\psi$ defined near $0$ and independent of $j$ (of course it depends on $z_0$).
The derivative of $\psi$ is bounded.
We have $d(z^\alpha) = \alpha z^\alpha dz/z$ and $|z^\alpha| \approx |z|^{\Re \alpha}$ in the sense that their quotient remains bounded when $z$ remains in the sector indexed by $j$.
So on each sector, using polar coordinates and $dx\wedge dy = r dr \wedge d\theta$, we have
\[\ds \int \left|\frac{d(z^\alpha)}{dz}\right|^2 \approx \int_0^R  r^{(2\,\Re \alpha -1)}dr,\]
which is convergent since $\Re \alpha >0$.

%
% ----------------------------------------------------------
%

\clearpage

\part{A connection as a limit of the similarity surfaces}\label{part:cnx}

In \Cref{part:beltrami} we introduced a series of approximations of the normalized solution of the Beltrami equation associated to a Beltrami form $\mu$.
They straighten a Beltrami form $\mu_n$ that is constant on small squares and equal to the average of $\mu$ on these squares.
Together with each of these approximations came a similarity surface, with many puncture type singularities, that are erasable for the underlying Riemann surface.
These surfaces were conformally mapped to the Riemann sphere, yielding a global Riemann chart of the similarity surface, on which the similarity charts are recovered via a Christoffel symbol $\zeta_n$.
The function $\zeta_n$ is rational with simple poles, and the number of poles is of order $n^4$. As $n$ tends to infinity, the poles get close to each other and the residues get small.

In this part we make regularity assumptions on $\mu$ (in particular it is $C^2$) and prove that the sum of Dirac masses at the poles weighted by the residues converges weakly to a limit complex measure $m$ that we relate to $\mu$ and its first two derivatives (\Cref{prop:lim-mes}). 
We prove that $\zeta_n$ converges weakly to a complex valued (but not holomorphic) function $\zeta$ that is a convolution product involving $m$ (\Cref{prop:wlz}).
This function $\zeta$ defines a symmetric and conformal affine connection that is not flat. To an affine connection is associated a notion of parallel transport.
\Cref{part:cnx} culminates with the proof that the parallel transport associated to $\zeta_n$ converges in some sense to the parallel transport associated to $\zeta$ (\Cref{lim:parallel-transport}).

\Cref{sec:affco} introduces the necessary notions from the domain of \emph{affine connections} in what we hope is a gentle way.

In \Cref{sec:sslim} we prove the theorems mentionned above.

In \Cref{sec:var-avg} we make comments on the influence of changing the averaging process in the definition of $\mu_n$ from $\mu$.

In this more investigative part, we insist less on being self-contained and on using low level methods. In particular we allow ourselves to use results proved by other authors using the Ahlfors-Beurling operator.

%
% ----------------------------------------------------------
%

\section{Affine connections}\label{sec:affco}

This section starts with a quick introduction to affine connections, presented in coordinates: it covers directional derivatives $\partial_X Y$, connections in the form of the $\nabla$ operator, their Christoffel symbol $\Gamma$ in charts, parallel transport, symmetric connections, holonomy, and the geodesic equation.
Then we specialize to conformal connections in Riemann surfaces, define their curvature form $\omega$, relate it to holonomy, and discuss the case when $\omega$ vanishes.

\subsection{A quick introduction}\label{sub:qi}

Consider a (real) $d$-dimensional differentiable manifold $M$. An \emph{affine connection} (abbreviated as a \emph{connection} in most of the remainder of the article) is a differential object on $M$ closely related to a notion of \emph{parallel transport}.
It is difficult to honestly motivate the precise formulation of affine connections (equivalently of the parallel transport) so as to make them appear natural notions, so we will not try to do that, but instead define them with their expressions in charts.
Let us just mention that they generalize notions that occur in Riemannian manifolds and generalize themselves in the setting of vector bundles over manifolds.

\begin{remark}
A short introduction to affine connections can be found in the encyclopedic article \cite{Ha}.
The book \cite{CCL} is a classic treating differential geometry, and contains a presentation of affine connections that rapidly focuses on coordinates. Similar remarks hold for \cite{Spi}, Vol~II, Chapter~6.
A standard and thorough reference for the subject is \cite{KN}, which bases the approach on the notion of connection on principal bundles, so reading the definitions and statements require absorbing part of this theory and going back and forth between the chapters.
A nice, progressive and motivated presentation can be found in \cite{Sh}, however the theory is formulated in the framework of Cartan geometries and the link with the material explained here is even less direct.
There are many other books on the subject.
We end this paragraph by pointing out the fact that what we today call affine connections were initially called linear connections, and affine connections used to be a related but slightly different notion.
\end{remark}

\medskip

\noindent\textsf{Directional derivative.} Consider a chart of $M$, with image the open subset $U\subset \R^d$.
Consider a differentiable object $\bullet$ of a type $\cal T$ that is a function, a vector field, a form, or more generally any object type that is expressed in the chart as a function from $U$ to a fixed vector space $E$ (tensors, for instance, can be expressed in charts as maps from $U$ to the set $E$ of multilinear maps $(\R^d)^n\to\R$).
In $U$ the \emph{directional derivative} of $\bullet$ at a point $x\in U$, given a vector $v\in T_x U = \R^d$ is
\[\partial_{x,v} \bullet = \sum_{i=1}^d v^i \frac{\partial\, \bullet}{\partial x_i}\{x\}\]
where $\{x\}$ means to evaluate the quantity at point $x$.
Note that we use the exponent notation $v^i$ for the $i$-th coordinate of the vector $v$, a convention that is traditional in tensor calculus.

If $X$ is a vector field over $M$, we can define in the chart
\[\partial_X \bullet = \sum_{i=1}^d X^i \frac{\partial\, \bullet}{\partial x_i}\]
where the $X^i$ are the components of $X$, and every term is a function of $x$ with values in $\R$ or in $E$.
This is a function from $U$ to $E$.
It is to be stressed that, if $\bullet$ is not a function ($0$-form), then $\partial_X \bullet$ does not satisfy the same formula of change of coordinates as objects of type $\cal T$. Or if one insists anyway on defining an object of type $\cal T$ using $\partial_X\bullet$, then this object will depend on the chart.

Note: The classical Lie derivative of a vector field is recovered as
\[\cal L_X Y = \partial_X Y - \partial_Y X\]
and is independent of the chart as a vector field.

\medskip

\noindent\textsf{Affine connection.} Denote $(e_1,\ldots,e_d)$ the canonical basis of $\R^d$. In the chart, an affine connection $\nabla$ is expressed as follows: to vector fields $X$, $Y$ it associates
\[\nabla_X Y = \de_X Y + \Gamma(X,Y)\]
where $\Gamma$ is a bilinear endomorphism of $\R^d=T_x U$ for each $x\in U$:
\[\Gamma(X,Y) = \sum_{ijk} \Gamma^i_{jk} X^jY^k e_i\]
where the $d^3$ coefficients $\Gamma^i_{jk}$ are functions of $x$ and are called the \emph{Christoffel symbols} of the connection in the given chart.
Expressed in another chart, one can check that the connection takes the same form for another function $\Gamma$.

The \emph{canonical connection of $\R^d$} is defined by $\Gamma=0$, i.e.\ 
\[\nabla_X Y = \de_X Y.\]
A connection on $M$ that is canonical in one chart will fail to be in many other charts. A connection for which there are local charts where it is canonical is called \emph{locally trivializable}.

\medskip

\noindent\textsf{Parallel transport.} For a curve $t\mapsto\gamma(t)$ whose differential does not vanish, the parallel transport of a vector $v$ along $\gamma$ with respect to the connection $\nabla$ is a vector $v(t)$ attached to $\gamma(t)$ and that is locally a solution of 
\[\nabla_X Y\{x\} = 0\]
for all $x$ along the curve $\gamma$ and for any vector fields $X$, $Y$ such that, locally for a given $t$, $X(\gamma(t)) = \gamma'(t)$ and $Y(\gamma(t)) = v(t)$. This elaborate definition amounts the simple ODE
\[v'(t) = -\Gamma\{\gamma(t)\} \big(\gamma'(t),v(t)\big)\]
which, omitting $t$, reads as
\[v' = -\Gamma\{\gamma\} \big(\gamma',v\big).\]
This formula allows to generalize parallel transport to paths whose derivative vanishes for some values of $t$. 
The map that associates $v(t)$ to $v(0)$ can be seen as a map from the tangent space of $M$ at the point represented by $\gamma(0)$ to the tangent space of $M$ at the point represented by $\gamma(1)$. This map is linear and is also called \emph{parallel transport along $\gamma$}.

\medskip

\noindent\textsf{Symmetric connections (a.k.a.\ torsion free connections).} The connection is symmetric whenever the bilinear form $\Gamma$ is symmetric, i.e.\
\[\forall i,j,k,\ \Gamma^i_{jk} = \Gamma^i_{kj}.\]
This is independent of the chart and can be defined in a coordinate-independent form by the cancellation of an associated tensor called the torsion, but we will not use this here.\footnote{See in Chapter~III of \cite{KN}, Theorem~5.1 page~133 together with Proposition~7.6 page~145.}

\medskip

\noindent\textsf{Holonomy}. The parallel transport along a closed curve is called a \emph{holonomy} and is a self-map of the tangent space at $\gamma(0)$.

\medskip

A connection is called \emph{flat} if, locally, its holonomies are all the identity. By this we mean that every point has a neighbourhood $U$ such that the holonomy of every path contained in $U$ is the identity. This is equivalent to asking that the \emph{curvature tensor} vanishes everywhere.\footnote{See \cite{BG}, Theorem~5.10.3 page~236. On page~232/233 one finds a definition of the curvature tensor and equation~5.10.10 page~234 gives a formula for the expression in coordinates. The statements in this reference depend on some differentiable function $\mu$ between manifolds and here we just take $\mu: M\to M$ to be the identity.}

\medskip

A connection is locally trivializable if and only if it is symmetric and flat.\footnote{See \cite{BG}, the corollary on page~238.}

\medskip

Note that there exist connections that are flat but not symmetric.
For instance one can take on $\R^2$ the connection whose symbol $\Gamma^1_{12} =1$ and all other $\Gamma^i_{jk}=0$. The connection whose only non-vanishing symbol is $\Gamma^1_{21}=1$ is another example.

\medskip

\noindent\textsf{Geodesics}. The famous geodesic equation for a path $\gamma$ asks that $\gamma'(t)$ be a parallel transport along $\gamma$, i.e.
\[\gamma'' = -\Gamma\{\gamma\} \big(\gamma',\gamma').\]
It depends only on the quadratic vector form associated to the bilinear vector form $\Gamma$. It has the same geodesics as the symmetrized connection with $\Gamma^{\mathrm{sym}}(u,v)=\frac{1}2{}\left(\Gamma(u,v)+\Gamma(v,u)\right)$.
For a symmetric connection, the collection of geodesics, as parameterized curves, characterizes $\Gamma$, so characterizes the connection.
This it not true if one only considers the support of the geodesics: for instance a projective transformation of $\R^2$ sends straight lines to straight lines but the trivial connection is sent to one whose $\Gamma$ is not $0$.

\medskip

\noindent\textsf{Riemannian metrics.}
Given a $C^1$ riemannian metric $g$, there is a unique symmetric connection whose parallel transport preserves $g$.\footnote{This is called the fundamental lemma/theorem of Riemannian geometry. See for instance \cite{CCL} Chapter~5, \cite{Spi}, Chapter~6 or \cite{KN}, Chapter~IV.} It is called the \emph{Levi-Civita connection} of $g$. Its holonomies are isometries of each tangent space.

\subsection{Affine and conformal connections in two dimensions}\label{sub:affcoconf}

If $d=2$, a connection has 8 independent coefficients, and a symmetric connection has 6 independent coefficients.

Assume now that $M$ is a Riemann surface (not to be confused with Riemannian surface) and see it as a two dimensional manifold. Then the parallel transport associated to a connection $\nabla$ preserves the conformal structure (angle between vectors) iff in a conformal chart, each map $v\mapsto\Gamma\{x\}(u,v)$ is a similitude for all $u$, iff $v\mapsto\Gamma\{x\}(e_1,v)$ and $v\mapsto\Gamma\{x\}(e_2,v)$ are two similitudes for each $x$. There are thus 4 independent (real) coefficients.

Finally, a symmetric connection that preserves the conformal structure has two independent real coefficients and it takes the following nice expression in a conformal chart, identifying $\R^2$ with $\C$ and using complex multiplication:
\[\nabla_X Y =\partial_X Y +\zeta X Y\]
where $\zeta$ is a complex valued function of $x$. We call this a \emph{conformal symmetric connection} and $\zeta$ is called its \emph{Christoffel symbol}.

\begin{remark} This is an abuse of language since we also call Christoffel symbol the coefficients $\Gamma^i_{jk}$.
Context should make clear which one we mean.
We can motivate this abuse as follows.
If one takes the definition of connections in one dimension (they have only one coefficient $\Gamma^1_{11}$) and complexifies it, i.e.\ interprets it on a complex curve, then one obtains exactly the formula above in charts, with the difference that a priori $\nabla$ is only defined to operate on holomorphic vector fields $X$ and $Y$ (because they are assumed complex-differentiable).
In a complex setting it is also natural to ask $\zeta$ to be a holomorphic function but we will need to consider here cases where it is not.
We distinguish this by saying that the conformal symmetric connection may or may not be holomorphic. 
\end{remark}

In this case the parallel transport equation takes the following nice form:
\[v'=-\zeta(\gamma)\gamma'v.\]
Its solution is
\[v(t)=v(0)\exp \tau(t)\]
with
\[\tau(t) = -\int_0^t \zeta(\gamma(s))\gamma'(s) ds = -\int_0^t \zeta(\gamma(s))d\gamma(s).\]
As we did not assume $\zeta$ holomorphic, this path integral usually \emph{does not only depend} on the homotopy class of $\gamma$.
If $\zeta$ is at least $C^1$ then by Stokes' formula, the holonomy for a path equal to the oriented boundary of a simply connected zone $S\subset U\subset\R^2$ is $v\mapsto \exp(\tau) v$ with
\[\tau = -2i\!\iint \bar\partial\zeta\, dx\wedge dy\]
where $\bar\partial \zeta := \partial_{\bar z} \zeta := \partial \zeta/\partial \bar z$ is the anticonformal part $b\in\C$ of the decomposition of the differential $d\zeta = a\, dz + b\, \overline{dz}$. 
The \emph{curvature form} is the form appearing in the formula above:
\[\omega =-2i\,\bar\partial\zeta\, dx\wedge dy.\]
This is a 2-form, well defined on the Riemann surface $M$ and the holonomy formula extends to $C^1$ Jordan domains in $M$.
It can also be expressed as follows, using $dz\wedge d\bar z = (dx+idy) \wedge (dx-idy) = -2idx\wedge dy$:
\[\omega = \bar\partial\zeta\, dz\wedge d\bar z.\]

\medskip

\noindent\textsf{Change of variable.} Let $z=\phi(w)$ be a holomorphic change of variable and $\zeta^w$ denote the Christoffel symbol of the connection in the variable $w$. Then
\[\zeta^w = \phi' \times \zeta^z \circ\phi +\frac{\phi''}{\phi'}.\]
This is the same formula as for the symbol $\zeta$ appearing in \Cref{eq:zetachvar} on page~\pageref{eq:zetachvar},
with the difference that, here, $\zeta$ is not assumed holomorphic. In fact the symbol $\zeta$ here is the same as the symbol $\zeta$ there.

The effect of the change of variable on the curvature form is the same effect as on any 2-form on a 2-dimensional real manifold.

\medskip

\noindent\textsf{Vanishing curvature form.}

\begin{lemma}\label{lem:curv-nul}
If the curvature form of a symmetric conformal connection vanishes identically, then there are local coordinate systems that trivialize the connection: $\zeta=0$ in these charts.
\end{lemma}
\begin{proof} Its local holonomies are the identity, i.e.\ it is flat. We mentionned earlier that locally trivializable connections are those that are flat and symmetric.
\end{proof}

Actually we can say more:
\begin{lemma}\label{lem:curv-nul-2}
The curvature form of a symmetric conformal connection vanishes identically if and only if in charts, its Christoffel symbol $\zeta$ is holomorphic.
The trivializing coordinates can be taken to be holomorphic.
\end{lemma}
\begin{proof}
The first claim follows from the expression of the curvature form as $\omega = \bar\partial\zeta\, dz\wedge d\bar z$, which it vanishes iff $\bar\partial\zeta=0$, which is equivalent to the differential $d\zeta$ being $\C$-linear at every point, i.e.\ to $\zeta$ being holomorphic.
In this case, a holomorphic coordinate $\phi$ trivializes the connection if and only if 
$\zeta = 0+\frac{\phi''}{\phi'}$ by the previous change of variable formula, % geodesics $\gamma$ map to uniform motion on straigh lines: $\phi\circ \gamma(t) = at+b$, which amounts to $(\phi\circ \gamma)''=0$, i.e.\ $\gamma''\phi'\circ\gamma + (\gamma')^2\phi''\circ\gamma=0$, i.e.\ $\frac{\phi''}{\phi'}\circ\gamma = -\frac{\gamma''}{(\gamma')^2} = -\Gamma \circ \gamma$, i.e.\ $\phi''/\phi' = -\Gamma$, 
an ODE which has (holomorphic) solutions.
\end{proof}

Post-composing a trivializing coordinate by a real affine map of $\R^2$ also yields a trivializing coordinate, so not all trivializing coordinates are holomorphic.

\medskip
\noindent\textsf{Conformal metrics.}
A conformal metric on a Riemann surface is a Riemannian metric of the form
\[\rho(z) |dz|\]
in conformal charts, with $\rho>0$, in other words $g=\rho^2(dx^2+dy^2)$. 
It will be convenient to write
\[\rho^2 = e^h\]
for some function $h$. Assume $h$ is $C^1$. Then the Levi-Civita connection of $g$ is conformal and we have\footnote{The motivated readers can deduce it from \cite{KN2}, Chapter~IX, section~5. Or they can check it directly using the formula expressing the Levi-Civita connection in terms of the metric; for an expression in terms of coordinates and the coefficients of the Christoffel symbols, see for instance \cite{KN}, Corollary~2.4, Chapter~4, page~160 or \cite{CCL}, formula~(1.35), Section 5--1, page~139.}
\[\zeta = \partial_z h\]
where $\partial _z h := \partial h/\partial z$ is the conformal part $a\in\C$ of the decomposition of the differential $dh = a\, dz + b\, \overline{dz}$. Since $\partial_{\bar z} \partial_z = \frac{1}{4}\Delta$,  assuming $\rho$ is $C^2$
the curvature form is
\[ \frac{\Delta \log \rho}{i}dx\wedge dy.\]
Note that it is purely imaginary.
Holonomies are rotations, and the holonomy around the oriented boundary of a simply connected zone $S\subset M$ is the rotation whose angle, in radians, is given by the imaginary part of the integral of the curvature form over $S$: it is the integral of the form defined in charts by $-(\Delta\log\rho)\, dx \wedge dy$.

The \emph{curvature} is the quotient of $-i$ times the curvature form by the area form $\rho^2\, dx\wedge dy$ naturally associated to $g$, and takes the following expression
\[K = \frac{\Delta \log \rho}{\rho^2}\]
which is purely real.

%
% ----------------------------------------------------------
%

\section{A limit of the similarity surfaces used in the construction of a solution to the Beltrami equation}\label{sec:sslim}

We state here the main results of \Cref{part:cnx}.
  We start by a short introduction to the notion of Beltrami forms and ellipse fields.
  Then from \Cref{sub:setup} on, we restrict ourselves to the situation of a Beltrami form on $\C$ that is $C^2$ with compact support and recall the approximating scheme for the solutions of the Beltrami equation done in \Cref{part:beltrami}; we define in \Cref{sub:ato} an associated complex valued atomic measure $m_n$ supported by the position of the singularities in $\hat\C$ after uniformization, with value at these points the residues of the associated meromorphic connection $\zeta_n$.
  These $m_n$ can be thought of as the curvature form $\omega$, but taken in the sense of distribution, of the associated meromorphic connection on $\wh\C$.
  We prove the weak convergence of $m_n$ to an (absolutely continuous) complex measure $m$, which we identify.
  In \Cref{sub:lack} we comment on the naturality of the measure $m$.
  The convergence of $m_n$ to $m$ implies the weak convergence of the functions $\zeta_n$ to a function $\zeta$ which is expressed as a convolution (\Cref{sub:limSn}).
  This limit is the Christoffel symbol of a conformal (but not holomorphic) connection, which we interpret in
  \Cref{sub:limPT} as giving the limit of the parallel transports associated to the $\zeta_n$, and in \Cref{sub:interpretLim} as the unique conformal connection satisfying some property related to the pull-back of the horizontal and vertical vector fields by the map $f$ straightening $\mu$.

\subsection{Preliminary notions}\label{sub:prelim}

We shall use the notions of \emph{Beltrami form} and of \emph{ellipse field}, which are two equivalent ways of defining a differential object on differentiable manifolds of dimension $2$. In the first case the manifold is a Riemann surface.

On a Riemann surface $S$ it is designed to represent the Beltrami differential of a differentiable or quasiconformal mapping $f:S\to \C$ 
and takes form $B_z(f)=\partial_{\bar z} f/\partial_{z} f = \mu(z)$ in charts.
If $z=\phi(u)$ is a change of Riemann chart on $S$, then $\phi$ is holomorphic and
$B_u(f) = \bar\partial_u f/\partial_u f = \mu(\phi(u)) \overline{\phi'(u)}/\phi'(u)$, which defines the notion of pull-back and push-forward of Beltrami forms under \emph{conformal} maps and also justifies the notation
\[\mu(z)\frac{d\bar z}{dz}\]
in charts.

A second point of view is \emph{ellipse fields}.
This one works for any differentiable manifold $M$ of dimension $2$.
For each $p\in M$, we look a the space of all ellipses centred on $0$ in $T_p M$, modulo multiplication by a positive real.
This defines a bundle whose sections are the ellipse fields.
This notion transports well under the differential of a differentiable function between two dimensional manifolds as long as the differential remains invertible. This allows to define pull-backs and push-forwards of ellipse fields under diffeomorphisms.

On $\C$ consider the global circular ellipse field $\cal O$, whose ellipses are all circles.
For a given $f: S\to \C$ that is almost everywhere differentiable with invertible differential preserving the orientation, there is a one to one correspondence between the pull-back $f^*\cal O$ and the Beltrami differential of $f$.
Indeed, in a chart the differential of $f$ at a point $z$ takes the form
\[v\in\C \mapsto a\,v+ b\,\bar v\]
with $a=\partial_z f$ and $b=\partial_{\bar z} f$.
The preimage of circles have equation
$|a\, v+b\, \bar v| = \cst$ for various constants $\cst$, which is the same as the equation
\[|v +\mu\, \bar v| = \cst\] for another constant and where $\mu = ba^{-1} = B_z f$.
If $\mu=0$ these preimages are homothetic circles. Otherwise these are homothetic ellipses with:
\begin{itemize}
\item the ratio of the major and minor axes is equal to $\frac{1+|\mu|}{1-|\mu|}$,
\item the direction of the minor axis is $\frac{1}{2}\arg\mu$ modulo $\pi$.
\end{itemize}
One sees that this defines a bijection between the set of all ellipses centred on $0$ quotiented by the group of homotheties, and the set of all $\mu\in\D$.

For a function $z\mapsto\mu(z)$ defined on an open subset of $\C$ (more generally for a Beltrami form $\mu$ defined on a Riemann surface), the Beltrami equation $B f (z)=\mu$ concerning the differentiable or quasiconformal map $f$, is equivalent to $f^*\cal O = \cal E$ where $\cal E$ is the ellipse field defined by $\mu$ (the equality is required to hold almost everywhere in the case of quasiconformal maps).

We will use the following facts:
\begin{itemize}
\item
Quasiconformal maps are Lebesgue regular in that the image of a set of Lebesgue measure $0$ has Lebesgue measure $0$. \cite{Abis}, Theorem~3 page~22.% page~33.
\item
The inverse (reciprocal) of a quasiconformal map is quasiconformal. The composition of two quasiconformal maps is quasiconformal. \cite{Abis} ``trivial'' properties~3 and~4 page~15 %page~22 
together with the equivalence between definitions~A page~15 %page~21 
and~B' page~19.%page~29.
\item
A quasiconformal map is differentiable almost everywhere and its differential is invertible almost everywhere. \cite{Abis} Lemma~1 page~17 %page~24
and Corollary~3 page~22.%and the corollary page~34.
\end{itemize}
As a consequence, a coherent notion of pull-back and push-forward of ellipse fields by quasiconformal maps can be defined. It satisfies $g^*(f^*\cal E)=(f\circ g)^* \cal E$ for all $f$, $g$ quasiconformal and $\cal E$ ellipse field.

Since there is a perfect correspondence between Beltrami forms (with $\|\mu(z)\|<1$) and ellipse fields, one can extend the notion of pull-back and push forward of Beltrami forms to non-holomorphic orientation preserving differentiable or quasiconformal maps.
An interesting alternative point of view can be found in \cite{Hu}, Section~4.8.

\subsection{Setup}\label{sub:setup}

From here onwards, we assume that $\mu:\C\to\D$ is $C^2$ with compact support.

We recall the straightening method of \Cref{part:beltrami}, \Cref{sub:MRMT}: for $n>0$ we define $\mu_n$ by 
dividing the square defined by $|\Re z|<n$, $|\Im z|<n$ into small squares of side $1/n$  and we let $\mu_n$ be constant in the interior of each of those $4n^4$ squares, equal to the average of $\mu$ on this square, and elsewhere we let $\mu_n(z)=0$.
We saw in \Cref{app:beltrami:cv} that the normalized solution $f_n$ of the Beltrami equation associated to $\mu_n$ converges to the normalized solution $f$ of the Beltrami equation associated to $\mu$.

\subsection{A sequence of atomic complex measures}\label{sub:ato}

To $\mu_n$ we associated a similarity surface $\cal S'_n$ with singularities corresponding to corners of the squares (and infinity): it is obtained by mapping each square $S$ to a quadrilateral by any real affine map that straightens the constant Beltrami form $\mu_n$ on $S$, then gluing together all the obtained quadrilaterals, and gluing this to the complement of the big square of side $2n$.
The similarity surface $\cal S'_n$ was completed at its singularities to form a Riemann surface $\cal S_n$ homeomorphic to a sphere.
An isomorphism $\cal S_n\to\hat\C$ was chosen, sending $\infty$ to $\infty$. On the target $\hat \C$, a meromorphic Christoffel symbol $z\mapsto \zeta_n(z)$ was defined on $\hat\C$ minus $\infty$ and minus the images $s$ of the other singularities, with simple poles at $s$ of residue $\res=\res(s)$ such that in particular $\exp(2\pi i\res)$ is the monodromy factor of the similarity surface at $s$.
The precise determination of $\res$ depends on the sum of angles of the quadrilaterals at the singularity, see \Cref{sub:presc}.

Consider the complex-valued measure $m_n$ in the range $\C\equiv\R^2$ defined by the sum of Dirac masses at each square corner $c$ with complex weight $\res(s)$ where $s=f_n(c)$.
The measure $m_n$ can also be considered as complex-valued $2$-form on $\R^2$ that is singular (a.k.a.\ a current).

\begin{proposition}\label{prop:lim-mes}
Under the conditions of \Cref{sub:setup}, $m_n$ weakly tends to a complex valued finite measure $m$ on $\C$ with continuous density w.r.t.\ the Lebesgue measure given by
\[m= -\frac{1}{2\pi}\left(\frac{2 \frac{\partial^2\mu}{\partial x\partial y}}{1-\mu^2}+\frac{ 4\mu\frac{\partial \mu}{\partial x}\frac{\partial \mu}{\partial y}}{(1-\mu^2)^2} \right)\on{Leb}
\]
in the following sense:
for every continuous function $\tau$ on $\C$,
\[\int \tau\, m_n \tend \int \tau\, m.\]
Moreover, the total mass converges: 
\[|m_n|\tend |m|
.\]
\end{proposition}
\begin{proof}
Let $c\in \C$ and consider the following four squares:
\bEA
C_0 &=& ``\Re (z-c)\in[0,\eps] \text{ and }\Im (z-c)\in[0,\eps]\text{''} \\
C_1 &=& ``\Re (z-c)\in[-\eps,0] \text{ and }\Im (z-c)\in[0,\eps]\text{''} \\ 
C_2 &=& ``\Re (z-c)\in[-\eps,0] \text{ and }\Im (z-c)\in[-\eps,0]\text{''} \\ 
C_3 &=& ``\Re (z-c)\in[0,\eps] \text{ and }\Im (z-c)\in[-\eps,0]\text{''} \\
\eEA
Consider also the Taylor expansion
\[\mu(c+z) = \mu_0 + \mu_x x + \mu_y y + \mu_{xx} x^2 + 2\mu_{xy}xy + \mu_{yy} y^2 + o(z^2)\]
where $z=x+iy$ and denote
\bEA
\Delta & = & \mu_{xx}+\mu_{yy},\\
\rtimes & = & 2\mu_{xy}.
\eEA
Later in this proof we will use that the $o(z^2)$ is uniform, i.e.\ that 
\[\mu(c+z) = \mu_0(c) + \mu_x(c) x + \mu_y(c) y + \mu_{xx}(c) x^2 + 2\mu_{xy}(c)xy + \mu_{yy}(c) y^2 + |z|^2 r(c,z) \]
where $r(c_n,z_n)\tend 0$ as $z_n\tend 0$ and $c_n\in \C$ is any sequence.\footnote{Since $\mu$ has compact support and is $C^2$, there is a uniform modulus of continuity $\delta\mapsto M(\delta)$ for its second order partial derivatives.
By two successive integrations, one finds that $|r(c,z)|\leq M(|z|)$.}

Let $M_i$ be the average of $\mu(z)$ for $z\in C_i$. 
Then
\bEA
M_0 & = & \mu_0 + \frac{\mu_x+\mu_y}{2}\eps + \left(\frac{\Delta}{3}+\frac{\rtimes}{4}\right)\eps^2+o(\eps^2)\\
M_1 & = & \mu_0 + \frac{-\mu_x+\mu_y}{2}\eps + \left(\frac{\Delta}{3}-\frac{\rtimes}{4}\right)\eps^2+o(\eps^2)\\
M_2 & = & \mu_0 + \frac{-\mu_x-\mu_y}{2}\eps + \left(\frac{\Delta}{3}+\frac{\rtimes}{4}\right)\eps^2+o(\eps^2)\\
M_3 & = & \mu_0 + \frac{\mu_x-\mu_y}{2}\eps + \left(\frac{\Delta}{3}-\frac{\rtimes}{4}\right)\eps^2+o(\eps^2)
\eEA
where the $o(\eps^2)$ are uniform as above.
Consider the $\R$-affine map
\[A_j(z) = z+M_j \ov z,\]
which satisfies $\partial_{\bar z} A_j /\partial_z A_j= M_j$.
Let $\theta$ be the sum of the angles at $0$ of the four quadrilaterals $A_j(C_j)$. As $\eps\tend 0$, we have $\theta\tend 2\pi$.
Then the monodromy factor $\Lambda = \lambda e^{i\theta}$ around $0$, which depends on $c$ and $\epsilon$, satisfies
\[\frac{1}{\Lambda} = \frac{A_1(i)}{A_0(i)} \times \frac{A_2(-1)}{A_1(-1)} \times \frac{A_3(-i)}{A_2(-i)} \times \frac{A_0(1)}{A_3(1)}\]
which, after simplification reads
\[\Lambda = \frac{1-M_0}{1+M_0}\times\frac{1+M_1}{1-M_1}\times\frac{1-M_2}{1+M_2}\times\frac{1+M_3}{1-M_3}\]
i.e.\ 
\[\Lambda = \frac{L_1\times L_3}{L_0 \times L_2}\]
with
\[L_j = \frac{1+M_j}{1-M_j}.\]
We have $\res = \frac{\log (\lambda) +i\theta}{2\pi i} - 1$
As $\eps\tend 0$, the monodromy factor tends to $1$.
Since $\theta\tend 2\pi$, we have that $\log (\lambda)+i\theta -2\pi i = \log_p \Lambda $ where $\log_p$ denotes the principal branch of the complex logarithm:
\begin{equation}\label{eq:rnb0}
\res = \frac{\log_p \Lambda}{2\pi i}
\end{equation}
As soon as $\eps$ is small enough (independently of $c$).
After a moderately complicated calculation involving the expansions of the $M_j$, one finds that
\begin{equation}\label{eq:lambda-expand}
\log_p \Lambda = -\left(\frac{2\rtimes}{1-\mu_0^2}+\frac{4\mu_0\mu_x\mu_y}{(1-\mu_0^2)^2}\right)\eps^2+o(\eps^2)
\end{equation}
where the $o(\eps^2)$ is uniform in the sense explained earlier.

We denote $\res=\res(c,\eps)$ and $\Lambda_p=\Lambda_p(c,\eps)$ to emphasize that they depend on $c$ and $\eps$.
By the above analysis, the convergence of $(\log_p\Lambda(c,\eps))/\eps^2$ is uniform on $c\in \C$ as $\eps\tend 0$.
Note that for $\eps = 1/n$, the denominator $\eps^2$ is the area of each small square. 

Let for $c\in \C$
\[h(c) = \frac{1}{2\pi i}\lim_{\eps\to 0} \frac{\log_p\Lambda(c,\eps)}{\eps^2} = \lim_{\eps\to 0} \frac{\res(c,\eps)}{\eps^2}\]
and let $m=h\times\Leb$.

The rest of the proof is standard but we include it anyway.

Let $\cal C$ denote the (finite) set of squares $S$ of horizontal and vertical sides of length $1/n$, centred on the points $s$ of coordinates $(i/n,j/n)$, $i,\,j\in\Z$, and $\cal C'$ the (finite) subset of those for which $2S$ has non-empty intersection with the support of $\mu$, where $2S$ denotes the square of the same centre but double side length.
Note that $\cal C$ and $\cal C'$ depend on $n$.
The square $2S$ is formed of the four little squares that have $s$ as a vertex and the condition ensures that every point in the support of $m_n$ is the centre of some $S\in\cal C$.
In fact if $s\in\cal C-\cal C'$ then $m_n$ and $m$ are zero on $S$ so in particular $\int_S m_n\tau=0=\int_S m\tau$.
Let $R>0$ big enough so that $\ov B(0,R)$ contains the support of $\mu$.
Note that if $S\in\cal C$ then $S\subset \ov B(0,R+\sqrt{2}/n)\subset \ov B(0,R+\sqrt{2})$.

Let $\tau$ be a continuous function on $\C$.
Let $\eta>0$. Then by uniform continuity of $h$ (it is continuous with compact support) and of the restriction of $\tau$ to $\overline B(0,R+\sqrt 2)$, for all $n$ big enough, for all square $S\in\cal C'$, denoting $s$ its centre, we have:
\[\sup_{z\in S} |h(z)-h(s)|<\eta,\ \ \sup_{z\in S} |\tau(z)-\tau(s)|<\eta,\ \ \left|\frac{\res(s)}{\Leb S}-h(s)\right|<\eta.\]
The third inequality is a consequence of the uniform convergence mentioned above.
For all $s\in\cal C$ we have $\int_S m_n\tau = \res(s)\tau(s)$ and thus for $S\in \cal C'$ we get
\[ \left|\int_S m_n\tau-\int_S m\tau \right| 
\leq \int_{z\in S} \left|\frac{\res(s)}{\Leb S}\tau(s) - h(z)\tau(z)\right|\times\Leb\]
\[\leq \int_{z\in S} \left(\left|\frac{\res(s)}{\Leb S}\tau(s) - h(s)\tau(s)\right| + |h(s)\tau(s) -h(s)\tau(z)|+|h(s)\tau(z)-h(z)\tau(z)|\right)\]
\[\leq \int_{z\in S} \left(\eta \|\tau\|_\infty + \eta \|h\|_\infty + \eta \|\tau\|_\infty \right)\times\Leb
\leq K \eta \on{Leb}(S) \]
for some $K>0$ independent of $n$ and of $\eta$, i.e.
%\begin{equation}\label{eq:rnb}
\[
\left|\int_S m_n\tau-\int_S m\tau \right| \leq K \eta \on{Leb}(S).
\]
%\end{equation}

We saw that for $S\in\cal C-\cal C'$ then $m_n$ and $m$ are zero on $S$.
Hence $\left|\int_\C m_n\tau-\int_\C m\tau \right|$ is bounded from above by the sum of $\left|\int_S m_n\tau-\int_S m\tau \right|$ over the squares $S\in\cal C'$.
Summing the upper bound of the previous paragraph, over this finite collection of squares, we get
\[ \left|\int_\C m_n\tau-\int_\C m\tau \right| \leq K\eta\on{Leb} W
\]
where $W = \ov B(0,R+\sqrt{2})$.

The proof for the total mass is based on an almost identical computation. We take the constant $1$ instead of $\tau$ and use uniform continuity of $|h|$ and uniform convergence of $|\res s|/\eps^2$ to $|h|$.
\end{proof}

We are also interested in the images (push-forward) of the measures $m_n$ by the straightening maps $f_n$.

\begin{corollary}\label{cor:lim-mes-2}
The measure $(f_n)_* m_n$ weakly tends to $f_* m$ in the same sense.
\end{corollary}
\begin{proof}
This follows from $f_n$ tending locally uniformly to $f$. Let $\tau$ be continuous over $\C$. The proof is easy but we detail it here.
\[\int (f_n)_* m_n \times\tau- \int f_* m \times\tau =
\int m_n \times \tau\circ f_n -\int m \times \tau\circ f\]
\[=\int m_n \times (\tau\circ f_n - \tau\circ f)
 + \int (m_n-m)\times \tau \circ f.\]
The second term tends to $0$ by a direct application of the second part of \Cref{prop:lim-mes}.
The first term is bounded by
\[|m_n| \times \sup_{B(0,R)}|\tau\circ f_n - \tau\circ f|\]
where $R$ is independent of $n$ and is taken big enough so that all $m_n$ and $m$ have support in $\ov B(0,R)$.
We have seen that $|m_n|$ is bounded.
By uniform continuity of $\tau$ on compact sets, we have $\sup_{B(0,R)}|\tau\circ f_n - \tau\circ f|\tend 0$.
\end{proof}

\medskip

\noindent\textbf{Note}: In \Cref{sec:var-avg} we show that if we average $\mu$ through a (reasonable) function of $\mu$ instead of directly taking $\mu_n = \int_S\mu /\int_S 1$ on each square $S$, we still get the same limit.

\subsection{On the (lack of) invariance of the limit measure $m$}\label{sub:lack}

We denote $m[f]$ to emphasize the dependence of $m$ on the $C^2$ function $f$.
\begin{enumerate}
\item The limit $m$ is invariant under a rescaling of the following form $g(x+iy) = a x + b i y + c$ with $a>0$, $b>0$ and $c\in\C$, i.e.\ $m[g\circ f] = g_* m[f]$ where $g_*$ refers to the push forward of measures;
\item it is invariant under rotation by a fourth of a turn: for all linear map $g(Z) = i z + c$ with $c\in\C$, $m[g\circ f] = g_* m[f]$;
\item it is \emph{not} generally invariant under a rotation by an angle $\alpha$ such that $\alpha\neq 0 \bmod \pi/2$;
\item and it is \emph{not} generally invariant under a non-conformal orientation-preserving $\R$-linear map.
\end{enumerate}

These claims are easy to check from the formula given in \Cref{prop:lim-mes}.

Another way to state them is that if instead of small squares directed by the main axes of $\R^2$ we had chosen (identical) rectangles with the same directions, it would still give the same limit $m$. But if we had chosen other (identical) quadrilaterals or rectangles directed by different axes, we would obtain a different limit $m$ for most choices of $\mu$, even though the corresponding maps $f_n$ converge to the same solution $f$ of $\mu$.

Hence, the weak limit $m$ of the residues $m_n$ attached to our approximation scheme for the solution of the Beltrami equation is less canonical than one would think from first examination, because not only it depends on $\mu$, but it also depends on the choice of a horizontal/vertical direction in $\R^2$.

\begin{remark}
We chose to cut the plane in small squares. If we had chosen for instance small hexagons, we would get a different formula for the limit $m$ of $m_n$. The exploration of all the possibilities goes beyond the scope of this article.
\end{remark}

\subsection{A limit for \texorpdfstring{$\cal S_n$}{Sn}}\label{sub:limSn}

Do the similarity surfaces $\cal S_n$ converge in some sense?

Consider the meromorphic Christoffel symbol $\zeta_n$ on $\hat\C$ associated to the construction. It is the symbol of a locally flat symmetric conformal connection with singularities at the poles of $\zeta_n$.

Consider also the push-forward by $f$ of the measure $m$: $f_* m$.

\medskip

It is quite easy to deduce from \Cref{cor:lim-mes-2} that $\zeta_n$ has a weak limit:

\begin{proposition}\label{prop:wlz}
The sequence $\zeta_n$ weakly tends on $\C$ to the following convolution product:
\[\zeta_n \wtends \zeta := (f_* m)*\frac{1}{z} \]
the weak limit is understood against continuous test functions, compactly supported in $\C$.
\end{proposition}
\begin{proof}
According to \Cref{thm:sc}, denoting $\kappa (z) = 1/z$ and $\on{neg}(z)=-z$ we have
\[\zeta_n = ((f_n)_* m_n) * \kappa.\]
Now
\[\int \tau\times(\zeta_n-\zeta) = \int \tau \times \big(((f_n)_* m_n-f_* m) * \kappa\big)\] 
\[= \int ((f_n)_* m_n - f_* m) (\tau * (\kappa \circ \on{neg})) \]
The last identity comes from the Fubini-Tonelli theorem for the product measure $\Leb \times (|m_n|+|m|)$, remarking that $\kappa$ is locally $L^1$ and that both integration variables remain bounded because $\tau$, $m_n$ and $m$ have compact support.
The function $\tau * (\kappa\circ\on{neg})$ is continuous (because $\tau$ is continuous and $\kappa$ is locally $L^1$).
By \Cref{cor:lim-mes-2} the integral tends to $0$.
\end{proof}

\begin{lemma}\label{lem:c1}
  If $\mu$ is a $C^1$ Beltrami form on $\C$ and has compact support, then its straightening is a $C^1$ diffeomorphism.
\end{lemma}
\begin{proof}
  The condition is too strong, but this statement is enough for our purposes.
  One can check that the conditions of Theorem~7.2, page~235 of \cite{LV} hold, and its conclusion directly give that $f$ is a $C^1$-diffeomorphism.\footnote{
      According to [Bo], in \cite{V2}
      one finds a proof of the $C^1$ character of the solution if $\mu$ is Hölder.
    } 
\end{proof}

In particular, $f$ is at least $C^1$ (actually, it is better\footnote{According to Theorem~15.6.2 of \cite{AIM}, since our $\mu$ is $C^2$, the map $f$ is $C^{2+\alpha}$ for all $\alpha\in(0,1)$.}).
Recall that $m = h\Leb$ for some continuous function $h$.
It follows that $f_* m = g \Leb$ where $g$ is a continuous function: $g(z) = h(f^{-1}(z))/\det Df(f^{-1}(z))$.
As a consequence, the convolution product $\zeta = (f_* m)*\frac{1}{z}$, defines a continuous function $\zeta$. We sum this up in the following statement.

\begin{corollary}
The function $f$ is $C^1$ and the function $\zeta$ is continuous.
\end{corollary}

\medskip

We now propose an interpretation of the limit $\zeta$ with \Cref{lim:parallel-transport,thm:interpret} below.

\bigskip

\subsection{Taking a limit of the parallel transport}\label{sub:limPT}

We saw above that even though the singularities of $\zeta_n$ form a more and more dense set, their decreasing residues have a sufficiently small influence so that the sequence $\zeta_n$ weakly converges. In the theorem below, we prove that the influence is small enough to allow the convergence of the parallel transport, provided we make a small modification near the beginning and the end of the path.

We will prove in \Cref{lim:parallel-transport} a result about parallel transport along sequences of paths that are $C^1$ by parts. Note that parallel transport of a vector is well-defined for such paths when they avoid the singularities, and we have the formula
\[v(t) = v(0) \exp\left(\int_{u=0}^t \zeta_n(\gamma(u)) \gamma'(u) du\right).\]

In \Cref{lim:parallel-transport} we will assume that $\gamma'$ stays away from $0$: $\inf |\gamma'|>0$ and we will use a sequence of paths $\gamma_n$ tending to $\gamma$, to allow for more flexibility (for instance to dodge the singularities in the case the limit $\gamma$ hits some singularities of $\zeta_n$ for some $n$).
The formula above reads
\[v(t) = v(0)\exp(\tau)\]
with
\[\tau = \int_{0}^t \zeta_n(\gamma_n(u)) \gamma_n'(u) du.\]
The convergence of parallel transport can thus be stated purely in terms of this quantity.
If $\gamma_n$ does not run through a singularity, $\tau$ is a path-integral in the sense used in holomorphic functions theory:
\[\tau = \int \zeta_n(\gamma_n) d\gamma_n.\]
Path integrals can also be defined as follows for non-holomorphic functions like $\zeta$ and paths that are $C^1$ by parts:
\[\int \zeta(\gamma) d\gamma = \int_{0}^t \zeta(\gamma(u)) \gamma'(u) du\]
but note that, unless $\zeta$ is holomorphic, it does not depend only on the homotopy class of $\gamma$. 
It would be natural to conjecture the convergence of the path integrals against $\zeta_n$ to the one against $\zeta$. 
However, there is a subtlety.
We have a problem when the path starts or ends too close to a singularity $s$: first note that if the path starts or ends on a singularity then integral is the sum of a converging term and of $\res(s) \log (\gamma -s)$, and the imaginary part of $\log(\gamma-s)$ has a limit if $\gamma'$ does not vanish but not the real part.
If the path ends or starts close to a singularity for which $\res s\neq 0$, the integral will take big values, which may prevent convergence.
We will deal with this case, and also with the case where the ending or starting points are too close to a singularity, by ignoring the singularity, i.e.\ subtracting $\frac{\res s}{z-s}$ from $\zeta_n$.
Before giving the criterion for ``too close'' above, let us gauge the influence of the closest singularity.
Assuming $n$ large and $\gamma'_n$ bounded, and because $\mu$ is $C^2$, we expect the singularities to sit at the vertices of a grid that is not too distorted and whose sides have lengths of order $1/n$ (see \Cref{lem:lip}).
Recall that the residues are a $\cal O(1/n^2)$. Hence $\int \frac{\res s}{z-s}dz =  - \res s \log (z-s)$ is a $\cal O((\log n)/n^2)$ at the midpoints between two singularities $s=f(c)$ and $s'=f(c')$ corresponding to adjacent points $c$ and $c'$ in the grid of small squares.
Now if we take a point $z$ and let $s$ be the closest singularity to $z$ and call $r=|z-s|$, we get $-\res s \log (z-s) = \cal O((\log r)/n^2)$.
If we let
\[r_n = \exp(-n),\] we have that $r_n$ is much smaller than the grid size, and that $|\log r_n|/n^2 = 1/n\tend 0$.
We thus make the following definition.
\begin{definition}\label{def:zetatilden}
Let $\wt \zeta_n$, which depends on $\gamma_n$, be the rational map $\zeta_n$, which is expressed as the sum $\zeta_n(z) = \sum_s \frac{\res s}{z-s}$, from which we subtract all the terms for which $s$ is at distance $<r_n$ from $\gamma_n(0)$ or $\gamma_n(1)$.
\end{definition}
For $n$ big enough, since $r_n$ is much smaller than the distance between the residues by the forthcoming \Cref{lem:lip}, at most two terms are subtracted above, whatever the path $\gamma_n$ is.
Moreover, typically none will be subtracted and we will have $\wt\zeta_n = \zeta_n$ in this case.

\begin{theorem}\label{lim:parallel-transport}
The parallel transport has a limit in the following sense.
Consider a sequence of $C^1$ by parts paths $\gamma_n : [0,1]\to U$
that converges to a $C^1$ path $\gamma$ in the sense that both $\|\gamma_n-\gamma\|_\infty\tend 0$ and $\|\gamma'_n-\gamma'\|_\infty\tend 0$ as $n\to\infty$.
Assume that for all $n$, $\gamma_n$ avoids the singularities of $\zeta_n$.
Assume that $\forall t$, $\gamma'(t)\neq 0$.
Consider the meromorphic function $\wt \zeta_n$ as above.
Then
\[ \int \wt \zeta_n(\gamma_n) d\gamma_n \tend \int \zeta(\gamma(t)) \gamma'(t)dt .\]
\end{theorem}

The proof of this theorem comes after preparatory statements.

\begin{proposition}\label{prop:cpt}
If we endow the set of homeomorphisms of the Riemann sphere $\hat \C$ with the metric $d(f,g) = \max d(f(z),g(z))+\max d(f^{-1}(z),g^{-1}(z))$ where $d$ is the spherical metric, then the set of $K$-quasiconformal homeomorphisms of the Riemann sphere $\hat \C$ fixing $0$, $1$ and $\infty$ is compact.
\end{proposition}
This follows from Theorem~2 in Chapter~III of \cite{A} (see \cite{Abis} page~33), remarking that $f^{-1}$ is also $K$-quasiconformal.

Let us come back to our situation, with a grid of step $1/n$ and the quasiconformal maps $f_n$ that tend to $f$.
\begin{lemma}\label{lem:lip}
For any compact subset $K$ of $\C$, there exists $C,C'$ such that for all $n$,
for any two grid corners $a$, $b$ that sit in $K$, we have
\[C'|b-a| \leq |f_n(a)-f_n(b)| \leq C|b-a|.\]
\end{lemma}
\begin{proof}
Let us write
\[ f_n = \hat f_n \circ f \]
with
\[\hat f_n=f_n\circ f^{-1}
.\]
Like $f$ and $f_n$, the map $\hat f_n$ is normalized.
Since $\mu$ is $C^2$, hence $C^1$, it follows that $f$ is also a $C^1$ diffeomorphism by \Cref{lem:c1}.
 Since it is holomorphic near infinity, and since $f$ is a homeomorphism from $\C$ to itself, it follows by classical arguments that $f'(z)$ has a limit $a\in\C^*$ as $z\to \infty$. Hence $f$ is actually globally bi-Lipschitz for the Euclidean metric on $\C$.
In particular there exists $\kappa>0$ such that for any $n$ and for any two distinct grid corners $u'$ and $v'$, we have
\[ |f(u')-f(v')| > \kappa/n .\]

The map $\hat f_n$ straightens a Beltrami form
\[\hat\mu_n = f_* \mu_n\]
which satisfies $\sup_\C|\hat\mu_n| \leq M/n$ for some $M>0$.
Indeed, $\mu_n$ is on each little square the average of $\mu = f^* 0$.
Since $\mu$ is $C^2$ with compact support, hence $C^1$ with compact support, we have $\|\mu_n-\mu\|_{\infty}\leq\cst/n$.
The value of $|\hat\mu_n|(f(z))$ is given by the smooth formula $\left|\frac{a-b}{1-a\ov b}\right|$ that depends only on $a=\mu_n(z)$ and $b=\mu(z)$ and defined for $(a,b)\in\D^2$.
Since $|\mu_n|$ and $|\mu|$ are uniformly bounded away from $1$,
the variables $a,b$ remain in a compact set and hence
\[\sup_\C|\hat\mu_n| \leq M/n\]
for some $M$, as claimed.

Let $\kappa>0$ and consider two points $u$, $v\in f(K)$ with
\[|u-v|>\on{\kappa}/n\]
and denote $u_n=\hat f_n(u)$ and $v_n=\hat f_n(v)$. 
We claim that 
\[\frac{1}{a_n}\leq \frac{|u_n-v_n|}{|u-v|}\leq a_n\]
where $a_n$ depends on $\kappa$ but is independent of $u$ and $v$ and $a_n\tend 1$ as $n\to +\infty$. Once the claim is proved the lemma follows, using the bi-Lipschitz character of $f$.

To prove the claim, let us introduce a holomorphic motion parametrized by \[\lambda\in B(0,R_n),\ R_n = 1/\|\hat\mu_n\|_\infty \geq n/M,\]
 the normalized straightening $\hat f_n(\lambda,z)$ of $\lambda \hat\mu_n$.
For a fixed $z$, $\lambda\mapsto f_n(\lambda,z)$ is holomorphic.
Let us restrict to $\lambda\in B(0,R_n/2)$, then $z\mapsto\hat f_n(\lambda,z)$ is a normalized $3$-quasiconformal map.
In particular there exists $R>0$ such that for any $z\in K$ and $\lambda\in B(0,R_n/2)$,
$\hat f_n(\lambda,z) \in B(0,R)$ (see \Cref{prop:cpt}).
Consider the function $\Phi:\lambda \in B(0,R_n/2) \mapsto \hat f_n(\lambda,v)-\hat f_n(\lambda,u)$.
It takes its values in $B'=B(0,2R)-\{0\}$.
We can use then the hyperbolic metric of $B'$ whose expression is
\[ \frac{|dz|}{2|z|\log|2R/z|}
\]
because a universal cover thereof is given by $\Re z<0$, $z\mapsto 2R\exp(z)$.
The distance in $B'$ between two concentric circles $C(0,r)$ and $C(0,r')$ is
\begin{equation}\label{eq:rnb2}
\frac{1}{2}\left|\log\frac{\log(2R/r')}{\log(2R/r)}\right|
.\end{equation}
Here and below, each occurrence of $\cst$ refers to a (possibly) different constant, which is independent of $u$, $v$ and $n$.
Recall that holomorphic maps between hyperbolic Riemann surfaces are distance non-increasing.
So the hyperbolic distance in $B'$ between $\Phi(1)=v_n-u_n$ and $\Phi(0)=v-u$ is at most the hyperbolic distance from $0$ to $1$ in $B(0,R_n/2)$ so it is at most $\cst/n$. 
So by \cref{eq:rnb2} it follows that
\[ e^{-\frac{\cst}{n}}\log |2R/\Phi(0)| \leq \log |2R/\Phi(1)| \leq e^{\frac{\cst}{n}} \log |2R/\Phi(0)|
.\]
Hence
\[ (e^{\frac{\cst}{n}}-1)\log |\Phi(0)/2R|\leq \log |\Phi(1)/\Phi(0)| \leq (e^{-\frac{\cst}{n}} -1)\log |\Phi(0)/2R|
.\]
Hence
\[ \frac{\cst}{n}\log \left|\frac{\Phi(0)}{2R}\right|\leq \log \left|\frac{\Phi(1)}{\Phi(0)}\right| \leq -\frac{\cst}{n}\log \left|\frac{\Phi(0)}{2R}\right|
.\]
Now $|\Phi(0)|\geq \frac{\cst}{n}$ hence
\[ \left|\log \left|\frac{\Phi(1)}{\Phi(0)}\right|\right| \leq \cst \frac{\cst+\log n}{n} \leq \cst.\]
The claim follows.
\end{proof}

\begin{proof}[Proof of \Cref{lim:parallel-transport}]
We have
\[\wt\zeta_n(z) = \sum_{s\in \mathrm{Sing}'} \frac{\res s}{z-s}\]
where $\mathrm{Sing}'$ is the set of all actual (with non-zero residue) singularities of $\wt\zeta_n$ (with $\infty$ omitted), which is the same as the set $\mathrm{Sing}$ of actual singularities of $\zeta_n$ with $0$, $1$ or $2$ points removed and $\infty$ omitted.
Let us write
\[I_n:=\int \wt\zeta_n(\gamma_n)d\gamma_n = \sum_{s\in\mathrm{Sing}'} \res (s) \int \frac{d\gamma_n}{\gamma_n-s} .\]
Then,
\[I_n = \sum_{s\in\mathrm{Sing}'} \res (s) \Psi_n(s) = \langle \tilde m_n',\Psi_n \rangle\]
where $\Psi_n(s) = \int \frac{d\gamma_n}{\gamma_n-s}$, and is
a determination of $\log\left(\frac{\gamma_n(1)-s}{\gamma_n(0)-s}\right)$ and where $\tilde m'_n$ is the sum of Dirac masses at points $s\in \mathrm{Sing}'$ ponderated by $\res s$; note that $\tilde m'_n$ is also equal to $(f_n)_* m_n$ with at most two residudes removed.

We have for $n$ big enough
\[|\Im\Psi_n|\leq C_0\]
where $C_0$ is independent from $\eps$ and $n$.
\begin{proof}
There exists a real $m>0$ such that for all $n$, $\min |\gamma'_n| \geq m$.
There exists $\eta>0$ such that for all $n$ big enough, $|t-s|\leq\eta$ $\implies$ $|\gamma'_n(t)-\gamma'_n(s)|\leq m/10$ (by uniform convergence of $\gamma'_n$ to a continuous function).
For any $z\in\C$, on an interval $[t_0,t_1]$ of values of $t$ of size $\eta$ (or less), the argument of $\gamma'_n(t)$ deviates less than $\asin(1/10)$ from its value at $t=t_0$ and it follows that the variation of argument of $\gamma_n(t)-s$ is at most $\pi +\asin(1/10)$ on $[t_0,t_1]$. 
Slicing $[0,1]$ into at most $\left\lceil\frac{1}{\eta}\right\rceil$ intervals of length at most $\eta$, we get $|\Im\Psi_n|\leq \left\lceil\frac{1}{\eta}\right\rceil (\pi + \asin(1/10))$.
\end{proof}

\medskip

For $\eps>0$ consider the set $G_n = G_n(\eps) = V_{\eps}(\gamma_n)$ defined as the $\eps$-neighbourhood of $\gamma_n$.
By \Cref{prop:cpt}, the sets $f_n^{-1}(G_n)$ are all contained in some common closed ball $K=\ov B(0,R)$.
By increasing $R$, we also assume that $K$ contains for all $n$ the corners of the litte squares on which the average of $\mu$ is not $0$. In particular, $f_n(K)$ contains $\on{Sing}$ and $G_n$.

Let us write
\[ I_n = J_n + J'_n\]
where
\[ J'_n = \sum_{s\in G_n\cap\mathrm{Sing}'} \res(s) \Psi_n(s) \]
\[ J_n = \sum_{s\in T_n} \res(s) \Psi_n(s) \]
where $T_n = \mathrm{Sing}'-G_n$.

\smallskip
\noindent \textbf{Note:} In the rest of this proof, inequalities involving constants $C_k$ for some $k\in\N$ will appear. They will be valid for all $\eps$ and $n$ subjected to the conditions $\eps<\eps_k$ and $n>n_k$, where $\eps_k$ and $n_k$ are constants that will be implied.

\medskip \noindent Bound on $J'_n$:\smallskip 

We have
\[\on{Leb} G_n \leq C_1\eps\]
for some $C_1$,
and for all $s\in\mathrm{Sing}$, by \cref{eq:rnb0,eq:lambda-expand}:
\begin{equation}\label{eq:rnb3}
|\res s| \leq C_2/n^2
\end{equation}
for some $C_2$. %where $C_1$ and $C_2$ are independent from $n$ big enough and $\eps$ small enough as explained above.

Because of \Cref{lem:lip}, the disks of radius $C_3/2n$ centered on the actual singularities are disjoint, because their preimage by $f_n$ lie in $K$.
So this is the case for the $s\in G_n\cap \mathrm{Sing}'$.
Moreover, as soon as $C_3/2n < \eps$ then these disks are contained in $G_n(2\eps)$, so the sum of the areas of these disks is at most $2C_1\eps$.
It follows that
\[\#(G_n\cap \mathrm{Sing}') \leq \frac{2C_1\eps}{\pi(C_3/2n)^2} = C_4 \eps n^2\]
for some $C_4$.

Case 1: $\gamma_n(0)=\gamma_n(1)$. Then $\Psi_n$ is purely imaginary and hence is bounded, according to one of the remarks above. Hence
\[|J'_n| \leq \sum_{s\in\mathrm{Sing'}\cap G_n} |\res(s) \Psi_n(s)| \leq \#(\mathrm{Sing'}\cap G_n) \times \max |\res s| \times \max |\Psi_n(s)|\]
\[ \leq  C_4 \eps n^2 \times \frac{C_2}{n^2} \times C_0 = C_4 C_2 C_0 \eps.\]

Case 2: $\gamma_n(0)\neq\gamma_n(1)$.

We have
\[|\Psi_n(s)|< \log \frac{1}{|\gamma_n(0)-s|} + \log \frac{1}{|\gamma_n(1)-s|} + C_5\]
for some $C_5$.
\begin{proof}
We have already seen that $|\Im\Psi_n(s)|\leq C_0$.
The real part is easier to deal with, since 
\[\Re\Psi_n(s) = \log \frac{|\gamma_n(1)-s|}{|\gamma_n(0)-s|} = \log \frac{1}{|\gamma_n(0)-s|} - \log \frac{1}{|\gamma_n(1)-s|}
.\]
For any $t$, $\left|\log \frac{1}{|\gamma_n(t)-s|}\right| \leq \log \frac{1}{|\gamma_n(t)-s|} + C$ for $C=\max(0,2\log \on{diam} B)$ where $B$ is a compact set containing all the curves $\gamma_n$ and the singularities of non-vanishing residue of $\zeta_n$ for all $n$. The constant $C$ is independent of $\eps$ and $n$.
\end{proof}

For any singularity $s$ at distance $\geq\eps$ from $\gamma_n(0)$ we have
\[\log\frac{1}{|\gamma_n(0)-s|} \leq \log \frac{1}{\eps}.\]

For the singularities $s\in \mathrm{Sing}'$ at distance $<\eps$ from $\gamma_n(0)$ we have $r_n \leq |s-\gamma_n(0)|<\eps$ with $r_n = \exp(-n)$.
Recall that the disks $B(s,C_3/2n)$ are disjoint.
There is hence at most one $s_0\in\mathrm{Sing}'$ with $|s_0-\gamma_n(0)|<C_3/2n$ and since $|s_0-\gamma_n(0)|>e^{-n}$ we have
\[\log\frac{1}{|\gamma_n(0)-s_0|}\leq n.\]
For the other ones, let $k$ be the smallest integer such that $2^k \frac{C_3}{2n} > \eps$. 
Recall that we assume that $n$ is big enough so that $C_3/2n<\eps$, hence $k>0$. We have
\[2^k\leq \frac{4 n \eps}{C_3}.\]
Let $j\in\N\cap [0,k)$. The number of $s\in\mathrm{Sing}'$ such that $|\gamma_n(0)-s|$ belongs to $[2^j C_3/2n, 2^{j+1}C_3/2n)$ is at most $\big((2^{j+1}+1)^2-(2^j-1)^2\big) = 3\times 4^j +6\times 2^j \leq 9\times 4^j$, because the union of the disjoint disks $B(s,C_3/2n)$ is contained in the round annulus of $z\in\C$ such that $(2^{j+1}+1)C_3/2n<|z-s|<(2^j-1)C_3/2n$.
For such an $s$, we have
\[\log\frac{1}{|\gamma_n(0)-s|}\leq \log \frac{2n}{2^j C_3}.\]

Hence
\[\sum_{s\in\mathrm{Sing'}\cap B(\gamma_n(0),\eps)} \log\frac{1}{|\gamma_n(0)-s|} \leq n+\sum_{j=0}^{k-1} 9\times 4^j \times\log \frac{2n}{2^jC_3}\]

Now
\[\sum_{j=0}^{k-1} 4^j \times \log \frac{2n}{2^jC_3}
= 4^k\sum_{j=0}^{k-1} \frac{1}{4^{k-j}} \log 2^{k-j}\frac{2n}{2^kC_3}
=4^k \sum_{p=1}^{k} \frac{1}{4^p} \left(p\log 2 + \log\frac{2n}{2^kC_3}\right)
\]
\[\leq 4^k \sum_{p=1}^{k} \frac{1}{4^p} \left(p\log 2 + \log\frac{1}{\eps}\right) \leq 4^k \left(\frac{4}{9}\log 2+\frac{1}{3}\log \frac{1}{\eps}\right)\]
\[\leq \left(\frac{4n\eps}{C_3}\right)^2 \times \left(\frac{4}{9}\log 2+\frac{1}{3}\log \frac{1}{\eps}\right)
\]
Hence
\[\sum_{s\in\mathrm{Sing'}\cap B(\gamma_n(0),\eps)} \log\frac{1}{|\gamma_n(0)-s|} \leq n+C_6 n^2 \eps^2 + C_7 n^2\eps^2 \log\frac{1}{\eps}.
\]
A similar identity holds near $\gamma_n(1)$ in place of $\gamma_n(0)$.

Let us put it all together (assuming $\eps< 1$): using $\max |\res s| \leq C_2/n^2$ we have
\[ |J'_n| \leq \frac{C_2}{n^2} \left(\sum \log\frac{1}{|\gamma_n(0)-s|} +
\sum \log\frac{1}{|\gamma_n(1)-s|} + \sum C_5\right);
\]
separating the sum of $\log\frac{1}{|\gamma_n(0)-s|}$ according to whether or not $|\gamma_n(0)-s|\geq\eps$,
doing the same with the terms involving $\gamma_n(1)$,
and then regrouping the terms and using that there are at most $C_4\eps n^2$ values of $s$ as we saw earlier:
\[ |J'_n| \leq \frac{C_2}{n^2} \left(
C_4\eps n^2\times \Big(2\log\frac{1}{\eps}+C_5\Big)
+2\Big(n+C_6 n^2 \eps^2 + C_7 n^2\eps^2 \log\frac{1}{\eps} \Big)\right)
\]
\[=C_2\times \left(
C_4\eps \times \Big(2\log\frac{1}{\eps}+C_5\Big)
+2\Big(\frac{1}{n}+C_6 \eps^2 + C_7 \eps^2 \log\frac{1}{\eps}\Big)
\right).\]
Hence
\[\limsup |J'_n| \leq C_8\eps\log \frac{1}{\eps}\]
for some $C_8$.
The important fact is that this bound tends to $0$ as $\eps\tend 0$.

\medskip \noindent Bound on $J_n$:\smallskip 

Similarly to the definition of $G_n$ as the $\eps$-neighbourhood of $\gamma_n$, we let $G$ be the $\eps$-neighbourhood of $\gamma$.
For $s\notin \gamma([0,1])$, let $\Psi(s) = \int \frac{d\gamma}{\gamma-s}$, which is also the limit of the functions $\Psi_n$.
Similarly as for $G_n$ we have that $|\Im \Psi|\leq C_0$ and
\[\Re\Psi(s) = \log \frac{|\gamma(1)-s|}{|\gamma(0)-s|} = \log \frac{1}{|\gamma(0)-s|} - \log \frac{1}{|\gamma(1)-s|}
.\]

Note that $|\Psi(z)|$ and $|\Psi_n(z)|$ tend to $0$ as $|z|\tend+\infty$.
There exists a sequence $\tilde\Psi_n$ of continuous extensions of the restriction of $\Psi_n$ to the closed set $\C-G_n$ and a continuous extension $\tilde\Psi$ of the restriction of $\Psi$ to the closed set $\C-G$, such that moreover $\tilde\Psi_n$ tends to $\tilde\Psi$ uniformly on compact subsets of $\C$ and such that $\|\tilde\Psi\|_{\infty} = \|\Psi|_{\C-G}\|_{\infty}$.
This follows from the Tietze extension theorem (Theorem~35.1 page~219 of \cite{Mun}). Apply it first to extend $\Psi|_{\C-G}$ into $\tilde\Psi$.
Then apply it again to the metric space $X = (\{0\}\cup\setof{1/n}{n\in\N^*})\times \C$, the continuous function that is the disjoint union of $\tilde\Psi$ and the restrictions of the $\Psi_n$ (the disjoint union of the domains of these maps is indeed a closed subset of $X$).

Now
\[J_n - \int \tilde m'_n \tilde\Psi_n = -\hskip-.5cm\sum_{s\in G_n\cap \mathrm{Sing}'}\hskip-.5cm\res(s) \tilde\Psi_n(s)\]
Hence
\[\left|J_n - \int \tilde m'_n \tilde\Psi_n\right| \leq \max_{s\in  \mathrm{Sing}'} |\res(s)| \times \#(G_n\cap \mathrm{Sing}')\times\max_{G_n} \big|\tilde\Psi_n\big| \]
\[\leq \frac{C_2}{n^2}\times C_4 n^2 \eps^2 \times \max_n\max_{G_n} |\tilde\Psi_n| \leq C_7 \eps^2\]
for some $C_7$, using \cref{eq:rnb3} to bound $|\res s|$.
Also,
\[\left|\int \tilde m'_n \tilde\Psi_n- \int \tilde m'_n \tilde\Psi\right| \leq \|\tilde m'_n\| \|\tilde\Psi_n-\tilde\Psi\|_{\infty}
\]
and as $n\tend+\infty$, the left factor stays bounded and the right one tends to $0$.
Let $m' = f_* m$.
By \Cref{cor:lim-mes-2} (removing two singularities does not change the convergence, since the total mass removed is $\cal O(1/n^2)$) we have
\[\int \tilde m'_n \tilde\Psi \tend \int m' \tilde\Psi\]
when $n\tend+\infty$.
Recall that $\|\tilde\Psi\|_{\infty} = \|\Psi|_{\C-G}\|_{\infty}$.
Since $\Psi$ is a branch of $\log\left(\frac{\gamma(1)-s}{\gamma(0)-s}\right)$, and since on $\C-G$ the distance to $\gamma$ is at least $\eps$, we get for all $s\in\C-G$:
$|\Re\Psi(s)| \leq C_8+\log \frac{1}{\eps}$
for some $C_8>0$. Since, as we saw, $|\Im \Psi|\leq C_0$ we get
\[\|\tilde\Psi\|_{\infty} \leq C_0+C_8+\log \frac{1}{\eps}.\]
Now $m=h\Leb$ where $h$ is a continuous function so $m'=f_* m =g\Leb$ with $g(f(z)) = h(z) / \det Df$ is continuous too.
In particular $g$ is bounded, say by $C_9$.
We have $\Leb G\leq C_1 \eps$, like for $G_n$.
It follows that
\[\left|\int_G m' \tilde\Psi\right|\leq C_9 \eps \log \frac{1}{\eps}\]
for some $C_9$, and this quantity tends to $0$ as $\eps\tend 0$. From the expression of $\Re \Psi$ and the bound on $\Im\Psi$, we get
\[\left|\int_G m' \Psi\right|\leq C_{10} \eps + C_9\eps^2\log \frac{1}{\eps}.\]

If we collect all the results, we know that for all $\eta>0$, by taking $\eps$ small enough, then for all $n$ big enough, the quantity
$I_n$ to be evaluated, differs from $\int_\C m'\Psi$ by at most $\eta$, i.e.\ 
\[\left|I_n-\int_\C m'\Psi\right|<\eta.\] 
This integral $\int_\C m'\Psi$ is equal to $\int \zeta(\gamma(t)) \gamma'(t)dt$:
indeed, recall that $\zeta=m'*\frac{1}{z}$,
whence
\[\int \zeta(\gamma(t)) \gamma'(t)dt
= \int_{t\in[0,1]}\left(\int_{s\in \C} \frac{1}{\gamma(t)-s} dm'(s)\right)\gamma'(t)dt
\]
\[ = \int_{s\in \C}\left(\int_{t\in[0,1]} \frac{d\gamma(t)}{\gamma(t)-s} \right) dm'(s) = \int_\C m' \Psi.\]
The permutation in the integral is justified because the integral of the absolute value
\[\int_{s\in \C} \frac{1}{|z-s|} dm'(s) = \int_{s\in \C} \frac{g(s)}{|z-s|} d\Leb(s)\]
is convergent ($g$ is bounded with compact support, $z\mapsto 1/|z|$ is $L^1$ on any compact set) and bounded for $z\in \gamma([0,1])$, so we are in the conditions of the Fubini Tonelli theorem.
\end{proof}

\begin{remark}
In the statement of \Cref{lim:parallel-transport}, we decided to restrict to pathes avoiding the singularities and to remove the 0, 1 or 2 singularities that are too close to the extremities of the path.
We could have proceeded differently.
For instance we could have allowed for every paths by either:
\begin{itemize}
\item removing the 0, 1 or 2 aforementioned singularities, and for every other singularity $s$ crossed by $\gamma_n$, give a meaning to $\res (s)\int \frac{d\gamma_n}{\gamma_n-s} $ if $\gamma_n'$ is continuous at all $t$ for which $\gamma_n(t)=s$, or if the jump of $\arg \gamma_n'$ at $t$ is smaller than, say, $\pi/2$;
\item chopping off $\zeta_n$ close to some/all poles, for instance by setting $\zeta_n(z)=0$ when $\exists s$ singularity such that $|z-s|<r'_n$, for an appropriate sequence $r'_n$;
\item there are several choice in $r'_n$ and the way $\zeta_n$ is chopped off; for instance one can arrange so that $\zeta_n$ remains continuous;
\item one could instead subtract all the poles of $\zeta$ that are at distance $<r'_n$ to $\gamma$, for an appropriate sequence $r'_n$.
\end{itemize}
Other approaches are possible.
\end{remark}

\subsection{Interpreting the limit}\label{sub:interpretLim}

Let us sum up what we have done up to now. Recall the similarity surfaces $\cal S_n$ associated to our approximating Beltrami forms $\mu_n$, and recall their completions $\cal S'_n$ into a Riemann surface isomorphic to $\hat\C$.
To this global chart $\hat \C$ we associated a singular connection, of symbol $\zeta_n$ that is a meromorphic function\footnote{\ldots\ A rational map in fact, since for each $n$ there are only finitely many non-erasable poles.} all whose poles are all simple, but get more and more densely packed as $n$ tends to infinity.
We proved that $\zeta_n$ has a weak limit\footnote{In the previous section, \Cref{lim:parallel-transport}, we improved the weak limit statement into a statement in terms of limit of associated parallel transport.} and that this weak limit is the function $\zeta$ introduced in \Cref{prop:wlz} via a convolution product.

\begin{theorem}\label{thm:interpret}
Let $(e_1,e_2)$ be the canonical basis of $\R^2$.
Let $A$ be the push-forward by $f$ of the constant vector field $e_1$ of $\R^2$, and $B$ the push-forward of the constant vector field $e_2$.
Note that $A$ and $B$ commute.
Consider the conformal symmetric connection of symbol\footnote{The notation $\partial_A B$ refers to the directional derivative, see \Cref{sub:qi}.}
\[ \tilde\zeta(z) = -\frac{\partial_A B}{AB}
\]
defined on $U$.
Then $\tilde\zeta$ is equal to the symbol $\zeta$ introduced in \Cref{prop:wlz}:
\[\tilde\zeta=\zeta.\]
\end{theorem}
The connection of symbol $\tilde\zeta$ expresses as
\[\nabla_X Y = \partial_X Y - \partial_A B \frac{XY}{AB}.\]

Before we prove the theorem, let us introduce an easy lemma.
\begin{lemma}\label{lem:AB}
Let $A$ and $B$ be two vector fields in an open subset $U\subset \R^2\equiv \C$. Assume that they vanish nowhere. 
Then there is a unique conformal symmetric connection $\nabla$ such that the parallel transport under integral lines of $A$ leaves $B$ invariant.
The symbol of $\nabla$ in the chart $U$ is
\[-\frac{\partial_A B }{AB}
.\]
\end{lemma}
\begin{proof}
A conformal symmetric connection $\nabla$ takes the form
\[\nabla_X Y = \partial_X Y + \zeta XY\]
where $\zeta$ is a continuous function from $U$ to $\C$.
The condition expresses as the pair the equation
\[\nabla_A B=0,\]
i.e.\ 
\[\partial_A B + \zeta AB = 0.\]
\end{proof}

Incidentally, if $A$ and $B$ commute, then
\[-\frac{\partial_A B }{AB} = -\frac{\partial_B A }{AB},\]
i.e.\ the unique conformal symmetric connection such that the parallel transport under integral lines of $A$ leaves $B$ invariant coincides with the unique conformal symmetric connection such that the parallel transport under integral lines of $B$ leaves $A$ invariant.

We will also need the following, harder to get, result.
Let, as earlier,
\[r_n =\exp(-n)\]
and recall that $f$ is $C^1$ (\Cref{lem:c1}).
\begin{proposition}\label{prop:supDD}
\[\sup_{\on{dist}(z,\frac{1}{n}\Z)>r_n} \hskip-.5cm \|D_z f_n-D_z f\| \underset{n\to\infty}\tend 0
\]
where the norm is any operator norm on the set of $\R$-linear self-maps of $\C$.
\end{proposition}
\begin{proof}
In \Cref{app:pf:supDD}.
\end{proof}

\begin{proof}[Proof of \Cref{thm:interpret}]
By \Cref{lem:AB}, to prove the claim, it is enough to prove that for $\zeta$, the associated parallel transport along integral lines of $A$ leave $B$ invariant, where $A = f_* e_1$ and $B = f_* e_2$.

We will use the fact that this parallel transport is the limit of the parallel transport associated to $\zeta_n$ (\Cref{lim:parallel-transport}).
Let $\gamma : [0,1]\to\C$ defined by
\[ \gamma(z) = f(x_0+(x_1-x_0)t+iy_0) \]
for some $\ell>0$ and let for $n$ big enough
\[ \gamma_n(t) = f_n(x_0+(x_1-x_0)t +i y_n) \]
where $y_n$ is chosen so that $y_n\tend y_0$ and $d(y_n,\frac{1}{n}\Z)>r_n$: for instance $y_n=(0.5+\lfloor n y_0\rfloor)/n$.

The path $\gamma$ is $C^1$ and the path $\gamma_1$ is $C^1$ by parts.
By \Cref{app:beltrami:cv}, $\gamma_n\tend\gamma$ uniformly,
and by \Cref{prop:supDD} $\gamma_n'\tend\gamma'$ uniformly.

\begin{figure}
\begin{tikzpicture}
\node at (0,0) {\includegraphics[width=12cm]{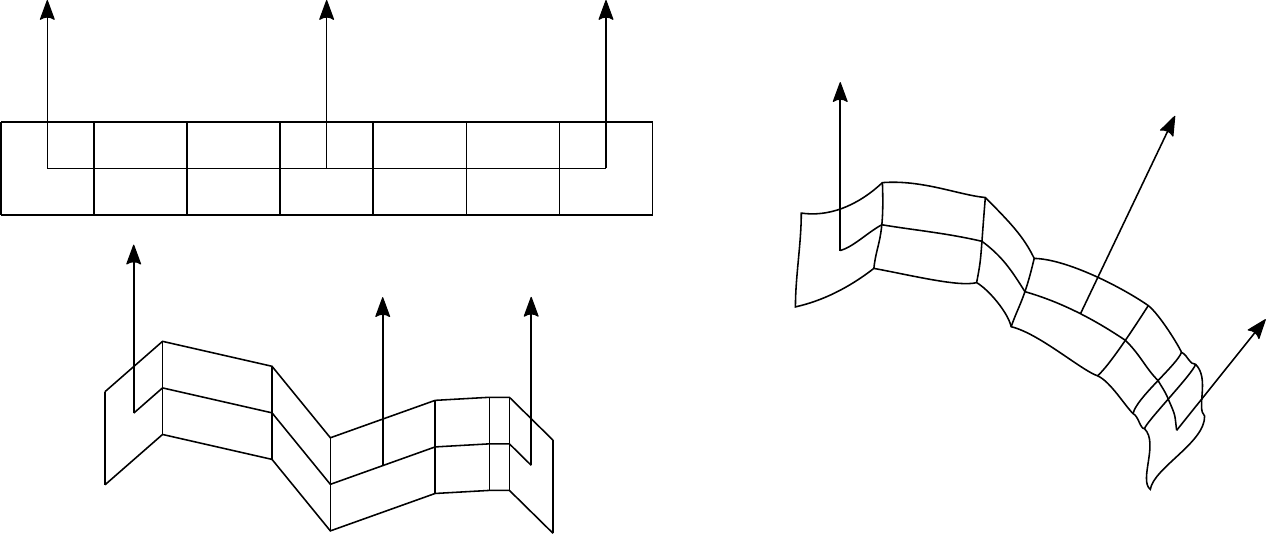}};
\draw[->] (0.5,0.8) --node[above right]{$f_n$} (1.3,0.2);
\end{tikzpicture}
\caption{Illustration of an argument in the proof of \Cref{thm:interpret}.}\label{fig:AB}
\end{figure}

Now, recall that the parallel transport for $\zeta_n$ is also the parallel transport associated to the similarity surface $\cal S_n$ built from gluing together the images of the small squares by maps of the form $a z + b \bar z$ where $a$ and $b$ are any complex constants such that $b/a\in\D$ is the value of $\mu_n$ on this square. The image of such a square is a parallelogram
The horizontal line $x_0+(x_1-x_0)t +i y_n$ runs through a row of those squares. Consider the corresponding parallelograms.
We can place the initial parallelogram so that the image of the vertical edge of the square is a vertical edge with the same orientation and same size $1/n$.
Then its opposite side is parallel and has the same size.
If we scale rotate and place the other parallelograms in this row so as to match this side, all the images of the vertical sides remain vertical with the same length $1/n$. See \Cref{fig:AB}.
One deduces that if we define the vector $v_n$ as the image by $f_n$ of the vertical unit vector $e_2$ attached to $x_0 +i y_n$, then
the parallel transport by $\nabla_{\zeta_n}$ of $v_n$ along $\gamma_n$ remains for all $t$ the image $v_n(t)$ by $f_n$ of the vertical vector $e_2$ attached to $x_0+(x_1-x_0)t+i y_n$.
Now by \Cref{prop:supDD} and continuity of $Df$, for each fixed $t$, the vector $v_n(t)$ tends as $n\tend+\infty$ to the push-forward $v(t)$ by $f$ of $e_2$ attached to $x_0+(x_1-x_0)t+iy_0$. Note that $v(t) = B(\gamma(t))$.

We can now invoke \Cref{lim:parallel-transport}: at the limit, the parallel transport along $\gamma$ of $v(0)=B(\gamma(0))$ is the vector $v(t)=B(\gamma(t))$.
\end{proof}

Unlike the case of $\tilde\zeta$, we have not been able to give an interesting expression of the curvature form $\omega=\bar\partial \tilde\zeta$ in terms of $A$ and $B$: a direct application of $\bar\partial$ to the expression $(\partial_B A)/AB$ seems to give rise to a complicated expression.

%
% ----------------------------------------------------------
%

\section{On averaging the Beltrami form}\label{sec:var-avg}

The approximation scheme that we introduced to prove the measurable Riemann mapping theorem and its version with holomorphic dependence, involves taking averages on squares of Beltrami forms coefficients $\mu$.
We may wonder, if we change the averaging method (in some way specified below), whether the results of this article still hold or break.
\Cref{sub:chmeth} contains positive results (convergence) and \Cref{sub:non-hol} negative ones (holomorphy).

\subsection{Changing the averaging method}\label{sub:chmeth}

In this article, up to now we used the average of the complex-valued function $\mu(z)$, assumed $L^\infty$ in \Cref{part:beltrami} or $C^2$ in \Cref{sec:sslim}, over the Lebesgue measure on squares.
But the Beltrami differential is just one way to represent ellipse fields. We could have taken another quantity $\nu=\nu(\mu)$, like $(1+\mu)/(1-\mu)$, or even functions of $\mu$ that are not holomorphic like $\frac{1+|\mu|}{1-|\mu|}e^{i\arg\mu}$.

So, what if we averaged $\nu$ instead of $\mu$? In other words, given $n$, we take $\mu_n$ on each small square $S$ in the construction to be constant equal to $\nu^{-1}\left(\frac{\int \nu\circ \mu}{\int 1}\right)$, where the integrals are still against the Lebesgue measure on $S$.
Note that the quantity $\frac{\int \nu\circ \mu}{\int 1}$ belongs to the convex hull of the image of $\nu$.

We assume that $\mu\mapsto\nu(\mu)$ is a $C^1$-diffeomorphism from the unit disk to a convex  open subset of $\R^2$ and ask two questions:
\begin{enumerate}
\item\label{item:avg:1} Does the sequence $f_n$ of normalized straightenings of $\mu_n$ still converge to the straightening of $\mu$?
\item\label{item:avg:2} Assuming that $\nu$ is $C^2$, and that $\mu$ is $C^2$ with compact support (the setting of \Cref{sec:sslim}), are the limit complex measure $m$, curvature form $\omega$ and Christoffel symbol $\zeta$ the same, and does the parallel transport still converge?
\end{enumerate}
The answer to both questions is: yes.

For point \eqref{item:avg:1}, the proof of convergence (\Cref{app:beltrami:cv}) relied on uniqueness of the straightening of $\mu$ (\Cref{thm:mainuniq}) and on the sole facts that $\|\mu_n\|_\infty$ is bounded away from $1$ and that locally, the $L^1$ norm of $\mu_n-\mu$ tends to $0$ (\Cref{lem:mnm}), which still hold: indeed $\nu\circ\mu_n$ is the average of the bounded measurable function $\nu\circ\mu$ on the little squares of the $n$-the generation; using the same proof as in \Cref{lem:mnm} we get that for any compact subset $S$ of $\C$, $\|\nu\circ\mu_n-\nu\circ\mu\|_{L^1(S)}\tend 0$ as $n\tend +\infty$.
The set $\nu\circ\mu(C)$ in $\C$ is compactly contained in the image of $\nu$ hence its convex hull $C$ too.
It follows that $\|\mu_n\|_\infty$ is bounded away from $1$. It also follows that the map $\nu^{-1}$ is uniformly Lipschitz on $C$, hence $\|\mu_n-\mu\|_{L^1(S)} = \|\nu^{-1}\circ\nu\circ\mu_n-\nu^{-1}\circ\nu\circ\mu\|_{L^1(S)} \leq c \|\nu\circ\mu_n-\nu\circ\mu\|_{L^1(S)}$ for some constant $c>0$ independent of $n$.

For point \eqref{item:avg:2}, recall that $\mu$ is in this case also $C^2$, hence continuous.
So in four squares around a corner, $\mu$ takes values that remain close to the value $\mu_0$ of $\mu$ at the corner.
Below, all the $o(\cdot)$ are uniform in $|m|$ small enough, $\mu_0$ subject to $|\mu_0|<\kappa<1$ and $z_0\in\C$.
Let 
\[\nu(\mu_0+m) = \nu_0 + Q_1(m)+Q_2(m)+o(|m|^2)\]
be an expansion where $Q_k$ is a degree-$k$ $\R$-homogeneous polynomial in $(\Re m,\Im m)$ taking values in $\C$ seen as a dimension 2 $\R$-vector space.
Let
\[\mu(z_0+z) = \mu_0 + P_1(z)+P_2(z)+o(|z|^2)\]
with a similar convention.
Note that $P_1$ and $Q_1$ are $\R$-linear endomorphisms of $\C$ and that $Q_1$ is invertible by hypothesis that $\nu$ is a diffeomorphism.
Then
\[ \nu^{-1}(\nu_0+v) = \mu_0 + Q_1^{-1}(v) - Q_1^{-1}\circ Q_2\circ Q_1^{-1}(v)+o(|v|^2),\]
and
\[ \nu\circ\mu(z_0+z) = \nu_0 + Q_1\circ P_1(z)+Q_1\circ P_2(z)+Q_2\circ P_1(z) + o(|z|)^2. \]
Let $M_j$ be the average of $\mu$ on the squares $C_j$ of side length $\eps$ appearing in the proof of \Cref{prop:lim-mes}.
We have
\[M_j = \mu_0 + P_1\big(\avgx{C_j} z\big) + \avgx{C_j}\big( P_2(z)\big) + o(\eps^2)\]
Denote $N_j$ the average of $\nu\circ\mu$ on $C_j$.
We get
\[N_j = \nu_0 + Q_1\circ P_1\big(\avgx{C_j} z\big) + Q_1\Big(\avgx{C_j} P_2(z)\Big)+  \avgx{C_j} \big(Q_2\circ P_1(z)\big)+ o(\eps^2)\]
Let $M_j' = \nu^{-1}(N_j)$.
Composing the previous formula with the expansion of $\nu^{-1}$ and comparing to $M_j$ we get
\[ M'_j= M_j + Q_1^{-1}\Big(\avgx{C_j} \big( Q_2\circ P_1(z)\big)-Q_2\circ P_1\big(\avgx{C_j} z\big)\Big) + o(\eps^2)\]
We have $Q_2\circ P_1(x+iy) = a_{xx} x^2+2a_{xy}xy+a_{yy} y^2$ for some coefficients $a_{\ldots}\in\C = \R+i\R$.
Then \[\avgx{C_j} \big( Q_2\circ P_1(z)\big)-Q_2\circ P_1\big(\avgx{C_j} z\big) = (a_{xx}+a_{yy})\left(\frac13-\frac14\right)\eps^2 = B \eps^2\]
where $B\in\C$ depends on $z_0$.
So
\[M'_j = M_j + Q_1^{-1}(B)\eps^2 + o(\eps^2).\]
Note that the term $Q_1^{-1}(B)\eps^2$ depends on $z_0$ but is independent of $j$.
Then the monodromy factor of the singularity for the new scheme is
\[\Lambda' = \frac{L'_1\times L'_3}{L'_0\times L'_2}\]
where $L'_j = (1+M'_j)(1-M'_j)$,
whereas for the original one it is
\[\Lambda = \frac{L_1\times L_3}{L_0\times L_2}\]
where $L_j = (1+M_j)(1-M_j)$.
Recall that $1+M_j$ is close to $1+\mu_0$ when $\eps$ is close to $0$.
A computation gives
\[\frac{L'_j}{L_j} = 1+Q_1^{-1}(B)\left(\frac{1}{1+\mu_0}+\frac{1}{1-\mu_0}\right)\eps^2+o(\eps^2) = 1 + B' \eps^2 + o(\eps^2)\]
Where $B'$ depends on $z_0$ but not on $j$.
As a consequence,
\[\frac{\Lambda'}{\Lambda} = 1 + o(\eps^2).\]
Hence the limit of $\log_p\Lambda'/\eps^2$ is the same as the limit of $\log_p\Lambda/\eps^2$.
One can check that the convergence is still uniform.
The proof of \Cref{prop:wlz} adapts verbatim to the new situation so we still have 
\[\zeta_n \wtends \zeta := (f_* m)*\frac{1}{z} \]
for the new $\zeta_n$, where $\zeta$ and $m$ stay the same.
In \Cref{lem:lip} one argument is that $\sup_\C|\hat\mu_n| \leq M/n$ for some $M>0$.
It is based on the fact that $\|\mu_n-\mu\|_\infty = \cal O(1/n)$, which is still true here ($\nu\circ\mu$ is $C^2$, hence $C^1$, with compact support, so its variation on a square of side $1/n$ is $\cal O(1/n)$; hence its average will differ from the central value by at most $\cal O(1/n)$ (in fact we have better than that); the preimage of this average by $\nu$ also differs by the same order since everything takes place in a compact subset of the domain of $\nu^{-1}$ which is at least $C^1$).
Several of the arguments of \Cref{lim:parallel-transport}, are based on the fact that that the residues are $\cal O(1/n^2)$, and this is still the case here. Also, we have seen that $f_n$ still tends to $f$.
So Point~\eqref{item:avg:2} holds.

\subsection{Non-holomorphy}\label{sub:non-hol}

On the other hand, in the case of holomorphic families of Beltrami forms, if $\mu\mapsto\nu(\mu)$ is not assumed holomorphic, we lose the holomorphic dependence, with respect to the parameter, of $\mu_n$. For instance let $\nu(x+iy) = x+x^3+iy$ and $\mu(x+iy)=\tau x$ for $|x+iy|<0.1$ where $\tau=u+iv$ is a complex parameter in $\D$ (choose any $C^2$ extension for $|x+iy|>0.1$).
Then when $n>100$, near $0$ we have on the square $\Re z \in [0,1/n]$, $\Im z\in [0,1/n]$ that
$\mu_n = P^{-1}(\frac{u}{2n}+\frac{u^3}{4n^3}) + \frac{v}{2n}i$
where $P(x)=x+x^3: \R\to \R$. This quantity does not depend holomorphically on $u+iv$.

This has the consequence that we also loose holomorphic dependence in the parameter of the normalized straightening $f_n$ of $\mu_n$. Indeed, if $f_n$ depends holomorphically on $\tau$, then we prove below that $\mu_n =\bar \partial f_n / \partial f_n$ also depends holomorphically on $\tau$, leading to a contradiction.

In fact we will prove a slightly more general lemma, probably already known:
\begin{lemma}\label{lem:fhmh}
Let $\tau\in B(0,\eps) \to f_\tau$ be a family of normalized $K$-quasiconformal homemomorphisms of $\C$ for some $K>1$.
Assume that for all $z$, the map $\tau\mapsto f_\tau(z)$ is holomorphic.
Then the Beltrami differential $\mu_\tau$ of $f_\tau$, as an element of the Banach space $\cal B = L^\infty(\C)$, depends analytically on $\tau$ in the sense of \Cref{prop:muhol}.
\end{lemma}
\begin{proof}
It is enough to prove holomorphic dependence near $\tau=0$.
We fix some ball $B$ in $\C$ that we will allow to grow near the end of the proof.
For all test function $\phi$ we have $\langle f_\tau,\phi\rangle$ that depends holomorphically too ($f_\tau$ satisfies local bounds that are uniform w.r.t.\ $\tau$ varying in a compact set).
Applying this to a test function that is a partial derivative of another test function, we get that criterion~\eqref{bw5} of \Cref{prop:holos} is satisfied for the Banach space $L^1(B)$ (a bound on $D f_\tau$ in $L^1(B)$ is needed: it follows for instancee from the $L^2$ bound of \Cref{lem:L2}).
It follows that the distribution derivatives of $f_\tau$, seen as $L^1$ functions on $B$, are holomorphic w.r.t.\ $\tau$ in the sense of criterion \eqref{bw1} of \Cref{prop:holos}, and so are $\partial f_\tau$ and $\bar\partial f_\tau$: there are $a_n$, $b_n$ in $L^1(B)$ and some $r>0$ such that for all $\tau$ with $|\tau|<r$, $\partial f_\tau = \sum a_n \tau^n$ and $\bar\partial f_\tau = \sum b_n \tau^n$ in the sense of distributions, with $\sum \left(\int_B|a_n|\right) r^n<+\infty$ and $\sum \left(\int_B|b_n|\right) r^n<+\infty$.
Consider representatives $A_n$ and $B_n$ of $a_n$ and $b_n$.
There is a subset $E$ of $B$ of Lebesgue measure $0$ for which, for all $z\in B-E$ the sums $\sum |A_n(z)| r^n$ and $\sum |A_n(z)| r^n$ converge.
For those $z$, the functions $\tau\mapsto A(z,\tau) = \sum A_n(z) \tau^n$ and $B(z,\tau) = \sum B_n(z) r^n$ are holomorphic.
We define $A(z,\tau)=0$ and $B(z,\tau)=0$ when $z\in E$.
The functions $A$ and $B$ are measurable.
The partial sums of the series $\sum a_n \tau^n$ tend to 
the function $z\mapsto A(z,\tau)$ in $L^1(B)$ so the latter is a representative of $\sum a_n \tau^n$. A similar statement holds for $\sum b_n \tau^n$.
Then for all $\tau$, $z\mapsto A(z,\tau)$ is a representative of $\partial f_\tau$ and $z\mapsto B(z,\tau)$ is a representative of $\bar\partial f_\tau$. 
Moreover we know by properties of quasiconformal maps (see \cite{A} or Corollary~3 page~22 of \cite{Abis}) that for all $\tau$, $\partial f_\tau$ is non-zero on the complement of a set (that depends a priori on $\tau$) of measure $0$.
So the set of $z$ such that the function $\tau\mapsto A(z,\tau)$ is identically equal to $0$ can only have measure $0$.
So the function $B/A$ is meromorphic in $\tau$ for almost every $z$.
For all $\tau$ with $|\tau|<r$, the function $z\mapsto B(z,\tau)/A(z,\tau)$ is the quotient of representatives of respectively $\bar\partial f_\tau$ and $\partial f_\tau$, so it is almost everywhere equal to $\mu_\tau$, which is bounded by some $\kappa<1$.
By measurability, $|B/A|>\kappa$ defines a measurable set, hence
for almost every $z$, the inequality $|B/A|\leq \kappa$ holds for almost every $\tau$. In particular the meromorphic function $\tau\mapsto B(z,\tau)/A(z,\tau)$ has all its poles removable. Let us extend this function into a holomorphic function of $\tau\in B(0,r)$.
Then the extended function is measurable in $(z,\tau)$ and can serve as the function $\mu_\tau(z)$ in criterion~\eqref{aw1} of \Cref{prop:muhol}.
\end{proof}

%
% ----------------------------------------------------------
%

\section{Appendix}

\Cref{sub:interpretLim} contained a technical statement on the convergence of the differential of $f_n$ to that of $f$, which we prove here.

\subsection{Proof of \Cref{prop:supDD}}\label{app:pf:supDD}

Outside a uniform ball $\ov B(0,R')$, the maps $f$ and $f_n$ are all holomorphic, injective, and tend to infinity at infinity.
In particular their derivatives $f_n'(z)$ and $f'(z)$ all have a limit as $z\to\infty$.
Moreover $f_n$ tends to $f$ uniformly on compact subsets of $\C$, according to \Cref{app:beltrami:cv}. It follows, by properties of analytic maps, that $f_n'$ tends to $f_n$ uniformly on $\C-\ov B(0,R'+1)$. So the claim is proved near infinity and there remains to prove it on $\ov B(0,R'+1)$.

Consider as in \Cref{lem:lip}
\[ \hat f_n = f_n\circ f^{-1} .\]
Since $f$ is a $C^1$ diffeomorphism, it is enough to prove the following claim: for all $\eps>0$, and all compact subset $K$ of $\C$,
\[\sup_{\scriptstyle z\in K,\,\on{dist}(z,\on{Sing})>\eps r_n} \hskip-1cm \|D_z\hat f_n-\on{Id}\| \underset{n\to\infty}\tend 0,
\]
where we recall that $\on{Sing}$ is the set of actual singularities of $\zeta_n$ (with $\infty$ omitted).
We have the following stronger version of \Cref{lem:lip}:
\begin{lemma}\label{lem:close}
For any compact subset $K$ of $\C$ and any sequence $\eps_n$ satisfying
\[\log 1/\eps_n = o(n)\] 
there exists $b_n\tend 0$ such that for all $n$ big enough,
for any two points $u$, $v$ in $K$, satisfying $|u-v|\geq\eps_n$ we have
\[\left|\frac{\hat f_n(v)-\hat f_n(u)}{v-u}-1\right| \leq b_n.\]
\end{lemma}
\begin{proof}
We use the notations of \Cref{lem:lip} and its proof. In particular, we set $u_n=\hat f_n(u)$ and $v_n=\hat f_n(v)$. Consider again the function $\Phi:\lambda \in B(0,R_n/2) \mapsto \hat f_n(\lambda,v)-\hat f_n(\lambda,u)\in B'=B(0,2R)-\{0\}$ and recall that $R_n\geq n/M$.
We still have that the hyperbolic distance in $B(0,2R)-\{0\}$ from $\Phi(0)=v-u$ to $\Phi(1)=v_n-u_n$ is at most the hyperbolic distance in $B(0,R_n/2)$ from $0$ to $1$, which is at most $\cst/n$.
Let $\Phi(0) = r e^{i\theta}$ and $\tau = \log (2R/r) >0$.
We now use the fact that the map $z\in\H_+ \mapsto 2R\exp (i\theta-\tau z) \in B'$ preserves the respective infinitesimal hyperbolic metrics and sends $1$ to $\Phi(0)$.
Consider the lift $\delta$ of the path $t\in[0,1]\mapsto\Phi(t)$ starting from $\delta(0)=1$.
We have 
\[\frac{\Phi(1)}{\Phi(0)} = \exp(\tau\delta(0)-\tau\delta(1)).\]
The point $\delta(1)$ is, in $\H_+$, at a hyperbolic distance from $\delta(0)=1$ that is at most $\cst/n$.
So  $|\delta(1)-\delta(0)| \leq \frac{\cst}{n}$ for another constant, provided $n$ is big enough (independently of $u$ and $v$).
Multiplying by $\tau$ (which depends on $u$ and $v$), we get \[\left|\log_p \frac{\Phi(1)}{\Phi(0)}\right| \leq |\tau\delta(1)-\tau\delta(0)| \leq \frac{\cst}{n}\tau .\]
Now 
\[ \frac{\tau}{n} = \frac{1}{n}\log \left|\frac{2R}{\Phi(0)}\right| \leq \frac{\cst+\log 1/\eps_n}{n}\underset{n\to\infty}\tend 0\]
since $|\Phi(0)|\geq \eps_n$ and $\log 1/\eps_n = o(n)$.
The claim follows.
\end{proof}
We assume that $n$ is big enough so that $K$ is contained in the array formed by the small squares, which we recall is the big square given by $|\Re z|\leq n$ and $|\Im z|\leq n$.
To get some margin we in fact take $n$ be big enough so that $K$ is contained in $|\Re z|\leq n-1/n$ and $|\Im z|\leq n-1/n$.

Consider now a set of four small squares around a square corner $c$. The Beltrami form $\mu_n$ takes four values on these squares.
Consider the four quadrants defined by $\arg(z-c)$ belonging to the intervals $(k\pi/2,(k+1)\pi/2)$, $k\in\{0,1,2,3\}$.
Consider the Beltrami form $\check\mu$ that is constant on each of these quadrants and takes the same value $\check\mu_k$ as the one of the four small squares with corner $c$ that the quadrant contains (we omit $n$ in these notations).
We know an explicit straightening of $\check\mu$, obtained as follows:
first apply
\[z\mapsto z-c ,\]
then on each quadrant based on $0$, indexed by $k\in\{0,1,2,3\}$, consider the map
\[ z\mapsto z+\check\mu_k \bar z .\]
Glue the four sectors thus obtained together along three of the four boundary pairs, using $\C$-affine maps.
For instance use the identity for $k=0$, then
\[ z\mapsto \frac{i+\check\mu_0 \bar \imath}{i+\check\mu_1 \bar \imath}\,z \]
for $k=1$, then
\[ z\mapsto \frac{-1-\check\mu_1}{-1-\check\mu_2}\cdot\frac{i+\check\mu_0 \bar\imath}{i+\check\mu_1 \bar\imath}\,z \]
for $k=2$ and
\[ z\mapsto \frac{-i-\check\mu_2 \bar \imath}{-i-\check\mu_3 \bar \imath}\cdot\frac{-1-\check\mu_1}{-1-\check\mu_2}\cdot\frac{i+\check\mu_0 \bar\imath}{i+\check\mu_1 \bar\imath}\,z \]
for $k=3$.
Note that each of these maps is close to the identity because the $\check \mu_k$ are close to each other, as averages of the continuous function $\mu$ on nearby small squares, and their modulus is never close to $1$ since $\|\mu\|_\infty<1$.
We obtain a sector of angle close to $2\pi$ (it may be bigger than $2\pi$). The sector is then closed appropriately as explained in \Cref{sec:simsurfpolygon} by applying
\[ z\mapsto z^\alpha =\exp(\alpha \log z) \]
for the branch of $\log z$ whose imaginary part belongs to $[0,2\pi)$, where
\[\alpha = \frac{2\pi i}{2\pi i+\log\tau}\]
and
$\tau = \frac{1+\check \mu_3}{1+\check \mu_0}\cdot\frac{-i-\check\mu_2 \bar \imath}{-i-\check\mu_3 \bar \imath}\cdot\frac{-1-\check\mu_1}{-1-\check\mu_2}\cdot\frac{i+\check\mu_0 \bar\imath}{i+\check\mu_1 \bar\imath}$, i.e.\ 
\[\tau = \frac{1+\check \mu_3}{1-\check\mu_3}\cdot\frac{1-\check\mu_2}{1+\check\mu_2}\cdot\frac{1+\check\mu_1}{1-\check\mu_1}\cdot\frac{1-\check\mu_0}{1+\check \mu_0}\]
which is close to $1$,
and where the determination of $\log \tau$ is the principal one.
In the proof of \Cref{prop:lim-mes} we have in fact evaluated $\tau$ (it corresponds to $\Lambda$ in that proposition) and found that
$\log\tau= \cal O(1/n^2)$,
uniformly: see \cref{eq:lambda-expand}.
As a consequence, $\alpha$ is close to $1$ and
\[|\alpha-1| = \cal O(1/n^2),\]
uniformly.
Last we apply
\[ z\mapsto az+b \]
where $a$ and $b$ will be chosen later.
Denote $\Phi_n$ the composition of the maps above, starting with $z\mapsto z-c$ and ending with $z\mapsto az+b$:
\[\Phi_n: F\to \C \]
where $F=\setof{z\in\C}{|\Re (z-c)<1/n,\ |\Im(z-c)|<1/n}$ is the interior of the union of the four squares (on each square, $\Phi_n$ is hence defined as the composition of five explicit simple functions).
The purpose of $\Phi_n$ is to serve as a model of the map $f_n$ near $c$.
Note that $\Phi_n$ is quasiconformal, and is $\R$-differentiable inside each square.
\begin{lemma}\label{lem:derclose1}
Let $r_n=\exp(-n)$.
There exists $\eta_n\tend 0$ such that for all $n$ big enough, for all $c$ as above, the differential of the corresponding $\Phi_n$ on the set $\setof{z\in F}{|z-c|>r_n}$ varies less than $\eta_n$, in the sense that for any two points $z$, $z'$ in this set, the corresponding differentials $L_z$ and $L_{z'}$ satisfy\footnote{Where $\|\cdot\|$ is any operator norm on chosen in advance on the set of $\R$-linear endomorphisms of $\C$.}
$\|L_{z'}\circ L_{z}^{-1}-\on{Id}\|<\eta_n$.
\end{lemma}
\begin{proof}
Of course the differential of $z\mapsto z-c$ is the identity.
Each $z\mapsto z+\check \mu_k \bar z$ is $\R$-linear and close to $z\mapsto z+\mu(c) \bar z$, because $|\check\mu_k-\mu(c)|<\cst /n$ where $\cst$ is independent of $c$ and $n$.
The differential of $z\mapsto z^\alpha$ is $\alpha z^{\alpha-1} dz$. Recall that $\alpha$ is close to $1$ and $z^{\alpha-1} = \exp((\alpha-1) \log z)$ is close to $1$ provided $(\alpha-1)\log|z|$ is close to $0$. Since $(\alpha-1)=\cal O(1/n^2)$ uniformly, it is enough that $|z|>r_n=\exp(-n)$ (we see that we could even have taken a much smaller $r_n$).
\end{proof}

Though we do not formally need it, it helps, for a clearer mental picture of the situation, to realize that, as we can see from the proof, $\Phi_n$ is close to the similitude $z\mapsto a(z-c)+b$ in some sense that we do not need to make explicit.

\begin{figure}
\begin{tikzpicture}
\node at (0,0) {\includegraphics[scale=0.75]{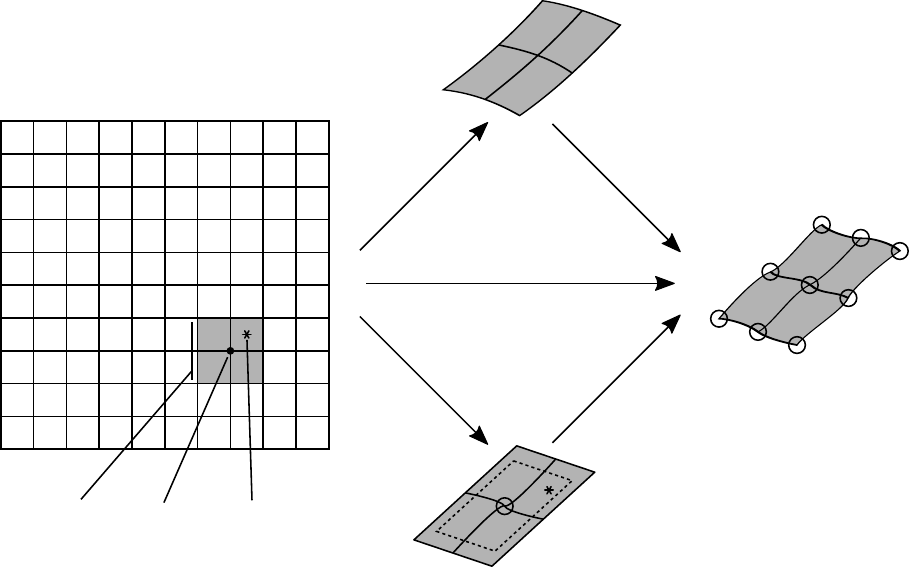}};

\node at (-0.2,1.0) {$f$};
\node at (1.8,1.0) {$\hat f_n$};
\node at (0.9,0.35) {$f_n$};
\node at (-0.15,-1.0) {$\Phi_n$};
\node at (1.8,-1.0) {$\check f_n$};

\node at (-4.9,-3) {$F$};
\node at (-3.8,-3.05) {$c$};
\node at (-2.5,-3) {$z^*$};

\node at (-0.5,-2.8) {$F'$};
\end{tikzpicture}
\caption{This commutative diagram shows the objects involved in the end of the proof of  \Cref{prop:supDD}. The model map $\Phi_n$ is (voluntarily) only defined on the four small squares depicted in gray. Circled are the places where the derivative of $\check f_n$ is not controlled, and the image of the places where the derivative of $f_n$ and $\hat f_n$ is not controlled. These places are of size of order $r_n = \exp(-n)$, so in fact much smaller than the grid step $1/n$, so one should imagine the circles nearly invisible. Dashed is the image by $\Phi_n$ of the gray square scaled about its centre by a factor of $0.75$.}
\label{fig:diag-1b}
\end{figure}
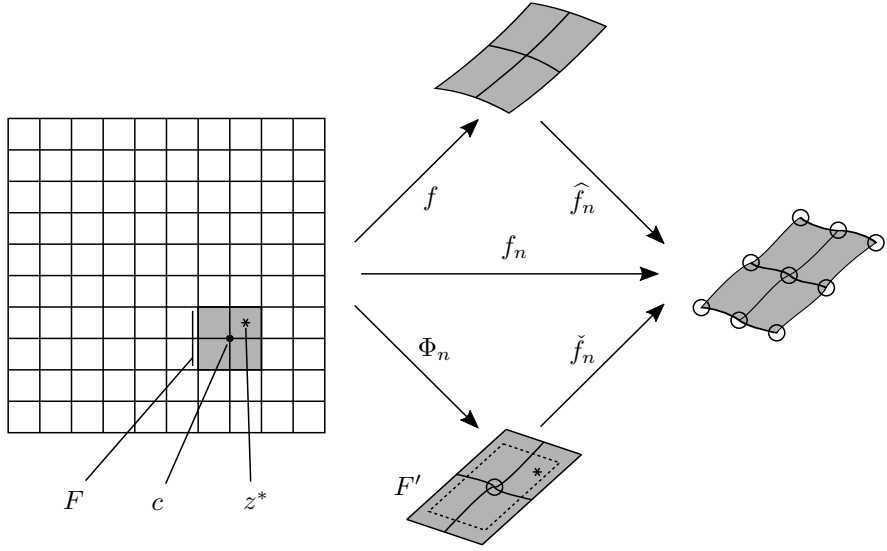

It follows from \Cref{lem:derclose1} that the differential of $\Phi_n^{-1}$ also varies less than $\eta_n$ on the image by $\Phi_n$ of $\setof{z\in F}{|z-c|>r_n}$.
The image $\Phi_n(F)$ of $\Phi_n$ is close to be a parallelogram and the part that must be removed is contained in a very small disk near its centre, of size of order $r_n$ times $\|L_z\|$, which is of the order of $|a|$ where $a$ is the constant yet to be chosen in the definition of $\Phi_n$, whereas the parallelogram has a diameter of order $|a|/n$, and bounded geometry (recall that $\|\mu\|_\infty<1$).
In the sequel we denote
\[F' = \Phi_n(\setof{z\in F}{|z-c|>r_n}) .\]

Let
$\check f_n = f_n\circ\Phi_n^{-1}$
so that 
$f_n = \check f_n\circ\Phi_n$.
Then $\check f_n$ is an injective holomorphic map because it is a composition of quasiconformal maps that sends the circle ellipse field to itself.
Moreover,
\[\check f_n = \hat f_n\circ(f\circ\Phi_n^{-1}).\]

In the definition of $\Phi_n$ there was to choose constants $a$ and $b$. We choose $b=f(c)$, so that $\Phi_n(c)=f(c)$. Consider the point $z^* =c+\frac{1+i}{2n}$, which is the centre of one of the four squares.
We choose $a$ so that $\check f_n'(\Phi_n(z^*)) = 1$.

Since $f$ is $C^1$ (actually $C^2$) and $F$ has small diameter, the differential of $f$ is nearly constant on $F$. It follows that the map $f\circ\Phi_n^{-1}$ has a nearly constant differential on $F'$.
In particular,\footnote{The set $F'$ is not convex but it is nearly a parallelogram with a removed hole that is nearly an ellipse.} pairs of points $u$, $v$ in $F'$ have images $u'$, $v'$ by $f\circ\Phi_n^{-1}$ such that $(u'-v')/L(u-v)$ is close to $1$, where $L$ is some $\R$-linear map that depends on $n$ and $c$ but not on $u$ and $v$.

We have also seen that $\hat f_n$ maps pairs of points $u,v \in K$ at distance $>r_n$ to pairs of points $u'$, $v'$ whose associated vector is nearly the same: $(u'-v')/(u-v)$ is close to $1$, uniformly.

It follows that pairs of points $u$, $v$ in $F'$ with distance at least $|a|r_n/M'$ for some $M'>0$ independent of $n$ and $c$, have images $u'$, $v'$ by $\check f_n$ such that $(u'-v')/L(u-v)$ is close to $1$.
But recall that $\check f_n$ is holomorphic and that its derivative at $\Phi_n(z^*)$ is equal to $1$. It follow by properties of univalent maps that $\check f_n$ is close to the identity on the compact subset of $\Phi_n(F'')$ of $\Phi_n(F)$ where $F''$ is defined by the equations $|\Re z-c|<0.75/n$ and $|\Im z-c|<0.75/n$, and that $\check f_n'$ is close to $1$ on that set.
It also follows that $L$ is actually close to the identity.

Recall that $f\circ\Phi_n^{-1}$ has a differential nearly constant and close to $L$ on $F'$. Hence $f\circ\Phi_n^{-1}$ has actually a differential close to the identity on $F'$.

Since, $\hat f_n = \check f_n \circ (f\circ\Phi_n^{-1})^{-1}$ we get \Cref{prop:supDD} by the chain rule.

\bibliographystyle{alpha}
\bibliography{bib} 

\begin{thebibliography}{dSLDS03}

\bibitem[AB60]{AhB}
Lars~V. Ahlfors and Lipman Bers.
\newblock Riemann's mapping theorem for variable metrics.
\newblock {\em Ann. Math. (2)}, 72:385--404, 1960.

\bibitem[AB18]{AB}
Fredric~D. Ancel and David~P. Bellamy.
\newblock On homogeneous locally conical spaces.
\newblock {\em Fundam. Math.}, 241(1):1--15, 2018.

\bibitem[AF03]{AF}
Robert~A. Adams and John J.~F. Fournier.
\newblock {\em Sobolev spaces}.
\newblock New York, NY: Academic Press, 2nd ed. edition, 2003.

\bibitem[Ahl54]{A2}
Lars~V. Ahlfors.
\newblock On quasiconformal mappings.
\newblock {\em J. Anal. Math.}, 3:1--58, 1954.

\bibitem[{Ahl}66]{A}
Lars~V. {Ahlfors}.
\newblock {Lectures on quasiconformal mappings.}
\newblock {Princeton, N.J.-Toronto-New York-London: D. Van Nostrand Company},
  1966.
\newblock First edition.

\bibitem[Ahl06]{Abis}
Lars~V. Ahlfors.
\newblock {\em Lectures on quasiconformal mappings. {With} additional chapters
  by {C}. {J}. {Earle} and {I}. {Kra}, {M}. {Shishikura} and {J}. {H}.
  {Hubbard}}, volume~38 of {\em Univ. Lect. Ser.}
\newblock Providence, RI: American Mathematical Society (AMS), 2nd enlarged and
  revised ed. edition, 2006.

\bibitem[Ahl10]{Ah4}
Lars~V. Ahlfors.
\newblock {\em Conformal invariants}.
\newblock AMS Chelsea Publishing, Providence, RI, 2010.
\newblock Topics in geometric function theory, Reprint of the 1973 original,
  With a foreword by Peter Duren, F. W. Gehring and Brad Osgood.

\bibitem[AIM09]{AIM}
Kari Astala, Tadeusz Iwaniec, and Gaven Martin.
\newblock {\em Elliptic partial differential equations and quasiconformal
  mappings in the plane}, volume~48 of {\em Princeton Math. Ser.}
\newblock Princeton, NJ: Princeton University Press, 2009.

\bibitem[Ale24]{Al}
J.~W. Alexander.
\newblock On the deformation of an {{\(n\)}}-cell.
\newblock {\em Bull. Am. Math. Soc.}, 30:10, 1924.

\bibitem[AN00]{AN}
W.~Arendt and N.~Nikolski.
\newblock Vector-valued holomorphic functions revisited.
\newblock {\em Math. Z.}, 234(4):777--805, 2000.

\bibitem[BG68]{BG}
R.~L. Bishop and S.~I. Goldberg.
\newblock Tensor analysis on manifolds.
\newblock New {York}: {The} {Macmillan} {Company}; {London}:
  {Collier}-{Macmillan} {Limited}. {VIII}, 280 p. 112 s. (1968)., 1968.

\bibitem[Boj58]{Bo1}
B.~V. Bojarskij.
\newblock Verallgemeinerte {L{\"o}sungen} eines {Systems} von
  {Differentialgleichungen} erster {Ordnung} vom elliptischen {Typus} mit
  unstetigen {Koeffizienten}.
\newblock {\em Mat. Sb., Nov. Ser.}, 43:451--503, 1958.

\bibitem[Boj09]{Bo2}
B.~V. Bojarski.
\newblock {\em Generalized solutions of a system of differential equations of
  the first order and elliptic type with discontinuous coefficients.
  {Translated} from the 1957 {Russian} original. {With} a foreword by {Eero}
  {Saksman}}, volume 118 of {\em Rep., Univ. Jyv{\"a}skyl{\"a}, Dep. Math.
  Stat.}
\newblock Jyv{\"a}skyl{\"a}: University of Jyv{\"a}skyl{\"a}, Department of
  Mathematics {and} Statistics, 2009.

\bibitem[CCL99]{CCL}
S.~S. Chern, W.~H. Chen, and K.~S. Lam.
\newblock {\em Lectures on differential geometry}, volume~1 of {\em Ser. Univ.
  Math.}
\newblock Singapore: World Scientific, 1999.

\bibitem[Ch{\'e}08]{C}
Arnaud Ch{\'e}ritat.
\newblock {B}eltrami forms, affine surfaces and the {S}chwarz-{C}hristoffel
  formula: a worked out example of straightening, 2008.
\newblock Preprint, arXiv:0811.2601v1 [math.CV].

\bibitem[dSLDS03]{dSLDS}
B.~de~Smit, H.~W. Lenstra, Jr., Douglas Dunham, and Reza Sarhangi.
\newblock Artful mathematics: the heritage of {M}. {C}. {E}scher.
\newblock {\em Notices Amer. Math. Soc.}, 50(4):446--457, 2003.

\bibitem[EK71]{EK}
R.~D. Edwards and R.~C. Kirby.
\newblock Deformations of spaces of imbeddings.
\newblock {\em Ann. Math. (2)}, 93:63--88, 1971.

\bibitem[GL00]{GL}
Frederick~P. Gardiner and Nikola Lakic.
\newblock {\em Quasiconformal {Teichm{\"u}ller} theory}, volume~76 of {\em
  Math. Surv. Monogr.}
\newblock Providence, RI: American Mathematical Society, 2000.

\bibitem[Glu08]{Glu}
Alexey~A. Glutsyuk.
\newblock Simple proofs of uniformization theorems.
\newblock In {\em Holomorphic dynamics and renormalization}, volume~53 of {\em
  Fields Inst. Commun.}, pages 125--143. Amer. Math. Soc., Providence, RI,
  2008.

\bibitem[GS68]{G}
Israel~M. Gel'fand and G.~E. Shilov.
\newblock {\em Generalized {F}unctions, {V}olume 2: {S}paces of {F}undamental
  and {G}eneralized {F}unctions}.
\newblock Academic Press, 2nd edition, 1968.

\bibitem[Gun90]{Gu}
Robert~C. Gunning.
\newblock {\em Introduction to holomorphic functions of several variables.
  {Volume} {I}: {Function} theory. {Volume} {II}: {Local} theory. {Volume}
  {III}: {Homological} theory. ({Rev}. version and complete rewriting of:
  {Analytic} functions of several complex variables by {Hugo} {Rossi} and the
  author)}.
\newblock Florence, KY: Wadsworth \&| Brooks/Cole Advanced Books \&| Software,
  1990.

\bibitem[Haz13]{Ha}
Michiel Hazewinkel, editor.
\newblock {\em Encyclopedia of mathematics. {Springer} online reference works}.
\newblock Berlin: Springer, 2013.

\bibitem[Hub06]{Hu}
John~Hamal Hubbard.
\newblock {\em Teichm{\"u}ller theory and applications to geometry, topology,
  and dynamics. {Volume} 1: {Teichm{\"u}ller} theory. {With} contributions by
  {Adrien} {Douady}, {William} {Dunbar}, and {Roland} {Roeder}, {Sylvain}
  {Bonnot}, {David} {Brown}, {Allen} {Hatcher}, {Chris} {Hruska}, {Sudeb}
  {Mitra}}.
\newblock Ithaca, NY: Matrix Editions, 2006.

\bibitem[Kat95]{Ka}
Tosio Kato.
\newblock {\em Perturbation theory for linear operators.}
\newblock Class. Math. Berlin: Springer-Verlag, reprint of the corr. print. of
  the 2nd ed. 1980 edition, 1995.

\bibitem[KN63]{KN}
Sh. Kobayashi and K.~Nomizu.
\newblock {\em Foundations of differential geometry. {I}}, volume~15 of {\em
  Intersci. Tracts Pure Appl. Math.}
\newblock Interscience Publishers, New York, NY, 1963.

\bibitem[KN69]{KN2}
Shoshichi Kobayashi and Katsumi Nomizu.
\newblock {\em Foundations of differential geometry. {Vol}. {II}}, volume~15 of
  {\em Intersci. Tracts Pure Appl. Math.}
\newblock Interscience Publishers, New York, NY, 1969.

\bibitem[Kne26]{K}
H.~Kneser.
\newblock Die {Deformationss{\"a}tze} der einfach zusammenh{\"a}ngenden
  {Fl{\"a}chen}.
\newblock {\em Math. Z.}, 25:362--372, 1926.

\bibitem[Lav35]{Lav1}
M.~A. Lavrentieff.
\newblock Sur une classe de repr{\'e}sentations continues.
\newblock {\em Rec. Math. Moscou}, 42:407--424, 1935.

\bibitem[Lav20]{Lav2}
Mikha{\"{\i}}l Lavrentieff.
\newblock On a class of continuous representations.
\newblock In {\em Handbook of Teichm\"uller theory. Volume VII}, pages
  417--439. Berlin: European Mathematical Society (EMS), 2020.

\bibitem[LL01]{LL}
Elliott~H. Lieb and Michael Loss.
\newblock {\em Analysis}, volume~14 of {\em Graduate studies in mathematics}.
\newblock American Mathematical Society, 2nd edition, 2001.

\bibitem[LV73]{LV}
O.~Lehto and K.~I. Virtanen.
\newblock {\em Quasiconformal mappings in the plane. {Translated} from the
  {German} by {K}. {W}. {Lucas}. 2nd ed}, volume 126 of {\em Grundlehren Math.
  Wiss.}
\newblock Springer, Cham, 1973.

\bibitem[Man72]{M}
Richard Mandelbaum.
\newblock Branched structures on {R}iemann surfaces.
\newblock {\em Trans. Amer. Math. Soc.}, 163:261--275, 1972.

\bibitem[Maz85]{Maz}
Vladimir~G. Maz'ya.
\newblock {\em Sobolev spaces. {Transl}. from the {Russian} by {T}. {O}.
  {Shaposhnikova}}.
\newblock Berlin etc.: Springer-Verlag, 1985.

\bibitem[Mor38]{Mo}
Charles B.~jun. Morrey.
\newblock On the solutions of quasi-linear elliptic partial differential
  equations.
\newblock {\em Trans. Am. Math. Soc.}, 43:126--166, 1938.

\bibitem[Muj86]{Muj}
Jorge Mujica.
\newblock {\em Complex analysis in {Banach} spaces. {Holomorphic} functions and
  domains of holomorphy in finite and infinite dimensions}, volume 120 of {\em
  North-Holland Math. Stud.}
\newblock Elsevier, Amsterdam, 1986.

\bibitem[Mun00]{Mun}
James~R. Munkres.
\newblock {\em Topology.}
\newblock Upper Saddle River, NJ: Prentice Hall, 2nd ed. edition, 2000.

\bibitem[Neh52]{Ne}
Zeev Nehari.
\newblock Conformal mapping.
\newblock ({International} {Series} in {Pure} and {Applied} {Mathematics}).
  {New} {York}-{Toronto}- {London}: {McGraw}-{Hill} {Book} {Company}, {Inc}.
  {VIII}, 396 p. (1952)., 1952.

\bibitem[Ric69]{Ri}
S.~Rickman.
\newblock Removability theorems for quasiconformal mappings.
\newblock {\em Ann. Acad. Sci. Fenn., Ser. A I}, 499:8, 1969.

\bibitem[Rob55]{Ro}
J.~H. Roberts.
\newblock {L}ocal arcwise connectivity in the space ${H}^n$ of homeomorphism of
  ${S}^n$ onto itself.
\newblock In {\em Summer Institute on Set Theoretic Topology, Summary of
  Lectures}, page 100. American Mathematical Society, Madison, Wisconsin, 1955.

\bibitem[Rud87]{R}
Walter Rudin.
\newblock {\em Real and Complex Analysis}.
\newblock Mathematics series. McGraw-Hill International Editions, 3rd edition,
  1987.

\bibitem[Sch66]{Sc}
Laurent Schwartz.
\newblock {\em Th{\'e}orie des distributions}.
\newblock Hermann, 1966.

\bibitem[Sha97]{Sh}
R.~W. Sharpe.
\newblock {\em Differential geometry: {Cartan}'s generalization of {Klein}'s
  {Erlangen} program. {Foreword} by {S}. {S}. {Chern}}, volume 166 of {\em
  Grad. Texts Math.}
\newblock Berlin: Springer, 1997.

\bibitem[Spi99]{Spi}
Michael Spivak.
\newblock {\em A comprehensive introduction to differential geometry. {Vol}.
  1-5}.
\newblock Houston, TX: Publish or Perish, 3rd ed. with corrections edition,
  1999.

\bibitem[Ste83]{Ste}
Shlomo Sternberg.
\newblock Lectures on differential geometry. 2nd ed.
\newblock New {York}: {Chelsea} {Publishing} {Company}. xv, 442 pp. (1983).,
  1983.

\bibitem[Str94]{St}
Robert Strichartz.
\newblock {\em A guide to distribution theory and {F}ourier transforms}.
\newblock Studies in Advanced Mathematics. CRC-Press, 1994.

\bibitem[Thu97]{T}
William~P. Thurston.
\newblock {\em Three-dimensional geometry and topology. {V}ol. 1}, volume~35 of
  {\em Princeton Mathematical Series}.
\newblock Princeton University Press, Princeton, NJ, 1997.
\newblock Edited by Silvio Levy.

\bibitem[Vee93]{Ve}
William~A. Veech.
\newblock Flat surfaces.
\newblock {\em Am. J. Math.}, 115(3):589--689, 1993.

\bibitem[Vek53]{V2}
I.~N. Vekua.
\newblock A boundary problem with oblique derivative for an equation of
  elliptic type.
\newblock {\em Doklady Akad. Nauk SSSR (N.S.)}, 92:1113--1116, 1953.

\bibitem[Vek55]{V}
I.~N. Vekua.
\newblock The problem of reduction to canonical form of differential forms of
  elliptic type and the generalized {C}auchy-{R}iemann system.
\newblock {\em Dokl. Akad. Nauk SSSR (N.S.)}, 100:197--200, 1955.

\bibitem[Zor06]{Z}
Anton Zorich.
\newblock Flat surfaces.
\newblock In {\em Frontiers in number theory, physics, and geometry I. On
  random matrices, zeta functions, and dynamical systems. Papers from the
  meeting, Les Houches, France, March 9--21, 2003}, pages 437--583. Berlin:
  Springer, 2nd printing edition, 2006.

\end{thebibliography}

\end{document}